\documentclass[a4paper,11pt]{article}
\usepackage{stmaryrd}
\usepackage{amsfonts}
\usepackage{bbm}
\usepackage{float}
\usepackage{amscd}
\usepackage{mathrsfs}
\usepackage{latexsym,amssymb,amsmath,amscd,amscd,amsthm,amsxtra,xypic}
\usepackage[dvips]{graphicx}
\usepackage[utf8]{inputenc}
\usepackage[T1]{fontenc}
\usepackage{lmodern}
\usepackage{amssymb}
\usepackage[all]{xy}
\usepackage{nicefrac,mathtools,enumitem}
\usepackage{microtype}
\usepackage{lipsum}
\usepackage{mathtools}
\usepackage{xfrac}

\usepackage{geometry}
\usepackage{hyperref}

\usepackage{times}
\usepackage{amsthm}
\usepackage{amsmath}
\usepackage{amsfonts}
\usepackage{amssymb}
\usepackage{mathdots}
\usepackage{color}
\usepackage{mathrsfs}
\usepackage{stmaryrd}
\SetSymbolFont{stmry}{bold}{U}{stmry}{m}{n}

\usepackage{bm}

\usepackage{tikz-cd}
\usepackage{tikz}
\usetikzlibrary{matrix,arrows,decorations.pathmorphing}
\usepackage[all]{xy}
\usepackage{graphicx}
\usepackage{subfigure}
\usepackage[toc,page]{appendix}
\usepackage[usestackEOL]{stackengine}

\usepackage{theoremref}
\newtheorem{thm}{Theorem}[section]
\newtheorem{lem}[thm]{Lemma}
\newtheorem{cor}[thm]{Corollary}
\newtheorem{pro}[thm]{Proposition}

\newtheorem{conj}[thm]{Conjecture} 

\theoremstyle{definition}
\newtheorem{ex}[thm]{Example}
\newtheorem{rmk}[thm]{Remark}
\newtheorem{defi}[thm]{Definition}

\newcommand{\be }{\begin{equation}}
\newcommand{\ee }{\end{equation}}

\newcommand{\pf}{\noindent{\bf Proof.}\ }

\newcommand{\CG}{{\mathcal{G}}}

\newcommand{\h}{\mathfrak h}
\newcommand{\frkt}{\mathfrak t}

\newcommand{\U}{{\rm U}}

\def\qed{\hfill ~\vrule height6pt width6pt depth0pt}

\newcommand{\br}[1]{   [ \cdot,    \cdot  ]   }

\newcommand{\g}{\mathfrak g}

\newcommand{\Ad}{\mathrm{Ad}}

\newcommand{\gl}{\mathfrak {gl}}

\newcommand{\I}{{\mathrm{i}}}

\newcommand{\D}{\mathbb{D}}

\newcommand {\A}{\mathcal A}

\newcommand{\Herm}{{\rm Herm}}

\title{Regularized limits of Stokes matrices, isomonodromy deformation and crystal basis}
\author{Xiaomeng Xu}
\date{}

\newcommand{\Addresses}{{
  \bigskip
  \footnotesize
\noindent \textsc{
School of Mathematical Sciences \& Beijing International Center
for Mathematical Research, Peking University, Beijing 100871, China}\par\nopagebreak
  \textit{E-mail address}: \texttt{xxu@bicmr.pku.edu.cn}
}}
\begin{document}

\maketitle
    
\begin{abstract}
In the first part of the paper, we solve the boundary and monodromy problems for the isomonodromy equation of the $n\times n$ meromorphic linear system of ordinary differential equations with Poncar\'{e} rank $1$. In particular, we derive an explicit expression of the Stokes matrices of the linear system, via the boundary value of the solutions of the isomonodromy equation at a critical point. Motivated by this result, we then describe the regularized limits of Stokes matrices as the irregular data $u={\rm diag}(u_1,...,u_n)$ in the linear system degenerates, i.e., as some $u_i, u_j,...,u_k$ collapse. 
The prescription of the regularized limit is controlled by the geometry of the De Concini-Procesi wonderful compactification space. As applications, many analysis problems about higher rank Painlev\'e transcendents can be solved. 

In the second part of the paper, we show some important applications of the above analysis results in representation theory and Poisson geometry: we obtain the first transcendental realization of crystals in representations of $\frak{gl}_n$ via the Stokes phenomenon in the WKB approximation; we develop a wall-crossing formula that characterizes the discontinuous jump of the regularized limits of Stokes matrices as crossing walls in the compactification space, and interpret the known cactus group actions on crystals arising from representation theory as a wall-crossing phenomenon; and we find the first explicit linearization of the standard dual Poisson Lie group for $U(n)$. 
\end{abstract}

 \tableofcontents 
\section{Introduction and main results}
In recent years there has been considerable interests in analyzing the Stokes matrices of a meromorphic linear systems of ordinary differential equations at a second order pole, and the associated isomonodromy deformation equations. The sources of these interests are
quite diverse, which include many subjects in mathematics and physics, like the theory of Gromov-Witten and Frobenius manifolds \cite{Dubrovin, Givental, GGI}, stability conditions \cite{Bridgeland}, higher Painlev\'e transcendents \cite{JMMS, Miwa}, Poisson groups and simple models for the wild non-abelian Hodge correspondence \cite{Boalch1, BB}, representation theory of quantum groups \cite{Xu2, Xu4}, and so on.

Yet some basic properties of the Stokes matrices and the solutions of associated isomonodromic deformation equations, as special functions, still remain almost unexplored, like their behavior at singularities, their WKB approximation, and their explicit expressions and so on. These unresolved problems are obstacles to the substantial role of the Stokes phenomenon in the related areas. In this paper, these problems are studied. Most importantly, an explicit expression of the Stokes matrices is derived via the isomonodromy approach. It provides us a manipulable analysis model for the study of the Stokes phenomenon itself and the relations with other subjects. Some applications of the model in representation theory and Poisson geometry are given in this paper: among them, the most noteworthy ones are the transcendental realization of crystal basis and the cactus group action. More applications of the model are presented in a series of follow-up works, some of them can be found in \cite{ANXZ, LX, TangXu, WXX, Xu3, Xu6, Xu7}. 

Other than various new transcendental realizations of the known algebraic structures, including Gelfand-Tsetlin basis, quantum groups, crystals and cactus groups given in this paper, we are more interested in the role of the algebraic structures in the study of the Stokes phenomenon itself. For example, the representation theoretic data can be used to characterize the Stokes phenomenon in the WKB approximation of the meromorphic ODEs, which is still an open problem. Additionally, they help to study some basic problems of the isomonodromy equations. So another main idea conveyed in this article is the importance of algebraic structures in understanding the Stokes phenomenon.

\subsection{Explicit expression of Stokes matrices via the boundary value}\label{firstsec}
Let $\h_{\rm reg}(\mathbb{R})$ denote the space of $n\times n$ diagonal matrices $u={\rm diag}(u_1,...,u_n)$ with distinct real eigenvalues. Let $\Herm(n)$ denote the space of $n\times n$ Hermitian matrices. 
Let us consider the $n\times n$ linear system of partial differential equations for a function $F(z,u_1,...,u_n)\in {\rm GL}(n)$ 
\begin{align}\label{introisoStokeseq1}
\frac{\partial F}{\partial z}&=\left(\I u-\frac{1}{2\pi\I }\frac{\Phi(u)}{z}\right)\cdot F,\\
\label{introisoStokeseq2}
\frac{\partial F}{\partial u_k}&=\left(\I E_kz-\frac{1}{2\pi\I } {\rm ad}^{-1}_u{\rm ad}_{E_{k}}\Phi(u)\right)\cdot F, \ \text{for all} \ k=1,...,n,
\end{align}
where the residue $\Phi(u)=\Phi(u_1,...,u_n)\in {\rm Herm}(n)$ is a solution of the {\it isomonodromy differential equation}
\begin{eqnarray}\label{introisoeq}
\frac{\partial \Phi}{\partial u_k} =\frac{1}{2\pi\I }[\Phi,{\rm ad}^{-1}_u{\rm ad}_{E_{k}}\Phi], \ \text{for all} \ k=1,...,n.\end{eqnarray} 
Here $\I=\sqrt{-1}$, $E_k$ is the $n\times n$ diagonal matrix whose $(k,k)$-entry is $1$ and other entries are $0$. And $u$ in the coefficients of the linear system represents the diagonal matrix ${\rm diag}(u_1,\ldots,u_n)$.
Note that ${\rm ad}_{E_{k}}\Phi$ takes values in the space ${\frak {gl}}_n^{od}$ of off diagonal matrices and that the adjoint operator ${\rm ad}_u$
is invertible when restricted to ${\frak {gl}}_n^{od}$. One checks that \eqref{introisoeq} is the compatibility condition of the linear system. Our first result solves a Riemann-Hilbert problem of the linear system. In particular, the following Theorem \ref{isomonopro} gives a natural parameterization of the linear systems by the boundary value of $\Phi(u)$, and Theorem \ref{mainthm} expresses the monodromy data explicitly via the parameterization. Let us now give more details.

For any fixed $u\in\h_{\rm reg}(\mathbb{R})$, the ordinary differential equation \eqref{introisoStokeseq1} has a unique formal solution $\hat{F}(z)$ around
$z = \infty$. Then the standard theory of resummation states that there exist certain  sectorial regions around $z=\infty$, such that on each of these sectors there is a unique (therefore canonical) holomorphic solution with the prescribed asymptotics $\hat{F}(z)$. These solutions are in general different (that reflects the Stokes phenomenon), and the transition between them can be measured by a pair of Stokes matrices $S_\pm(u,\Phi(u))\in{\rm GL}(n)$. The Stokes matrices $S_+$ and $S_-$ are upper and lower triangular matrices and, due to the real condition $u\in \frak t_{\rm reg}(\mathbb{R})$ and $A\in\Herm(n)$, are complex conjugate to each other. See Section \ref{beginsection} for more details. Varying $u$, the Stokes matrices $S_\pm(u,\Phi(u))\in {\rm GL}(n)$ of the system are constant (independent of $u$), and this is why the equation \eqref{introisoeq} is called isomonodromy. 

Following Miwa \cite{Miwa}, the $\frak{gl}_n$-valued solutions $\Phi(u)$ of the equation \eqref{introisoeq} with $u_1,...,u_n\in \mathbb{C}$ have the strong Painlev\'{e} property: they are multi-valued meromorphic functions of $u_1,...,u_n$ and the branching occurs when $u$ moves along a loop around the fat diagonal
\[\Delta=\{(u_1,...,u_n)\in \mathbb{C}^n~|~u_i = u_j, \text{for some } i\ne j \}.\]
Thus, according to the original idea of Painlev\'{e}, they can be a new class of special functions. The problem of determining their behavior at the fixed
critical singularities is left open. The following theorem treats this problem for the $\Herm(n)$-valued solutions. First according to Boalch \cite{Boalch1}, the $\Herm(n)$-valued solutions $\Phi(u)$ of \eqref{introisoeq} are real analytic on each connected component of $u\in \h_{\rm reg}(\mathbb{R})$. In Section \ref{asyisoeq} and \ref{proof11}, we prove
\begin{thm}\label{isomonopro}
For any solution $\Phi(u)$ of the isomonodromy equation \eqref{introisoeq} on the connected component $U_{\rm id}:=\{u\in \h_{\rm reg}(\mathbb{R})~|~u_1<\cdots <u_n\}$, there exists a unique constant $\Phi_0\in\Herm(n)$ such that as the real numbers $\frac{u_{k+1}-u_{k}}{u_{k}-u_{k-1}}\rightarrow +\infty$ for all $k=2,...,n-1$,
\begin{align}\label{firstasy}
\Phi(u)=&{\rm Ad}{\left((u_2-u_1)^{\frac{\delta_1(\Phi_0)}{2\pi\I }}\cdot \overrightarrow{\underset{k=2,...,n-1}{\prod} }\left(\frac{u_{k+1}-u_{k}}{u_{k}-u_{k-1}}\right)^{\frac{\delta_k(\Phi_0)}{2\pi\I }}\right)}\Phi_0 +\sum_{k=2}^{n-1}\mathcal{O}\left(\frac{u_{k}-u_{k-1}}{u_{k+1}-u_{k}}\right),
\end{align}
where ${\rm Ad}(g)X=gXg^{-1}$ for any $g\in U(n)$ and $X\in \Herm(n)$, the product $\overrightarrow{\prod}$ is taken with the index $i$ to the right of $j$ if $i>j$. And $\delta_k(\Phi)$ is the Hermitian matrix with entries
\[\delta_k(\Phi)_{ij}=\left\{
          \begin{array}{lr}
             \Phi_{ij},   & \text{if} \ \ 1\le i, j\le k, \ \text{or} \ i=j  \\
           0, & \text{otherwise}.
             \end{array}
\right.\] (Here we use the big $\mathcal{O}$ notation, in particular each $\mathcal{O}\left(\frac{u_{k}-u_{k-1}}{u_{k+1}-u_{k}}\right)$ stands for a remainder whose norm is less than $M \times\left(\frac{u_{k+1}-u_{k}}{u_{k}-u_{k-1}}\right)^{-1}$ for a positive real number $M$ as $\frac{u_{k+1}-u_{k}}{u_{k}-u_{k-1}}$ big enough.) Furthermore, given any $\Phi_0\in\Herm(n)$ there exists a unique real analytic solution $\Phi(u)$ of \eqref{introisoeq} with the prescribed asymptotics \eqref{firstasy}.
\end{thm}
From the isomonodromy equation, we see that any solution $\Phi(u)$ is invariant under the translation action on $\frak{h}_{\rm reg}(\mathbb{R})$, that is $\Phi(u_1,...,u_n)=\Phi(u_1+c,...,u_n+c)$. Let $\mathbb{R}$ act on $\frak{h}_{\rm reg}(\mathbb{R})$ by translation, then $\Phi(u)$ is defined on $\frak{h}_{\rm reg}(\mathbb{R})/\mathbb{R}$. Let $u_{\rm cat}$ denote the following limit \begin{equation}\label{limit}
    u_2-u_1\rightarrow 0+ \text{ and } \frac{u_j-u_{j-1}}{u_{j-1}-u_{j-2}}\rightarrow +\infty \text{ for all } j=3,...,n.
\end{equation}
The limit is a point, called a caterpillar point in the literature \cite{Sp}, in the De Concini-Procesi space that is a certain compactification of
$\frak{h}_{\rm reg}(\mathbb{R})/\mathbb{R}$, see Section \ref{closuresec}. In this paper, we do not distinguish the limit \eqref{limit} and the point in the De Concini-Procesi space. 
We call the regularized limit $\Phi_0\in\Herm(n)$ the boundary value at $u_{\rm cat}$, and denote by $\Phi(u;\Phi_0)$ the unique solution of \eqref{introisoeq} on $U_{\rm id}$ with the prescribed boundary value $\Phi_0$.
Following Theorem \ref{isomonopro}, the boundary value $\Phi_0\in\Herm(n)$ describes the leading asymptotics of $\Phi(u;\Phi_0)$ as $u\rightarrow u_{\rm cat}$. 

In Section \ref{exStokesviaiso}, we solve the monodromy problem of the linear system \eqref{introisoStokeseq1}-\eqref{introisoStokeseq2} with $\Phi(u)=\Phi(u;\Phi_0)$. 
We denote by $\{\lambda^{(k)}_i\}_{i=1,...,k}$ the eigenvalues of the left-top $k\times k$ submatrix of $\Phi_0\in\Herm(n)$, and $(\Phi_0)_{kk}$ the $k$-th diagonal element. Then
\begin{thm}\label{mainthm}
The sub-diagonal entries of the Stokes matrices $S_\pm\left(u,\Phi(u;\Phi_0)\right)$ of the linear system \eqref{introisoStokeseq1} are given by
\begin{align*}
(S_+)_{k,k+1}&=2\pi\I\cdot {\rm exp}{\left(\frac{\small{(\Phi_0)_{kk}+(\Phi_0)_{k+1,k+1}}}{4}\right)} \\
\times &\sum_{i=1}^k\frac{\prod_{l=1,l\ne i}^{k}\Gamma\left(1+\frac{\lambda^{(k)}_l-\lambda^{(k)}_i}{2\pi \I }\right)}{\prod_{l=1}^{k+1}\Gamma\left(1+\frac{\lambda^{(k+1)}_l-\lambda^{(k)}_i}{2\pi \I }\right)}\frac{\prod_{l=1,l\ne i}^{k}\Gamma\left(\frac{\lambda^{(k)}_l-\lambda^{(k)}_i}{2\pi \I }\right)}{\prod_{l=1}^{k-1}\Gamma\left(1+\frac{\lambda^{(k-1)}_l-\lambda^{(k)}_i}{2\pi \I }\right)}\cdot {\Delta^{1,...,k-1,k}_{1,...,k-1,k+1}\left(\frac{\Phi_0-\lambda^{(k)}_i}{2\pi\I }\right)},\\
(S_-)_{k+1,k}&=-2\pi\I \cdot {\rm exp}{\left(\frac{\small{(\Phi_0)_{kk}+(\Phi_0)_{k+1,k+1}}}{4}\right)}\\
\times &
\sum_{i=1}^k \frac{\prod_{l=1,l\ne i}^{k}\Gamma\left(1-\frac{\lambda^{(k)}_l-\lambda^{(k)}_i}{2\pi \I }\right)}{\prod_{l=1}^{k+1}\Gamma\left(1-\frac{\lambda^{(k+1)}_l-\lambda^{(k)}_i}{2\pi \I }\right)}\frac{\prod_{l=1,l\ne i}^{k}\Gamma\left(-\frac{\lambda^{(k)}_l-\lambda^{(k)}_i}{2\pi \I }\right)}{\prod_{l=1}^{k-1}\Gamma\left(1-\frac{\lambda^{(k-1)}_l-\lambda^{(k)}_i}{2\pi \I }\right)}\cdot {\Delta^{1,...,k-1,k+1}_{1,...,k-1,k}\left(\frac{\lambda^{(k)}_i-\Phi_0}{2\pi\I }\right)}.
\end{align*}
where $k=1,...,n-1$ and $\Delta^{1,...,k-1,k}_{1,...,k-1,k+1}(\frac{\Phi_0-{\lambda^{(k)}_i}}{2\pi\I })$ is the $k$ by $k$ minor of the matrix $\frac{1}{2\pi\I }(\Phi_0-{\lambda^{(k)}_i}\cdot {\rm Id}_{n})$ formed by the first $k$ rows and $1,...,k-1,k+1$ columns (here ${\rm Id}_n$ is the rank $n$ identity matrix). 
Furthermore, the other entries are determined by the sub-diagonal ones in a systematic way, and thus are also given by explicit expressions, see Section \ref{classicalRLL}.
\end{thm}
\begin{rmk}
For the special case $n=3$, Theorem \ref{isomonopro} and \ref{mainthm} recover Jimbo's asymptotic and monodromy formula for the Painlev\'e VI transcendents, see \cite{WXX} for more details. From this viewpoint, they provide the asymptotic and monodromy formula for the higher Painlev\'e transcendents.

We also remark that the asymptotics of $\Phi(u)$ in Theorem \ref{isomonopro} is only a local analysis, while many global properties of $\Phi(u;\Phi_0)$ can be obtained from Theorem \ref{mainthm}. That is the power of the Riemann-Hilbert approach to the study of the nonlinear differential equations. 
\end{rmk}

Theorem \ref{isomonopro} gives a parameterization of the Hermitian matrix valued solutions of the isomonodromy equation \eqref{introisoeq}, and Theorem \ref{mainthm} computes explicitly the Stokes matrices of the corresponding linear equation \ref{introisoStokeseq1} via the parameterization. Therefore, we obtain an explicit Riemann-Hilbert-Birkhoff map (a diffeomorphism) from the space $\Herm(n)$ of meromorphic linear systems to the space of Stokes matrices via the equivalences
\begin{eqnarray*}
    \Big\{\Phi_0\in\Herm(n)\Big\}&\Longleftrightarrow &\Big\{\text{solutions $\Phi(u;\Phi_0)\in\Herm(n)$ of the isomonodromy equation \eqref{introisoeq} on $U_{\rm id}$}\Big\}\\ &\Longleftrightarrow&  \Big\{\text{linear systems of PDEs \eqref{introisoStokeseq1} and \eqref{introisoStokeseq2}}\Big\} \\  & \Longleftrightarrow & \Big\{\text{space of Stokes matrices $S_\pm\left(u,\Phi(u;\Phi_0)\right)$}\Big\}. 
\end{eqnarray*}
The Poisson geometric nature of the diffeomorphism is unveiled in Theorem \ref{eBoalchthm}.

\subsection{Regularized limits of Stokes matrices at a caterpillar point}

Consider a linear system of meromorphic ordinary differential equations for a function $F(z)\in {\rm GL}(u,\mathbb{C})$
\begin{eqnarray}\label{introeq}
\frac{dF}{dz} = \left(\I u-\frac{1}{2\pi\I }\frac{A}{z}\right)\cdot F,\end{eqnarray}
where $u\in \h_{\rm reg}(\mathbb{R})$ and $A\in\Herm(n)$. For a generic $A\in\Herm(n)$, the limit of the Stokes matrices $S_\pm(u,A)$, as some components $u_i$ of $u$ collapse, do not exist. See the rank 2 example in Example \ref{rank2ex}. Motivated by Theorem \ref{isomonopro} and \ref{mainthm}, the behaviour of the Stokes matrices at the singularities is studied in this paper. For example, a manipulation of Theorem \ref{isomonopro} and \ref{mainthm} shows that (see Section \ref{leadterm} for a proof)

\begin{pro}\label{StokesinGZ}
For any $A\in\Herm(n)$, the sub-diagonal entries of $S_+(u,A)$, as $\frac{u_{k+1}-u_{k}}{u_{k}-u_{k-1}}\rightarrow +\infty$ for all $k=2,...,n-1$, are
\begin{align*}
&(S_+)_{k,k+1}=2\pi\I e^{\frac{(A)_{kk}+(A)_{k+1,k+1}}{4}}\\ & \times \sum_{i=1}^k\frac{\prod_{l=1,l\ne
i}^{k}\Gamma\left(1+\frac{\lambda^{(k)}_l-\lambda^{(k)}_i}{2\pi \I }\right)}{\prod_{l=1}^{k+1}\Gamma\left(1+\frac{\lambda^{(k+1)}_l-\lambda^{(k)}_i}{2\pi \I }\right)}\frac{\prod_{l=1,l\ne i}^{k}\Gamma\left(\frac{\lambda^{(k)}_l-\lambda^{(k)}_i}{2\pi \I }\right)}{\prod_{l=1}^{k-1}\Gamma\left(1+\frac{\lambda^{(k-1)}_l-\lambda^{(k)}_i}{2\pi \I }\right)}\Delta^{1,...,k-1,k}_{1,...,k-1,k+1}\left(\frac{A-{\lambda^{(k)}_i}}{2\pi\I }\right)
\cdot f^{(k)}_i(A)\\
& \ \ \ +\sum_{k=2}^{n-1}\mathcal{O}\left(\left(\frac{u_{k+1}-u_{k}}{u_{k}-u_{k-1}}\right)^{-1}\right),
\end{align*}
where the singular part
\begin{equation}\label{fastspin}
f^{(k)}_i(A)=\left\{
          \begin{array}{lr}
            (u_k-u_{k-1})^{\frac{\lambda^{(k)}_i(A)-A_{kk}}{2\pi\I }} (u_{k+1}-u_{k})^{\frac{A_{k+1,k+1}-\lambda^{(k)}_i(A)}{2\pi\I }},   & \text{if} \ \ k\ge 2, \\
           (u_{2}-u_{1})^{\frac{A_{22}-\lambda^{(1)}_1(A)}{2\pi\I }}, & \text{if} \ \ k=1.
             \end{array}
\right.
\end{equation}
Here recall $\{\lambda^{(k)}_i\}_{i=1,...,k}$ denote the eigenvalues of the left-top $k\times k$ submatrix of $A$.
Similarly, the leading terms of all entries of $S_\pm(u,A)$ can be given explicitly.
\end{pro}
Following Proposition \ref{StokesinGZ}, the fast spin terms \eqref{fastspin}
(whose norm is $1$ since $\lambda^{(k)}_i$ are real numbers) in the expression prevent the Stokes matrix $S_+(u,A)$ from having a limit. 
However, the upshot is that $S_\pm(u,A)$ have properly regularized limit at $u_{\rm cat}$, i.e., as all $u_1,...,u_n$ collapse in the speed that $u_i$ approaches to $u_{i-1}$ much faster than $u_{i-1}$ approaches to $u_{i-2}$
. To be more precise, for any $u$ and $A$, let us introduce the unitary matrix
\begin{equation}\label{unig}
G(u,A)=({u_{2}-u_{1}})^{\frac{{\rm log}\left(\delta_1(S_-(u,A))\delta_1(S_+(u,A))\right)}{2\pi\I }}\cdot\overrightarrow{\underset{k=2,...,n-1}{\prod} }\left(\frac{u_{k+1}-u_{k}}{u_{k}-u_{k-1}}\right)^{\frac{{\rm log}\left(\delta_k(S_-(u,A))\delta_k(S_+(u,A))\right)}{2\pi\I }},
\end{equation}
where ${\rm log}(\delta_k(S_-)\delta_k(S_+))$ is the logarithm of the positive definite Hermitian matrix $\delta_k(S_-(u,A))\delta_k(S_+(u,A))$ (see Section \ref{beginsection} for the positive definiteness). Then
we have (see Section \ref{pf:introthm3} for a proof)
\begin{thm}\label{introcor}
For $A\in\Herm(n)$, the limit of the matrix valued function
\begin{equation}\label{gaugeRHB}
\lim_{u\rightarrow u_{\rm cat}} G(u,A)\left(S_-(u,A)S_+(u,A)\right) G(u,A)^{-1}=S_{-}\left(u, \Phi(u;A)\right)S_{+}\left(u, \Phi(u;A)\right).
\end{equation}
Here by definition, $\Phi(u;A)$ is the solution of isomonodromy equation on $U_{\rm id}$ with the given boundary value $A$ at $u_{\rm cat}$. Therefore, by Theorem \ref{mainthm} the limit, i.e., the right hand side of \eqref{gaugeRHB}, as a function of $A$, has a closed formula (provided replacing $\Phi_0$ by $A$ in Theorem \ref{mainthm}).
\end{thm}
The regularization of $S_-(u,A)S_+(u,A)$ by the unitary matrix given in \eqref{unig} has a geometric interpretation as a fast spin on Liouville tori of the Gelfand-Testlin integrable system, see Section \ref{pf:introthm3}. The theorem states that $S_\pm(u,A)$ have regularized limits $S_{\pm}(u, \Phi(u;A))$ (independent of $u$) at $u_{\rm cat}$.
Based on Theorem \ref{introcor}, we introduce an important notion (the meaning of a caterpillar point will be clear in the next subsection).
\begin{defi}\label{caterpillar}
For any $A\in \Herm(n)$, we call the regularized limits $S_{\pm}(u_{\rm cat}, A):=S_{\pm}\left(u, \Phi(u;A)\right)$ the Stokes matrices at the caterpillar point $u_{\rm cat}$ with respect to the connected component $U_{\rm id}$.
\end{defi}

\subsection{Regularized limit of Stokes matrices in the De Concini-Procesi space and wall-crossing formula}\label{closuresec}
The results in this subsection are not necessary for the second main result of this paper, i.e., the following Application II. 

On the one hand, by definition, the Stokes matrices $S_\pm(u,A)$ are invariant under the translation action on $\frak{h}_{\rm reg}(\mathbb{R})$. Thus let $\mathbb{R}$ act on $\frak{h}_{\rm reg}(\mathbb{R})$ by translation, then for any fixed $A$, $S_\pm(u,A)$ are parameterized by $\frak t_{\rm reg}(\mathbb{R})\cong \frak{h}_{\rm reg}(\mathbb{R})/\mathbb{R}$. Here $\frak t_{\rm reg}(\mathbb{R})$ is the space of $n\times n$ diagonal matrices $u={\rm diag}(u_1,...,u_n)$ with distinct real eigenvalues and $\sum_{i=1}^n u_i=0$. On the other hand, any solution $\Phi(u)$ is also invariant under the translation action on $\frak{h}_{\rm reg}(\mathbb{R})$. Therefore, we can actually assume the irregular data $u$ is in the subspace $\frak t_{\rm reg}(\mathbb{R})\subset \h_{\rm reg}(\mathbb{R})$. In this paper, for simplicity, we use $U_{\sigma}=\{u_{\sigma(1)}<\cdots <u_{\sigma(n)}\}$ for some $\sigma\in S_n$ to denote the connected components both in $\h_{\rm reg}(\mathbb{R})$ and $\frak t_{\rm reg}(\mathbb{R})$ according to the context.

In Section \ref{closureStokes} we study of the regularized limit of the Stokes matrices $S_\pm(u,A)$ of equation \eqref{introeq}, as some components $u_i$ of $u={\rm diag}(u_1,...,u_n)\in \frak t_{\rm reg}(\mathbb{R})$ collapse in a comparable speed. It completely describes the asymptotic behaviour of the Stokes matrices at the singularities $u_i=u_j$.
And the prescription of the regularized limits is controlled by the geometry of the De Concini-Procesi wonderful space $\widetilde{\frak t_{\rm reg}}(\mathbb{R})$. Let us give a brief introduction. Here given a finite set of subspaces of a vector space, the De Concini-Procesi space \cite{dCP} replaces the set of subspaces by a divisor with normal crossings, and leaves the complement of these subspaces unchanged. As for the root hyperplanes of type $A$ Lie algebra, 
the associated De Concini-Procesi space $\widetilde{\frak t_{\rm reg}}(\mathbb{R})$ contains $\frak t_{\rm reg}(\mathbb{R})$ as an open part, and 
roughly speaking, a point in the boundary $\widetilde{\frak t_{\rm reg}}(\mathbb{R})\setminus{\frak t_{\rm reg}}(\mathbb{R})$ is a limit point $u={\rm diag}(u_1,...,u_n)$, where some $u_i$ collapse in a comparable speed. See Section \ref{closureStokes} for more details. 
In particular, the limit of $u={\rm diag}(u_1,...,u_n)$, as $\frac{u_{k+1}-u_{k}}{u_{k}-u_{k-1}}\rightarrow +\infty$ for all $k=2,...,n-1$ and $u_2-u_1\rightarrow 0$, is a point $u_{\rm cat}$ in the $0$-dimensional stratum of $\widetilde{\frak t_{\rm reg}}(\mathbb{R})$, called a caterpillar (see \cite{Sp}, page 16). There are many caterpillar points, and $u_{\rm cat}$ is one of them.

Note that $u$ can approach to a same point in the boundary $\widetilde{\frak t_{\rm reg}}(\mathbb{R})\setminus{\frak t_{\rm reg}}(\mathbb{R})$, for example $u_{\rm cat}$, from different connected components of ${\frak t_{\rm reg}}(\mathbb{R})$.  As $u$ approaches to a fixed boundary point from two different components, the regularized limits can be different, see Proposition \ref{introwall}. Therefore, when study the regularized limits of Stokes matrices at a boundary point, we should specialize the connected component from which we take the limit. For example, in this paper by saying $u\rightarrow u_{\rm cat}$ from $U_{\rm id}$ we mean that $u\in U_{\rm id}$ and in the meanwhile $\frac{u_{k+1}-u_{k}}{u_{k}-u_{k-1}}\rightarrow +\infty$ for all $k=2,...,n-1$, and $u_2-u_1\rightarrow 0+$.

We generalize the previous three theorems and Definition \ref{caterpillar} from a caterpillar point to a general boundary point (see Section \ref{Limitiso} for more precise statements with full details).
\begin{thm}\label{mainthm2}
Let us denote by $U_{\sigma}=\{u\in \frak t_{\rm reg}(\mathbb{R})~|~u_{\sigma(1)}<\cdots <u_{\sigma(n)}\}$ the connected component of $\frak t_{\rm reg}(\mathbb{R})$ associated to an element $\sigma\in S_n$ in the permutation group. For any $A\in\Herm(n)$,  
\begin{itemize} 
\item[(a).] as $u\in {\frak t_{\rm reg}}(\mathbb{R})$ approaches to a boundary point $u_0\in \widetilde{\frak t_{\rm reg}}(\mathbb{R})\setminus {\frak t_{\rm reg}}(\mathbb{R})$ from $U_\sigma$, the Stokes matrices $S_\pm(u,A)$ of the linear system \eqref{introeq} have regularized limits, denoted by $S_\pm^\sigma(u_0,A)$ (that encodes the first order approximation of  $S_\pm(u,A)$ as $u\rightarrow u_0$ from $U_\sigma$);

\item[(b).] the regularized limits can be expressed by
\begin{eqnarray}\label{limitviaiso}
S_\pm^\sigma(u_0,A)=S(u,\Phi_{u_0}(u;A)), \ \ \text{for all } \ u\in U_\sigma,
\end{eqnarray}
where $\Phi_{u_0}(u;A)$ is the solution of isomonodromy equation \eqref{introisoeq} on $U_\sigma$ with the prescribed asymptotics parameterized by $A$ as $u\rightarrow u_0$ (i.e., the prescribed boundary value $A$ at $u_0$);
\item[(c).] more importantly, associated to the boundary point $u_0$ there exists a collection of linear systems of differential equations with lower ranks and simpler forms, such that the regularized limits $S^\sigma_\pm(u_0,A)$ are explicitly expressed in terms of the Stokes/connection matrices of the collection of linear systems.
\end{itemize}
\end{thm}
When $u_0=u_{\rm cat}$ the caterpillar point and $U_\sigma=U_{\rm id}$, the $(a), (b)$ parts of the theorem recover Theorem \ref{introcor} and Theorem \ref{isomonopro}. And in this case, the collection of simpler linear systems given in part $(c)$ are exact solvable, see Section \ref{explicitev}. The explicit evaluation of the regularized limits gives rise to the formula in Theorem \ref{mainthm}, see Section \ref{eleabove}. Thus in this special case, part $(c)$ recovers Theorem \ref{mainthm}.

Like Definition \ref{caterpillar}, we call the regularized limits $S^\sigma_\pm(u_0,A)$ as the Stokes matrices at a boundary point $u_0$ with respect to the choice of $U_{\sigma}$. There are different choices of $\sigma\in S_n$ such that $u_0$ is in the closure of the connected component $U_\sigma$ in $\widetilde{\frak t_{\rm reg}}(\mathbb{R})$, therefore we use the upper index $\sigma$ in $S^\sigma_\pm$ to stress that the limit is taken as $u$ approaches to $u_0$ from inside of $U_{\sigma}\subset \widetilde{\frak t_{\rm reg}}(\mathbb{R})$. 

Actually, as $u$ approaches to a fixed boundary point $u_0$ from two different components $U_\sigma$ and $U_{\sigma'}$, the regularized limits are in general different. The wall-crossing formula of the Stokes matrices, as $u$ crosses the common face of $\overline{U_\sigma}$ and $\overline{U_{\sigma'}}$ in $\widetilde{\frak t_{\rm reg}}(\mathbb{R})$, is derived, see Section \ref{connembedd} for more details.

In particular, for each $1\le i\le n$, let $\tau_i\in S_n$ be the permutation reversing the order of the subset $[1,...,i]$ of $[1,...,n]$. Then the caterpillar point $u_{\rm cat}$ is in the intersection of all the $\overline{U_{\tau_i}}$ in $\widetilde{\frak t_{\rm reg}}(\mathbb{R})$, and 
\begin{pro}[Wall-crossing formula at $u_{\rm cat}$]\label{introwall}
Denote by $S_+(u_{\rm cat},A)$ and $S_+^{\tau_i}(u_{\rm cat},A)$ the regularized limits of the Stokes matrices $S_+(u,A)$ as $u\rightarrow u_{\rm cat}$ from $U_{\rm id}$ and $U_{\tau_i}$ respectively. Then we have \begin{eqnarray}
S_+(u_{\rm cat},A)=\left(\begin{array}{cc}
S_{i+} & B \\
0 & C
\end{array} \right)\rightarrow S_+^{\tau_i}(u_{\rm cat},A)=\left(\begin{array}{cc}
P_{i}S_{i+}^\dagger P_i^{-1} & P_i(S_{i+} S_{i+}^\dagger)^{-\frac{1}{2}}S_{i+} B \\
0 & C
\end{array} \right).
\end{eqnarray}
where $S_{i+}$ are the left-top $i$-th principal submatrices of $S_+(u_{\rm cat},A)$, and $S_{i+}^\dagger$ is the complex conjugate of $S_{i+}$, and $P_i$ is the $i\times i$ permutation matrix associated to $\tau_i$.
\end{pro}

The wall-crossing formula is used to realize the known cactus group action on $\frak{gl}_n$-crystals constructed in representation theory, see Theorem \ref{cactusact}. Here the cactus group is the ($S_n$-equivariant) fundamental group of $\widetilde{\frak t_{\rm reg}}(\mathbb{R})$, see \cite{DJS}.
\begin{rmk}
One should compare the new wall-crossing phenomenon in $\widetilde{\frak t_{\rm reg}}(\mathbb{R})$ to the known wall-crossing phenomenon in $\frak t_{\rm reg}(\mathbb{C})$: Stokes matrices $S(u,A)$ parameterized by $u\in\h_{\rm reg}(\mathbb{C})$ and $A\in\frak{gl}_n$ depend on the discrete choices of initial Stokes sectors, and as $u$ varies in $\frak t_{\rm reg}(\mathbb{C})$, the configuration of the Stokes sectors and anti-Stokes rays vary accordingly. As $u$ crosses some particular real codimensional one wall in $\frak t_{\rm reg}(\mathbb{C})$, the chosen initial Stokes sector first collapses into a line and then gives birth to a new Stokes sector. This wall-crossing phenomenon amounts to the braid group actions on the space of Stokes matrices. See \cite[Appendix F]{Dubrovin} and \cite{BoalchG}.
\end{rmk}

\begin{rmk}
Theorem \ref{mainthm2} can be seen as an analytic branching rule of the linear system \eqref{introeq}, as they "decouple" the linear system \eqref{introeq} into the multiple linear systems of lower ranks, according to the branches of the planar tree representing the way in which the components $u_i$ of $u$ collapse (in a comparable speed). The analytic branching rule here is related to the branching rule in representation theory, i.e., the rules decomposing the restriction of an irreducible representation into irreducible representations of the subgroup. See Sections \ref{sec:qStokesGZ} and \ref{GZcrystal} for a special case, where, in the WKB approximation, the analytic branching rule at the caterpillar point $u_{\rm cat}$ is shown to recover the Gelfand-Tsetlin basis, arsing from the branching rules from $U(k)$ to $U(k-1)$ for all $2\le k\le n$. See Conjecture \ref{conj2} for other boundary point cases and a relation with the eigenbasis of the shift of argument subalgebras. In particular, the analytic branching rule as the two sets $u_1,...,u_k$ and $u_{k+1},...,u_n$ collapse is related to the branching rule for the reduction from ${\gl}_n$ to ${\gl}_k\times {\gl}_{n-k}$.
\end{rmk}

\begin{rmk}
The De Concini-Procesi space was introduced to describe the asymptotic behavior of solutions of the Khniznik-Zamolodchikov (KZ) equation, see \cite{dCP}. Since the isomonodromy equation \eqref{introisoeq} is a non-linear differential equation closely related to the KZ equation, what we get can be seen a nonlinear analog. However, the nonlinear analysis contains new phenomenon and is much harder, see \cite{TangXu}.
\end{rmk}

\begin{rmk}
By the expressions in Theorem \ref{reglimitbound} of Section \ref{closureStokes}, if the Stokes matrices $S_\pm(u,A)$ take very special forms, the limits of the Stokes matrices $S_\pm(u,A)$ exist as some components $u_i$ of $u$ collapse (independent of the relative speed) and coincide with the corresponding one on the boundary $\frak t(\mathbb{R})\setminus \frak t_{\rm reg}(\mathbb{R})$. It should be related to the work \cite{CDG} of Cotti, Dubrovin and Guzzetti, where continuous deformation of $S_\pm(u,A)$ from $\frak t_{\rm reg}$ to $\frak t$ was studied.
\end{rmk}

The Stokes matrices of the system \eqref{introisoStokeseq1}-\eqref{introisoStokeseq2} or the confluent hypergeometric system \eqref{introeq} have appeared in many fields of mathematics and physics, but the difficult in practical application is that in general they don't have explicit expression. However, since the problems, to which Stokes phenomenon applies, usually depend on the parameter $u$ in an isomonodromy/isospectral way, and the structures in question are preserved under taking the regularized limit of Stokes matrices in $\widetilde{\frak t_{\rm reg}}(\mathbb{R})$, thus one can use the closed formula at the caterpillar point $u_{\rm cat}$ and then to study the problems via the isomonodromy approach. See the applications in Poisson geometry and representation theory given in the following. Therefore, Theorem \ref{mainthm} provides a framework that brings manipulable analysis tools to many seemingly disparate fields, and more importantly, helps to solve some open analysis problems for the confluent hypergeometric system \eqref{introeq} and the associated isomonodromy deformation equations. In particular, it enables us to give/understand 
\begin{itemize}
    \item the first explicit linearization of dual Poisson Lie groups for $U(n)$, see Theorem \ref{eBoalchthm}; 
    \item a realization of crystals in the representation of $\frak{gl}_n$ via the WKB approximation of the quantum Stokes matrices, see Theorem \ref{WKBthm}, and more importantly an algebraic characterization of the WKB approximation by explicit $\frak{gl}_n$-crystals, see Conjecture \ref{conj2};
    \item a realization of the cactus group actions on the $\frak{gl}_n$-crystals as a wall-crossing phenomenon of the regularized limit of Stokes matrices in the De Concini-Procesi space $\widetilde{\frak t_{\rm reg}}(\mathbb{R})$, see Theorem \ref{cactusact};
    \item the geometry of spectral networks in the WKB approximation of equation \eqref{introeq} and its relation with cluster algebras and Cauchy interlacing inequality \cite{ANXZ};
    \item the boundary conditions and the asymptotic expansion of generic $\frak{gl}_n(\mathbb{C})$ valued solutions of the nonlinear isomonodromy deformation equation, which is a generalization of Jimbo's formula \cite{Jimbo} from the Painlev\'{e} VI to higher rank case \cite{TangXu}. As applications, we construct some algebraic solutions and obtain a connection formula of the nonlinear isomonodromy deformation equation \cite{Xu3, TangXu}.
\end{itemize}

The first three applications are contributions of this paper, and others are given in follow-up works. In the rest of the introduction, let us give more details on them, with a stress on the fact that the involved structures (either Poisson geometric nature of Riemann-Hilbert-Birkhoff maps, the quantum groups, or the crystals) we are interested in are preserved under taking the regularized limit of Stokes matrices in $\widetilde{\frak t_{\rm reg}}(\mathbb{R})$. 

\subsection{Application I: the first explicit Ginzburg-Weinstein diffeomorphism}\label{introPoisson}
In this subsection, we concern the Poisson geometric nature of the Riemann-Hilbert-Birkhoff map of the linear system \eqref{introeq} initiated by Boalch \cite{Boalch1}. In particular, we prove that the Riemann-Hilbert-Birkhoff map appeared in Section \ref{firstsec} is a Poisson diffeomorphism. As a corollary, we derived the first explicit expression of the Ginzburg-Weinstein linearization.

Following the irregular Atiyah-Bott construction \cite{Boalch1, Boalch2}, the theory of Stokes matrices can be placed into the context of Poisson/symplectic
geometry. For that, let us consider the Lie algebra ${ \frak u}(n)$ of the unitary group ${\rm U}(n)$, consisting of skew-Hermitian matrices, and identify $\Herm(n)\cong { \frak u}(n)^*$ via the pairing $\langle A,\xi\rangle=2{\rm Im}({\rm tr}A\xi)$. Thus $\Herm(n)$ inherits a Poisson structure from the canonical linear (Kostant-Kirillov-Souriau) Poisson structure on ${\frak u}(n)^*$. Furthermore, the
unitary group ${\rm U}(n)$ carries a standard structure as a Poisson Lie
group (see e.g. \cite{LW}). The dual Poisson Lie group ${\rm U}(n)^*$, which is the group of complex
upper triangular matrices with strictly positive diagonal entries,
is identified with the space $\Herm^+(n)$ of positive definite Hermitian $n\times n$-matrices, by
taking the upper triangular matrix $X\in {\rm U}(n)^*$ to the positive Hermitian matrix $(X^*X)^{1/2}\in
\Herm^+(n)$. The Ginzburg-Weinstein linearization theorem \cite{GW} states that the dual Poisson Lie group ${\rm U}(n)^*\cong \Herm^+(n)$ is Poisson isomorphic to the dual of the Lie algebra $\frak u(n)^*\cong\Herm(n)$. 
We remark that the linearization theorem works for any compact Lie group $K$ with its standard Poisson structure, and there are various proofs and generalizations of Ginzburg-Weinstein diffeomorphism, from the different perspectives of Moser's trick in symplectic geometry, Stokes phenomenon, the theory of quantum algebras and so on, see e.g., \cite{Anton,AM, Boalch2,BoalchG,EEM}.

Although there are many proofs of the existence of Ginzburg-Weinstein diffeomorphisms, the explicit expression of such a diffeomorphism was not known before (except $n=2$ case). In Section \ref{rRHvsGZ}, we prove
\begin{thm}\label{eBoalchthm}
The map \[\nu(u_{\rm cat})\colon\Herm(n)\cong\frak u(n)^*\rightarrow \Herm^+(n)\cong{\rm U}(n)^*; \ A\mapsto S_{-}(u_{\rm cat},A) S_+(u_{\rm cat},A),\] is a Poisson isomorphism (here since the Stokes matrices satisfy $S_-=S_+^\dagger$, the conjugate transpose of $S_+$, the product $S_-S_+$ is a positive definite Hermitian matrix).
\end{thm}
As an immediate corollary of Theorem \ref{eBoalchthm} and Theorem \ref{mainthm}, we find the first explicit Ginzburg-Weinstein diffeomorphism.

Theorem \ref{eBoalchthm} can be understood as an extension of the following theorem of Boalch to a caterpillar point $u_{\rm cat}$.
\begin{thm}\label{Boalchthm}\cite{Boalch1}
For any fixed $u\in\h_{\rm reg}(\mathbb{R})$, the Riemann-Hilbert-Birkhoff map (also known as the dual exponential map)
\begin{eqnarray}\label{introRH}
\nu(u)\colon\Herm(n)\cong\frak u(n)^*\rightarrow \Herm^+(n)\cong{\rm U}(n)^*; \ A\mapsto S_-(u,A) S_+(u,A),
\end{eqnarray}
is a Poisson isomorphism. 
\end{thm}
In Section \ref{endsection}, we give a "linear algebra" proof of Theorem \eqref{Boalchthm}.
The above maps were first studied in \cite{Boalch1}, whose surjectivity was new (i.e had not been used in the context of Riemann-Hilbert-Birkhoff before). Using Theorem \ref{eBoalchthm} and the Hamiltonian formulation of the isomonodromy deformation equations, in Section \ref{endsection} we give a new proof of the Poisson map part of Theorem \ref{Boalchthm}.
\begin{rmk}
In \cite{AM}, Alekseev and Meinrenken constructed a distinguished Ginzburg-Weinstein linearization via the Gelfand-Tsetlin integrable systems. As an application of Theorem \ref{isomonopro} and \ref{mainthm}, in a follow up work \cite{Xu6}, we clarify the relation between $\nu(u)$ and the Alekseev-Meinrenken diffeomorphism. Furthermore, we find the explicit expression of the Alekseev-Meinrenken diffeomorphism. See \cite{Xu6} for more details.
\end{rmk}

\subsection{The quantum case: representations of quantum groups arising from the quantum Stokes matrices}

From now on, we switch to the quantum analog of the differential equation \eqref{introeq}, with a focus on its Stokes phenomenon and WKB approximation. And in the following, we unveil the quantum analogs of all the above results. We first recall that how the representation of quantum groups naturally arises from the Stokes matrices of the quantum confluent hypergeometric equation \eqref{introqeq}, and then deepen the relation between the Stokes phenomenon and representation theory, by proposing a conjecture (and prove an important special case) that the Stokes phenomenon in the WKB approximation of quantum confluent hypergeometric equation is characterized by explicit crystal structures. Here on the one hand, the crystal structures in representation theory were introduced by Kashiwara \cite{Kashiwara1,Kashiwara2} and Lusztig \cite{Lu} in the 1990’s. Since then, they have become ubiquitous in algebra and geometry. On the other hand, the characterization of the WKB approximation, as a singular perturbation problem, is still open. So it is rather striking that the crystal structures enable us to study and understand open analysis problems. We mention that a classical analog, i.e., a geometric (instead of a representation theoretic) characterization of the WKB approximation of the (classical) confluent hypergeometric equation \eqref{introeq} in terms of integral periods over the associated spectral curves, is studied in \cite{ANXZ}.

Let us take the Lie algebra ${\frak {gl}_n}$ over the field of complex numbers, and its universal enveloping algebra $U({\frak {gl}}_n)$ generated by $\{e_{ij}\}_{1\le i,j\le n}$ subject to the relation $[e_{ij},e_{kl}]=\delta_{jk}e_{il}-\delta_{li} e_{kj}$. Let us take the $n\times n$ matrix $T=(T_{ij})$ with entries valued in $U({\frak {gl}}_n)$
\begin{eqnarray*}
T_{ij}=e_{ij}, \ \ \ \ \ \text{for} \ 1\le i,j\le n.
\end{eqnarray*}
Given any finite-dimensional irreducible representation $L(\lambda)$ of ${\gl}_n$ with a highest weight $\lambda$, let us consider the quantum confluent hypergeometric system
\begin{eqnarray}\label{introqeq}
\frac{dF_h}{dz}=h\left(\I u+\frac{1}{2\pi \I }\frac{T}{z}\right)\cdot F_h,
\end{eqnarray}
for $F_h(z)\in {\rm End}(L(\lambda))\otimes {\rm End}(\mathbb{C}^n)$ an $n\times n$ matrix function with entries in ${\rm End}(L(\lambda))$. Here $\I =\sqrt{-1}$, $h$ is a complex parameter, $u\in\frak t_{\rm reg}(\mathbb{C})$ is seen as a $n\times n$ matrix with scalar entries in $U({\frak {gl}}_n)$, and the action of the coefficient matrix on $F_h(z)$ is given by matrix multiplication and the representation of ${\frak {gl}}_n$. 

Let us assume $h\notin \I \mathbb{Q}$. The equation \eqref{introqeq} is then nonresonant and thus has a unique formal solution $\hat{F}_h(z)$ around
$z = \infty$ (note that the equation can be seen as a block linear system of rank $n\times {\rm dim}(L(\lambda))$, that is a special case of \eqref{introeq}). Again the standard theory of resummation states that there exist certain  sectorial regions around $z=\infty$, such that on each of these sectors there is a unique (therefore canonical) holomorphic solution with the prescribed asymptotics $\hat{F}$. These solutions are in general different (that reflects the Stokes phenomenon), and the transition between them can be measured by a pair of Stokes matrices $S_{h\pm}(u)\in{\rm End}(L(\lambda))\otimes{\rm End}(\mathbb{C}^n)$. See \cite{Xu2} for more details. 

\begin{thm}\cite{Xu2}\label{introthm1}
For any fixed $h\notin \I \mathbb{Q}$ and $u\in\frak t_{\rm reg}(\mathbb{C})$, the map (with $q=e^{h/2}$)
\begin{equation}\label{mapbasis}
\begin{split}
\mathcal{S}_q(u): U_q(\frak{gl}_n)&\rightarrow {\rm End}(L(\lambda))~;\\
e_i&\mapsto \frac{S_+(u)_{i,i}^{-1}\cdot S_+(u)_{i,i+1}}{q^{-1}-q}, \\
f_i&\mapsto \frac{S_-(u)_{i+1,i}\cdot S_-(u)_{i,i}^{-1}}{q-q^{-1}}, \\
q^{h_i}&\mapsto S_+(u)_{i,i}
\end{split}
\end{equation}
defines a representation of the Drinfeld-Jimbo quantum group $U_q(\frak{gl}_n)$ on the vector space $L(\lambda)$. Here recall that $U_q(\frak{gl}_n)$ is a unital associative algebra with generators $q^{\pm h_i}, e_j, f_j,$ $1\le j\le n-1, 1\le i\le n$ and relations:
\begin{itemize}
    \item for each $1\le i\le n$, $1\le j\le n-1$,
\[q^{h_i}q^{-h_i}=q^{-h_i}q^{h_i}=1, \
q^{h_i}e_jq^{-h_i}=q^{\delta_{ij}}q^{-\delta_{i,j+1}}e_j, \ q^{h_i}f_jq^{-h_i}=q^{-\delta_{ij}}q^{\delta_{i,j+1}}f_j;
\]
\item for each $1\le i,j\le n-1$,
\[
[e_i,f_j] = \delta_{ij} \frac{q^{h_i-h_{i+1}}-q^{-h_i+h_{i+1}}}{q-q^{-1}};
\]
\item for $|i-j|=1$, 
\[
e_i^2e_j - (q+q^{-1})e_ie_je_i + e_je_i^2=0, 
\]
\[
f_i^2f_j - (q+q^{-1})f_if_jf_i + f_jf_i^2=0,
\]
and for $|i-j|\ne 1$, $[e_i,e_j]=0=[f_i,f_j]$.
\end{itemize}
\end{thm}

The theorem associates to any representation $L(\lambda)$ of $U(\frak{gl}_n)$ a representation $\mathcal{S}_q(u)$ of $U_q(\frak{gl}_n)$ on the same vector space $L(\lambda)$. In the following, we will call $S_{h\pm}(u)$ the quantum Stokes matrices.
In \cite{Xu2}, the (formal solution and quantum Stokes matrices of) linear system \eqref{introqeq} is interpreted as a quantization of the (ones of) linear system \eqref{introeq} in the framework of deformation quantization. In particular, Theorem \ref{introthm1} can be seen as a quantum analog of Theorem \ref{Boalchthm}.
\begin{rmk}
The study of the quantum Stokes matrices is generalized from the second order pole case to arbitrary order pole cases: in \cite{Xu7}, a quantum analog of meromorphic linear systems of ODEs with pole of order $k$, as well as its quantum Stokes matrices, is introduced. As for $k=2$, it becomes the equation \eqref{introqeq}. The quantum Stokes matrices at pole of order $k$ is then interpreted as a quantization of the space of the classical Stokes matrices.
\end{rmk}

\subsection{The WKB approximation in the Stokes matrices is the crystal limit in the quantum groups}
Now let us introduce the most important application of Theorem \eqref{mainthm}, i.e., a transcendental realization of $\frak{gl}_n$-crystals in representation theory. The following observation, on a correspondence between the WKB approximation in the differential equation and the crystal limit in quantum groups, is the starting point. 

On the one hand, the WKB method, named after Wentzel, Kramers, and Brillouin, is for approximating solutions of a differential equation whose highest derivative is multiplied by a small parameter (other
names, including Liouville, Green, and Jeffreys are sometimes attached to this method). Accordingly, we would like to study directly the asymptotics of $S_{h\pm}(u)$ as $h\rightarrow +\infty$ along the positive real axis, which describe the Stokes phenomenon in the WKB approximation of \eqref{introqeq}. For that, we fix a lowest vector $\xi_0$, and an inner product on $L(\lambda)$ given by the conditions
$\langle\xi_0,\xi_0\rangle = 1$ and $\langle e_{ij}v_1, v_2\rangle=\langle v_1,  e_{ji}v_2\rangle$
for any $v_1,v_2\in L(\lambda)$. 

On the other hand, a ${\gl}_n$-crystal (see Definition \ref{WKBdata}) is a combinatorial object, that is a finite set along with some operators called crystal operators satisfying certain conditions, where the finite set models a weight basis for a representation of ${\gl}_n$, and crystal operators indicate the leading order behaviour of the simple root vectors on the basis under the crystal limit $q\rightarrow \infty$ in quantum group $U_q({\gl}_n)$.

Note that in the realization of quantum $U_q(\frak{gl}_n)$ via the quantum Stokes matrices $S_{h\pm}(u)$ given in Theorem \ref{introthm1}, the parameters $q$ and $h$ are related by $q=e^{h/2}$. Therefore, the WKB leading asymptotics as $h\rightarrow+\infty$ (along the positive real axis) corresponds to the crystal limit $q\rightarrow \infty$,
\[\fbox{WKB approximation of the equation \eqref{introqeq}} \longleftrightarrow \fbox{crystal limit of the quantum group $U_q({\gl}_n)$}\] 
In the following, let us make the above correspondence precise. 

Following the standard process of WKB analysis, the leading asymptotics of solutions encodes the eigenvalues and eigenbasis of the coefficient $( \I u+\frac{1}{2\pi\I }\frac{T}{z})$ of equation \eqref{introqeq}. In our case, let us take the shift of argument subalgebra $\mathcal{A}(u)$ of $U({\gl}_n)$, which is a maximal commutative subalgebras parameterized by $u\in \frak t_{\rm reg}(\mathbb{C})$. The action of $\mathcal{A}(u)$ on the representation $L(\lambda)$ has simple spectral for the real $u\in\frak t_{\rm reg}(\mathbb{R})$. We denote by $E(u;\lambda)$ an eigenbasis of the action of $\mathcal{A}(u)\subset U(\frak{gl}_n)$ on $L(\lambda)$. See e.g., \cite{FFR,HKRW} for more details. In the discussion below, let us assume $h\in\mathbb{R}_{>0}$, and $u\in\frak t_{\rm reg}(\mathbb{R})$. The assumption guarantees the existence of the $E(u;\lambda)$. 

The action of the off-diagonal entry $S_{h+}(u)_{k,k+1}$ of the quantum Stokes matrix on the eigenbasis vectors $\{v_i(u)\}_{i\in I}$ of $E(u;\lambda)$ should have the WKB type asymptotic behaviour as $h\rightarrow +\infty$,
\begin{equation}\label{WKBasy}
S_{h+}(u)_{k,k+1} \cdot v_i(u)=\sum_{j\in I}e^{h\phi^{(k)}_{ij}(u)+\I  g_{ij}^{(k)}(u,h)} \left(v_j(u)+\mathcal{O}(h^{-1})\right), 
\end{equation}
where $\phi^{(k)}_{ij}(u)$ are real valued functions independent of $h$, and $g_{ij}^{(k)}(u,h)$ are real valued functions for all $1\le i,j\le k\le n-1.$
An element $v_i(u)$ of $E(u;\lambda)$ is called generic if there exists only one index $j\in I$ such that $\phi^{(k)}_{ij}(u)$ is the biggest in the collection $\{\phi^{(k)}_{il}(u)\}_{l\in I}$ of real numbers.
Thus, the WKB approximation of $S_{h+}(u)_{k,k+1}$ naturally defines an operator $\widetilde{e_k}$ on the generic elements of $E(u;\lambda)$ by picking the unique leading term in \eqref{WKBasy}, i.e.,
\begin{eqnarray}
\widetilde{e_k}(v_i(u)):= v_j(u), \ \ \ \text{if} \ \ \ \phi^{(k)}_{ij}(u)={\rm max}\{\phi^{(k)}_{il}(u)~|~{l\in I}\}.
\end{eqnarray}
Similarly, by considering the WKB approximation of $S_{h-}(u)_{k+1,k}$, one defines an operator $\widetilde{f_k}$ on (some other) generic elements of $E(u;\lambda)$. In a universal sense, the operators $\{\widetilde{e_k}, \widetilde{f_k}\}_{k=1,...,n-1}$ uniquely extend to the whole set $E(u;\lambda)$ of eigenbasis. See Section \ref{sec:WKBdatum} for more details on the extension for the case of caterpillar point.

In this heuristic spirit, the correspondence between the WKB approximation and the crystal limit predicts that the finite set $E(u;\lambda)$ equipped with the operators $\{\widetilde{e_k}(u), \widetilde{f_k}(u)\}_{k=1,...,n-1}$ is a $\frak{gl}_n$-crystal, i.e., the WKB approximation of the Stokes matrices is characterized by a crystal structure. To be more precise, our conjecture states that (the conjecture is proved in an important special case, see Section \ref{GZcrystal})
\begin{conj}\label{conj2}
For any $u\in\frak t_{\rm reg}(\mathbb{R})$ and each $k=1,...,n-1$, there exists canonical operators $\widetilde{e_k}(u)$ and $\widetilde{f_k}(u)$ acting on the finite set $E(u;\lambda)$, and real valued functions $c_{kj}(\xi(u))$, $\theta_{kj}(h,u,\xi(u))$ with $j=1,2$ such that for any generic element $\xi(u)\in E(u;\lambda)$,
\begin{eqnarray} \label{eku}
\mathop{\rm lim}\limits_{h\rightarrow +\infty}\left(S_{h+}(u)_{k,k+1}\cdot e^{c_{k1}(\xi) h+\I  \theta_{k1}(h,u,\xi)} \xi(u)\right)=\widetilde{e_k}(\xi(u)),\\ \label{fku}
\mathop{\rm lim}\limits_{h\rightarrow +\infty}\left(S_{h-}(u)_{k+1,k}\cdot e^{c_{k2}(\xi) h+\I  \theta_{k2}(h,u,\xi)} \xi(u)\right)=\widetilde{f_k}(\xi(u)).
\end{eqnarray}
Furthermore, the WKB datum $(E(u;\lambda), \widetilde{e_k}(u), \widetilde{f_k}(u))$ is a $\frak{gl}_n-$crystal.
\end{conj}
\begin{rmk}
The functions $c_{k1}(\xi)$ and $c_{k2}(\xi)$ are also determined by the representation theoretic data in the ${\gl}_n$-crystal. But we do not need it in this paper. 
\end{rmk}
\begin{rmk}
In the literature of KZ equations associated to a simple Lie algebra $\g$, see e.g., \cite{FFR0}, the parameter $h$ equals to $k+h^\vee$, where $k$ is the level of the representation of the affine Lie algebra $\hat{\g}$ and $h^\vee$ the dual Coxeter number of $\g$. Then the Stokes phenomenon in the WKB approximation of the equation \eqref{introqeq} is expected to be related to the theory of representation of $\hat{\g}$ at the critical level $k=-h^\vee$. 
\end{rmk}

In one of our next papers, an isomonodromy deformation approach to this conjecture is proposed. It decomposes the proof of the conjecture into a problem of quantitative analysis and a problem of qualitative analysis. The quantitative analysis problem is then solved in Theorem \ref{WKBthm}, and the qualitative one is a pure analysis problem. 

\subsection{Leading asymptotics of quantum Stokes matrices in terms of Gelfand-Tsetlin basis}\label{sec:qStokesGZ}
This subsection gives the quantum analog of Proposition \ref{StokesinGZ}.  
\begin{thm}\label{explicitS0}
The leading asymptotics of the off-diagonal entries $S_{h+}(u)_{k,k+1}\in {\rm End}(L(\lambda))$ of $S_{h\pm}(u)$, as $\frac{u_{k+1}-u_{k}}{u_{k}-u_{k-1}}\rightarrow +\infty$ for all $k=2,...,n-1$, are given by
\begin{align*}
&S_{h+}(u)_{k,k+1}\sim 2\pi\I  h^{\frac{h\small{(e_{kk}-e_{k+1,k+1}-1})}{2\pi\I }}e^{\frac{-h(e_{kk}+e_{k+1,k+1}+1)}{4}}  \\
&\times \sum_{i=1}^k\frac{\prod_{l=1,l\ne i}^{k}\Gamma(h\frac{\zeta^{(k)}_i-\zeta^{(k)}_l}{2\pi \I })}{\prod_{l=1}^{k+1}\Gamma(1+h\frac{\zeta^{(k)}_i-\zeta^{(k+1)}_l-1}{2\pi \I })}\frac{\prod_{l=1,l\ne i}^{k}\Gamma(1+h\frac{\zeta^{(k)}_i-\zeta^{(k)}_l-1}{2\pi \I })}{\prod_{l=1}^{k-1}\Gamma(1+h\frac{\zeta^{(k)}_i-\zeta^{(k-1)}_l}{2\pi \I })}\cdot f^{(k)}_{hi}(\Delta_L)^{1,...,k}_{1,...,k-1,k+1}\left(\frac{hT(\zeta^{(k)}_i)}{2\pi\I }\right).
\end{align*}
Here the singular part
\begin{equation}\label{fg}
f^{(k)}_{hi}=\left\{
          \begin{array}{lr}
            (u_k-u_{k-1})^{\frac{he_{kk}-h\zeta^{(k)}_i}{2\pi\I }}\left({u_{k+1}-u_k}\right)^{\frac{h\zeta^{(k)}_i-he_{k+1,k+1}}{2\pi\I }},   & \text{if} \ \ k\ge 2, \\
           (u_{2}-u_{1})^{\frac{h\zeta^{(1)}_1-he_{22}}{2\pi\I }}, & \text{if} \ \ k=1,
             \end{array}
\right.
\end{equation}
and $(\Delta_L)^{1,...,k}_{1,...,k-1,k+1}\left(\frac{hT(\zeta^{(k)}_i)}{2\pi\I }\right)\in {\rm End}(L(\lambda))$ and $\zeta^{(k)}_i\in {\rm End}(L(\lambda))$ are the quantum analog of the minors and eigenvalues defined as in Definition \ref{lrminor} and \ref{qroots}.
\end{thm} 
In practice, the above formula can be computed under the orthonormal Gelfand-Tsetlin basis. Denote by ${\gl}_{k}$ the subalgebra of  ${\gl}_{n}$ spanned by the elements  $\{e_{ij}\}_{i,j=1,...,k}$, and denote by the $n$-tuples of numbers $(\lambda^{(n)}_1,...,\lambda^{(n)}_n)$ parameterizing the highest weight $\lambda$.
Then the orthonormal Gelfand-Tsetlin basis $E_{GT}(\lambda)=\{\xi_\Lambda(u_{\rm cat})\}$ in $L(\lambda)$, associated to the chain of subalgebras
\[ {\gl}_1\subset\cdots\subset {\gl}_{n-1}\subset {\gl}_n\] is parameterized by the Gelfand-Tsetlin patterns $\Lambda$. Such a pattern $\Lambda$ is a collection of numbers $\{\lambda^{(i)}_j(\Lambda)\}_{1\le j\le i\le n}$ with the fixed $\{\lambda^{(n)}_k\}_{k=1,...,n}$ satisfying the interlacing conditions
\begin{eqnarray}\label{interrel}
\lambda_j^{(i)}(\Lambda)-\lambda^{(i-1)}_j(\Lambda)\in\mathbb{Z}_{\ge 0}, \hspace{5mm} \lambda^{(i-1)}_j(\Lambda)-\lambda^{(i)}_{j+1}(\Lambda)\in\mathbb{Z}_{\ge 0}.
\end{eqnarray}
The action of the quantum minors and the elements $\{\zeta^{(k)}_i\}_{1\le i\le k\le n}$ on the basis ${\xi_\Lambda}(u_{\rm cat})$ of $L(\lambda)$ are given in Proposition \ref{GZaction}. And we refer the reader to \cite{Molev} for a general theory of Gelfand-Tsetlin basis. 

Just as in the classical case, there exist regularized limits $S_{h\pm}(u_{\rm cat})$ of $S_{h\pm}(u)$ as $\frac{u_{k+1}-u_{k}}{u_{k}-u_{k-1}}\rightarrow +\infty$, called the quantum Stokes matrices at $u_{\rm cat}$ (with respect to the connected component $U_{\rm id}$). In Section \ref{sec:qStokescat} we give the explicit expression of $S_{h\pm}(u_{\rm cat})$ at the caterpillar point $u_{\rm cat}$ (see Theorem \ref{explicitS}), as a quantum version of the formula in Theorem \ref{mainthm}. (Indeed, the formula in Theorem \ref{mainthm} is given in terms of the action-angle variables of the classical Gelfand-Tsetlin integrable systems, while the formula in Theorem \ref{mainthm} is given in terms of the quantum variables).

More generally, the quantum Stokes matrices $S_{h\pm}(u)$ have canonically regularized limits as some components $u_i$ of $u={\rm diag}(u_1,...,u_n)$ collapse in a comparable speed, that is controlled by the geometry of the De Concini-Procesi space $\widetilde{\frak t_{\rm reg}}(\mathbb{R})$. Furthermore, Theorem \ref{introthm1} can be generalized from $u\in\frak t_{\rm reg}(\mathbb{R})$ to the boundary points of the De Concini-Procesi space. For example, we refer to Theorem \ref{quantumgpucat}, that gives the quantum analog of Theorem \eqref{eBoalchthm}.

\subsection{Application II: a realization of $\frak{gl}_n$-crystals via the WKB approximation in Stokes phenomenon}\label{GZcrystal}
Theorem \ref{explicitS0} enables us to compute explicitly the WKB approximation of quantum Stokes matrices in the limit $u\rightarrow u_{\rm cat}$. In particular, a straightfoward computation verifies the expression \eqref{WKBasy} at the infinite point, and the result leads to a realization of the $\frak{gl}_n$-crystal via the Stokes phenomenon. 
\begin{thm}\label{WKBthm}
For each $k=1,...,n-1$, there exists canonical operators $\widetilde{e_k}(u_{\rm cat})$ and $\widetilde{f_k}(u_{\rm cat})$ acting on the finite set $E_{GT}(\lambda)$, and real valued functions $c_{ki}(\xi(u_{\rm cat}))$, $\theta_{ki}(h,u,\xi(u_{\rm cat}))$ with $i=1,2$ such that such that for any generic element $\xi(u_{\rm cat})\in E_{GT}(\lambda)$,
\begin{eqnarray*}
\mathop{\rm lim}\limits_{h\rightarrow +\infty}\left(\mathop{\rm lim}\limits_{u\rightarrow u_{\rm cat} \text{ from } U_{\rm id}} S_{h+}(u)_{k,k+1}\cdot e^{c_{k1}(\xi) h+\I  \theta_{k1}(h,u,\xi)} \xi(u_{\rm cat})\right)=\widetilde{e_k}(\xi(u_{\rm cat})),\\
\mathop{\rm lim}\limits_{h\rightarrow +\infty}\left(\mathop{\rm lim}\limits_{u\rightarrow u_{\rm cat} \text{ from } U_{\rm id}}  S_{h-}(u)_{k+1,k}\cdot e^{c_{k2}(\xi) h+\I  \theta_{k2}(h,u,\xi)} \xi(u_{\rm cat})\right)=\widetilde{f_k}(\xi(u_{\rm cat})).
\end{eqnarray*}
Furthermore, the set $E_{GT}(\lambda)$ equipped with the operators $\widetilde{e_k}(u_{\rm cat})$ and $\widetilde{f_k}(u_{\rm cat})$ is a ${\gl}_n$-crystal, that (under the natural bijection between semistandard Young tableaux and Geland-Testlin patterns) coincides with the known $\frak{gl}_n$-crystal structure on semistandard Young tableaux.
\end{thm}

The shift of argument subalgebras $\mathcal{A}(u)$ of $U({\gl}_n)$ extend from $u\in {\h_{\rm reg}}(\mathbb{R})$ to the de Concini-Procesi space $u\in\widetilde{\frak t_{\rm reg}}(\mathbb{R})$. In particular, the subalgebra $\mathcal{A}(u)$ at $u_{\rm cat}$ becomes the Gelfand-Tsetlin subalgebra, and the eigenbasis $E(u_{\rm cat};\lambda)$ at $u_{\rm cat}$ becomes the Gelfand-Tsetlin basis $E_{GT}(\lambda)$ (this is why we denote the basis vector in $E_{GT}(\lambda)$, corresponding to a pattern $\Lambda$, by $\xi_\Lambda(u_{\rm cat})$). See e.g., \cite{HKRW} for more details. Therefore, Theorem \ref{WKBthm} proves a limit case of Conjecture \ref{conj2}.

The $\frak{gl}_n$-crystals are unique, see e.g., \cite[6.4.21]{Jo}, and there are a number of ways to construct them: combinatorially using Littelmann’s path
model \cite{Li}, representation theoretically using crystal bases of a quantum group representation \cite{Kashiwara2},
and geometrically using the affine Grassmannian \cite{BG}. See also \cite{HKbook} for the $\frak{gl}_n$-crystal structure on semistandard Young tableaux. As far as we know, Theorem \ref{WKBthm} gives the first transcendental construction of them.

\subsection{Application III: cactus group actions on $\frak{gl}_n$-crystals arising from the wall-crossing phenomenon}

Recall that the regularized limits $S_{h\pm}(u_{\rm cat})$ simply encode the leading terms of $S_{h\pm}(u)$ as $u\rightarrow u_{\rm cat}$ from the connected component $U_{\rm id}$ (see Section \ref{qGZsystem} for more details). However, as $u\rightarrow u_{\rm cat}$ from a different connected component $U_{\sigma}=\{u\in \frak t_{\rm reg}(\mathbb{R})~|~u_{\sigma(1)}<\cdots <u_{\sigma(n)}\}$, the regularized limits are given by different $S_{h\pm}^\sigma(u_{\rm cat})\in{\rm End}(L(\lambda))\otimes {\rm End}(\mathbb{C}^n)$. And the different regularized limits $S_{h\pm}(u_{\rm cat})$ and $S_{h\pm}^{\sigma}(u_{\rm cat})$ are related by an explicit wall-crossing formula (with respect to the codimension one common face of (the closure of) $U_{\rm id}$ and $U_{\sigma}$ in $\widetilde{\frak t_{\rm reg}}(\mathbb{R})$).

In particular, for each $1\le i\le n$, let $\tau_i\in S_n$ be the permutation reversing the segment $[1,...,i]$. Then the caterpillar point $u_{\rm cat}$ sets in the codimension one common face of the closure $\overline{U_{\rm id}}$ and $\overline{U_{\tau_i}}$. Similar to Theorem \ref{WKBthm}, the WKB approximation of $S_{h\pm}^{\tau_i}(u_{\rm cat})$ also leads to a $\frak{gl}_n$-crystal \[\left(E_{GT}(\lambda),\{\widetilde{e_k}^{\tau_i}(u_{\rm cat})\}_{k}, \{\widetilde{f_k}^{\tau_i}(u_{\rm cat})\}_k\right),\] but with a different set of crystal operators. Then a straightforward computation shows that the wall-crossing formula (relating $S_{h\pm}(u_{\rm cat})$ and $S_{h\pm}^\sigma(u_{\rm cat})$) in the WKB approximation generate the cactus group $Cact_n$ action on the Gelfand-Tsetlin basis, first introduced by Berenstein-Kirillov \cite{BK}, see also \cite{HKRW}. Here 
the cactus group $Cact_n$ appeared in the work of Davis-Januszkiewicz-Scott \cite{DJS} as the $S_n$-equivariant fundamental group of $\widetilde{\frak t_{\rm reg}}(\mathbb{R})$. That is (see Section \ref{cactuswall})

\begin{thm}\label{cactusact0}
The operators $\{\rho_{i}\}_{i=1,...,n-1}$ on the finite set of Gelfand-Tsetlin basis $E_{GT}(\lambda)$, uniquely determined by the identities
\begin{eqnarray}
\rho_k\circ \widetilde{e_k}=  \widetilde{e_k}^{\tau_i}\circ \rho_k,  \text{ for all } k=1,...,n-1,
\end{eqnarray}
generate the known cactus group $Cact_n$ action.
\end{thm}
In the theory of crystals,
the cactus group plays a role analogous to that of the braid group in representations of the quantum group. The famous Drinfeld-Kohno theorem shows that the action of braid group can
be realized as the monodromy of the Knizhnik-Zamalodchikov equation, see \cite{TL}. Theorem \ref{cactusact0} can be seen as a Drinfeld-Kohno type theorem in the WKB/crystal approximation. 

\subsection{Witham dynamics and Halacheva-Kamnitzer-Rybnikov-Weekes covers} 
Although it is only a limit case, the (computation of) WKB datum at $u_{\rm cat}$, given in Theorem \ref{explicitS0} and Theorem \ref{WKBthm}, provides an applicable model for general $u$ via the isomonodromy deformation, as explained in the following. The discussion below is based on Conjecture \ref{conj2}.

In Section \ref{reglimit}, the system of PDEs of $S_{h\pm}(u)$ with respect to $u$, controlling the variation of the irregular data $u$, is given. It is equivalent to the isomonodromic deformation equation of the linear system \eqref{introqeq}. Following a general principal, an isomonodromic deformation equation in the WKB approximation degenerates into the combination of a fast (isospectral deformation) dynamics and a slow (Whitham) dynamics. See e.g., \cite{Takasaki} and \cite[Section 7]{Krichever}. In the classical case, the spectral data is encoded by the underlying spectral curves. And in our quantum case, the spectral data is replaced by the eigenbasis $E(u;\lambda)$. 

By \eqref{WKBasy} the $h\rightarrow+\infty$ limit of the dynamics of $S_{h\pm}(u)$ with respect to $u$ degenerates into a slow variation of the eigenbasis $v_j(u)'s$ (the Whitham dynamics that changes the spectral data) and a fast spin $\I g_{ij}^{(k)}(u,h)$ on the basis $v_j(u)'s$ (the fast isospectral spin that preserves the spectral data). For example, at the caterpillar point, the fast spin terms are given in \eqref{fg} in the explicit expression of the quantum Stokes matrices.  

Therefore, the slow/Whitham dynamics part of the variation of $u$, in the WKB approximation of the quantum Stokes matrices, should give rise to a variation of the WKB datum $(E(u;\lambda),\widetilde{e_k}(u), \widetilde{f_k}(u))$ over any connected component $U_\sigma$ of $\frak t_{\rm reg}(\mathbb{R})$. And the discrete combinatorial structure encoded in the WKB datum should be locally independent of $u\in U_\sigma$, i.e., if $u_{\rm cat}\in \overline{U_\sigma}$ the closure, then \[\left(E(u;\lambda),\widetilde{e_k}(u), \widetilde{f_k}(u)\right)\cong \left(E_{GT}(\lambda),\widetilde{e_k}(u_{\rm cat}), \widetilde{f_k}(u_{\rm cat})\right).\] The quantum Stokes matrices $S_{h\pm}(u)$ have regularized limits as $u\in U_\sigma\subset \h_{\rm reg}(\mathbb{R})$ approaches to a boundary point in the closure $\overline{U_\sigma}\subset \widetilde{\frak t_{\rm reg}}(\mathbb{R})$, and the regularization terms are fast spin on the eigenbasis $E(u;\lambda)$. For example, following Theorem \ref{explicitS0} as $u\in U_{\rm id}$ approaches to the caterpillar point $u_{\rm cat}$, the blow up terms \eqref{fg} are just fast spin on the Gelfand-Tsetlin basis. In general, the WKB datum $(E(u;\lambda),\widetilde{e_k}(u), \widetilde{f_k}(u))$ extends from any connected component $U_\sigma$ to its closure $\overline{U_\sigma}$ in the space $\widetilde{\frak t_{\rm reg}}(\mathbb{R})$.

Therefore, on each connected component $U_\sigma$ we have a cover whose fibre at $u$ is given by the finite set $E(u;\lambda)$. And along the common face of (the closure of) two different connected components $U_{\sigma}$ and $U_{\sigma'}$, the finite set $E(u;\lambda)$, treated as two different extensions from either $u\in U_{\sigma}$ or $u\in U_{\sigma'}$, glue according to the wall-crossing formula of the Stokes matrices. For example, let us consider the caterpillar point $u_{\rm cat}$ setting in the codimension one common face of $U_{\rm id}$ and $U_{\tau_i}$, where $\tau_i\in S_n$ is the permutation reversing the segment $[1,...,i]$. Then by Theorem \ref{cactusact0} we should have the following commutative diagram, that shows the two covers by the finite set $E(u;\lambda)\cong E_{GT}(\lambda)$ over $\overline{U_{\rm id}}$ and $\overline{U_{\tau_i}}$ glue together along the common face by the operator $\rho_i:E_{GT}(\lambda)\rightarrow E_{GT}(\lambda)$ 
{\scriptsize\[
\begin{CD}
S_{h\pm}(u) @> \text{as $u\rightarrow u_{\rm cat}$  from $U_{\rm id}$}>> S_{h\pm}(u_{\rm cat})\overset{\text{wall-crossing formula}}{\Longleftrightarrow}  S^\sigma_{h\pm}(u_{\rm cat}) @< \text{as $u\rightarrow u_{\rm cat}$  from $U_{\sigma}$}<< S_{h\pm}(u) \\
@V \text{WKB datum} VV @V \text{WKB datum} VV @V \text{WKB datum} VV \\
(E(u;\lambda),\widetilde{e_k}(u))@> \text{as $u\rightarrow u_{\rm cat}$ from $U_{\rm id}$}>> (E(u_{\rm cat};\lambda)),\widetilde{e_k}(u_{\rm cat}))\overset{\text{cactus group action}}{\Longleftrightarrow} (E(u_{\rm cat};\lambda)),\tilde{e}^\sigma_k(u_{\rm cat}))  @< \text{as $u\rightarrow u_{\rm cat}$  from $U_{\sigma}$}<<(E(u;\lambda),\widetilde{e_k}(u))
\end{CD}\]}
In this way, the Whitham dynamics, that describes the variation of the spectral data $v_j(u)'s$ with respect to $u$ in the $h\rightarrow+\infty$/WKB approximation of the isomonodromy deformation, leads to a cover of the space $\widetilde{\frak t_{\rm reg}}(\mathbb{R})$ whose fibre at any $u$ is $E(u;\lambda)$ and is equipped with a $\frak{gl}_n$-crystal structure. Then by Theorem \ref{cactusact0} and Proposition \ref{knowncactus}, the monodromy of the cover is given by the known cactus group action (see Section \ref{cactuswall}) arising from representation theory. 

In \cite{HKRW}, Halacheva, Kamnitzer, Rybnikov and Weekes (HKRW) constructed crystal structures on the set of eigenlines of shift of argument subalgebras $\mathcal{A}(u)$, and defined a covering of the space $\widetilde{\frak t_{\rm reg}}(\mathbb{R})$ by the eigenlines. They proved that the monodromy representation of the cover, with respect to the base point $u_{\rm cat}$, coincides with the Berenstein-Kirillov cactus group action on $E_{GT}(\lambda)$. Therefore, this subsection interprets the HKRW cover over $\widetilde{\frak t_{\rm reg}}(\mathbb{R})$ as a local system arsing from the Whitham dynamics in the WKB approximation of the differential equation \eqref{introqeq}. 

The construction of \cite{HKRW} works for any simple Lie algebra $\g$. One can ask the similar relation between the $\g$-crystals and the WKB approximation of quantum Stokes matrices associated to $\g$. Furthermore,  using the same idea, we expect that various wall-crossing type formula in representation theory can be interpreted via Stokes phenomenon. For example, motivated by \cite[Conjecture 1.14]{HKRW}, we expect that the cactus group action on the Weyl group constructed by Losev \cite{Lo}, via perverse equivalences coming from wall-crossing functors in category $\mathcal{O}$ (a more elementary definition was given shortly
afterwards by Bonnaf$\rm\acute{e}$ \cite{Bo}), can be realized via the stokes phenomenon of the affine KZ equations introduced by Cherednik \cite{Ch}.

\subsection{A continuous path connecting the canonical basis and the Gelfand-Tsetlin basis}
Let $\mathbb{U}^+$ be the subalgebra of $U_q(\frak{gl}_n)$ generated by the elements $\{e_i\}_{i=1,...,n-1}$, and let $\mathbb{B}$
be the canonical basis in $\mathbb{U}^+$. We refer the reader to \cite{Lu} for the construction of $\mathbb{B}$. Upon acting on the lowest weight vector, the image of the canonical basis $\mathbb{B}$, under the map $\mathcal{S}_q(u)$ given in Theorem \ref{introthm1}, defines a set $B_q(u;\lambda)$ of vectors in $L(\lambda)$. 
\begin{conj}\label{conj1}
The set $B_q(u;\lambda)$ is a basis of $L(\lambda)$ for all $q\in (0,\infty)$, whose leading asymptotics as $q\rightarrow \infty$ correspond to an eigenbasis $E(u;\lambda)$ of the action of the shift of argument subalgebra $\mathcal{A}(u)\subset U(\frak{gl}_n)$ on $L(\lambda)$. 
\end{conj}
In the case $n=2$, the conjecture can be verified directly using the closed formula of the Stokes matrices. In the case $n = 3$ the canonical basis is given explicitly in \cite[Example 3.4]{Lu}. Then using the method of isomonodromy deformation, Conjecture \ref{conj1} for $\frak{gl}_3$ can be verified using the properties of the Painlev\'e VI function. 

On the one hand, as $q=1$, $U_q(\frak{gl}_n)$ becomes the undeformed $U(\frak{gl}_n)$, and the representation in \eqref{mapbasis} coincides with the given representation of $U(\frak{gl}_n)$ on $L(\lambda)$. On the other hand, by specializing $q=1$ and acting on $L(\lambda)$, the canonical basis $\mathbb{B}$ of $\mathbb{U}^+$ recovers the canonical basis $B(\lambda)$ in $L(\lambda)$, i.e., $B_{q=1}(u;\lambda)$ coincides with $B(\lambda)$. Thus, the conjecture enables us to get a one parameter family of basis $B_q(u;\lambda)$ connecting the canonical basis $B(\lambda)$ and eigenbasis $E(u;\lambda)$ in $L(\lambda)$, by varying $q=e^{\pi\I h}$ from $1$ to $\infty$ (i.e., varying $h$ from $0$ to $+\infty$ along the real axis). In particular, as for the caterpillar point $u_{\rm cat}$ we expect to get a $q$-continuous family of basis $B_q(u_{\rm cat};\lambda)$ connecting the canonical basis $B(\lambda)$ and the Gelfand-Testlin basis $E_{GT}(\lambda)$.

\subsection{Extension to $\frak{gl}_n(\mathbb{C})$ valued solutions and higher rank Painlev\'e transcendents}
The purpose of this subsection is to explain the importance of the results in Section \ref{firstsec} from the viewpoint of the theory of Painlev\'e transcendents. It outlines some results in the follow-up papers \cite{TangXu, WXX, Xu2}.

Following Miwa \cite{Miwa}, the $\frak{gl}_n$-valued solutions $\Phi(u)$ of the equation \eqref{introisoeq} with $u_1,...,u_n\in \mathbb{C}$ have the strong Painlev\'{e} property. According to the original idea of Painlev\'{e}, they can be a new class of special functions. Indeed, they arise from and play important roles in various branches of mathematics and physics. In the following, let us interpret $\Phi(u)$ as higher rank Painlev\'e transcendents, and outline how Theorem \ref{isomonopro} and Theorem \ref{mainthm} can be used to unveil many basic properties of the higher rank Painlev\'e transcendents. 

First, note that the expression of the Stokes matrices given in Theorem \ref{mainthm} is an analytic function of $\Phi_0$ for all $\Phi_0\in \frak{gl}_n(\mathbb{C})$ (not just $\Herm(n))$ satisfying the condition 
\begin{equation}
\label{boundcondtion}
|{\rm Im} (
\lambda^{(k)}_{i}(\Phi_0) - 
\lambda^{(k)}_{j}(\Phi_0)
) | < 2\pi,
\quad
\text{for every $1\le i,j\le k\le n$}.
\end{equation}
Furthermore, one checks the Riemann-Hilbert map 
\[\nu:~\Big\{\Phi_0\in \frak{gl}_n(\mathbb{C})~|~\Phi_0 \text{ satisfies \eqref{boundcondtion}}\Big\}\rightarrow \Big\{S_-(u,\Phi(u;\Phi_0))S_+(u,\Phi(u;\Phi_0))\in {\rm GL}_n(\mathbb{C})\Big\}\]
is one to one onto an open dense subset of ${\rm GL}_n(\mathbb{C})$, therefore of the space of all possible Stokes matrices. For example, it can be seen from compatibility of 
the map $\nu$ with the (complex) Gelfand-Tsetlin systems, i.e., (the complex version of) Proposition \ref{ConnGZ}.
Then by "analytical continuation", Theorem \ref{isomonopro} can be generalized to $\frak{gl}_n$ case: that is if $\Phi(u)$ is a generic $\frak{gl}_n$ valued solutions of \eqref{introisoeq}, then its asymptotics can be parameterized by a constant boundary value $\Phi_0$ satisfying the boundary condition \eqref{boundcondtion}. A precise such statement is proved in our follow up work \cite{TangXu}, by studying the complete series expansion of generic $\frak{gl}_n(\mathbb{C})$ valued solution $\Phi(u)$. And Theorem \ref{mainthm} is true simply by analytical continuation.

Then, let us give a brief introduction to the Painlev\'e transcendents. The six classical Painlev\'{e} equations were introduced at the turn of the twentieth century by Painlev\'{e} \cite{Pa} and Gambier \cite{Ga}, in a specific
classification problem for second order ODEs. Since then, they have appeared in the integrable nonlinear PDEs, 2D Ising models, random matrices, topological field theory and so on. 
We refer the reader to the book of Fokas, Its, Kapaev and Novokshenov \cite{FIKN} for a thorough introduction to the history and developments of the study of Painlev\'{e} equations. In particular, the Painlev\'{e} VI equation is the nonlinear differential equation
\begin{eqnarray*}
\frac{d^2y}{dx^2}&=&\frac{1}{2}\Big[\frac{1}{y}+\frac{1}{y-1}+\frac{1}{y-x}\Big](\frac{dy}{dx})^2-\Big[\frac{1}{x}+\frac{1}{x-1}+\frac{1}{y-x}\Big]\frac{dy}{dx}\\
&+&\frac{y(y-1)(y-x)}{x^2(x-1)^2}\Big[\alpha+\beta\frac{x}{y^2}+\gamma\frac{x-1}{(y-1)^2}+\delta\frac{x(x-1)}{(y-x)^2}\Big], \ \alpha,\beta,\gamma,\delta\in\mathbb{C}.
\end{eqnarray*}
A solution $y(x)$ has $0,1,\infty$ as critical points, and can be analytically
continued to a meromorphic function on the universal covering of $\mathbb{P}^1\setminus\{0, 1,\infty\}$. 
As stressed in \cite{FIKN,IN}, the solutions of Painlev\'{e} equations (called Painlev\'{e} transcendents) are seen as nonlinear special functions, because they play the same role in nonlinear mathematical physics as that of classical special functions, like Airy functions, Bessel functions, etc., in linear physics. And it is the answers of the following questions that make Painlev\'{e} transcendents as efficient in applications as linear special functions (here we list some of them, see \cite{FIKN, IN} and the references therein for more details):
\begin{itemize}
    \item[(a).] The parametrization of Painlev\'{e} transcendents by their asymptotic behaviour at critical points;
     
\item[(b).] The explicit expression of the monodromy of the associated linear problem via the parametrization at critical points;

\item[(c).] The construction of the connection 
formula from one critical point to another.

\end{itemize}
In particular, the problems were solved in generic case by Jimbo \cite{Jimbo}. As shown by Harnad \cite{Harnad} (see also \cite{Ma2}, \cite[Section 3]{Boalch4} for a detailed way to do the Harnad duality), that Painlev\'{e} VI is equivalent to the equation \eqref{introisoeq} for $n=3$ with suitable matrices $\Phi(u)$. In this way, the generic solutions $\Phi(u)$ of the isomonodromy equation \eqref{introisoeq} can be seen as higher rank Painlev\'e transcendents. 

Note that (the complex version of) Theorem \ref{isomonopro} and \ref{mainthm} already give answers to problem $(a)$ and $(b)$ for the transcendents $\Phi(u)$. Actually, for $n=3$, Theorem \ref{isomonopro} and \ref{mainthm} exactly recover Jimbo's formula for Painlev\'e VI, see \cite{WXX} for details.
And just like the expression of Stokes matrices of linear differential equations with small ranks $2$ or $3$ in terms of the asymptotics of the solutions of the associated nonlinear isomonodromy equations has been a major tool in the analysis of Painlev\'e transcendents, see e.g., \cite{Boalch4, Ma2} and the book \cite{FIKN}, the Theorem \ref{isomonopro} and \ref{mainthm} lay a foundation for our study of the higher rank Painlev\'e transcendents $\Phi(u)$. For example, we have used Theorem \ref{isomonopro} and \ref{mainthm} to give partial answers to the above problem $(c)$ and boundary conditions for $\Phi(u)$, see e.g, \cite{TangXu, Xu3}. 

As another example, let us outline how to use Theorem \ref{mainthm} to find algebraic solutions of the isomonodromy equation \eqref{introisoeq}. Let us first recall the case $n=3$: based on Jimbo's formula, the algorithm in \cite{DM, Boalch2, Boalch4} derives various algebraic solutions of the Painlev\'e VI. See \cite[Section 5]{Boalch2} for an important and detailed example of the algorithm. These in turn give algebraic solutions $\Phi(u)$ of \eqref{isoeq} in the case of rank $n=3$. 

Paralelly, let us consider a general $n$ case. We know from Miwa’s theorem that the matrix function $\Phi(u)$ is meromorphic on the universal covering of $\mathbb{C}^n\setminus \Delta$. Continuation along closed paths in the deformation space interchanges the branches of $\Phi(u)$. And such monodromy of the nonlinear isomonodromy equation is explicitly given in terms of the geometric terms, i.e., by an explicit braid group action on the corresponding Stokes matrices. See \cite{Dubrovin,BoalchG}.
In particular, if a solution $\Phi(u)$ of the isomonodromy equation \eqref{isoeq} with the given Stokes matrix $S_\pm(u,\Phi(u))$ is an algebraic function with branching along the diagonals $u_i=u_j$ only if $S_\pm$ belong to a finite orbit of the action of
the pure braid group. See e.g., \cite[Appendix F]{Dubrovin}. 

Therefore, we have the following systematic way to find the algebraic solutions $\Phi(u)$,
\begin{equation}
    \Big\{\text{Stokes matrices $S_\pm(u,\Phi(u;\Phi_0))$}\Big\}\Longrightarrow \Big\{\text{boundary value} \ \Phi_0\Big\}\Longrightarrow \Big\{\text{solution} \ \Phi(u;\Phi_0) \Big\}.
\end{equation}
The algorithm is as follows: starting from such a pair of explicit $S_\pm$ that belong to a finite orbit, by Theorem \ref{mainthm} we can get the leading term (boundary value $\Phi_0$) in the multivariable Puiseux expansion of the solution $\Phi(u;\Phi_0)$. Following \cite{TangXu}, the boundary value $\Phi_0$ determines explicitly the series expansion of the corresponding solution. Substituting back the leading term into the series expansion of the solutions of isomonodromy equation would determine, algebraically, any desired term in the multivariable Puiseux expansion, and then the solution itself. See our next work on the algebraic solutions $\Phi(u)$ for more details. 

In summary, given the many known applications and connections to other subjects, we believe that the higher Painlev\'e transcendents $\Phi(u)$ have richer structures and applications remained to be found. And just like Painlev\'e VI case, we expect that the answers to the above problems $(a)- (c)$ for $\Phi(u)$ will play crucial
roles in other problems from mathematical physics.

\vspace{3mm}
The organization of the paper is as follows. Section \ref{beginsection} gives the preliminaries of Stokes data of meromorphic linear systems. Section \ref{exStokesviaiso} studies the boundary values/asymptotics of the solutions of the isomonodromy equation, and then derives the expression of Stokes matrices of the associated linear system via the boundary values, i.e., Theorem \ref{mainthm}. That is the first main result of this paper. Section \ref{appPoisson} shows some applications of the analysis results in Poisson geometry, including the explicit Ginzburg-Weinstein diffeomorphism and a new proof of Theorem \ref{Boalchthm}. Section \ref{closureStokes} introduces the De Concini-Procesi space, and studies the
regularized limit of Stokes matrices $S_\pm(u,A)$ as the irregular data $u$ approaches to a boundary point of the De Concini-Procesi space. Section \ref{qStokes} introduces the quantum Stokes matrices of the quantum confluent hypergeometric differential equation, and obtains the expression of their regularized limits at caterpillar points, i.e., Theorem \ref{explicitS0}. In the end, Section \ref{tropicalisomo} gives the second main result, that is a transcendental realization of the crystals and the cactus group actions via the WKB approximation.

\subsection*{Acknowledgements}
\noindent
We would like to thank Anton Alekseev, Roman Bezrukavnikov, Philip Boalch, Pavel Etingof, Giovanni Felder, Michio Jimbo, Alexander Its, Yanpeng Li, Leonid Rybnikov, Qian Tang and Valerio Toledano Laredo for their interests and valuable feedback on the paper. The author is supported by the National Key Research and Development Program of China (No. 2021YFA1002000) and by the National Natural Science Foundation of China (No. 12171006).

\section{Stokes phenomenon and monodromy data}\label{beginsection}
In Section \ref{Canonicalsol}, we construct the canonical solutions of the meromorphic linear systems of differential equations \eqref{introeq}, and prove a  uniform property of their asymptotics. In Section \ref{defiStokes}, we introduce the Stokes matrices and connection matrices of the linear systems, as well as the monodromy relation relating connection matrices to Stokes matrices.
\subsection{Canonical solutions}\label{Canonicalsol}
Let $\h(\mathbb{R})$ (resp. $\h_{\rm reg}(\mathbb{R})$) denote the set of diagonal matrices with (resp. distinct) real eigenvalues. Let us consider the meromorphic linear system
\begin{eqnarray}\label{Stokeseq}
\frac{dF}{dz}=\left(\I u-\frac{1}{2\pi\I }\frac{A}{z}\right)\cdot F,
\end{eqnarray}
where $F(z)$ is valued in $\mathbb{C}^n$, $u\in\h(\mathbb{R})$ and $A\in{\Herm}(n)$. The system has an order two pole at $\infty$ and (if $A\neq 0$) a first order pole at $0$. 

\begin{pro}
The system has a unique formal fundamental solution taking the form \begin{eqnarray}\label{formalsol}
\hat{F}(z;u,A)=({\rm Id}_n+\hat{H}(z^{-1};u,A))\cdot e^{\I uz}z^{-\frac{[A]}{2\pi\I }},
\end{eqnarray}
where $\hat{H}(z;u,A)=\sum_{k\ge 1} H_k(u,A)z^{k}$ is a $n\times n$ matrix-valued formal power series. (Here for the convenience, we transfer the power series at $z=\infty$ to the series at $z=0$ after the change of variable $z\mapsto 1/z$.)
\end{pro}
\begin{proof}Actually, by plugging the expression \eqref{formalsol} in the equation, one checks that the coefficients $H_k$ are determined by the recursive relation
\begin{eqnarray}\label{recuH}
[\I u,H_{k+1}]=\frac{[A]}{2\pi\I }\cdot H_{k}-H_{k}\cdot \frac{A}{2\pi\I }+ kH_{k}, \ \text{for all} \ k\ge 0,
\end{eqnarray}
here we assume $H_0={\rm Id}_n$. Actually, the relation \eqref{recuH} can be rewritten, in terms of the components of $A=(a_{ij})_{i,j=1}^n$ and $H_k=(H_{k,ij})_{i,j=1}^n$, as: for $i\ne j$ 
\begin{equation}\label{Hknel}
(\I u_i-\I u_j) H_{m+1, ij}=kH_{m,ij}+\sum_{l=1}^n \frac{a_{il}}{2\pi\I } H_{m, lk}- H_{m, ij} \frac{a_{jj}}{2\pi\I },
\end{equation}
and for $i=j$ (replacing $k$ by $k+1$ in \eqref{recuH}), 
\begin{eqnarray}\label{Hkel}
0=\sum_{l=1,l\ne i}^n \frac{a_{il}}{2\pi\I } H_{k+1, li}+(k+1) H_{k+1, ii}.
\end{eqnarray}
Due to the assumption that the diagonal elements of $u$ are distinct, we see that the off-diagonal part of the matrix $H_{k+1}$ is uniquely determined by $H_k$ from \eqref{Hknel}, and then the diagonal part of $H_{k+1}$ is uniquely determined by \eqref{Hkel}.
\end{proof}

\begin{defi}\label{Stokesrays}
The {\it Stokes sectors} of the system
are the right/left half planes ${\rm Sect}_\pm=\{z\in\mathbb{C}~|~ \pm{\rm Re}(z)>0\}$.
\end{defi}

For any two real numbers $a, b$, an open sector and a closed sector with opening
angle $b-a>0$ are respectively denoted by
\[S(a,b):=\{z\in\mathbb{C}~|~a<{\rm arg}(z)<b\}, \hspace{5mm} \bar{S}(a,b):=\{z\in\mathbb{C}~|~a\le {\rm arg}(z)\le b\}. \]
For any $d>0$, let 
\[\h_{\rm reg}^d(\mathbb{R}):=\{u\in\h_{\rm reg}(\mathbb{R})~|~|u_i-u_j|>d, \ \text{for all} \ i\ne j, \ i,j=1,...,n\}\] denote the subset of $\h_{\rm reg}(\mathbb{R})$, consisting of all $u$ whose components keep a fixed positive distance $d$ from each other. 

Let us choose the branch of ${\rm log}(z)$, which is real on the positive real axis, with a cut along the nonnegative imaginary axis $\I \mathbb{R}_{\ge 0}$. Then by convention, ${\rm log}(z)$ has imaginary part $-\pi$ on the negative real axis in ${\rm Sect}_-$.
\begin{thm}\label{uniformresum}
For any $u\in\h_{\rm reg}(\mathbb{R})$, on ${\rm Sect}_\pm$ there
is a unique fundamental solution $F_\pm:{\rm Sect}_\pm\to {\rm GL}(n,\mathbb{C})$ of equation \eqref{Stokeseq} such that $F_+(z;u,A)\cdot e^{-\I uz}\cdot z^{\frac{[A]}{2\pi\I }}$ and $F_-(z;u,A)\cdot e^{-\I uz}\cdot z^{\frac{[A]}{2\pi\I }}$ can be analytically continued to $S(-\pi,\pi)$ and $S(-2\pi,0)$ respectively, and for any small $\varepsilon_0>0$,
\begin{eqnarray*}
\lim_{z\rightarrow\infty}F_+(z;u,A)\cdot e^{-\I uz}\cdot z^{\frac{[A]}{2\pi\I }}&=&{\rm Id}_n, \ \ \ \text{as} \ \ \ z\in \bar{S}(-\pi+\varepsilon_0,\pi-\varepsilon_0),
\\
\lim_{z\rightarrow\infty}F_-(z;u,A)\cdot e^{-\I uz}\cdot z^{\frac{[A]}{2\pi\I }}&=&{\rm Id}_n, \ \ \ \text{as} \ \ \ z\in \bar{S}(-2\pi+\varepsilon_0,-\varepsilon_0),
\end{eqnarray*}
Here ${\rm Id}_n$ is the rank n identity matrix, and $[A]$ is the diagonal part of $A$.
Furthermore, the above limits are uniform for $u\in \h_{\rm reg}^d(\mathbb{R})$.
\end{thm}
\pf The construction of the canonical solutions via the Laplace-Borel transforms is standard, see e.g., \cite{Balser0, BJL,LR, MR}. The new phenomenon in our case is the uniform property with respect to $u$, which relies on the fact that all components $u_i$ of $u$ lie in the same line. Since the uniform property will be used in the proof of Proposition \ref{inductionstep}, in the following, we will review the construction of the canonical solutions with a stress on the uniform property with respect to $u$. We will divide the construction into three parts, i.e., the analytic property of Borel transform of formal solutions, the analytic property of the Laplace-Borel transform, then the construction of canonical solutions.

\vspace{2mm}
{\bf Borel transform.}
Denote by $\tilde{H}=\mathcal{B}(\hat{H})$ the formal Borel transform of $\hat{H}-1$ (in the Borel plane with complex variable $\xi$), i.e., 
\[\tilde{H}(\xi;u,A):=\sum_{k\ge 1}\frac{H_k(u,A)}{\Gamma(k)}\xi^{k-1}.\]
In the rest of the proof, let us fix a positive $d>0$. Let us also take a matrix norm that is sub-multiplicative, i.e., $|BC|\le |B|\cdot |C|$ for two matrices $B$ and $C$. Then it follows from the formula \eqref{recuH} that
\begin{lem}\label{constant1}
There exists a constant $\tilde{K}>0$ such that $|H_k|/\Gamma(k)\le \tilde{K}^{k-1}$ for $u\in \h_{\rm reg}^d(\mathbb{R})$.
\end{lem}
\pf 
Let us denote by $H_k^{d}$ and $H_k^{od}$ the diagonal and the off-diagonal part of the $n\times n$ matrix $H_k$.
On the one hand, note that the entries of $[\I u,H_k]$ takes the form $(\I u_i-\I u_j)H_{k,ij}$. Since $u\in \h_{\rm reg}^d(\mathbb{R})$, we have $|\I u_i-\I u_j|>d$. Thus there exists a constant $\tilde{c}>0$ such that 
\begin{eqnarray}\label{dodineq}
|H_k^{od}|\le \tilde{c}\cdot |[u,M]|, \ \text{for all} \ u\in \h_{\rm reg}^d(\mathbb{R}).
\end{eqnarray}
On the other hand, assume $|A/2\pi \I |=a$, then \eqref{Hkel} implies that
\begin{eqnarray}\label{normineq}
|H^d_k|\le \frac{a}{k+1}|H^{od}_{k}|<|H^{od}_{k}|, \ \ \text{for k large enough}.
\end{eqnarray}

By \eqref{dodineq} and \eqref{normineq} we have 
\[|H_k|\le |H_k^d|+|H_k^{od}|\le 2|H_k^{od}|\le 2\tilde{c}\cdot |[u,H_k]|, \ \text{for $k$ large and} \ u\in \h_{\rm reg}^d(\mathbb{R}).\]
Set $t_k=H_k/\Gamma(k)$ and racall that $a=|A/2\pi\I |$, it follows from \eqref{recuH} that
\[t_k\le 2\tilde{c}(2a+1)t_{k-1}.\]
Since $a$ and $\tilde{c}$ are independent of the choices of $k$ and $u\in \h_{\rm reg}^d(\mathbb{R})$, we see that taking $\tilde{K}=2\tilde{c}(2a+1)$ verifies the Lemma. 
\qed


\vspace{2mm}
Thus the Borel series $\tilde{H}(\xi)$ is convergent in a small neighborhood $B_\rho(0):=\{\xi\in\mathbb{C}~|~|\xi|<\rho\}$ of $\xi=0$. Furthermore,
\begin{lem}\label{exgrow}
Given any number $\theta\in (-\frac{\pi}{2}, \frac{\pi}{2})\cup (-\frac{3\pi}{2}, -\frac{\pi}{2})$, there exists a sufficiently small $\varepsilon>0$, such that the Borel sum $\tilde{H}(\xi)$ can be analytically continued to $S(\theta-\varepsilon,\theta+\varepsilon)$ and there exist constants $\alpha, \beta>0$ such that
\[|\tilde{H}(\xi;u,A)|\le \alpha e^{\beta|\xi|}, \ \text{ for all } \ \xi\in  S(\theta-\varepsilon,\theta+\varepsilon) \ \ \text{and} \ u\in \h_{\rm reg}^d(\mathbb{R}).\]
\end{lem}
\pf 
Recall that $\tilde{H}(\xi;u,A)=\sum_{k\ge 1}\frac{H_k(u)}{\Gamma(k)}\xi^{k-1}$. The identities \eqref{recuH} are formally equivalent to the integral equation
\begin{equation}\label{inteq}
\tilde{H}\cdot \I u-(\I u-\xi{\rm Id}_n)\cdot \tilde{H}=\frac{[A]}{2\pi\I }-\frac{A}{2\pi\I }+\int_{t=0}^\xi\left(\tilde{H}(t)\cdot \frac{[A]}{2\pi\I }-\frac{A}{2\pi\I }\cdot \tilde{H}(t)\right)dt.
\end{equation}
To study the integral equation, we employ an iteration, by beginning $\tilde{H}^{(0)}(\xi;u,A)\equiv 0$, and plugging $\tilde{H}^{(m)}(\xi;u,A)$ into the right hand side of \eqref{inteq} and determining $\tilde{H}^{(m+1)}(\xi;u,A)$ from the left hand side. The sequences so obtained
are holomorphic near $\xi=0$, and the limit of $\tilde{H}^{(m)}(\xi;u,A)$ as $m\rightarrow \infty$ coincides with the convergent Borel series $\tilde{H}(\xi)$ for all $\xi\in B_\rho$. 

Now let $\varepsilon>0$ be small enough such that $S(\theta-\varepsilon,\theta+\varepsilon)$ do not overlap with the purely imaginary axis. 
Since the entries of $M\cdot \I u-(\I u-\xi{\rm Id}_n)\cdot M$ of any $n\times n$ matrix $M$ takes the form $(\I u_j-\I u_i+\xi)M_{ij}$, we can choose a big enough constant $c>0$ such that the off-diagonal part $M^{od}$ of $M$ satisfies
\begin{eqnarray}\label{xidod}
|M^{od}|\le c\cdot |M\cdot \I u-(\I u-\xi{\rm Id}_n)\cdot M|,
\end{eqnarray}
for all $n\times n$ matrices $M$, $\xi\in S_{>\rho'}(\theta-\varepsilon,\theta+\varepsilon)$, and  $u\in \h_{\rm reg}^d(\mathbb{R})$.
Here $c$ can be chosen, because $u_1,...,u_n$ lie on the same line, and the differences $\I u_j-\I u_i$ for all $u\in \h_{\rm reg}^d(\mathbb{R})$ keep a fixed positive distance from the sector $S(\theta-\varepsilon,\theta+\varepsilon)$.

Set $a=|A/2\pi\I |$. Let $\rho'>0$ be a fixed small enough positive number such that $\rho>\rho'>0$ and $\frac{2c}{a^2}>\rho'$. Let us introduce the region $S_{\ge \rho'}(\theta-\varepsilon,\theta+\varepsilon):=\{\xi\in S(\theta-\varepsilon,\theta+\varepsilon)~|~|\xi|\ge \rho'\}.$ 

In terms of the components $\tilde{H}^{(m)}(\xi)_{ij}$, the integral equation \eqref{inteq} can be written as

\begin{align*}
  (\I u_j-\I u_i+\xi)  \tilde{H}^{(m)}(\xi)_{ij}&=-\frac{a_{ij}}{2\pi\I }+\int_{t=0}^\xi\left(\tilde{H}^{(m-1)}(t)\cdot \frac{[A]}{2\pi\I }-\frac{A}{2\pi\I }\cdot \tilde{H}^{(m-1)}(t)\right)_{ij}dt \ \ \text{ for } i\ne j, \\
  \tilde{H}^{(m)}(\xi)_{ii}&=-\frac{1}{\xi}\int_{t=0}^\xi\sum_{l=1,l\ne i}^n\frac{a_{il}}{2\pi\I }\cdot \tilde{H}^{(m)}(t)_{lj} dt.
\end{align*}
Then we prove inductively that for all $\xi$ in the small ball $ B_{\rho'}(0)$ and for all integer $m\ge 0$,
\[|\tilde{H}^{(m)}(\xi)^{od}|\le \frac{ac}{1-2ca\rho'(a+1)}, \ \ \text{ and } \ \  |\tilde{H}^{(m)}(\xi)^{d}|\le \frac{a^2c}{1-2ca\rho'(a+1)}.\]
Set $\alpha={\rm Max}(\frac{a}{\rho'}, \frac{ac+a^2c}{1-2ca\rho'(a+1)})$. Then we have for all $\xi\in B_{\rho'}(0)$ and $m\ge 0$, 
\begin{equation}\label{allm}
    |\tilde{H}^{(m)}(\xi)|\le \alpha.
\end{equation}
We are now ready to show that each $\tilde{H}^{(m)}(\xi)$ has exponential growth of order $1$ in $S(\theta-\varepsilon,\theta+\varepsilon)$. 

Set $c'={\rm Max}(c,\frac{1}{\rho'})$. 
Let us introduce the sequence of positive real numbers $t^{(m)}_k$ with $m,k\in \mathbb{Z}_{\ge 0}$ determined by the recursive relation
\begin{align}\label{trecu}
t^{(0)}_k&=0, \ \text{ for all } k\in\mathbb{Z}_{\ge 0},\\
t^{(m)}_1&=\alpha, \ \text{ for all } m\in\mathbb{Z}_{> 0},\\
t^{(m+1)}_k&=2ac'\cdot t^{(m)}_{k-1}, \ \text{ for all } m\in\mathbb{Z}_{\ge 0}.
\end{align}
Then let us inductively show
estimates of the form
\begin{equation}\label{Hmest}
    |\tilde{H}^{(m)}(\xi;u,A)|\le \sum_{k\ge 1}{t^{(m)}_k}|\xi|^{k-1}/\Gamma(k).
\end{equation}
For $m=0$, since $\tilde{H}^{(0)}(\xi;u,A)\equiv 0$, the estimate $\eqref{Hmest}$ is true. Suppose the estimate is true for all $\tilde{H}^{(l)}$ with $l\le m$. First, the diagonal part

\[\tilde{H}^{(m+1)}(\xi)_{ii}=-\frac{1}{\xi}\int_{t=0}^\xi\sum_{l=1,l\ne i}^n\frac{a_{il}}{2\pi\I }\cdot \tilde{H}^{(m)}(t)_{lj} dt\]

Applying \eqref{xidod} to the case $M=\tilde{H}^{(m+1)}(\xi;u,A)$ and using the integral equation, we get the estimate 
\begin{itemize}
    \item for $\xi\in S_{>\rho'}(\theta-\varepsilon,\theta+\varepsilon)$,
\begin{align*}
|\tilde{H}^{(m+1)}(\xi)|&\le c'\cdot |\tilde{H}^{(m+1)}\cdot \I u-(\I u-\xi{\rm Id}_n)\cdot \tilde{H}^{(m+1)}|\\
    &\le c'\cdot \Big|\frac{[A]}{2\pi\I }-\frac{A}{2\pi\I }\Big|+c'\cdot\Big|\int_{t=0}^\xi\left(\tilde{H}^{(m)}(t)\cdot \frac{[A]}{2\pi\I }-\frac{A}{2\pi\I }\cdot \tilde{H}^{(m)}(t)\right)dt\Big|.
    \\
    &\le ac'+2ac'\int_{t=0}^\xi(
\sum_{k\ge 1}{t^{(m)}_k}|t|^{k-1}/\Gamma(k))dt
, \\
&= ac'+2ac'
\sum_{k\ge 1}{t^{(m)}_k}|\xi|^{k}/\Gamma(k+1)\\
&\le {t^{(m+1)}_{1}}+\sum_{k\ge 2}{{t^{(m+1)}_{k}}}|\xi|^{k-1}/\Gamma(k).
\end{align*}
which proves the estimate \eqref{Hmest} for $m+1$. 
Here in the last inequlity, we use the identities $t^{(m+1)}_{1}=\alpha={\rm Max}(\frac{a}{\rho'}, \frac{ac+a^2c}{1-2ca\rho'(a+1)})$ and $t^{(m+1)}_k=2ac'\cdot t^{(m)}_{k-1}$ in the defining relation \eqref{trecu}.
\item for $\xi\in S_{\le \rho'}(\theta-\varepsilon,\theta+\varepsilon):= S(\theta-\varepsilon,\theta+\varepsilon)\cap \overline{B_{\rho'}(0)}$, we have
\[|\tilde{H}^{(m+1)}(\xi)|\le \alpha=t^{(m+1)}_{1}\le {t^{(m+1)}_{1}}+\sum_{k\ge 2}{{t^{(m+1)}_{k}}}|\xi|^{k-1}/\Gamma(k).\]
\end{itemize}

Set $\alpha=a$ and $\beta=2ac$, then by \eqref{Hmest} and \eqref{trecu}, for $\xi\in S_{>\rho'}(\theta-\varepsilon,\theta+\varepsilon)$,
\[|\tilde{H}^{(m)}(\xi;u,A)|\le \sum_{k\ge 1}{t^{(m)}_k}|\xi|^{k-1}/\Gamma(k)\le \alpha \sum_{k\ge 1}\frac{(\beta|\xi|)^{k-1}}{\Gamma(k)}=\alpha e^{\beta|\xi|}. \]
Thus, each $H^{(m)}(\xi;u,A)$ is estimated by $\alpha e^{\beta|\xi|}$ with suitable $\alpha,\beta$ independent of $m$ and $u\in \h_{\rm reg}^d(\mathbb{R})$. Here we remark that the constant $c$ (therefore the other constants) depends on the choice of the small $\varepsilon$.

Therefore, the proof will be completed provided that we show the convergence of $\tilde{H}^{(m)}(\xi;u,A)$ as $m\rightarrow\infty$ locally uniform for all $\xi\in S_{>\rho'}(\theta-\varepsilon,\theta+\varepsilon)$. For that, let us denote by $W^{(m)}(\xi):=\tilde{H}^{(m)}(\xi;u,A)-\tilde{H}^{(m-1)}(\xi;u,A)$ the difference, then by \eqref{inteq} we have
\begin{equation}\label{Wm}
W^{(m)}(\xi)\cdot \I u-(\I u-\xi{\rm Id}_n)\cdot W^{(m)}(\xi)=\int_{t=0}^\xi\left(W^{(m-1)}(t)\cdot \frac{[A]}{2\pi\I }-\frac{A}{2\pi\I }\cdot W^{(m-1)}(t)\right)dt.
\end{equation}
Applying \eqref{xidod} to the case $M=W^{(m)}(\xi)$, from \eqref{Wm}, we deduce that
\[|W^{(m)}(\xi)|\le 2ca\cdot \int_{t=0}^\xi |W^{(m-1)}(t)|dt. \]
Similar to the estimate \eqref{Hmest}, by induction one proves that there exists a sequence of constants $b_m$ such that $b_m=2ac\cdot b_{m-1}$ and 
\[|W^{(m)}(\xi)|\le b_m \xi^{m-1}/\Gamma(m).\]
Therefore, there exists a large enough real number $b$ such that $b_m\le b^k$ for all $m,k$. That is $|W^{(m)}(\xi)|\le b^m \xi^{m-1}/\Gamma(m)$ and
\[\sum_{m=1}^\infty|W^{(m)}(\xi)|\le \sum_{m=1}^{\infty} b^m \xi^{k-1}/\Gamma(k).\]
Since $\tilde{H}^{(m)}(\xi)=\sum_{k=1}^m W^{(m)}(\xi)$ is the partial sum of the above absolute convergence series, we get that $\tilde{H}^{(m)}(\xi;u,A)$ converge uniformly on every compact subset of the domain $\xi\in  S_{>\rho'}(\theta-\varepsilon,\theta+\varepsilon)$. The ${\rm lim}_{m\rightarrow \infty} \tilde{H}^{(m)}(\xi)$ coincides with the Borel series $\tilde{H}(\xi)$ in the common domain $B_\rho\cap S_{>\rho'}(\theta-\varepsilon,\theta+\varepsilon)$.
\qed

\vspace{2mm}
{\bf Laplace transform.} The Laplace transform $H_\theta=\mathcal{L}_\theta(\tilde{H})$ of the function $\tilde{H}(\xi)$ in the direction $e^{\I \theta}$ is a function (in the Laplace plane of the initial
variable $z$) defined by
\[ H_\theta(z;u,A)=1+\int_{\xi=0}^{\infty(e^{\I \theta})}e^{-\frac{\xi}{z}}\tilde{H}(\xi;u,A)d\xi.\]
Lemma \ref{exgrow} ensures that for any $\theta$ with $e^{\I \theta}\in {\rm Sect}_\pm$, the integrand is indeed defined on the integral path, and that for any fixed $l>\beta$ the integral exists for all $z\in R(\theta,l)$, where the domain
\[R(\theta,l):=\{z\in\mathbb{C}~|~{\rm Re}(e^{\I \theta}/z)>l>\beta\}.\] The function \[H_\theta=\mathcal{L}_\theta(\tilde{H})=\mathcal{L}_\theta(\mathcal{B}(\hat{H})) \hspace{3mm} \text{on} \ \ R(\theta,l)\] is called the Laplace-Borel transform of $\hat{H}$ in the direction $e^{\I \theta}$. The following proposition shows that the original form power series $\hat{H}(z)$ is an asymptotics expansion of $H_\theta(z)$ as $z\rightarrow 0$ within $R(\theta,l)$. 
\begin{lem}\label{BLest}
For any $u\in \h_{\rm reg}^d(\mathbb{R})$ and $\theta\in (-\frac{\pi}{2}, \frac{\pi}{2})\cup (-\frac{3\pi}{2}, -\frac{\pi}{2})$, there exist constants $l, C,D>0$ such that 
\begin{eqnarray}\label{BLestineq}
|H_\theta(z;u,A)-1-\sum_{k=1}^{N-1}H_k(u)z^k|\le CN^Ne^{-N}\frac{|z|^N}{D^N}, \ \text{for all} \ z\in R(\theta,l), \ \ N\in\mathbb{N}_+.
\end{eqnarray}
Furthermore, the above constants $l, C, D$ can be chosen independent of $u\in \h_{\rm reg}^d(\mathbb{R})$.
\end{lem}
\begin{rmk}
For a fixed $u$, the proof of the inequality \eqref{BLestineq} is standard. In the following, we will go through the proof given in \cite[Theorem 5.3.9]{LR}, and show that in our case the involved constants can be chosen independent of $u\in \h_{\rm reg}^d(\mathbb{R})$.
\end{rmk}
\pf Without lose of generality, let us assume that $\theta=0$. By Lemma \ref{constant1} and Lemma \ref{exgrow}, there exist constants $K$, $\varepsilon$, $\alpha$ and $\beta>0$, such that $\tilde{H}(\xi;u,A)$ is holomorphic in the union of the disk $\{\xi:|\xi|\le 1/K\}$ and the sector $S(-\varepsilon,\varepsilon)$, and satisfies
\[|\tilde{H}(\xi;u,A)|\le \alpha e^{\beta |\xi|}, \ \text{for all} \ \ \xi\in  S(-\varepsilon,\varepsilon) \ \ \text{and} \ u\in \h_{\rm reg}^d(\mathbb{R}).\]

Just as in \cite[Theorem 5.3.9]{LR}, let us take a point $b$ with argument $\pi/4$ and small enough norm $|b|<1/K$ such that the path, following a straight line from $0$ to $b$ and
continues along a horizontal line from $b$ to $+\infty$, lies in the domain $\{\xi:|\xi|\le 1/K\}\cup S(-\varepsilon,\varepsilon)$ of the $\xi$ plane. Since the path is homotopy to $[0,+\infty)$, 
by the Cauchy’s theorem, the Laplace integral $H_{\theta=0}(z;u,A)$ decomposes to
$H_{\theta=0}(z;u,A)=1+H^b(z;u,A)+Y^b(z;u,A)$ along the path, where
\[H^b(z;u,A)=\int_{\xi=0}^{b}\tilde{H}(\xi;u,A)e^{-\frac{\xi}{z}}d\xi, \hspace{5mm} Y^b(z;u,A)=\int_{\xi=b}^{+\infty}\tilde{H}(\xi;u,A)e^{-\frac{\xi}{z}}d\xi. \]
On the one hand, given $0<\delta<\pi/2$, following \cite[Lemma 1.3.2]{LR}, we have
\begin{eqnarray}
|H^b(z;u,A)-\sum_{k=1}^{N-1}H_k(u)z^k|\le C'N^Ne^{-N}\frac{|z|^N}{{D'}^N}, \ \text{for} \ \ z\in S(-\pi/4+\delta,3\pi/4-\delta), \ \ N\in\mathbb{N}_+,
\end{eqnarray}
where the constants $C', D'$ are
\begin{eqnarray}
C':=\sum_{k\ge 1}\frac{|H_k(u)|}{\Gamma(k)}b^k, \hspace{5mm} D':=b\cdot {\rm sin}(\delta).
\end{eqnarray}
On the other hand, let us take a constant $l>\beta$, then following the proof of \cite[Theorem 5.3.9]{LR}, we have 
\begin{eqnarray}\label{exgrowest}
|Y^b(z;u,A)|\le he^{-\frac{c}{z}}, \ \ \text{for} \ z\in S(-\pi/4+\delta,3\pi/4-\delta)\cap R(\theta=0,l),
\end{eqnarray} with the constants given by
\begin{eqnarray}
h=\frac{\alpha e^{\beta|b|}}{l-\beta}, \hspace{5mm} c=|b|{\rm cos}(\pi/2-\delta).
\end{eqnarray}
The estimation \eqref{exgrowest} further implies, see e.g,. \cite[Proposition 1.2.17]{LR}, for $z\in S(-\pi/4+\delta,3\pi/4-\delta)\cap R(\theta=0,l)$, $N\in\mathbb{N}_+$,
\begin{eqnarray*}
|Y^b(z;u,A)-\sum_{k=1}^{N-1}H_k(u)z^k|\le C''N^Ne^{-N}\frac{|z|^N}{{D''}^N}, 
\end{eqnarray*}
with the constants $C''$ and $D''$ determined by $h$ and $c$. In conclusion, if we take $C = {\rm max}(C', C'')$ and $D = {\rm max}(D', D'')$, then $H_\theta(z;u,A)$ satisfies the inequality \eqref{BLestineq} on the domain
\[S(-\pi/4+\delta,3\pi/4-\delta)\cap R(\theta=0,l).\]
According to \cite[Theorem 5.3.9]{LR}, an argument using the symmetry with respect to the real axis, i.e., by choosing $\bar{b}$ instead of $b$ and the corresponding path, shows that $H_{\theta=0}(z)$ satisfies the inequality \eqref{BLestineq} on the symmetric domain 
\[S(-3\pi/4+\delta,\pi/4-\delta)\cap R(\theta=0,l)\]
with respect to the real axis. Since the union of the above two domains cover $R(\theta=0,l)$, we get $H_{\theta=0}(z)$ satisfies the inequality \eqref{BLestineq} on $R(\theta=0,l)$. 

In the end, let us check the independence of constants $C = {\rm max}(C', C'')$ and $D = {\rm max}(D', D'')$ on $u_n$. First, following Lemma \ref{constant1}, we have $C'\le \frac{Kb}{1-Kb}$. Thus we can set $C'=\frac{Kb}{1-Kb}$, and then the constants $C,D$ are determined by $\varepsilon, K,\alpha,\beta,b$ and $\delta$. By Lemma \ref{constant1} and Lemma \ref{exgrow}, as a sufficiently small $\varepsilon$ fixed, those constants can be chosen independent of $u\in \h_{\rm reg}^d(\mathbb{R})$. It therefore proves the proposition. \qed

\vspace{2mm}
{\bf Canonical solutions on Stokes sectors.} The functions $H_\theta$ and $H_{\theta'}$, with $e^{\I \theta}, e^{\I \theta'}\in {\rm Sect}_+$ in the same Stokes sector, coincide in the overlapping of their defining domains. See e.g., \cite[Proposition 5.3.7]{LR} or \cite[Section 6.2]{Balser0}. (While if $\theta$ and $\theta'$ is not in the same Stokes sector,  $H_\theta(z;u,A)$ is in general not equal to $H_{\theta'}(z;u,A)$ at the points $z$ where both functions are defined.)  Thus the functions $H_\theta$ for all $e^{\I \theta}\in {\rm Sect}_+$ glue together into a holomorphic function $H_+(z;u,A)$ defined on the domain $S(-\pi,\pi)$. Furthermore, since the Laplace-Borel transform is a morphism of differential algebras from the algebra of power series to the algebra of holomorphic functions, $H_+(z;u,A)$ satisfies the same equation obeyed by $\hat{H}(z;u,A)$. In this way, we see that the function \[F_+(z;u,A)=({\rm Id}_n+H_+(z^{-1};u))e^{\I uz}z^{-\frac{[A]}{2\pi\I }}\] is a solution of \eqref{Stokeseq} with the prescribed asymptotics at $z=\infty$ within $S(-\pi,\pi)$. It follows from the Waston's Lemma (see \cite[Section 4.7]{Balser0}) that once an actual solution of \eqref{Stokeseq} is obtained which has the formal fundamental solution matrix \eqref{formalsol} as its asymptotic expansion as $z\rightarrow\infty$ in a sector whose opening is larger than $\pi$, then such an actual solution is unique. Similarly, we can construct the unique (therefore canonical) solution $F_-(z;u,A)$ in $S(-2\pi,0)$.

As for the uniform property with respect to $u$, we just observe that by the proof of Lemma \ref{BLest}, the constants in the inequality \eqref{BLestineq} can be chosen independent of $u\in \h_{\rm reg}^d(\mathbb{R})$, as long as $\theta$ keeps a fixed positive distance $2\varepsilon>0$ with $\pi$ and $-\pi$. Therefore, in each proper closed subsector $\bar{S}(-\pi+\varepsilon_0,\pi-\varepsilon_0)$ of $S(-\pi,\pi)$ for a small $\varepsilon_0>0$, by \eqref{BLestineq} the limit
\[\lim_{z\rightarrow \infty}|H_+(z^{-1};u)|=\lim_{z\rightarrow 0}|H_+(z;u,A)|=0\]
is uniform for all $u\in \h_{\rm reg}^d(\mathbb{R})$. The same statement is true for $H_-(z;u,A)$. It finishes the proof of Theorem \ref{uniformresum}. \qed

\begin{rmk}
The above argument implies stronger results (but the present form is enough for the use in the proof of Proposition \ref{inductionstep}): the limit
\[\lim_{z\rightarrow 0}z^{-m}(|H_\theta(z;u,A)-\sum_{k=0}^{N-1}H_k(u,A)z^k|)=0, \ \ \ z\in R(\theta,l)\]
for all $m>0$, is uniform for $u\in \h_{\rm reg}^d(\mathbb{R}).$
\end{rmk}
Similar argument shows that Theorem \ref{uniformresum} generalizes to the case $u\in\h(\mathbb{R})$, the space of $n\times n$ real diagonal matrices $u={\rm diag}(u_1,...,u_n)$. To be more precise, let us take any partition of the set \[\{1,...,n\}=I_{1}\cup I_2\cup\cdots\cup I_k.\] 
Then for any $u={\rm diag}(u_1,...,u_n)$, satisfying $u_i=u_j$ if $i,j \in I_s$ for some $s$; $u_i\ne u_j$ otherwise, there exist canonical solutions $F_\pm(z;u,A)$ in ${\rm Sect}_\pm$ such that $F_+\cdot e^{-\I uz}\cdot z^{\frac{\delta_u(A)]}{2\pi\I }}$ and $F_-\cdot e^{-\I uz}\cdot z^{\frac{\delta_u(A)}{2\pi\I }}$ can be analytically continued to $S(-\pi,\pi)$ and $S(0,2\pi)$ respectively, and for every small $\varepsilon>0$,
\begin{eqnarray*}
\lim_{z\rightarrow\infty}F_+(z;u,A)\cdot e^{-\I uz}\cdot z^{\frac{\delta_u(A)}{2\pi\I }}&=&{\rm Id}_n, \ \ \ as \ \ \ z\in \bar{S}(-\pi+\varepsilon,\pi-\varepsilon),
\\
\lim_{z\rightarrow\infty}F_-(z;u,A)\cdot e^{-\I uz}\cdot z^{\frac{\delta_u(A)}{2\pi\I }}&=&{\rm Id}_n, \ \ \ as \ \ \ z\in \bar{S}(-2\pi+\varepsilon,-\varepsilon),
\end{eqnarray*}
Here $\delta_u(A)$ is the projection of $A$ on the centralizer of $u$ in $\frak{gl}_n$, i.e., 
\[\delta_u(A)_{ij}=\left\{
          \begin{array}{lr}
            A_{ij},   & \text{if} \ \ i, \ j \in  \ I_s \ \text{for some} \ s\\
           0, & \text{otherwise},
             \end{array}
\right. \]
Furthermore, the above limits are uniform for all $u$ such that if $u_i$ is not identically equal to $u_j$, then $|u_i-u_j|$ is bigger than a fixed positive number $d$.


\subsection{Stokes matrices and connection matrices}\label{secSC}
For any $\sigma\in S_n$, let us denote by $U_{\sigma}$ the component $\{ u_{\sigma(1)}<\cdot\cdot\cdot < u_{\sigma(n)}\}$ of $\h_{\rm reg}(\mathbb{R})$, and denote by $P_\sigma\in{\rm GL}_n$ the corresponding permutation matrix.

\begin{defi}\label{defiStokes}
For any $u\in U_\sigma$, the {\it Stokes matrices} of the system \eqref{Stokeseq} (with respect
to ${\rm Sect}_+$ and the chosen branch of ${\rm log}(z)$) are the elements $S_\pm(u,A)\in {\rm GL}(n)$ determined by
\begin{align*}
F_+(z;u,A)&=F_-(z;u,A)\cdot e^{-\frac{[A]}{2}}P_\sigma S_+(u,A)P_\sigma^{-1},\\
F_-(ze^{-2\pi \I };u,A)&=F_+(z;u,A)\cdot P_\sigma S_-(u,A)P_\sigma^{-1}e^{\frac{[A]}{2}},
\end{align*}
where the first (resp. second) identity is understood to hold in ${\rm Sect}_-$
(resp. ${\rm Sect}_+$) after $ F_+$ (resp. $F_-$)
has been analytically continued anticlockwise around $z=\infty$. 
\end{defi} 
The prescribed asymptotics of $F_\pm(z;u,A)$ at $z=\infty$, as well as the identities in Definition \ref{defiStokes}, ensures that the Stokes matrices $S_+(u,A)$ and $S_-(u,A)$ are upper and lower triangular matrices respectively. see e.g., \cite[Chapter 9.1]{Balser0}. Furthermore, the following lemma follows from the fact that if $F(z;u,A)$ is a solution, so is $F(\bar{z};u,A)^\dagger$, see \cite{Boalch1}.
\begin{lem}
Let $S_+(u,A)^\dagger$ denote the conjugation transpose of $S_+(u,A)$, then $S_-(u,A)=S_+(u,A)^\dagger$.
\end{lem}

Since the system \eqref{Stokeseq} is non-resonant, i.e., no two eigenvalues of $\frac{A}{2\pi\I }$ for $A\in\Herm(n)$ are differed by a positive integer, we have (see e.g \cite[Chapter 2]{Balser0}).
\begin{lem}\label{le:nr dkz}
There is a unique holomorphic fundamental solution
$F_0(z;u,A)\in {\rm GL}(n)$ of the system \eqref{Stokeseq} on a neighbourhood of $\infty$ slit along $\I \mathbb{R}_{\ge 0}$, such that $F_0\cdot z^{\frac{A}{2\pi\I }}\rightarrow {\rm Id}_n$ as $z\rightarrow 0$.
\end{lem}

\begin{defi}\label{connectionmatrix}
The {\it connection 
matrix} $C(u,A)\in {\rm GL}_n(\mathbb{C})$ of the system \eqref{Stokeseq} (with respect to ${\rm Sect}_+$) is determined by 
\[F_0(z;u,A)=F_+(z;u,A)\cdot C(u,A), \]
as $F_0(z;u,A)$ is
extended to the domain of definition of $F_+(z;u,A)$.
\end{defi}
The connection matrix $C(u,A)$ is valued in ${\rm U}(n)$ (see e.g., \cite[Lemma 29]{Boalch1}).
Thus for any fixed $u$, by varying $A\in{\Herm}(n)$ we obtain the connection map
\begin{eqnarray}\label{Cmap}
C(\cdot, u)\colon{\rm Herm}(n)\rightarrow {\rm U}(n).\end{eqnarray}

In a global picture, the connection matrix is related to the Stokes matrices
by the following monodromy relation, which follows from the fact that a simple negative loop (i.e., in clockwise direction) around $0$ is a simple positive loop (i.e., in anticlockwise direction) around $\infty$: for any $u\in U_\sigma\subset \h_{\rm reg}(\mathbb{R})$,
\begin{eqnarray}\label{monodromyrelation}
C(u,A)e^{A}C(u,A)^{-1}=P_\sigma S_-(u,A)S_+(u,A) P_\sigma^{-1}.
\end{eqnarray}

\section{Expression of Stokes matrices via the boundary values of solutions of isomonodromy equations}\label{exStokesviaiso}
In Section \ref{sec:iso}, we recall the isomonodromy equations of the meromorphic linear systems of differential equation \eqref{Stokeseq}. In Section \ref{asyisoeq} and \ref{proof11}, we study the asymptotics of the solutions of the isomonodromy equation at a critial point, and then give a proof of Theorem \ref{isomonopro}. In Sections \ref{RHtorelRH} and \ref{exviaasy}, we prove an analytic branching rule of Stokes matrices with respect to the collapse of the components of the irregular parameter $u$. The branching rule enables us to express the Stokes matrices via the asymptotics of solutions of isomonodromy equation at an "infinite" point. Then in Section \ref{explicitnurel}, we prove our first main result, i.e, Theorem \ref{mainthm}, using the known global behavior of the solutions of confluent hypergeometric systems. In Sections  \ref{diffPhi}-\ref{pf:introthm3}, we give a proof of Theorem \ref{introcor}. In Section \ref{leadterm}, we get the explicit leading order terms of $S_\pm(u,A)$ as $u$ approaches to the "infinite" point. In the end, Section \ref{dirreg} studies the analog of the results in Section \ref{asyisoeq}-\ref{leadterm} for a degenerate irregular parameter $u$.

\subsection{Isomonodromy deformation}\label{sec:iso}
In this subsection, we recall some facts about the theory of isomonodromy deformation. In general, the Stokes matrices $S_\pm(u,A)$ of the system \eqref{Stokeseq} will depend on the irregular term $u$.
The isomonodromy deformation (also known as monodromy preserving) problem is to find the matrix valued function $\Phi(u)$ such that the Stokes matrices $S_\pm(u,\Phi(u))$ are (locally) constant. The following definition and proposition are well known. See more detailed discussions in e.g., \cite{BoalchG,JMMS, JMU}, \cite[Chapter 3]{Dubrovin}. In particular, the equation \eqref{isoeq} is the Jimbo-Miwa-M\^ori-Sato equation \cite{JMMS} with one irregular singularity and one regular singularity. 
\begin{defi}\label{isomonodromyequation}
The isomonodromy equation is the differential equation for a matrix valued function $\Phi(u): \h_{\rm reg}(\mathbb{R})\rightarrow {\rm Herm}(n)$
\begin{equation}\label{isoeq}
\frac{\partial \Phi}{\partial u_k} =\frac{1}{2\pi\I }[\Phi,{\rm ad}^{-1}_u{\rm ad}_{E_{k}}\Phi], \ \text{for all} \ k=1,...,n.\end{equation} 
Here $E_k$ is the $n\times n$ diagonal matrix whose $(k,k)$-entry is $1$ and other entries are $0$.
Note that ${\rm ad}_{E_{k}}\Phi$ takes values in the space ${\frak {gl}}_n^{od}$ of off diagonal matrices and that ${\rm ad}_u$
is invertible when restricted to ${\frak {gl}}_n^{od}$.
\end{defi} 
\begin{rmk}
The isomonodromy equation with respect to the derivation of $u_j$ is generated by the time dependent quadratic Hamiltonian $H_j:=(-\frac{1}{2\pi\I })\sum_{k\ne j}\frac{\phi_{kj}\phi_{jk}}{u_k-u_j}$, where $\phi_{ij}$'s are the entry functions on $\Herm(n)$, see e.g., \cite{Dubrovin, BoalchG}.
\end{rmk}
Set $\Phi(u)=(\phi_{ij}(u))$, then in terms of the components, the equation \eqref{isoeq} becomes
\begin{eqnarray*}
  \frac{\partial}{\partial u_k} \phi_{i j} (u)
  & = & 
  \frac{1}{2\pi\I }\left(
  \frac{1}{u_k - u_j}-\frac{1}{u_k - u_i} \right) 
  \phi_{i k}(u) \phi_{k j}(u), \quad i, j \neq k,
  \label{uisoeqij}
  \\
  \frac{\partial}{\partial u_k} \phi_{i k} (u)
  & = & \frac{1}{2\pi\I }
  \sum_{j \neq k} 
  \frac{\phi_{i j}(u)\phi_{j k}(u)-\delta_{i j}\phi_{k k}(u) \phi_{j k}(u)}{u_k - u_j} ,
  \quad i \neq k ,
  \label{uisoeqik}
  \\
  \frac{\partial}{\partial u_k} \phi_{k j} (u)
  & = & \frac{1}{2\pi\I }
  \sum_{i \neq k} 
  \frac{\delta_{i j}
  \phi_{k k}(u)\phi_{k i}(u)-\phi_{k i}(u)\phi_{i j}(u)}{u_k - u_i},
  \quad j \neq k ,
  \label{uisoeqkj}\\
  \frac{\partial}{\partial u_k} \phi_{k k}(u) 
  & = & 
  0.
  \label{uisoeqkk}
\end{eqnarray*}

\begin{pro}\label{isomonodef}
For any solution $\Phi(u)$ of the isomonodromy equation, the canonical solutions $F_\pm$ of the system
\begin{eqnarray}\label{isoStokeseq1}
\frac{dF}{dz}=\left(\I u-\frac{1}{2\pi\I }\frac{\Phi(u)}{z}\right)\cdot F,
\end{eqnarray}
satisfy (as a function of $u={\rm diag}(u_1,...,u_n)$)
\begin{eqnarray}\label{isoStokeseq2}
\frac{\partial F}{\partial u_k}=\left(\I E_kz-\frac{1}{2\pi\I } {\rm ad}^{-1}_u{\rm ad}_{E_{k}}\Phi(u)\right)\cdot F.
\end{eqnarray}
In particular, the Stokes matrices $S_\pm(u,\Phi(u))$ of \eqref{isoStokeseq1}
are locally constants (independent of $u$). Furthermore the isomonodromy equation \eqref{isoeq} is the compatibility condition of the systems \eqref{isoStokeseq1} and \eqref{isoStokeseq2}.
\end{pro}

\subsection{Boundary values of the solutions of the isomonodromy equation}\label{asyisoeq}
In this subsection, we will prove Theorem \ref{isomonopro}. Recall that we denote by $U_{\rm id}$ the connected component $\{u\in \h_{\rm reg}(\mathbb{R})~|~u_1<\cdots <u_n\}$ of $\h_{\rm reg}(\mathbb{R})$, and by $\delta_k(A)$ the matrix
 \begin{eqnarray}\label{delta1}\delta_k(A)_{ij}=\left\{
          \begin{array}{lr}
             A_{ij},   & \text{if} \ \ 1\le i, j\le k, \ \text{or} \ i=j  \\
           0, & \text{otherwise}.
             \end{array}
\right.\end{eqnarray}
We first prove a stronger statement, from which Theorem \ref{isomonopro} follows. We introduce the new coordinates
\begin{align}\label{newcoor}
z_0=u_1+u_2+\cdots +u_n, \ z_1=u_2-u_1, \ z_2=\frac{u_3-u_2}{u_2-u_1}, \ ...... , \ z_{n-1}=\frac{u_n-u_{n-1}}{u_{n-1}-u_{n-2}}.
\end{align}
Then $\sum_{i=1}^n u_i=z_0$, and
\begin{align}
 u_2=u_1+z_1, \ u_3=u_1+z_1+z_1z_2, \ ......, \ u_n=u_1+z_1+z_1z_2+z_1z_2z_3+...+z_1\cdots z_{n-1}.
\end{align}
In terms of the new coordinates, equation \eqref{isoeq} becomes
\begin{align*}
\frac{\partial\Phi}{\partial z_{j}}&=\frac{1}{2\pi\I }\Big[\Phi,{\rm ad}_{u}^{-1}{\rm ad}_{\frac{\partial u}{\partial z_j}}\Phi\Big], \ \text{for} \ j=1,...,n-1,\\
\frac{\partial\Phi}{\partial z_{0}}&=0,
\end{align*}
where the diagonal matrix
\[\frac{\partial u}{\partial z_j}=z_j^{-1}\left((u_{j+1}-u_{j})E_{j+1}+(u_{j+2}-u_{j})E_{j+2}+\cdots+(u_n-u_{j})E_{n}\right).\]
Therefore, any solution $\Phi$ is independent of $z_0$, and is thus a function of $z_1,...,z_{n-1}$.
\begin{pro}\label{strongerpro}
For any Hermitian matrix valued solution $\Phi(z_1,...,z_{n-1})$ of \eqref{isoeq} on $U_{\rm id}$, there exists a chain of functions $\Phi_{k}(z_1,...,z_{k})\in\Herm(n)$ for all $k=0,...,n-1$ such that $\Phi_{n-1}=\Phi$ and $\Phi_0$ is a constant, and for any $k=2,...,n-1$,
    \begin{align}
      \label{Asym:dk-1_Phik}
      \delta_{k}(\Phi_k(z_1,...,z_k)) 
      & =  
      \delta_{k} (\Phi_{k-1}(z_1,...,z_{k-1})) + \mathcal{O}(z_k^{-1}),\\
      \label{Asym:uPhik-1_Phik}
      z_k^{\frac{-{\rm ad}\delta_{k} (\Phi_{k-1})}{2\pi\I }} (\Phi_k)
      & = 
      \Phi_{k-1} + \mathcal{O}(z_k^{-1}).
    \end{align}
as $z_k\rightarrow +\infty$, uniformly with respect to $z_1,...,z_{k-1}$. 
Here $z_k^{\frac{{\rm ad}X}{2\pi\I }}Y:=z_k^{\frac{X}{2\pi\I }} Y z_k^{-\frac{X}{2\pi\I }}$ for any $X,Y\in\Herm(n)$.
Furthermore, the function $\Phi_k$ satisfy the differential equations
\begin{eqnarray}\label{kdiffeq}
\frac{\partial\Phi_{k}}{\partial z_{j}}=\frac{1}{2\pi\I }\Big[\Phi_{k},({\rm ad}_{u}^{-1}{\rm ad}_{\frac{\partial u}{\partial z_j}}-z_j^{-1})\delta_{k+1}(\Phi_{k})\Big], \ \text{ for } \ j=1,...,k.
\end{eqnarray}
In the end, for $k=1$,
\begin{equation}\label{k=1}
    z_1^{\frac{-{\rm ad}\delta_{1} (\Phi_{0})}{2\pi\I }} (\Phi_1)= \Phi_{0}.
\end{equation}
We call $\Phi_{k-1}$ the boundary value of $\Phi_k$ at $z_k=+\infty$.
\end{pro}
\begin{proof} Case I: assume that an integer $k\in\{2,...,n-1\}$, and a solution $\Phi_k(z_1,...,z_k)$ of the equation \eqref{kdiffeq}, are given, let us construct a function $\Phi_{k-1}(z_1,...,z_{k-1})$ such that \eqref{Asym:dk-1_Phik} and \eqref{Asym:uPhik-1_Phik} hold. For this purpose, let us consider the differential equation of $\Phi_k$ with respect to the parameter $z_k$.

Firstly, we have

\begin{lem}\label{bound}
The solution $\Phi_k$ satisfies \[\frac{d}{dz_j}{\rm Tr}(\Phi_k^\dagger\Phi_k)=0, \text{ for } j=1,...,k.\]
Therefore, the norm of $\Phi_k$ is uniformly bounded for all $z_j$.
\end{lem}
\begin{proof} We have
\begin{eqnarray*}
\frac{\partial}{\partial z_j}({\rm Tr}(\Phi_k^\dagger\Phi_k))={\rm Tr}(\frac{\partial\Phi_k^\dagger}{\partial z_j} \Phi_k+\Phi_k^\dagger \frac{\partial\Phi_k}{\partial z_j})=\frac{1}{\pi \I }{\rm Im}\Big[\Phi_{k},({\rm ad}_{u}^{-1}{\rm ad}_{\frac{\partial u}{\partial z_j}}-z_j^{-1})\delta_{k+1}(\Phi_{k})\Big]=0.
\end{eqnarray*}
Here the second identity follows from $\Phi_k=\Phi_k^\dagger\in\Herm(n)$, and the last identity follows from the fact that ${\rm Tr}(X\cdot [X,Y]])=0$ for any matrices $X,Y$.
\end{proof}

Secondly, we have 
\begin{lem}
There exist functions $J_{k-1}(z_1,...,z_{k-1})$ and $\Phi_{k-1}(z_1,...,z_{n-1})\in \Herm(n)$ such that
\begin{align}\label{dbound}
    \delta_{k}(\Phi_k(z_1,...,z_k))&=J_{k-1}(z_1,...,z_{k-1})+\mathcal{O}(z_k^{-1}),\\ \label{est1}
z_k^{-\frac{{\rm ad}J_{k-1}}{2\pi\I }}\Phi_k(z_1,...,z_k)&= \Phi_{k-1}(z_1,...,z_{k-1})+\mathcal{O}(z_k^{-1}),
\end{align}
as $z_k\rightarrow\infty$, uniformly with respect to $z_1,...,z_{n-1}$. 
\end{lem} 
\begin{proof} Note that
\begin{eqnarray}
  \label{uzchange}
  {\rm ad}_u^{- 1} {\rm ad}_{\partial u / \partial z_k} \delta_{k+1}(\Phi_k) =
  - {\rm ad}_{D^{(k)}_k} {\rm ad}_{E_{k+1}}
  \delta_{k+1}(\Phi_k),
\end{eqnarray}
where the $n \times n$ diagonal matrix with $k+1$ none zero elements
\begin{eqnarray}
\label{Dkk:Expression}
D^{(k)}_k ={\rm diag} \left( 
\frac{1}{z_k + \frac{u_{k} - u_1}{u_{k} - u_{k - 1}}}
, \ldots, 
\frac{1}{z_k + \frac{u_{k} - u_{k - 2}}{u_{k} - u_{k - 1}}}, 
\frac{1}{z_k + 1},
\frac{1}{z_k}
, 0, \ldots, 0 \right).
\qquad\quad
\end{eqnarray}
Then we have (in terms of the block matrix form)
\begin{eqnarray*}
\frac{\partial\delta_{k}(\Phi_k)}{\partial z_k}&=&\frac{1}{2\pi\I }\delta_{k}\left(\Big[\Phi_{k}, (-{\rm ad}_{D^{(k)}_k} {\rm ad}_{E_{k+1}}-z_k^{-1})
  \delta_{k+1}(\Phi_{k})\Big]\right)\\
&=&\frac{1}{2\pi\I }\left(\begin{array}{cc}
    \left(((D^{(k)}_k)_{jj}-(D^{(k)}_k)_{ii})(\Phi_k)_{i,k+1}(\Phi_k)_{k+1,j}\right)_{i,j=1}^{k} & 0  \\
    0 & 0
  \end{array}\right).
\end{eqnarray*}
Since $0\le \frac{u_{k} - u_i}{u_{k} - u_{k - 1}}\le 1$ for all $i=1,...,k$, we have 
\[(D^{(k)}_k)_{jj}-(D^{(k)}_k)_{ii}=\frac{1}{z_k + \frac{u_{k} - u_j}{u_{k} - u_{k - 1}}}-\frac{1}{z_k + \frac{u_{k} - u_i}{u_{k} - u_{k - 1}}}=\mathcal{O}(z_k^{-2}) \]
as $z_k\rightarrow+\infty$ uniformly with respect to $z_1,...,z_{k-1}$. Furthermore, by Lemma \ref{bound},
$\Phi_k=\mathcal{O}(1)$ as $z_k\rightarrow\infty$. We thus get
\begin{align}\label{firstest}
    \frac{\partial \delta_{k}(\Phi_k)}{\partial z_k}=\mathcal{O}(z_k^{-2}),
\end{align}
uniformly with respect to $z_1,...,z_{k-1}$. Therefore, by integrating \eqref{firstest}, we see that for any $z_1,...,z_{n-1}$ there exists $J_{k-1}(z_1,...,z_{n-1})$ such that \eqref{dbound} holds. 

Now let us consider the function $z_k^{-\frac{{\rm ad}J_{k-1}}{2\pi\I }}\Phi_k:=z_k^{-\frac{J_{k-1}}{2\pi\I }}\cdot \Phi_k \cdot z_k^{\frac{J_{k-1}}{2\pi\I }}$. It satisfies
\begin{eqnarray}\label{nthstep}
\frac{d (z_k^{-\frac{{\rm ad}J_{k-1}}{2\pi\I }}\Phi_k)}{d z_k} =\frac{1}{2\pi\I }[z_k^{-\frac{{\rm ad}J_{k-1}}{2\pi\I }}\Phi_k,B(z_k)], 
\end{eqnarray}
with 
\begin{eqnarray*}
B(z_k)=z_k^{-\frac{J_{k-1}}{2\pi\I }}\left( 
 (-{\rm ad}_{D^{(k)}_k} {\rm ad}_{E_{k+1}}-z_k^{-1})
  \delta_{k+1}(\Phi_{k})
  -z_k^{-1}(\Phi_k-J_{k-1})\right) z_k^{\frac{J_{k-1}}{2\pi\I }}.
\end{eqnarray*}

By $\Phi_k=\mathcal{O}(1)$ and the expression \eqref{uzchange} of $D^{(k)}_k$, we get 
\begin{eqnarray}\label{est0}
(-{\rm ad}_{D^{(k)}_k} {\rm ad}_{E_{k+1}}-z_k^{-1})
  \delta_{k+1}(\Phi_{k})=z_k^{-1}(\delta_{k+1}(\Phi_k)-\delta_{k}(\Phi_k))+\mathcal{O}(z_k^{-2}).
\end{eqnarray}
Therefore, by Lemma \ref{dbound},
\[
 (-{\rm ad}_{D^{(k)}_k} {\rm ad}_{E_{k+1}}-z_k^{-1})
  \delta_{k+1}(\Phi_{k})
  -z_k^{-1}(\Phi_k-J_{k-1})=\mathcal{O}(z_k^{-2}). \]
Since $J_{k-1}/{2\pi\I }$ is skew-Hermitian, i.e., $|(z_k)^{\frac{J_{k-1}}{2\pi\I }}|=1$, we get $B(z_k)=\mathcal{O}(z_k^{-2})$ and
\[
\frac{d (z_k^{-\frac{{\rm ad}J_{k-1}}{2\pi\I }}\Phi_k)}{d z_k} =\mathcal{O}(z_k^{-2}). \]
Integrating the above identity, we find a function $\Phi_{k-1}(z_1,...,z_{n-1})\in \Herm(n)$ such that \eqref{est1} holds. \end{proof}

Thirdly, we show that $J_{k-1}$ is actually equal to $\delta_{k}(\Phi_{k})$. The identity \ref{dbound} implies that
\begin{eqnarray}\label{exest}
z_{k}^{-\frac{\delta_{k}(\Phi_{k})}{2\pi\I }}=z_{k}^{-\frac{J_{k-1}}{2\pi\I }}+\mathcal{O}({\rm In}(z_k)/z_k).
\end{eqnarray}
Here we use the fact that for any two skew-Hermitian matrices $X$ and $Y$, 
the inequality $|e^X-e^Y|\le |X-Y|$ holds (simply note that $e^X-e^Y=\int_0^1\frac{d}{dt}(e^{(1-t)Y}e^{tX})dt=\int_0^1(e^{(1-t)Y}(X-Y)e^{tX})dt,$
and $|e^{(1-t)Y}|=|e^{tX}|=1$).

Thus, identities \eqref{est1} and \eqref{exest} lead to 
\begin{equation}\label{k-1klimit}
z_k^{-\frac{{\rm ad}\delta_{k}(\Phi_k)}{2\pi\I }}\Phi_k(z_1,...,z_k)= \Phi_{k-1}(z_1,...,z_{k-1})+\mathcal{O}({\rm In}(z_k)/z_k).
\end{equation}
Taking the operator $\delta_{k}$ on both sides of \eqref{k-1klimit} leads to \[\lim_{z_k\rightarrow+\infty}\delta_{k}(\Phi_k)= \delta_{k}(\Phi_{k-1}),\]
which (together with \eqref{dbound}) implies 
\[J_{k-1}=\delta_{k}(\Phi_{k-1}).\]
Then \eqref{dbound} and \eqref{est1} becomes \eqref{Asym:dk-1_Phik} and \eqref{Asym:uPhik-1_Phik} respectively provided replacing $J_{k-1}$ by $\delta_{k}(\Phi_{k-1})$.

In the end, we have to prove that the new function $\Phi_{k-1}(z_1,...,z_{k-1})$ satisfies the equation \eqref{kdiffeq} provided replacing $k$ by $k-1$ in \eqref{kdiffeq}. Similar to the above discussion, we verify that as $z_k\rightarrow+\infty$ for all $j<k$,
\begin{align*}
\frac{\partial\delta_{k}(\Phi_k)}{\partial z_j}&=\frac{1}{2\pi\I }\Big[\delta_{k}(\Phi_k), ({\rm ad}_{u}^{-1}{\rm ad}_{\frac{\partial u}{\partial z_j}}-z_j^{-1})\delta_{k}(\Phi_{k})\Big]+\mathcal{O}(z_k^{-1}),\\
\frac{\partial\Phi_{k}}{\partial z_j}&=\frac{1}{2\pi\I }\Big[\Phi_{k}, ({\rm ad}_{u}^{-1}{\rm ad}_{\frac{\partial u}{\partial z_j}}-z_j^{-1})\delta_{k}(\Phi_{k})\Big]+\mathcal{O}(z_k^{-1}),
\end{align*}
uniformly for all $z_j$ with $j<k$. The above two identities give rise to
\begin{equation}\label{kjequation}
\frac{\partial}{\partial z_j}\left(z_k^{-\frac{{\rm ad}\delta_{k}(\Phi_k)}{2\pi\I }}(\Phi_k)\right)=\frac{1}{2\pi\I }\Big[z_k^{-\frac{{\rm ad}\delta_{k}(\Phi_k)}{2\pi\I }}(\Phi_k), \ ({\rm ad}_{u}^{-1}{\rm ad}_{\frac{\partial u}{\partial z_j}}-z_j^{-1})\delta_{k}(\Phi_{k})\Big]+\mathcal{O}(z_k^{-1}),
\end{equation}
as $z_k\rightarrow +\infty$. Based on \eqref{k-1klimit}, letting $z_k\rightarrow +\infty$ in \eqref{kjequation} leads to
\[\frac{\partial \Phi_{k-1}}{\partial z_j}=\frac{1}{2\pi\I }\Big[\Phi_{k-1}, \ ({\rm ad}_{u}^{-1}{\rm ad}_{\frac{\partial u}{\partial z_j}}-z_j^{-1})\delta_{k}(\Phi_{k-1})\Big],\]
which is just the equation \eqref{kdiffeq} for $\Phi_{k-1}$ (replacing $k$ by $k-1$ in \eqref{kdiffeq}). By induction, it finishes the proof of the proposition for any $k=2,...,n$.

\vspace{2mm}
Case II: $k=1$. Then $\Phi_1(z_1)$ is a solution of 
\[\frac{d \Phi_{1}}{d z_1}=\frac{1}{2\pi\I }\Big[\Phi_{1}, \ ({\rm ad}_{u}^{-1}{\rm ad}_{\frac{\partial u}{\partial z_1}}-z_j^{-1})\delta_{2}(\Phi_{1})\Big]. \]
This equation can be simply solved explicitly, and any solution takes the form of \eqref{k=1}, i.e., 
\[
    z_1^{\frac{-{\rm ad}\delta_{1} (\Phi_{0})}{2\pi\I }} (\Phi_1)= \Phi_{0}.
\]
for a unique constant $\Phi_0$. 

Therefore, given any solution $\Phi$ of \eqref{isoeq}, by induction it finishes the proof of existence of the chain of functions $\Phi_{n-1}=\Phi, \Phi_{n-2}, ..., \Phi_1,$ and the constant $\Phi_0$.
\end{proof}

\subsection{The proof of Theorem \ref{isomonopro}}\label{proof11}
\begin{proof}[The proof of Theorem \ref{isomonopro}]
On the one hand, by Proposition \ref{strongerpro}, we have
\begin{align} \nonumber
     &z_{k-1}^{\frac{-\delta_{k-1} (\Phi_{k-2})}{2\pi\I }}  z_k^{\frac{-\delta_{k} (\Phi_{k-1})}{2\pi\I }} (\Phi_k)z_k^{\frac{\delta_{k} (\Phi_{k-1})}{2\pi\I }} z_{k-1}^{\frac{\delta_{k-1} (\Phi_{k-2})}{2\pi\I }}  \\ \nonumber
      =& \label{k-2k-1}
      z_{k-1}^{\frac{-\delta_{k-1} (\Phi_{k-2})}{2\pi\I }}  \Phi_{k-1} z_{k-1}^{\frac{\delta_{k-1} (\Phi_{k-2})}{2\pi\I }} + \mathcal{O}(z_k^{-1})\\ 
      =&\Phi_{k-2}+\mathcal{O}(z_{k-1}^{-1})+\mathcal{O}(z_k^{-1}).
\end{align}

On the other hand, we have
\begin{equation}\label{Wdefi}
    z_{k-1}^{\frac{-\delta_{k-1} (\Phi_{k-2})}{2\pi\I }}  z_k^{\frac{-\delta_{k} (\Phi_{k-1})}{2\pi\I }} = z_k^{\frac{-W}{2\pi\I }} z_{k-1}^{\frac{-\delta_{k-1} (\Phi_{k-2})}{2\pi\I }} 
\end{equation}
where \[W= z_{k-1}^{\frac{-\delta_{k-1} (\Phi_{k-2})}{2\pi\I }} (\delta_{k} (\Phi_{k-1}))z_{k-1}^{\frac{\delta_{k-1} (\Phi_{k-2})}{2\pi\I }}\]
By Proposition \ref{strongerpro}, we have
\begin{equation}\label{k-1W}
    W=\delta_{k}(\Phi_{k-2})+\mathcal{O}(z_{k-1}^{-1})
\end{equation}
as $z_{k-1}\rightarrow +\infty$. Furthermore, it follows from the equation of $\Phi_{k-1}$ that the eigenvalues of $\delta_k(\Phi_{k-1})$ are constant (independent of $z_1,...,z_{k-1}$). In particular, the $n\times n$ Hermitian matrices $W$ and $\delta_{k}(\Phi_{k-2})$ have same eigenvalues. Therefore, for any fixed $z_1,...,z_{k-2}$ we can find a family of $n\times n$ unitary matrices $P(z_{k-1})$ such that 
\[W=P\delta_{k}(\Phi_{k-2})P^{-1}, \ \text{ and } \ \lVert P-{\rm Id}_{n}\rVert=\mathcal{O}(z_{k-1}^{-1}).\]
For example, just like Remark \ref{rmkrho}, the equation of $W$ respect to the variable $z_k$ is equivalent to a differential equation of $P(z_k)$. Then the required $P(z_k)$ can be determined as a solution of the corresponding equation with the prescribed asymptotics.

Thus, \[z_k^{\frac{-W}{2\pi\I }}=Pz_k^{\frac{-\delta_{k}(\Phi_{k-2})}{2\pi\I }} P^{-1}= z_k^{\frac{-\delta_{k}(\Phi_{k-2})}{2\pi\I }}+\mathcal{O}(z_{k-1}^{-1}),\]
which gives rise to
\begin{equation}\label{kk-1}
   z_k^{\frac{-W}{2\pi\I }} z_{k-1}^{\frac{-\delta_{k-1} (\Phi_{k-2})}{2\pi\I }}= z_{k}^{\frac{-\delta_{k} (\Phi_{k-2})}{2\pi\I }}  z_{k-1}^{\frac{-\delta_{k-1} (\Phi_{k-2})}{2\pi\I }}+\mathcal{O}(z_{k-1}^{-1}).
\end{equation}
The identities \eqref{Wdefi} and \eqref{kk-1} give
\begin{align} \nonumber
     &z_{k-1}^{\frac{-\delta_{k-1} (\Phi_{k-2})}{2\pi\I }}  z_k^{\frac{-\delta_{k} (\Phi_{k-1})}{2\pi\I }} (\Phi_k)z_k^{\frac{\delta_{k} (\Phi_{k-1})}{2\pi\I }} z_{k-1}^{\frac{\delta_{k-1} (\Phi_{k-2})}{2\pi\I }}  
      \\ \label{k-1k-2}
      =&
z_{k}^{\frac{-\delta_{k} (\Phi_{k-2})}{2\pi\I }}  z_{k-1}^{\frac{-\delta_{k-1} (\Phi_{k-2})}{2\pi\I }} (\Phi_k)z_{k-1}^{\frac{\delta_{k-1} (\Phi_{k-2})}{2\pi\I }} z_{k}^{\frac{\delta_{k} (\Phi_{k-2})}{2\pi\I }}
+\mathcal{O}(z_{k-1}^{-1})+ \mathcal{O}(z_k^{-1}).
\end{align}
Combining \eqref{k-2k-1} and \eqref{k-1k-2} gives rise to
\[z_{k}^{\frac{-\delta_{k} (\Phi_{k-2})}{2\pi\I }}  z_{k-1}^{\frac{-\delta_{k-1} (\Phi_{k-2})}{2\pi\I }} (\Phi_k)z_{k-1}^{\frac{\delta_{k-1} (\Phi_{k-2})}{2\pi\I }} z_{k}^{\frac{\delta_{k} (\Phi_{k-2})}{2\pi\I }}=\Phi_{k-2}+\mathcal{O}(z_{k-1}^{-1})+\mathcal{O}(z_k^{-1}). \]
Continuing the process step by step, we get (recall that $\Phi=\Phi_{n-1}(z_1,...,z_{n-1})$)
\begin{align}\nonumber
&\left(\overleftarrow{\underset{k=1,...,n-1}{\prod}} z_{k}^{\frac{-\delta_{k} (\Phi_{0})}{2\pi\I }} \right)\cdot \Phi \cdot \left(\overleftarrow{\underset{k=1,...,n-1}{\prod} }z_{k}^{\frac{-\delta_{k} (\Phi_{0})}{2\pi\I }}\right)^{-1}\\ \nonumber
=&z_{1}^{\frac{-\delta_{1} (\Phi_{0})}{2\pi\I }}\cdot \Phi_{1}(z_1)\cdot z_{1}^{\frac{\delta_{1} (\Phi_{0})}{2\pi\I }}+\sum_{k=2}^{n-1}\mathcal{O}(z_k^{-1})\\ \label{11stronger}
=&\Phi_0+\sum_{k=2}^{n-1}\mathcal{O}(z_k^{-1})
. 
\end{align}
Here the last identity uses \eqref{k=1}, and the product $\overleftarrow{\prod}$ is taken with the index $i$ to the right of $j$ if $i>j$. Note that the identity \eqref{11stronger} is stronger version of \eqref{firstasy}, i.e., it implies that \eqref{firstasy} holds for all $z_1>0$, not just in the limit $z_1\rightarrow 0$. (We add the redundant condition $z_1\rightarrow 0$ in Theorem \ref{isomonopro}, simply because the infinite point $u=(u_1,...,u_n)$ with $z_1\rightarrow 0$ and $z_k\rightarrow+\infty$ for $k=2,...,n-1$ is the caterpillar point $u_{\rm cat}$.) It proves the first part of Theorem \ref{isomonopro}.

The existence of a real analytic solution $\Phi(u)$ of \eqref{isoeq} with the boundary value $\Phi_0$ follows from Proposition \ref{inversesol}.
\end{proof}


\begin{rmk}
In a follow-up work \cite{TangXu}, we find the convergent series expansion of $\Phi(u)$, not just the leading expansion, in terms of $\Phi_0$, as $z_k\rightarrow \infty$ for $k=2,...,n-1$.
\end{rmk}

\begin{rmk}\label{rmkrho}
There exists a $U(n)$ valued function $\rho(z_1,...,z_{n-1})$ such that $\Phi(u)$ is related to its boundary value $\Phi_0$ by the conjuation of $\rho$, that is $\Phi(z_1,...,z_{n-1})=\rho \Phi_0 \rho^{-1}$. Actually, following from the equation and the asymptotics of $\Phi(u)$, we see that $\rho$ is the solution of the equation
\[
\frac{\partial\rho}{\partial z_{j}}=-\frac{1}{2\pi\I }{\rm ad}_{u}^{-1}{\rm ad}_{\frac{\partial u}{\partial z_j}}\left(\rho\Phi_0\rho^{-1}\right), \ \text{for} \ j=1,...,n-1,\]
with the prescribed asymptotics
\begin{equation*}
\left(\overleftarrow{\underset{j=1,...,n-1}{\prod}} z_{j}^{\frac{-\delta_{j} (\Phi^{(11)}_{0})}{2\pi\I }} \right)\cdot\rho(z_1,...,z_{n-1})={\rm Id}_{n}+\sum_{j=2}^{n-1}\mathcal{O}(z_j^{-1}).
\end{equation*}
\end{rmk}
\begin{rmk}\label{phikblock}
It follows from the equation \eqref{kdiffeq} that if we write the $n\times n$ matrix $\Phi_k$ as a block form 
\[\Phi_k(z_1,...,z_k)=\left(\begin{array}{cc}
    \Phi_k^{(11)} & \Phi_k^{(12)}  \\
    \Phi_k^{(21)} & \Phi_k^{(22)}
  \end{array}\right),\]
then the $(k+1)\times (k+1)$ block $\Phi_k^{(11)}(z_1,...,z_k)$ is a solution of the isomonodromy equation \eqref{isoeq} of rank $k+1$ (i.e., replacing $n$ by $k+1$ in \eqref{isoeq}). Let us explain how the whole $\Phi_k$ can be determined from the upper left submatrix $\Phi_k^{(11)}$. By the equation \eqref{kdiffeq}, the $(k+1)\times (n-k+1)$ block $\Phi_k^{(12)}$ and the $(n-k+1)\times (n-k+1)$ block $\Phi_k^{(22)}$ satisfy for $j=1,...,k$,
\begin{align*}
    \frac{\partial\Phi_k^{(12)}}{\partial z_{j}}&=\frac{1}{2\pi\I }\Phi_k^{(12)}\cdot \delta_0(\Phi_k^{(22)})- \frac{1}{2\pi\I }\left(({\rm ad}_{u}^{-1}{\rm ad}_{\frac{\partial u}{\partial z_j}}-z_j^{-1})\Phi_k^{(11)}\right)\cdot \Phi_k^{(12)} ,\\
     \frac{\partial\Phi_k^{(21)}}{\partial z_{j}}&=-\frac{1}{2\pi\I }\delta_0(\Phi_k^{(22)})\cdot \Phi_k^{(21)}+\frac{1}{2\pi\I }\Phi_k^{(21)}\cdot ({\rm ad}_{u}^{-1}{\rm ad}_{\frac{\partial u}{\partial z_j}}-z_j^{-1})\Phi_k^{(11)},\\
      \frac{\partial\Phi_k^{(22)}}{\partial z_{j}}&=\frac{1}{2\pi\I }\delta_0(\Phi_k^{(22)})\cdot \Phi_k^{(22)}-\frac{1}{2\pi\I }\Phi_k^{(22)}\cdot \delta_0(\Phi_k^{(22)}).
\end{align*}
Here recall $\delta_0(\Phi_k^{(22)})$ is the diagonal part of the matrix $\Phi_k^{(22)}$. Therefore, knowing the block $\Phi_k^{(11)}$ and the boundary value $\Phi_0$, the other blocks can be expressed in a simple form. In particular, there exists a $U(k+1)$ valued function $\rho_k(z_1,...,z_k)$ such that the upper left block 
\begin{equation*}
\Phi^{(11)}_k=\rho_k \cdot \Phi^{(11)}_0\cdot \rho_{k}^{-1}, 
\end{equation*}
and the other blocks are then determined simply by
\begin{equation*}
\Phi_k=\left(\begin{array}{cc}
    \rho_k & 0  \\
    0 & (z_1\cdots z_k)^{\frac{\delta_0(\Phi_0^{(22)})}{2\pi\I }}
  \end{array}\right) \cdot \Phi_0\cdot \left(\begin{array}{cc}
    \rho_k & 0  \\
    0 & (z_1\cdots z_k)^{\frac{\delta_0(\Phi_0^{(22)})}{2\pi\I }}\end{array}\right)^{-1}, 
\end{equation*}
where $\delta_0(\Phi_0^{(22)})$ is the $(n-k+1)\times (n-k+1)$ diagonal matrix with diagonal entries $(\Phi_0)_{k+2,k+2},...,(\Phi_0)_{n,n}$. It follows from the asymptotics of $\Phi_k$ that the unitary matrix $\rho_k$ satisfies 
\begin{equation*}
\left(\overleftarrow{\underset{j=1,...,k}{\prod}} z_{j}^{\frac{-\delta_{j} (\Phi^{(11)}_{0})}{2\pi\I }} \right)\cdot\rho_k(z_1,...,z_k)={\rm Id}_{k+1,k+1}+\mathcal{O}(z_j^{-1};j=2,...,k).
\end{equation*}
\end{rmk}

A similar proof as above shows that
\begin{pro}\label{reverseasym}
For any solution $\Phi(u)$ of the isomonodromy equation \eqref{introisoeq} on the connected component $U_{\rm id}$, there exists a unique constant $\Phi_0\in\Herm(n)$ such that as $z_k\rightarrow+\infty$ for all $k=2,...,n-1$,
\begin{align}\label{reverse}
&\left(\overrightarrow{\underset{k=1,...,n-1}{\prod}} z_{k}^{\frac{-\delta_{k} (\Phi)}{2\pi\I }} \right)\cdot \Phi \cdot \left(\overrightarrow{\underset{k=1,...,n-1}{\prod} }z_{k}^{\frac{-\delta_{k} (\Phi)}{2\pi\I }}\right)^{-1}=\Phi_0+\sum_{k=2}^{n-1}\mathcal{O}(z_k^{-1})
. 
\end{align}
\end{pro}
We remark that the ordered product in \eqref{reverse} is reverse to the one in \eqref{11stronger}, and Proposition \ref{reverseasym} will be used in the proof of Theorem \ref{introcor} in Section \ref{pf:introthm3}.
\begin{defi}\label{solasy}
We call $\Phi_0$ the boundary value of the solution $\Phi(u)$ of isomonodromy equation \eqref{isoeq} at $u_{\rm cat}$. And for any $\Phi_0\in\Herm(n)$, we denote by $\Phi(u;\Phi_0)$ the solution of \eqref{isoeq} with the boundary value $\Phi_0$.
\end{defi}
The meaning of the infinite point will become clear in Section \ref{closureStokes}. In the rest of this section, we will derive an explicit formula of the Stokes matrices $S_\pm(u,\Phi(u;\Phi_0))$ via the boundary value $\Phi_0$, i.e., Theorem \ref{mainthm}.

\subsection{A branching rule of the system \eqref{Stokeseq} via isomonodromy deformation: one recursive step}\label{RHtorelRH}
It is more convenient to work with the following coordinates 
\begin{align*}
\widetilde{z_0}=u_1, \ z_1=u_2-u_1, \ z_2=\frac{u_3-u_2}{u_2-u_1}, \ ...... , \ z_{n-1}=\frac{u_n-u_{n-1}}{u_{n-1}-u_{n-2}}.
\end{align*}

In this section, we will show that as $z_{n-1}\rightarrow +\infty$, the boundary value $\Phi_{n-2}(z_1,...,z_{n-1})$, in the sense of Proposition \ref{strongerpro}, of a given solution $\Phi_{n-1}=\Phi(z_1,...,z_{n-1})$ of the isomonodromy equation can be used to "decouple" the system \eqref{Stokeseq} of rank $n$ into two lower rank systems, with rank $n-1$ and rank $1$ (therefore trivial) respectively. To fix some notations, for any $u\in U_{\rm id}$ and $A\in\Herm(n)$, we denote by
\begin{itemize}
    \item $S_\pm(u,A)$ the Stokes matrices, and $C(u,A)$ the connection matrix of $\frac{dF}{dz}=(\I u-\frac{1}{2\pi\I }\frac{A}{z})F$;
    \item $C(E_n, A)\subset{\rm U}(n)$, the connection matrix of $\frac{dF}{dz}=(\I E_n-\frac{1}{2\pi\I }\frac{A}{z})F$, with $E_n={\rm diag}(0,...,0,1)$;
    \item $C(u^{(n-1)},\delta_{n-1}(A))\in{\rm U}(n)$, the connection matrix of the $n\times n$ system $\frac{dF}{dz}=(u^{(n-1)}-\frac{1}{2\pi\I }\frac{\delta_{n-1}(A)}{z})F$. Here $u^{(n-1)}:={\rm diag}(u_1,u_2,...,u_{n-1},0)$.
\end{itemize}

\begin{pro}\label{inductionstep}
Let $\Phi_{n-1}(z_{n-1};\Phi_{n-2})\in\Herm(n)$ denote the solution of the isomonodromy equation \eqref{isoeq} with the boundary value $\Phi_{n-2}(z_1,...,z_{n-2})$ at $z_{n-1}=\infty$ in the sense of Proposition \ref{strongerpro}. Then we have the identity
\begin{align}\nonumber
&C\left(u,\Phi_{n-1}(z_{n-1};\Phi_{n-2})\right)\cdot e^{\Phi_{n-1}}\cdot C\left(u,\Phi_{n-1}(z_{n-1};\Phi_{n-2})\right)^{-1}\\ \nonumber
=&S_-\left(u,\Phi_{n-1}(z_{n-1};\Phi_{n-2})\right) S_+\left(u,\Phi_{n-1}(z_{n-1};\Phi_{n-2})\right)\\ \label{codim1}
=&{\rm Ad}{\left(C\left(u^{(n-1)}, \delta_{n-1}(\Phi_{n-2})\right)\cdot \left({u_{n-1}-u_{n-2}}\right)^{-\frac{\delta_{n-1}(\Phi_{n-2})}{2\pi\I }}\cdot C\left(E_n, \Phi_{n-2}\right)\right)}e^{\Phi_{n-2}}.
\end{align}
\end{pro}
\pf The first identity simply follows from the monodromy relation \eqref{monorelation}. We prove the second identity by showing that the left and right hand sides compute respectively the monodromy of the linear system of equation \eqref{introisoeq} along two homotopy paths. Since the paths are homotopy, the monodromy are equal.

Let us assume $n>2$ (for the $n=2$ case the involved system can be solved exactly). Let us fix the $n-1$ variables $\widetilde{z_0}=u_1$ and $z_1,...,z_{n-2}$. Thus $\Phi_{n-2}(z_1,...,z_{n-2})\in\Herm(n)$ is constant, and $u_{n-1}-u_{n-2}=z_1z_2\cdots z_{n-2}$ is a constant real number. The compatible linear system of PDEs \eqref{introisoStokeseq1}-\eqref{introisoStokeseq2} reduces to the equation with respect to $z$ and $z_{n-1}$ 
\begin{align}\label{relequation1}
\frac{\partial F}{\partial z}&=\left(\I u-\frac{1}{2\pi\I }\frac{\Phi_{n-1}(z_{n-1};\Phi_{n-2})}{z}\right)F,\\
\label{relequation2}\frac{\partial F}{\partial z_{n-1}}&=\left(\I (u_{n-1}-u_{n-2})E_nz-\frac{1}{2\pi\I } {\rm ad}^{-1}_u{\rm ad}_{(u_{n-1}-u_{n-2})E_n}\Phi_{n-1}(z_{n-1};\Phi_{n-2})\right)F.
\end{align}
Let $F_\pm(z,z_{n-1})$ be the canonical fundamental solutions of the first equation with the prescribed asymptotics
\begin{align*}
\lim_{z\rightarrow\infty}F_+(z,z_{n-1})\cdot e^{-\I uz}\cdot z^{\frac{[\Phi_{n-1}]}{2\pi\I }}&={\rm Id}_n, \ \ \ \text{as} \ \ \ z\in S(-\pi,\pi),
\\
\lim_{z\rightarrow\infty}F_-(z,z_{n-1})\cdot e^{-\I uz}\cdot z^{\frac{[\Phi_{n-1}]}{2\pi\I }}&={\rm Id}_n, \ \ \ \text{as} \ \ \ z\in S(-2\pi,0).
\end{align*}
Then following Proposition \ref{isomonodef}, the canonical solutions $F_\pm(z,z_{n-1})$ of the first equation also satisfy the second equation, and therefore are solutions of the compatible system \eqref{relequation1}-\eqref{relequation2}.

By definition, for any fixed $z_{n-1}$, the monodromy of $F_+(z,z_{n-1})$, along a loop $\gamma_1(t)=(\frac{e^{-i\theta}}{\varepsilon},z_{n-1})$ around $(\infty,z_{n-1})$, is just the first row in \eqref{codim1}. 
Another loop homotopy to $\gamma_1$ is $\gamma_2\circ \gamma_3\circ \gamma_2^{-1}$, where $\gamma_2$ is a simple path from $(\frac{1}{\varepsilon},z_{n-1})$ to $(\varepsilon, \frac{1}{\varepsilon^2})$, and $\gamma_3(\theta)=(e^{-i\theta}\varepsilon,\frac{1}{\varepsilon^2})$ a loop around $(0, \frac{1}{\varepsilon^2})$. To show the monodromy along $\gamma_2\circ \gamma_3\circ \gamma_2^{-1}$ coincide with the last row in \eqref{codim1}, we need to compute respectively the monodromy along the loops $\gamma_2$ and $\gamma_3$. This can be done by choosing certain reference solution $Y_+(z,z_{n-1})$ around $(\varepsilon, \frac{1}{\varepsilon^2})$.

To introduce $Y_+(z,z_{n-1})$, first note that under the change of coordinates \[x=z \ \ \text{  and  } \ \ y=(u_{n-1}-u_{n-2})\cdot zz_{n-1},\] 
(here by assumption $u_{n-1}-u_{n-2}=z_1z_2\cdots z_{n-2}$ is a fixed positive real number) the systems \eqref{relequation1} and \eqref{relequation2} become (for simplicity, we write $\Phi_{n-1}$ for $\Phi_{n-1}(z_{n-1};\Phi_{n-2})$)
\begin{eqnarray}\label{changeequation1}
&&\frac{\partial F}{\partial x}=\left(\I u^{(n-1)}-\frac{1}{2\pi\I }\frac{\delta_{n-1}(\Phi_{n-1})}{x}-\frac{1}{2\pi\I }{\rm ad}_{D_1}{\rm ad}_{E_n} \Phi_{n-1}\right)F,\\
\label{changeequation2}&&\frac{\partial F}{\partial y}=\left(\I E_n-\frac{1}{2\pi\I }{\rm ad}_{D_2}{\rm ad}_{E_n} \Phi_{n-1}\right)F,
\end{eqnarray}
where \begin{align*}
    D_1&={\rm diag}\left(\frac{u_1}{y+x(u_{n-1}-u_1)},\frac{u_2}{y+x(u_{n-1}-u_2)},...,\frac{u_{n-1}}{y},0\right),\\D_2&={\rm diag}\left(\frac{1}{y+x(u_{n-1}-u_1)},\frac{1}{y+x(u_{n-1}-u_2)},...,\frac{1}{y},0\right)
\end{align*}are two $n\times n$ diagonal matrices. Here recall that $\delta_{n-1}(\Phi_{n-1})$ is defined in \eqref{delta1}.

Then for any fixed $y\ne 0$,
\begin{itemize}
\item let $G_0(x,y)$ be the solution of equation \eqref{changeequation1} with the asymptotics $G_0\cdot x^{\frac{\delta_{n-1}(\Phi_{n-2})}{2\pi\I }}\sim 1$ at $x=0$. The existence of such a solution can be seen by the estimation of the coefficient matrix of equation \eqref{changeequation1} at $x=0$: recall that $z_{n-1}=\frac{y}{(u_{n-1}-u_{n-2})x}$, by the proof of Proposition \ref{strongerpro}, we have that \[\delta_{n-1}(\Phi_{n-1})-\delta_{n-1}(\Phi_{n-2})=\mathcal{O}(z_{n-1}^{-1})=\mathcal{O}(x),\] and that $\frac{1}{2\pi\I }{\rm ad}_{D_1}{\rm ad}_{E_n}\Phi_{n-1}$ has a limit at $x=0$. 

\item Let $F^{(n)}_\pm(y)$ be the canonical solution of $\frac{dF}{dy}=\left(\I E_n-\frac{1}{2\pi\I }\frac{\Phi_{n-2}}{y}\right)F$ in the two Stokes sectors (right and left half planes). 

\end{itemize}

\begin{lem}\label{1satisfies2}
The functions $Y_\pm(x,y):=G_0(x,y)\cdot (\frac{y}{u_{n-1}-u_{n-2}})^{\frac{\delta_{n-1}(\Phi_{n-2})}{2\pi\I }}\cdot F^{(n)}_\pm(y)$ satisfy the equations \eqref{changeequation1} and \eqref{changeequation2}.
\end{lem}
\pf We denote the coefficients of \eqref{changeequation1} and \eqref{changeequation2} by $A_1$ and $A_2$ respectively. By the compatibility of the two equations, we have that $(d_x-A_1)(d_yG_0-A_2G_0)=0$. 
Since $(d_x-A_1)G_0=0$, we can set $d_yG_0-A_2G_0=G_0X(y)$ for a function $X(y)$ of $y$. That is $X(y)=G_0^{-1}d_yG_0-G_0^{-1}A_2G_0$. To get the expression of $X(y)$, we use the asymptotics of $G_0$ and $A_2$ to compute the asymptotics of $G_0^{-1}d_yG_0-G_0^{-1}A_2G_0$ as $x\rightarrow 0$ (while fixing $y$).

Firstly, it follows from the asymptotics of $\Phi_{n-1}$ that the asymptotics of $G_0^{-1}A_2G_0$ as $x\rightarrow 0$ along the real axis is
\begin{align*}
&G_0^{-1}\cdot\left(\I E_n-\frac{1}{2\pi\I }{\rm ad}_{D_2}{\rm ad}_{E_n}\Phi_{n-1}\right)\cdot G_0\\
\sim& \ x^{\frac{\delta_{n-1}(\Phi_{n-2})}{2\pi\I }}\cdot\left(\I E_n-\frac{1}{2\pi\I }{\rm ad}_{D_2}{\rm ad}_{E_n}\Phi_{n-1}\right)\cdot x^{-\frac{\delta_{n-1}(\Phi_{n-2})}{2\pi\I }} \\
\sim& \ (\frac{y}{u_{n-1}-u_{n-2}})^{\frac{\delta_{n-1}(\Phi_{n-2})}{2\pi\I }}\cdot\left(\I E_n-\frac{1}{2\pi\I }\frac{\Phi_{n-2}-\delta_{n-1}(\Phi_{n-2})}{y}\right)\cdot (\frac{y}{u_{n-1}-u_{n-2}})^{\frac{-\delta_{n-1}(\Phi_{n-2})}{2\pi\I }}.
\end{align*}
Here the norms $\lVert x^{\frac{\delta_{n-1}(\Phi_{n-2})}{2\pi\I }}\rVert=1$, $\lVert (\frac{y}{u_{n-1}-u_{n-2}})^{\frac{\delta_{n-1}(\Phi_{n-2})}{2\pi\I }}\rVert=1$ and in the last step, we use the asymptotics of $\Phi_{n-1}$ and $D_2$ as $x\rightarrow 0$ (recall $y=({u_{n-1}-u_{n-2}})xz_{n-1}$) 
\begin{eqnarray*} {z_{n-1}}^{\frac{-\delta_{n-1}(\Phi_{n-2})}{2\pi\I }}\cdot  \Phi_{n-1}\cdot {z_{n-1}}^{\frac{\delta_{n-1}(\Phi_{n-2})}{2\pi\I }} \sim  \Phi_{n-2}, \hspace{5mm} \delta_{n-1}(\Phi_{n-1})\sim \delta_{n-1}(\Phi_{n-2}), \hspace{5mm}
D_2\sim {\rm diag}\left(\frac{1}{y},...,\frac{1}{y},0\right),\end{eqnarray*}
to get
\begin{align*}
&    x^{\frac{\delta_{n-1}(\Phi_{n-2})}{2\pi\I }}\cdot({\rm ad}_{D_2}{\rm ad}_{E_n}\Phi_{n-1})\cdot x^{-\frac{\delta_{n-1}(\Phi_{n-2})}{2\pi\I }}\\ \sim & \ (\frac{y}{u_{n-1}-u_{n-2}})^{\frac{\delta_{n-1}(\Phi_{n-2})}{2\pi\I }}\cdot \frac{\Phi_{n-2}-\delta_{n-1}(\Phi_{n-2})}{y}\cdot (\frac{y}{u_{n-1}-u_{n-2}})^{-\frac{\delta_{n-1}(\Phi_{n-2})}{2\pi\I }}.
\end{align*}

Secondly, 
the limit of $G_0^{-1}d_yG_0$ is zero as $x\rightarrow 0$.
Therefore, we have computed the limit of  $G_0^{-1}d_yG_0-G_0^{-1}A_2G_0$ as $x\rightarrow 0$, and deduce that $d_yG_0-A_2G_0=G_0\cdot X(y)$, for \[X(y)=(\frac{y}{u_{n-1}-u_{n-2}})^{\frac{\delta_{n-1}(\Phi_{n-1})}{2\pi\I }}\cdot \left(\I E_n-\frac{1}{2\pi\I }\frac{\Phi_{n-2}-\delta_{n-1}(\Phi_{n-2})}{y}\right)\cdot (\frac{y}{u_{n-1}-u_{n-2}})^{-\frac{\delta_{n-1}(\Phi_{n-1})}{2\pi\I }}.\]

By the explicit formula of $X(y)$ and the defining equation of $F_\pm^{(n)}$, one verifies directly that
\[d_y\left((\frac{y}{u_{n-1}-u_{n-2}})^{\frac{\delta_{n-1}(\Phi_{n-2})}{2\pi\I }} F^{(n)}_\pm\right) =X(y)\cdot (\frac{y}{u_{n-1}-u_{n-2}})^{\frac{\delta_{n-1}(\Phi_{n-1})}{2\pi\I }} F^{(n)}_\pm.\]
Therefore $Y_\pm=G_0\cdot (\frac{y}{u_{n-1}-u_{n-2}})^{\frac{\delta_{n-1}(\Phi_{n-2})}{2\pi\I }} F^{(n)}_\pm$ satisfy the equation $d_yY_\pm=A_2Y_\pm$, i.e., the equation \eqref{changeequation2}.
\qed

\vspace{2mm}
Under the coordinates transformation, the function $Y_+(x,y)$ is viewed as solutions of \eqref{relequation1} and \eqref{relequation2} with prescribed asymptotics at $(z,z_{n-1})=(\varepsilon,\frac{1}{\varepsilon^2})$.

\begin{lem}
The monodromy of $Y_+$ along the loop $\gamma_3$ coincides with the monodromy of $F^{(n)}_+$ around $y=\infty$, i.e., the monodromy $C(E_n,\Phi_{n-2})e^{\Phi_{n-2}}C(E_n,\Phi_{n-2})^{-1}$ of $\frac{dF}{dy}=\left(\I E_n-\frac{1}{2\pi\I }\frac{\Phi_{n-2}}{y}\right)F$ around $y=\infty$.
\end{lem}
\pf On the one hand, in the $(z,z_{n-1})$ coordinates, $\gamma_3(\theta)=(z(\theta),z_{n-1}(\theta))=(e^{-i\theta}\varepsilon,\frac{1}{\varepsilon^2})$ is the loop around $(0, \frac{1}{\varepsilon^2})$. On the other hand, we have the coordinates transformation $y=(u_{n-1}-u_{n-2})zz_{n-1}$ and $x=z$. Thus in $(x,y)$ coordinates, $\gamma_3$ is a loop $\gamma_3(\theta)=(x(\theta),y(\theta))=(e^{-i\theta}\varepsilon,\frac{e^{-i\theta}}{\varepsilon})$, that is a loop around $x=0$ and $y=\infty$. Analytic continuation along such a loop gives
\begin{equation*}
G_0(e^{2\pi\I }x,e^{2\pi\I }y)=G_0(x,y)e^{-\delta_{n-1}(\Phi_{n-2})},
\end{equation*}
and
\begin{align*}
&(\frac{e^{2\pi\I }y}{u_{n-1}-u_{n-2}})^{\frac{\delta_{n-1}(\Phi_{n-2})}{2\pi\I }} \cdot F^{(n)}_\pm(e^{2\pi\I }x,e^{2\pi\I }y)\\= \ &e^{\delta_{n-1}(\Phi_{n-2})}(\frac{y}{u_{n-1}-u_{n-2}})^{\frac{\delta_{n-1}(\Phi_{n-2})}{2\pi\I }}\cdot F^{(n)}_\pm(x,y) C(E_n,\Phi_{n-2})e^{\Phi_{n-2}}C(E_n,\Phi_{n-2})^{-1}.
\end{align*}
The result follows from the above two identities and the expression $Y_+=G_0\cdot (\frac{y}{u_{n-1}-u_{n-2}})^{\frac{\delta_{n-1}(\Phi_{n-1})}{2\pi\I }} F^{(n)}_+$. \qed

\vspace{2mm}
The rest is to compare the two solutions $Y_+$ and $F_+$ of the same equation \eqref{relequation1} (to get the monodromy along the path $\gamma_2$). For that, we only need to study the asymptotics of $F_+^{-1}Y_+$ as $z_{n-1}\gg z\gg 0$. 

Let us first study the asymptotics of $Y_+$. For any $x=z$, as $z_{n-1}=y/(u_{n-1}-u_{n-2})x \rightarrow+\infty$ the equation \eqref{changeequation1} approaches to (here we use the fact $\delta_{n-1}(\Phi_{n-1})\rightarrow\delta_{n-1}(\Phi_{n-2})$ as $z_{n-1}\rightarrow\infty$ and $\Phi_{n-1}$ is bounded, see Proposition \ref{strongerpro}) \begin{eqnarray}\label{xequation}
\frac{dF}{dx}=\left(\I u^{(n-1)}-\frac{1}{2\pi\I }\frac{\delta_{n-1}(\Phi_{n-2})}{x}\right)F,
\end{eqnarray}
and $G_0(x,y)$ approaches to the fundamental solution $T_0(x)$ of \eqref{xequation} with the asymptotics $T_0x^{\frac{\delta_{n-1}(\Phi_{n-2})}{2\pi\I }}\sim 1$ at $x=0$. Let $T_+$ be the canonical solution of \eqref{xequation} around $x=\infty$ in the Stokes sector ${\rm Sect}_+$ (right half plane), then by the definition of connection matrix, we have \[T_+=T_0\cdot C(u^{(n-1)},\delta_{n-1}(\Phi_{n-2})).\]
Then from the expression $Y_+(x,y)=G_0(x,y)\cdot (\frac{y}{u_{n-1}-u_{n-2}})^{\frac{\delta_{n-1}(\Phi_{n-1})}{2\pi\I }} F^{(n)}_+(y)$, the asymptotics of $F^{(n)}_+(y)$ as $y\rightarrow+\infty$ and the fact that $G_0$ approaches to $T_0$, we have that for any fixed real number $x=z$, as $z_{n-1}\rightarrow +\infty$ or equivalently $y\rightarrow +\infty$,
\begin{align} \nonumber
&\lim_{y\rightarrow \infty}\left(Y_+(x,y) \cdot ({u_{n-1}-u_{n-2}})^{\frac{\delta_{n-1}(\Phi_{n-2})}{2\pi\I }}\cdot C(u^{(n-1)},\delta_{n-1}(\Phi_{n-2}))^{-1} e^{-\I u x}x^{\frac{[\Phi_{n-2}]}{2\pi\I }}\right)\\ \nonumber
=&\lim_{y\rightarrow \infty}\Big(T_0(x)(\frac{y}{u_{n-1}-u_{n-2}})^{\frac{\delta_{n-1}(\Phi_{n-2})}{2\pi\I }} (1+\mathcal{O}(y^{-1}))(\frac{y}{u_{n-1}-u_{n-2}})^{-\frac{\delta_{n-1}(\Phi_{n-2})}{2\pi\I }}\\ \nonumber
& \ \ \ \ \ \ \ \ \times C(u^{(n-1)},\delta_{n-1}(\Phi_{n-2}))^{-1} e^{\I y E_n}e^{-\I u x}x^{\frac{[\Phi_{n-2}]}{2\pi\I }} \Big)\\ \nonumber
=& T_0(x)C(u^{(n-1)},\delta_{n-1}(\Phi_{n-2}))^{-1} e^{-\I x u^{(n-1)}}x^{\frac{[\Phi_{n-2}]}{2\pi\I }} \\ \label{limitYC}
=& T_+(x)e^{-\I x u^{(n-1)}}x^{\frac{[\Phi_{n-2}]}{2\pi\I }}.
\end{align}
Here the first identity uses the asymptotics $F^{(n)}_+(y)=(1+\mathcal{O}(y^{-1}))y^{-\frac{\delta_{n-1}(\Phi_{n-1})}{2\pi\I }}e^{\I y E_n}$, and the second identity uses the fact that the norm $|y^{\frac{\delta_{n-1}(\Phi_{n-1})}{2\pi\I }}|=1$ for $y\in\mathbb{R}$. 
By the prescribed asymptotics of the canonical solution $T_+(x)$ at $x=+\infty$, the right hand side of \eqref{limitYC} approaches to the identity matrix as $x\rightarrow+\infty$. Thus, we get (after changing $(x,y)$ to $(z,z_{n-1})$ and restricting $z$ to be real)
\begin{eqnarray}\label{1asy}
\lim_{z_{n-1}\gg z\gg 0}\left(Y_+(z,u_n)\cdot ({u_{n-1}-u_{n-2}})^{\frac{\delta_{n-1}(\Phi_{n-2})}{2\pi\I }}\cdot C\left(u^{(n-1)},\delta_{n-1}(\Phi_{n-2})\right)^{-1} e^{-\I uz}z^{\frac{[\Phi_{n-1}]}{2\pi\I }}\right)={\rm Id}_n.
\end{eqnarray}
Here we use that the diagonal part $[\Phi_{n-1}]$ of $\Phi_{n-1}$ equals to the diagonal part $[\Phi_{n-2}]$ of $\Phi_{n-2}$.

Let us now study the asymptotics of $F_+$. By Theorem \ref{uniformresum} we have that (while fixing $u_1,...,u_{n-1}$) \begin{eqnarray}\label{2asy}
\lim_{z\rightarrow+\infty}\left(F_+(z,z_{n-1}) e^{-\I uz}z^{\frac{[\Phi_{n-1}]}{2\pi\I }}\right)={\rm Id}_n,  \hspace{3mm} \text{uniformly for } u_n>u_{n-1}+1.
\end{eqnarray} 
Therefore, the limits \eqref{1asy} and \eqref{2asy} imply that as $z_{n-1}\gg z\gg 0$ (for $z\in \mathbb{R}$)
\begin{eqnarray}\label{realest}e^{\I uz}z^{-\frac{[\Phi_{n-1}]}{2\pi\I }}F_+^{-1}Y_+ \cdot ({u_{n-1}-u_{n-2}})^{\frac{\delta_{n-1}(\Phi_{n-2})}{2\pi\I }}\cdot C\left(u^{(n-1)},\delta_{n-1}(\Phi_{n-1})\right)^{-1} e^{-\I uz}z^{\frac{[\Phi_{n}]}{2\pi\I }}\rightarrow {\rm Id}_n.
\end{eqnarray}
Since $e^{\I uz}z^{-\frac{[\Phi_{n}]}{2\pi\I }}$ is a diagonal matrix whose all diagonal elements have norm $1$ for all $z>0$, the limit \eqref{realest} gives rise to \[Y_+=F_+\cdot C(u^{(n-1)},\delta_{n-1}(\Phi_{n-2}))\cdot ({u_{n-1}-u_{n-2}})^{-\frac{\delta_{n-1}(\Phi_{n-2})}{2\pi\I }}.\]
That computes the monodromy along $\gamma_2$. 

In summary, the monodromy $M(\gamma_i)$ along $\gamma_i$ with respect to the reference solutions $F_+$ and $Y_+$ are
\begin{align*}
M\left({\gamma_1};F_+\mapsto F_+(ze^{-2\pi\I })\right) \colon & \  S_-\left(u, \Phi_{n-1}(u_n;\Phi_{n-1})\right) S_+\left(u, \Phi_{n-1}(u_n;\Phi_{n-2})\right);\\
M\left({\gamma_2};F_+\mapsto Y_+\right)\colon & \ C\left(u^{(n-1)},\delta_{n-1}(\Phi_{n-1})\right)\cdot ({u_{n-1}-u_{n-2}})^{-\frac{\delta_{n-1}(\Phi_{n-2})}{2\pi\I }}; \\
M\left({\gamma_3}; Y_+(x,y)\mapsto Y_+(e^{2\pi\I }x,e^{2\pi\I }y)\right)\colon &C(E_n,\Phi_{n-1})e^{\Phi_{n-1}}C(E_n,\Phi_{n-1})^{-1}.
\end{align*}
The proposition now follows from that the two loops $\gamma_1$ and $\gamma_2\circ \gamma_3\circ \gamma_2^{-1}$ are homotopy.
\qed

\vspace{2mm}
Similarly, we can prove a generalization of Proposition \ref{inductionstep}.
\begin{pro}\label{morestep}
Let $\Phi_{n-1}(z_1,...,z_{n-1})\in\Herm(n)$ be the solution of the isomonodromy equation \eqref{isoeq}, and $\Phi_k(z_1,...,z_k)$ for $k=1,...,n-1$ the corresponding chain of functions given as in Proposition \ref{strongerpro}. Then for each $k$ we have the identity   
\begin{align} \nonumber
   & C\left(u^{(k+1)},\delta_{k+1}(\Phi_{k})\right)\cdot \left({u_{k+1}-u_{k}}\right)^{\frac{\delta_{k+1}(\Phi_{k})}{2\pi\I }}\cdot  e^{\Phi_{k}}\cdot \left({u_{k+1}-u_{k}}\right)^{-\frac{\delta_{k+1}(\Phi_{k})}{2\pi\I }} C\left(u^{(k+1)},\delta_{k+1}(\Phi_{k})\right)^{-1}\\ \label{morecodim1}
=&{\rm Ad}{\left(C\left(u^{(k)}, \delta_{k}(\Phi_{k-1})\right)\cdot \left({u_{k}-u_{k-1}}\right)^{\frac{\delta_{k}(\Phi_{k-1})}{2\pi\I }}\cdot C\left(E_{k+1}, \delta_{k+1}(\Phi_{k-1})\right)\right)}e^{\Phi_{k-1}}.
\end{align}
Here recall that
\begin{itemize}
    \item $C(E_{k+1}, \Phi_{k-1})\subset{\rm U}(n)$, the connection matrix of $\frac{dF}{dz}=(\I E_{k+1}-\frac{1}{2\pi\I }\frac{\Phi_{k-1}}{z})F$;
    \item $C(u^{(k)},\delta_{k}(\Phi_{k-1}))\in{\rm U}(n)$, the connection matrix of the $n\times n$ system $\frac{dF}{dz}=(u^{(k)}-\frac{1}{2\pi\I }\frac{\delta_k(\Phi_{k-1})}{z})F$ with $u^{(k)}:={\rm diag}(u_1,...,u_{k},0,...0)$.
\end{itemize}
\end{pro}
The proof is omitted due to the similarity with Proposition \ref{inductionstep}. Actually Proposition \ref{morestep} is simply reduced to the case of Proposition \ref{inductionstep}, in aware of Remark \ref{phikblock}.

\subsection{Recursive branching rules and expression of Stokes matrices via the boundary value}\label{exviaasy}
In this subsection, let us consider the recursive branching of the system \eqref{Stokeseq} as $\frac{u_{k+1}-u_{k}}{u_{k}-u_{k-1}}\rightarrow+\infty$ (i.e., $z_k\rightarrow +\infty$) for all $k=2,...,n-1$. 
The recursive procedure in Section \ref{RHtorelRH} implies that for any $A\in{\Herm}(n)$ and $1\le k\le n$, we need to consider the system
\begin{eqnarray}\label{eq:relativeStokes}
\frac{dF}{dz}=\left(\I E_k-\frac{1}{2\pi \I }\frac{\delta_k(A)}{z}\right)F,
\end{eqnarray}
where $E_k={\rm diag}(0,...,0,1,0,..,0)$ with $1$ at the $k$-th position. 
For any $1\le k\le n$, we take the obvious inclusion of $U(k)$ as the upper left corner of $U(n)$, extended
by $1$'s along the diagonal. 
Since $U(k-1)\in U(n)$ is in the centralizer of the irregular term $\I E_k$ of the equation \eqref{eq:relativeStokes}, the connection matrix $C(E_k,\delta_k(A))$ of \eqref{eq:relativeStokes} has the following ${U}(k-1)$-equivariance.
\begin{lem}\label{commute}
For any $g\in {U}(k-1)\subset {U}(n)$, 
\[C\left(E_k, g\delta_k(A)g^{-1}\right)=g\cdot C\left(E_k,\delta_k(A)\right)\cdot g^{-1}.\]
\end{lem}

\begin{thm}\label{relautomorphism}
For any solution $\Phi(u;\Phi_0)$ of the isomonodromy equation \eqref{isoeq} on $u\in U_{\rm id}$ with the boundary value $\Phi_0\in\Herm(n)$ at $u_{\rm cat}$ (as in Theorem \ref{isomonopro}), we have
\[S_{-}(u, \Phi(u;\Phi_0)) S_{+}(u, \Phi(u;\Phi_0))=\left(\overrightarrow{\prod_{k=2,...,n}}C\left(E_k,\delta_k(\Phi_0)\right)\right)\cdot e^{\Phi_0}\cdot \left(\overrightarrow{\prod_{k=2,...,n}}C\left(E_k,\delta_k(\Phi_0)\right)\right)^{-1},\] 
where the product $\overrightarrow{\prod}$ is taken with the index $i$ to the right of $j$ if $j<i$.
\end{thm} 
\pf It can be proved by recursively using Proposition \ref{inductionstep} and Proposition \ref{morestep}.
For any $u\in\h_{\rm reg}(\mathbb{R})$ and $A\in\Herm(n)$, let us denote by $R_1(u^{(1)},A)=A$, and for each $k=2,...,n$,
\begin{eqnarray*}
R_{k}(u^{(k)},A):=C(u^{(k)},\delta_{k}(A))\cdot \left({u_{k}-u_{k-1}}\right)^{\frac{\delta_{k}(A)}{2\pi\I }} A\cdot\left({u_{k}-u_{k-1}}\right)^{-\frac{\delta_{k}(A)}{2\pi\I }}\cdot C(u^{(k)},\delta_{k}(A))^{-1},
\end{eqnarray*}
where $u^{(k)}={\rm diag}(u_1,u_2,..,u_k,0,...,0)$.
Then Proposition \ref{morestep} states that
\begin{eqnarray}\nonumber
&&R_k(u^{(k)},\Phi_{k-1})\\ \nonumber
&=&{\rm Ad}{\left(C\left(u^{(k-1)},\delta_{k-1}(\Phi_{k-2})\right)\left({u_{k-1}-u_{k-2}}\right)^{\frac{\delta_{k-1}(\Phi_{k-2})}{2\pi\I }}C\left(E_k,\delta_k(\Phi_{k-2})\right)\right)}{\Phi_{k-2}}\\ \label{recurC}
&=&{\rm Ad}\left({C\left(E_k, \delta_k\left(R_{k-1}(u^{(k-1)}, \Phi_{k-2})\right)\right)}\right) R_{k-1}(u^{(k-1)}, \Phi_{k-2}).
\end{eqnarray}
Here in the second identity we use the ${\rm U}(k-1)$-equivariant property in Lemma \ref{commute} (note that $C(u^{(k-1)},\delta_{k-1}(\Phi_{n-2}))\in U(k-1)\subset U(n)$)
\begin{eqnarray*}
&&{\rm Ad}{\left(C\left(u^{(k-1)},\delta_{k-1}(\Phi_{k-2})\right)\left({u_{k-1}-u_{k-2}}\right)^{\frac{\delta_{k-1}(\Phi_{k-2})}{2\pi\I }}\right)} C\left(E_k,\delta_k(\Phi_{k-2})\right)\\&=&C\left(E_k, \  {\rm Ad}{\left(C\left(u^{(k-1)},\delta_{k-1}(\Phi_{k-2})\right)\left({u_{k-1}-u_{k-2}}\right)^{\frac{\delta_{k-1}(\Phi_{k-2})}{2\pi\I }}\right)} \delta_k(\Phi_{k-2})\right).
\end{eqnarray*}
Replacing $k$ by $k-1$ in \eqref{recurC} we get
\begin{eqnarray*}\nonumber
R_{k-1}(u^{(k-1)},\Phi_{k-2})={\rm Ad}\left({C\left(E_{k-1}, \delta_{k-1}\left(R_{k-2}(u^{(k-2)}, \Phi_{k-3})\right)\right)}\right)R_{k-2}(u^{(k-2)}, \Phi_{k-3}).
\end{eqnarray*}
Then replacing $R_{k-1}$ by $R_{k-2}$ in the right hand side of the identity \eqref{recurC} and using $U(k-1)$- equivariant property
\begin{align*}
&{C\left(E_k, \ {\rm Ad}\left({C\left(E_{k-1}, \delta_{k-1}\left(R_{k-2}(u^{(k-2)}, \Phi_{k-3})\right)\right)}\right)\delta_k\left(R_{k-2}(u^{(k-2)}, \Phi_{k-3}) \right)\right)}\\
=&{{\rm Ad}\left({C\left(E_{k-1}, \delta_{k-1}\left(R_{k-2}(u^{(k-2)}, \Phi_{k-3})\right)\right)}\right)C\left(E_k, \delta_k\left(R_{k-2}(u^{(k-2)}, \Phi_{k-3}) \right)\right)},
\end{align*}
we derive
\begin{align*}\nonumber
&R_k(u^{(k)},\delta_k(\Phi_{k-1}))\\
=&{\rm Ad}\left({C\left(E_k, \delta_k(R_{k-1}(u^{(k-1)}, \Phi_{k-2}))\right)}\right)R_{k-1}(u^{(k-1)}, \Phi_{k-2})\\
=&{\rm Ad}\left({{C\left(E_{k-1}, \delta_{k-1}(R_{k-2}(u^{(k-2)}, \Phi_{k-3}))\right)} C\left(E_k, \delta_{k}(R_{k-2}(u^{(k-2)}, \Phi_{k-3})) \right)}\right) R_{k-2}(u^{(k-2)}, \Phi_{k-3}).
\end{align*}
Keep doing this for $k=n-1,n-2,...,2$ and note that $R_1(u^{(1)},\Phi_0)=\Phi_0$, we get
\[C\left(u, \Phi_{n-1}(u;\Phi_0)\right)\cdot e^{ \Phi_{n-1}}\cdot C\left(u, \Phi_{n-1}(u;\Phi_0)\right)^{-1}=\left(\overrightarrow{\prod_{k=2,...,n}}C\left(E_k,\delta_k(\Phi_0)\right)\right)\cdot e^{\Phi_0}\cdot \left(\overrightarrow{\prod_{k=2,...,n}}C\left(E_k,\delta_k(\Phi_0)\right)\right)^{-1}.\]
Therefore, by the monodromy relation, we get
\begin{align*}
    &S_-\left(u, \Phi_{n-1}(u;\Phi_0)\right) S_+\left(u, \Phi_{n-1}(u;\Phi_0)\right)\\ 
    &= C\left(u, \Phi_{n-1}(u;\Phi_0)\right)\cdot e^{ \Phi_{n-1}}\cdot C\left(u, \Phi_{n-1}(u;\Phi_0)\right)^{-1}\\ &=\left(\overrightarrow{\prod_{k=2,...,n}}C\left(E_k,\delta_k(\Phi_0)\right)\right)\cdot e^{\Phi_0}\cdot \left(\overrightarrow{\prod_{k=2,...,n}}C\left(E_k,\delta_k(\Phi_0)\right)\right)^{-1}.
\end{align*}
It finishes the proof. \qed

\vspace{2mm}
Theorem \ref{relautomorphism} motivates the following definitions. Let us denote by $u_{\rm cat}$ the infinite point $z_1\rightarrow 0+, z_k\rightarrow+\infty$ for all $k$. 
\begin{defi}\label{connmap}
For any $A\in\Herm(n)$, the connection matrix $C(u_{\rm cat},A)$, at the infinite point $u_{\rm cat}$ with respect to $U_{\rm id}$, is the pointwise ordered multiplication of all $C(E_k,\delta_k(A))$'s for $k=1,...,n$. That is 
\[C(u_{\rm cat},A):=C(E_1,\delta_1(A))C(E_2,\delta_2(A))\cdot\cdot\cdot C(E_n,\delta_n(A)), \ \ \ \text{for any} \ A\in\Herm_0(n).\]
\end{defi}

\begin{defi}\label{relStokes}
For any $A\in{\Herm}(n)$, the {\it Stokes matrices $S_{\pm}(u_{\rm cat},A)$} at $u_{\rm cat}$, with respect to the connected component $U_{\rm cat}$, are respectively the upper and lower triangular matrices uniquely determined by the identity (Gauss decomposition) 
\begin{eqnarray}\label{monorelation}
S_{-}(u_{\rm cat}, A) S_{+}(u_{\rm cat}, A)=C(u_{\rm cat},A)e^AC(u_{\rm cat},A)^{-1},
\end{eqnarray}
and by imposing the diagonal part $[S_+(u_{\rm cat},A)]=[S_-(u_{\rm cat},A)]=e^{\frac{[A]}{2}}$.
\end{defi}

It follows from Theorem \ref{relautomorphism} that 
Definition \ref{monorelation} coincides with the Definition \ref{caterpillar} given in the introduction. That is
\begin{cor}\label{uucat}
For any $u\in U_{\rm id}$, the Stokes matrices $S_{\pm}(u, \Phi(u;\Phi_0))$ coincides with $S_\pm(u_{\rm cat},\Phi_0)$.
\end{cor}

\subsection{Explicit expression of Stokes matrices in terms of Gelfand-Tsetlin coordinates}\label{explicitnurel}
In this subsection, we will give an explicit formula of the Stokes matrices $S_\pm(u_{\rm cat},A)$ for any $A\in\Herm(n)$, and prove Theorem \ref{mainthm}. To this end, we only need to compute explicitly the connection matrices of the system \eqref{eq:relativeStokes} for all $1\le k\le n$. The computational procedure is better to be understood in terms of the Gelfand-Tsetlin integrable systems \cite{GS}. So before computing the Stokes matrices, let us introduce the action and angle coordinates of the Gelfand-Tsetlin systems.

\subsubsection{The Gelfand-Tsetlin coordinates}\label{GTcoor}
{\bf Gelfand-Tsetlin maps.} For $k\le n$ let $A^{(k)}\in \Herm(k)$ denote the upper left $k$-th submatrix (upper left $k\times k$ corner) of a Hermitian matrix $A\in \Herm(n)$,
and $\lambda^{(k)}_i(A)$-its ordered set of
eigenvalues, $\lambda_1^{(k)}(A)\ge \cdots\ge \lambda_k^{(k)}(A)$. 
The map 
\begin{equation}\label{eq:momentmap}
\lambda\colon \Herm(n)\to \mathbb{R}^{\frac{n(n+1)}{2}},
\end{equation}
taking $A$ to the collection of numbers $\lambda_i^{(k)}(A)$ for $1\le i\le
k\le n$, is continuous and is called the Gelfand-Tsetlin map.
Its image $\mathcal{C}(n)$ is the Gelfand-Tsetlin cone, cut out by the following inequalities,
\begin{equation}\label{eq:cone}
\lambda_i^{(k+1)}\ge \lambda_i^{(k)}\ge \lambda_{i+1}^{(k+1)},\ \ 1\le
i\le k\le n-1.
\end{equation}

\vspace{2mm}
{\bf Thimm torus actions.}
Let $\mathcal{C}_0(n)\subset \mathcal{C}(n)$ denote the subset where all of the eigenvalue inequalities
\eqref{eq:cone} are strict. Let $\Herm_0(n):=\lambda^{-1}(\mathcal{C}_0(n))$ be the corresponding dense open subset of $\Herm(n)$. The $k$-torus 
$T(k)\subset {\rm U}(k)$ of diagonal matrices acts on 
$\Herm_0(n)$ as follows,  
\begin{equation}\label{eq:taction}
 t\bullet A=\Ad_{U^{-1} t U}A,\ \ \ \ t\in T(k),\ A\in \Herm_0(n).
\end{equation}
Here $U\in {\rm U}(k)\subset {\rm U}(n)$ is a unitary matrix such that $\Ad_{U}A^{(k)}$ is
diagonal, with entries $\lambda_1^{(k)},\ldots,\lambda_k^{(k)}$. The action
is well-defined since $U^{-1}t U$ does not depend on the choice of
$U$, and preserves the Gelfand-Tsetlin map \eqref{eq:momentmap}.  The actions
of the various $T(k)$'s commute, hence they define an action of the Gelfand-Tsetlin torus $T(1)\times \cdots \times T(n-1)\cong \U(1)^\frac{(n-1)n}{2}.$
Here the torus $T(n)$ is excluded, since the 
action \eqref{eq:taction} is trivial for $k=n$.  

\vspace{2mm}
{\bf Action-angle coordinates.}
If $A\in \rm Herm_0(n)$, then there exists a unique unitary matrix $P_k(A)\in {\rm U}(k)\subset {\rm U}(n)$, whose entries in the $k$-th row are positive and real, such that the upper left $k$-th submatrix of $A_k:=P_k(A)^{-1}AP_k(A)$ is the diagonal matrix ${\rm diag}(\lambda^{(k)}_1,...,\lambda^{(k)}_k)$, i.e., 
\begin{eqnarray}\label{zhuazi}
    A_k=P_k(A)^{-1}AP_k(A)=\begin{pmatrix}
    \lambda^{(k)}_1 & & & &a^{(k)}_1 & \cdots\\
    & \ddots & && \vdots & \cdots\\
    & &  \lambda^{(k)}_k && a^{(k)}_k & \cdots \\
    &&&&\\
    \overline{a_1^{(k)}}&\cdots &\overline{a_k^{(k)}}&& \lambda_{k+1}^{(k)} & \cdots\\
    \cdots & \cdots & \cdots &&\cdots & \cdots
    \end{pmatrix}.
\end{eqnarray}
The $(i,k+1)$ entries $a^{(k)}_i(A)$, for $1\le i\le k\le n-1$, are seen as functions on $\Herm_0(n)$.
Then the functions $\{\lambda^{(k)}_i\}_{1\le i\le k\le n}$ and $\{\psi^{(k)}_i={\rm Arg}(a^{(k)}_i)\}_{1\le i\le k\le n-1}$ on ${\rm Herm}_0(n)$ are called the Gelfand-Tsetlin action and angle coordinates.

\subsubsection{Diagonalization in stages}\label{diagonalization}
Now let us consider the system \eqref{eq:relativeStokes}. By the ${\rm U}(k-1)$-equivariant property in Lemma \ref{commute}, to simplify the computation, we can first diagonalize the upper left $k-1$-th submatrix of $A$. According to Section \ref{GTcoor}, there is a systematic way to do it for all $A\in\Herm_0(n)$ via the unique unitary matrix $P_k(A)$ in \eqref{zhuazi}.

These unitary matrices can also be inductively defined: suppose that we already have $P_k(A)\in {\rm U}(k)\subset {\rm U}(n)$ for $k<n$, such that $A_{{k}}=P_k(A)^{-1}AP_k(A)$ is the form in \eqref{zhuazi}. Let $L^{(k+1)}(A)\in {\rm U}(k+1)\subset {\rm U}(n)$ be the matrix given by \begin{eqnarray}\label{L}
&&L^{(k+1)}_{ij}(A):=\frac{a^{(k)}_i}{ N_j^{(k+1)}(\lambda^{(k)}_i-\lambda^{(k+1)}_j)}, \ \text{for} \ i\ne k+1, \ j=1,...,k+1,\\ \label{L1}
&&L^{(k+1)}_{k+1,j}(A):=\frac{1}{N_j^{(k+1)}}, \ \ \ \ \ \ \text{for} \ j=1,...,k+1,
\end{eqnarray}
where the normalizer \begin{equation}\label{normalizer}
N_j^{(k+1)}(A):=\sqrt{1+\sum_{l=1}^{k}\frac{|a_l^{(k)}|^2}{(\lambda^{(k)}_l-\lambda^{(k+1)}_j)^2}}=\sqrt{\frac{\prod_{v=1, v\ne j}^{k+1}(\lambda^{(k+1)}_j-\lambda^{(k+1)}_v) }{\prod_{v=1}^{k}(\lambda^{(k+1)}_j-\lambda^{(k)}_v)}}.
\end{equation}
(Here the second identity follows from the combinatorial identity given by computing the character polynomial of the upper left $k+1$-th submatrix of $A_k$ in two different ways.)
The upper left $k+1$-th submatrix of $L^{(k+1)}(A)^{-1}A_{k}L^{(k+1)}(A)$ is ${\rm diag}(\lambda^{(k+1)}_1,...,\lambda^{(k+1)}_{k+1})$ and the entries in $k+1$-th row of $L^{(k+1)}(A)$ is real and positive, thus we can simply define $P_{k+1}(A)$ by 
\begin{equation}\label{diagP}
P_{k+1}(A):=P_k(A)\cdot L^{(k+1)}(A).
\end{equation}

Furthermore, using the Laplace expansion, we see that for any $2\le k\le n$, the $n\times n$ matrix $P_k(A)$ has the explicit form
\begin{align}\label{clsP}
(P_k)_{ij}:=&{(-1)^{k+i}}\frac{\Delta_{1,...,\hat{i},...,k}^{1,...,k-1}\left(A-\lambda^{(k)}_j\right)}{\sqrt{\prod_{l=1,l\ne i}^k(\lambda^{(k)}_j-\lambda^{(k)}_l)\prod_{l=1}^{k-1}(\lambda^{(k)}_j-\lambda^{(k-1)}_l)}}, \ \text{if} \ 1\le i,j\le k\\ 
\nonumber (P_k)_{ii}:=&1, \ \ \text{if} \ i>k, \\ 
\nonumber (P_k)_{ij}:=&0, \ \ \text{otherwise},
\end{align}
and its inverse is,
\begin{align}\label{clsQ}
(P_k^{-1})_{ij}=&{(-1)^{k+j}}\frac{\Delta^{1,...,\hat{j},...,k}_{1,...,k-1}\left(A-\lambda^{(k)}_i\right)}{\sqrt{\prod_{l=1,l\ne i}^k(\lambda^{(k)}_i-\lambda^{(k)}_l)\prod_{l=1}^{k-1}(\lambda^{(k)}_i-\lambda^{(k-1)}_l)}}, \ \text{if} \ 1\le i,j\le k\\ 
\nonumber (P_k^{-1})_{ii}=&1, \ \ \text{if} \ i>k, \\ 
\nonumber (P_k^{-1})_{ij}=&0, \ \ \text{otherwise}.
\end{align}
Therefore, by definition the function
\begin{equation}\label{akiexp}
    a^{(k)}_i(A)= \sum_{v=1}^n(P_k(A)^{-1})_{iv} \cdot (A)_{v,k+1}=\frac{(-1)^{k+i}\Delta^{1,...,k}_{1,...,k-1,k+1}\left(A-\lambda^{(k)}_i\right)}{\sqrt{\prod_{l=1,l\ne i}^k(\lambda^{(k)}_i-\lambda^{(k)}_l)\prod_{l=1}^{k-1}(\lambda^{(k)}_i-\lambda^{(k-1)}_l)}}.
\end{equation}

Combing \eqref{normalizer} and \eqref{akiexp}, we have the following relation that will be used in Section \ref{eleabove}.
\begin{lem}\label{aNrelation}
For any $A\in\Herm(n)_0$, we have the identity
 \[ \frac{a^{(k)}_i(A)}{N_i^{(k)}(A)}=\frac{\Delta^{1,...,k}_{1,...,k-1,k+1}\left(A-\lambda^{(k)}_i\right)}{\prod_{l=1,l\ne i}^k(\lambda^{(k)}_i-\lambda^{(k)}_l)}\]
\end{lem}

Since $P_{k-1}(A)$ commutes with $E_k$, the conjugation by $P_{k-1}$ simplifies the equation \eqref{eq:relativeStokes} 
\[
\frac{dF}{dz}=\left(\I E_k-\frac{1}{2\pi \I }\frac{\delta_k(A)}{z}\right)F\]
to
\[
\frac{dF}{dz}=\left(\I E_{k}-\frac{1}{2\pi \I }\frac{\delta_{k}(A_{k-1})}{z}\right)\cdot F.\]
In particular, if $F(z)$ is a fundamental solution of the second equation, 
then ${P}_{k-1}(A) F(z) {P}_{k-1}(A)^{-1}$ is a fundamental solution of the first equation.
Moreover, the Stokes/connection matrices of above two equations are related by
\begin{align}\label{PCP}
C(E_k,\delta_k(A))&={P}_{k-1}(A) C(E_k,\delta_k(A_{k-1}) {P}_{k-1}(A)^{-1}\\
S_{\pm}(E_k, \delta_k(A))&={P}_{k-1}(A) S_{+}(E_k, \delta_k(A_{k-1})){P}_{k-1}(A)^{-1}.
\end{align}

\subsubsection{Normalized connection matrices}
Following \eqref{PCP}, for any $A$, the product appearing in the Definition \ref{connmap} is
\begin{align}\nonumber 
&C(E_1,\delta_1(A))\cdots C(E_{k-1},\delta_{k-1}(A))C(E_k,\delta_k(A))C(E_{k+1},\delta_{k+1}(A))\cdots C(E_n,\delta_n(A))\\ \label{normproduct}
    =&\cdots C(E_{k-1},\delta_{k-1}(A_{k-2})) {P}_{k-2}^{-1}{P}_{k-1} \cdot C(E_k,\delta_k(A_{k-1})) {P}_{k-1}^{-1}{P}_{k} \cdots
\end{align} 
It motivates to compute $C(E_k,\delta_k(A_{k-1}) {P}_{k-1}^{-1}{P}_{k}$ directly. Let us introduce (recall that $L^{(k)}(A)={P}_{k-1}^{-1}{P}_{k}$ was defined in \eqref{L})
\begin{defi}\label{defitildeC}
For any integer $1\le k\le n$ and any $A\in\Herm_0(n)$, we define the {\it normalized connection matrix}
\begin{eqnarray}\label{tildeC}
\widetilde{C}\left(E_k,\delta_k(A_{k-1})\right):=C\left(E_k, \delta_k(A_{k-1})\right)\cdot L^{(k)}(A).
\end{eqnarray}
Here recall from \eqref{zhuazi} that $A_{k-1}=P_{k-1}(A)^{-1}AP_{k-1}(A)$ is the diagonalization of the upper left $(k-1)$-th submatrix of $A$.
\end{defi}
\begin{pro}\label{mono2}
For any $A\in {\rm Herm}_0(n)$, if we define \[\widetilde{C}(u_{\rm cat},A)=\widetilde{C}(E_1,\delta_1(A)) \widetilde{C}(E_2,\delta_2(A))\cdot\cdot\cdot \widetilde{C}(E_n,\delta_n(A)) \in U(n),\]
as the pointwise multiplication, then we have \[\widetilde{C}(u_{\rm cat}, A)e^{A_{n}}\widetilde{C}(u_{\rm cat},A)^{-1}=C(u_{\rm cat},A)e^{A}C(u_{\rm cat},A)^{-1}.\] Here recall that $A_{n}={\rm diag}(\lambda_1^{(n)},...,\lambda^{(n)}_n)$. 
\end{pro}
\pf It follows from the identities \eqref{normproduct}, \eqref{diagP} and Definition \eqref{tildeC} that $\widetilde{C}(u_{\rm cat},A)=C(u_{\rm cat},A)\cdot P_n(A)$. Then the proposition follows from the identity $P_n(A)e^{A_{n}}P_n(A)^{-1}=e^A$, i.e., the definition $A_{n}=P_n(A)^{-1}AP_n(A)$ of $A_n$.
\qed

\subsubsection{Stokes matrices of special confluent hypergeometric equations}\label{explicitev}
For any $A\in{\Herm}_0(n)$, the Stokes matrix and normalized connection matrix of the system, for $k=1,...,n-1$,
\begin{eqnarray}\label{relequation}
\frac{dF}{dz}=\left(\I E_{k+1}-\frac{1}{2\pi \I }\frac{\delta_{k+1}(A_{k})}{z}\right)\cdot F,
\end{eqnarray}
are described by the following proposition. Here recall that $A_{k}$ is defined in \eqref{zhuazi}. \begin{rmk}
The equation \eqref{relequation} is a very special type of confluent hypergeometric equation, and can be solved exactly via the confluent hypergeometric functions $_kF_k$. Thus the expression of its Stokes matrices simply follows from the well known asymptotics of $_kF_k$, following the line of \cite[Proposition 8]{BJL} for $n=2$ case. The explicit expression has been derived in \cite{Balser}, see Remark \ref{Bal}, see also our paper \cite{LX}.
However, since our convention is different from theirs, and in order to illustrate how the computation can be generalized to the quantum case, we outline a proof using the asymptotics of $_kF_k$. 
\end{rmk}
\begin{pro}\label{explicitmCS}\hfill

$(1).$ The entries of the Stokes matrix $S_+\left(E_{k+1},\delta_{k+1}(A_{k})\right)$ of \eqref{relequation} are
\begin{align*}
(S_+)_{j,k+1}&= \frac{ e^{\frac{\lambda^{(k)}_j+(A)_{k+1,k+1}}{4}}\prod_{i=1}^{k}\Gamma(1+\frac{\lambda^{(k)}_i-\lambda^{(k)}_j}{2\pi \I })}{\prod_{i=1}^{k+1}\Gamma(1+\frac{\lambda^{(k+1)}_i-\lambda^{(k)}_j}{2\pi \I })}\cdot {a^{(k)}_j(A)},
 \ \ \ \text{for} \ j=1,...,k;\\
(S_+)_{ii}&=e^{\frac{A_{ii}}{2}}, \ \ \text{for} \ i=1,...,n, \\
(S_+)_{ij}&=0, \ \ \text{otherwise}.
\end{align*}

$(2).$ The entries of the normalized connection matrix $\widetilde{C}\left(E_{k+1},\delta_{k+1}(A_{k})\right)$ are given by\begin{equation*} \widetilde{C}_{ij}=\frac{ e^{\frac{\lambda^{(k)}_i-\lambda^{(k+1)}_j}{4}}}{(\lambda^{(k)}_i-\lambda^{(k+1)}_j)}\frac{\prod_{v=1}^{k+1}\Gamma(1+\frac{\lambda^{(k+1)}_v-\lambda^{(k+1)}_j}{2\pi \I })\prod_{v=1}^{k}\Gamma(1+\frac{\lambda^{(k)}_v-\lambda^{(k)}_i}{2\pi \I })}{\prod_{v=1,v\ne i}^{k}\Gamma(1+\frac{\lambda^{(k)}_v-\lambda^{(k+1)}_j}{2\pi \I })\prod_{v=1,v\ne j}^{k+1}\Gamma(1+\frac{\lambda^{(k+1)}_v-\lambda^{(k)}_i}{2\pi \I })}\cdot \frac{a^{(k)}_i(A)}{N_j^{(k+1)}(A)},
\end{equation*}
for $1\le j\le k+1, 1\le i\le k$, and
\begin{align*} \widetilde{C}_{k+1,j}&=\frac{e^{\frac{\lambda^{(k+1)}_j-(A)_{k+1,k+1}}{4}}\prod_{v=1}^{k+1}\Gamma(1+\frac{\lambda^{(k+1)}_v-\lambda^{(k+1)}_j}{2\pi \I })}{N_j^{(k+1)}(A)\cdot \prod_{v=1}^{k}\Gamma(1+\frac{\lambda^{(k)}_v-\lambda^{(k+1)}_j}{2\pi \I })}, \ \ \ \text{for} \ 1\le j\le k+1.\\
\widetilde{C}_{ii}&=1, \ \ \ \ \ \ \text{for} \ k+1<i\le n,\\
\widetilde{C}_{ij}&=0, \ \ \ \ \ \ \text{otherwise}.
\end{align*}
\end{pro}
\pf 
Under certain generic additional assumptions one can explicitly compute a Floquet solution of the equation \eqref{relequation}, using the generalized confluent hypergeometric functions. Recall that they are the functions, for any $m\ge 1$, $\alpha_j\in\mathbb{C}$, $\beta_j\in\mathbb{C}\setminus \{0, -1, -2, ...\}, 1\le j\le m$, are
\begin{eqnarray*}
_kF_k(\alpha_1,...,\alpha_m,\beta_1,...,\beta_m;z)=\sum_{n=0}^\infty \frac{(\alpha_1)_n\cdots(\alpha_m)_n}{(\beta_1)_n\cdots(\beta_m)_n}\frac{z^n}{n!},
\end{eqnarray*}
where $(\alpha)_0=1$, $(\alpha)_n=\alpha\cdots (\alpha+n-1),$ $n\ge 1$. 

\begin{lem}\label{hypersystem}
The equation \eqref{relequation} has a fundamental solution $F\left(z;E_{k+1},\delta_{k+1}(A_{k})\right)$ taking the form
\begin{eqnarray}
F\left(z;E_{k+1},\delta_{k+1}(A_{k})\right)=Y\cdot H\left(z;E_{k+1},\delta_{k+1}(A_{k})\right)\cdot z^{-\frac{1}{2\pi \I }\delta_{k+1}(A_{k+1})},
\end{eqnarray}
where $Y={\rm diag}(a^{(k)}_1, ..., a^{(k)}_k, 1,...,1)$ and $H\left(z;E_{k+1},\delta_{k+1}(A_{k})\right)$ is the $n\times n$ matrix given by
\begin{align*}
&H(z)_{ij}=\frac{1}{\lambda^{(k+1)}_j-\lambda^{(k)}_i}\cdot {_kF_k}(\alpha_{ij,1},...,\alpha_{ij,k},\beta_{ij,1},...,\widehat{\beta_{ij,j}},...,\beta_{ij,k+1};\I z), \ \ \ 1\le i\le k, 1\le j\le k+1,\\
&H(z)_{k+1,j}= {_kF_k}(\alpha_{k+1j,1},...,\alpha_{k+1j,k},\beta_{k+1j,1},...,\widehat{\beta_{k+1j,j}},...,\beta_{k+1j,k+1};\I z), \ \ i=k+1, \ 1\le j\le k+1,\\
&H(z)_{ii}=1, \ \ \ \ \ \ \text{for} \ k+1<i\le n,\\
&H(z)_{ij}=0, \ \ \ \ \ \ \text{otherwise},
\end{align*}
with the variables $\{\alpha_{ij,l}\}$ and $\{\beta_{ij,l}\}$ given by
\begin{eqnarray*}
&&\alpha_{ij,i}=\frac{1}{2\pi \I }(\lambda^{(k)}_i-\lambda^{(k+1)}_j), \ \ \ 1\le i\le k, \ 1\le j\le k+1,\\
&&\alpha_{k+1j, i}=1+\frac{1}{2\pi \I }(\lambda^{(k)}_{i}-\lambda^{(k+1)}_j), \ \ \ 1\le j\le k+1,\\
&&\alpha_{ij, l}=1+\frac{1}{2\pi \I }(\lambda^{(k)}_l-\lambda^{(k+1)}_j), \ \ \ \ l\ne i, \  1\le l\le k, \ 1\le i, j\le k+1, \\
&&\beta_{ij, l}=1+\frac{1}{2\pi \I }(\lambda^{(k+1)}_l-\lambda^{(k+1)}_j), \ \ \ \ l\ne j, \ 1\le l\le k+1, \ \ 1\le i, j\le k+1.
\end{eqnarray*}
\end{lem}
\pf The first $k$ rows of the matrix equation \eqref{relequation} follows from the definition of the functions $H(z)_{ij}$ and the special arguments $\alpha_{ij,l}$ and $\beta_{ij,l}$, that the functions $H\left(z;E_{k+1},\delta_{k+1}(A_{k})\right)_{ij}$ satisfy
\begin{eqnarray*}
z\frac{dH_{ij}}{dz}=\frac{{\small\lambda^{(k+1)}_j-\lambda^{(k)}_{i}}}{2\pi \I } H_{ij}-\frac{1}{2\pi \I }H_{k+1,j}, \ \ \ for \ 1\le i\le k.
\end{eqnarray*}
For the rest of the equation, one just needs the identity $a_i^{(k)}b_i^{(k)}=-\frac{\prod_{j=1}^{k+1}(\lambda^{(k+1)}_j-\lambda^{(k)}_i)}{\prod_{j\ne i}^{k}(\lambda^{(k)}_j-\lambda^{(k)}_i)}$, for any $1\le i\le k$, which follows from
the character polynomial of $A_{k}$. \qed

\vspace{2mm}
The asymptotics expansion of $_kF_k$ via gamma functions are (see \cite[Page 411]{OLBC}),
\begin{align}\nonumber
&\frac{\prod_{l=1}^k\Gamma(\alpha_l)}{\prod_{l=1}^k\Gamma(\beta_l)}\cdot  {_kF_k}(\alpha_1,...,\alpha_k,\beta_1,...,\beta_k;z)\\ \label{asyKummer}
&\sim \sum_{m=1}^k \Gamma(\alpha_m)\frac{\prod_{l=1,l\ne m}^k\Gamma(\alpha_l-\alpha_m)}{\prod_{l=1}^k\Gamma(\beta_l-\alpha_m)}(\mp z)^{-\alpha_m}(1+\mathcal{O}(z^{-1}))+e^z z^{\sum_{l=1}^k(\alpha_l-\beta_l)}(1+\mathcal{O}(z^{-1})).
\end{align}
where upper or lower signs are chosen as $z$ lies in the upper or lower half-plane. Using this, one can get explicitly the asymptotics of $F\left(z;E_{k+1},\delta_{k+1}(A_{k})\right)$ as $z\rightarrow \infty$ in the two different Stokes sectors, and its comparison with the unique formal solution of \eqref{relequation}\[\hat{F}\left(z;E_{k+1},\delta_{k+1}(A_{k})\right)=({\rm Id}_n+\mathcal{O}(z^{-1}))e^{\I zE_{k+1}}z^{-\frac{1}{2\pi \I }[\delta_{k+1}(A_{k})]}.\] In particular, we have 
\begin{eqnarray*}
&&F\left(z;E_{k+1},\delta_{k+1}(A_{k})\right)\sim \hat{F}\left(z;E_{k+1},\delta_{k+1}(A_{k})\right)\cdot Y  U_+, \ \ \ \text{as} \ z\rightarrow \ \infty \ \text{in} \  S(-\pi,\pi)\\ 
&&F\left(z;E_{k+1},\delta_{k+1}(A_{k})\right)\sim \hat{F}\left(z;E_{k+1},\delta_{k+1}(A_{k})\right)\cdot Y U_-, \ \ \ \text{as} \ z\rightarrow \ \infty \ \text{in} \ S(-2\pi,0),
\end{eqnarray*} 
where $Y={\rm diag}(a^{(k)}_1, ..., a^{(k)}_k, 1,...,1)$
and $U_\pm$ are certain explicit invertible matrices with entries given by gamma functions. Then by the uniqueness in Theorem \ref{uniformresum}, we know that the canonical solutions 
\begin{eqnarray*}
F\left(z;E_{k+1},\delta_{k+1}(A_{k})\right))=F_+\left(z;E_{k+1},\delta_{k+1}(A_{k})\right)Y U_+ \hspace{3mm} in \ \ {\rm Sect}_+,\\
F\left(z;E_{k+1},\delta_{k+1}(A_{k})\right)=F_-\left(z;E_{k+1},\delta_{k+1}(A_{k})\right)Y U_- \hspace{3mm} in \ \ {\rm Sect}_-.
\end{eqnarray*}
Then by definition, the Stokes matrices are given by 
(here to derive the second formula,
the change in choice of ${\rm log}(z)$ in $F(z)$ is accounted for)
\begin{eqnarray}\label{exlicitcom}
S_+\left(E_{k+1},\delta_{k+1}(A_{k})\right)= YU_-U_+^{-1}Y^{-1}, \ \ \ \ \ \ S_-\left(E_{k+1},\delta_{k+1}(A_{k})\right)= YU_-e^{\delta_{k+1}(A_{k+1})}U_+^{-1}Y^{-1}.
\end{eqnarray}
The explicit computation of $U_\pm$ using the asymptotics of $_kF_k$ is straight-forward but 
lengthy, and may be omitted here. For the explicit expression of $U_\pm$ and more details on the computation, see our paper \cite{LX}. Similarly, one can get the expression of the connection matrix by considering the asymptotics of $F(z)$ at $z=0$.
\qed

\begin{rmk}\label{Bal}
The proof in \cite{Balser} didn't make use of any known results on the global behaviour of
the functions $_kF_k$. In particular,  the expression of the central connection factors $\widetilde{\Omega_0}$ of the meromorphic linear system $\frac{dF}{dz}=(E_{k+1}+\frac{A}{z})\cdot F$ are given in \cite[Formula 7.3-7.4]{Balser}. Following \cite[Formula 6.3]{Balser} (where our $C(E_{k+1},A_{k})$ is denoted by $\Omega_0$ there), the connection factor $\widetilde{\Omega_0}$ is related to the connection matrix by
\[\widetilde{\Omega_0}=C(E_{k+1},A_{k})\cdot \widetilde{L}_0,\]
where (see the definition in \cite[Formula 5.2-5.3]{Balser}) the matrix $\widetilde{L}_0=L^{(k+1)}\cdot D$ with $D$ given by
\begin{eqnarray*}
&&D_{jj}=\frac{\prod_{t=1}^k(\lambda^{(k+1)}_j-\lambda^{(k)}_t)}{\prod_{t=1}^k\Gamma(1+\lambda^{(k+1)}_j-\lambda^{(k+1)}_t)}, \ \ \ 1\le j\le k,\\
&&D_{k+1,k+1}=\frac{1}{\prod_{t=1}^k\Gamma(1+\lambda^{(k+1)}_j-\lambda^{(k+1)}_t)}.
\end{eqnarray*}
Thus the normalized connection matrix $C^{(k+1)}(A_{k})L^{(k+1)}$, what we are computing, differs from the connection factor $\widetilde{\Omega_0}$ by $D$. One checks that multiplying the formula \cite[Formula 7.3-7.4]{Balser} for $\widetilde{\Omega_0}$ by $D$ gives rise to the formula in Proposition \ref{explicitmCS} (provided that the replacement of the matrix $A$ by $-A$ and the irregular term $E_{k+1}$ by $\I E_{k+1}$ are accounted for).
\end{rmk}
\begin{rmk}\label{clsdiag}
We note that the expression $P_k(A) F(z;A) P_k(A)^{-1}$, with $F(z;A)$ the solution of the diagonalized equation \eqref{relequation} given as in Lemma \ref{hypersystem}, smoothly extends from $A\in\Herm(n)_0$ to $\Herm(n)$ and gives the fundamental solution of the equation \ref{eq:relativeStokes} for all $A$. In the meanwhile, the Stokes matrices of the equation \ref{eq:relativeStokes} are just
\[S_+\left(E_{k+1},\delta_{k+1}(A)\right)=P_k(A)S_+\left(E_{k+1},\delta_{k+1}(A_{k})\right)P_k(A)^{-1}\]
i.e., the $P_k(A)$ conjugation of the ones of the diagonalized equation \eqref{relequation} given in Proposition \ref{explicitmCS}. The Stokes matrices $S_+\left(E_{k+1},\delta_{k+1}(A)\right)$ are real analytic function with respect to all $A\in\Herm(n)$, which also follows from the definition of the Stokes matrices of the nonresonant equation.
\end{rmk}

\subsubsection{Entries of Stokes matrices}\label{eleabove}
Recall that for any $1\le k\le n-1$, $A^{(k)}\in \Herm_0(k)$ denotes the upper left $k$-th submatrix (upper left $k\times k$ corner) of a Hermitian matrix $A\in \Herm_0(n)$. For any $k$, we define a map $\widetilde{C}_k\colon \Herm_0(k)\rightarrow U(k)$ by
\[\widetilde{C}_k(A^{(k)})=\widetilde{C}(E_1,\delta_1(A^{(k)}))\cdots \widetilde{C}(E_k,\delta_k(A^{(k)})), \hspace{3mm} \text{for} \ A^{(k)}\in\Herm_0(k).\] 
Then we define upper and lower $k\times k$ matrices $S_{k\pm}(A^{(k)})$ via 
\begin{eqnarray}\label{defSk}
\widetilde{C}_k(A^{(k)})e^{A^{(k)}_{k}}\widetilde{C}_k(A^{(k)})^{-1}=S_{k-}(A^{(k)}) S_{k+}(A^{(k)}),
\end{eqnarray}
with the diagonal part $[S_{k+}(A^{(k)})]=[S_{k-}(A^{(k)})]=e^{\frac{A^{(k)}_{k}}{2}}$. Here recall from the definition \eqref{zhuazi} that the lower index $k$ denotes the diagonalization of the upper left k-th submatrix of $A^{(k)}$, i.e., $A^{(k)}_k$ is the diagonal matrix ${\rm diag}(\lambda^{(k)}_1,...,\lambda^{(k)}_k)$ with the ordered eigenvalues. By definition $S_{k\pm}(A^{(k)})$ are just the $k\times k$ Stokes matrices at $u_{\rm cat}$.

Recall that the $n\times n$ Stokes matrix $S_{+}(u_{\rm cat},A)$ at $u_{\rm cat}$ is upper triangular. For any $1\le k\le n-1$, let us denote by $b_{k+1}(A)$ the column vector consisting of the first $k$ elements of the $k+1$-th column of $S_{+}(u_{\rm cat}, A)$, i.e., \[b_{k+1}(A)=\left(S_{+}(u_{\rm cat}, A)_{1,k+1},..., S_{+}(u_{\rm cat}, A)_{k,k+1}\right)^T.\]
\begin{lem}\label{mono3}
The column vector $b_{k+1}(A)$ is equal to the product of matrices \begin{eqnarray}
b_{k+1}(A)=S_{k+}(A^{(k)}_{{k-1}})\widetilde{C}_{k}(A^{(k)})e^{-\frac{A_{k}^{(k)}}{2}}b^{(k+1)}(A_{k}),
\end{eqnarray}
where by the definition \eqref{zhuazi} $A^{(k)}_{{k-1}}=P_{k-1}(A^{(k)})^{-1}A^{(k)}P_{k-1}(A^{(k)})$ is the matrix from the diagonalization of the upper left $(k-1)$-th submatrix of $A^{(k)}$, $b^{(k+1)}(A_{k})$ is the column vector consisting of the first $k$ entries of the $k+1$-th column of the Stokes matrix $S_+\left(E_{k+1},\delta_{k+1}(A_k)\right)$ of equation \eqref{relequation}, i.e., \[b^{(k+1)}(A)=\left(S_+\left(E_{k+1},\delta_{k+1}(A_k)\right)_{1,k+1},..., S_+\left(E_{k+1},\delta_{k+1}(A_k)\right)_{k,k+1}\right)^T.\]
\end{lem}
\pf It follows from the expression of the monodromy relation \eqref{monorelation} that the upper left $k+1$ submatrix of $S_+(u_{\rm cat},A)$ coincides with $S_{k+1}(A^{(k+1)})$. So it is enough to prove the case $k=n-1$. 
Using the monodromy relation 
\[\widetilde{C}(E_n,A)e^{A_{n}}\widetilde{C}(E_n,A)^{-1}=S_-(E_n,A_{n-1})S_+(E_n, A_{n-1}),\]
and the expression \[S_-(E_n, A_{n-1})=\left(\begin{array}{cc}
e^{\frac{1}{2}A_{n-1}^{(n-1)}} \ \ & 0\\
b^{(n)}(A_{n-1})^\dagger \ \ & \ast
\end{array} \right), \ \ \ S_+(E_n,A_{n-1})=\left(\begin{array}{cc}
e^{\frac{1}{2}A_{n-1}^{(n-1)}} & b^{(n)}(A_{n-1})\\
0& \ast
\end{array} \right),\]
we get 
\begin{align*}
&S_{n-}(A) S_{n+}(A)\\
=&\left(\begin{array}{cc}
\widetilde{C}_{n-1}(A^{(n-1)})&0\\
0& 1
\end{array} \right)\widetilde{C}(E_n,A)e^{A_{n}}\widetilde{C}(E_n,A)^{-1}\left(\begin{array}{cc}
\widetilde{C}_{n-1}(A^{(n-1)})^{-1}&0\\
0& 1
\end{array} \right)\\
=&\left(\begin{array}{cc}
\widetilde{C}_{n-1}(A^{(n-1)})&0\\
0& 1
\end{array} \right)\left(\begin{array}{cc}
e^{\frac{1}{2}A_{n-1}^{(n-1)}}&0\\
{b^{(n)}}^\dagger& \ast
\end{array} \right)\left(\begin{array}{cc}
e^{\frac{1}{2}A_{n-1}^{(n-1)}}& b^{(n)}\\
0& \ast
\end{array} \right)\left(\begin{array}{cc}
\widetilde{C}_{n-1}(A^{(n-1)})^{-1}&0\\
0& 1
\end{array} \right)\\
=& \left(\begin{array}{cc}
\widetilde{C}_{n-1}(A^{(n-1)})e^{A_{n-1}^{(n-1)}}\widetilde{C}_{n-1}(A^{(n-1)})^{-1} \ &\widetilde{C}_{n-1}(A^{(n-1)})e^{\frac{1}{2}A_{n-1}^{(n-1)}}b^{(n)}\\
{b^{(n)}}^\dagger e^{\frac{1}{2}A_{n-1}^{(n-1)}}\widetilde{C}_{n-1}(A^{(n-1)})^{-1} \ & \ast
\end{array} \right)\\
=& \left(\begin{array}{cc}
S_{n-1,-}&0\\
{b^{(n)}}^\dagger e^{\frac{1}{2}A_{n-1}^{(n-1)}}(S_{n-1,+}\widetilde{C}_{n-1})^{-1} & \ast
\end{array} \right)\left(\begin{array}{cc}
S_{n-1,+}& S_{n-1,+} \widetilde{C}_{n-1}e^{-\frac{1}{2}A_{n-1}^{(n-1)}}b^{(n)} \\
0& \ast
\end{array} \right).\end{align*}
Here in the last equality, we use again the monodromy relation
\[\widetilde{C}_{n-1}(A^{(n-1)}) e^{A_{n-1}^{(n-1)}} \widetilde{C}_{n-1}(A^{(n-1)})^{-1}=S_{n-1,-}\left(A_{n-2}^{(n-1)}\right) S_{n-1,+}\left(A_{n-2}^{(n-1)}\right).\]
Then the proof is finished for $k=n-1$. 
\qed

\subsubsection{The proof of Theorem \ref{mainthm}} 
\begin{proof}[The proof of Theorem \ref{mainthm}]
Note that $(S_+)_{k,k+1}$ is the $k$-th entry of the column vector $b_{k+1}(A)$. Lemma \ref{mono3} says that
\[b_{k+1}(A)=S_{k+}(A^{(k)}_{k-1})\widetilde{C}_{k}(A^{(k)})e^{-\frac{A_{k}^{(k)}}{2}}b^{(k+1)}(A_{k}).\]
On the one hand, the $k$-th row of $S_{k+}(A^{(k)}_{k-1})$ is simply $(0,...,0,e^{\frac{1}{2}A_{kk}})$. Here $A_{kk}$ is the $(k,k)$-entry of $A$. On the other hand, by the definition of $\widetilde{C}_{k}$ for $k=1,...,n$, we have \[\widetilde{C}_{k}(A^{(k)})=\left(\begin{array}{cc}
\widetilde{C}_{k-1}(A^{(k-1)})&0\\
0& 1
\end{array} \right)\cdot \widetilde{C}(E_k, A^{(k)}).\] 
Thus the $k$-th row of $S_{k+}(A^{(k)}_{k-1})\widetilde{C}_{k}(A^{(k)})$ is just the scalar multiplication of $e^{\frac{1}{2}A_{kk}}$ with the $k$-th row of the normalized connection matrix $\widetilde{C}(E_k, A^{(k)})$. 

Therefore, the $k$-th entry of the column vector $b_{k+1}(A)$ is given by multiplying the $k$-th row of (the scalar product with matrix) $e^{\frac{1}{2}A_{kk}}\cdot \widetilde{C}(E_k, A^{(k)})$ with the column vector $e^{\frac{-1}{2}A_{k}^{(k)}}b^{(k+1)}(A_{k})$. By Proposition \ref{explicitmCS}, that is
\begin{equation*}
(S_+)_{k,k+1}=e^{\frac{(A)_{k,k}+(A)_{k+1,k+1}}{4}}\sum_{i=1}^k\left(\frac{\prod_{l=1,l\ne i}^{k}\Gamma(1+\frac{\lambda^{(k)}_l-\lambda^{(k)}_i}{2\pi \I })}{\prod_{l=1}^{k+1}\Gamma(1+\frac{\lambda^{(k+1)}_l-\lambda^{(k)}_i}{2\pi \I })}\frac{\prod_{l=1,l\ne i}^{k}\Gamma(1+\frac{\lambda^{(k)}_l-\lambda^{(k)}_i}{2\pi \I })}{\prod_{l=1}^{k-1}\Gamma(1+\frac{\lambda^{(k-1)}_l-\lambda^{(k)}_i}{2\pi \I })}\right)\cdot \frac{a^{(k)}_i({A})}{N_i^{(k)}({A})}.
\end{equation*}
Then the expression of $(S_+)_{k,k+1}$ in the theorem follows immediately from the identity given in Lemma \ref{aNrelation}
 \[ \frac{a^{(k)}_i({A})}{N_i^{(k)}({A})}=2\pi\I \frac{\Delta^{1,...,k}_{1,...,k-1,k+1}\left(\frac{A-\lambda^{(k)}_i}{2\pi\I }\right)}{\prod_{l=1,l\ne i}^k(\frac{\lambda^{(k)}_i-\lambda^{(k)}_l}{2\pi \I })}\]
and the fact \[\prod_{l=1,l\ne i}^k\Gamma(1+\frac{\lambda^{(k)}_l-\lambda^{(k)}_i}{2\pi \I })=\prod_{l=1,l\ne i}^k\Gamma(\frac{\lambda^{(k)}_l-\lambda^{(k)}_i}{2\pi \I })\cdot \prod_{l=1,l\ne i}^k(\frac{\lambda^{(k)}_i-\lambda^{(k)}_l}{2\pi\I }).\]
The expression of $(S_-)_{k+1,k}$ follows in a similar way.
\end{proof}

\begin{rmk}\label{smoothex}
Since each $C(E_k,\delta_k(A))$ smoothly depend on $A$, the Stokes matrices $S_\pm(u_{\rm cat},A)$ are smoothly defined on $A\in \Herm(n)$. However, in Definition \ref{defitildeC} we introduce the normalized connection matrix $\widetilde{C}(E_k,\delta_k(A))$, associated to a systematically defined chain of matrices $P_k(A)$ diagonalizing $A^{(k)}$. The matrix $\widetilde{C}(E_k,\delta_k(A))$ is only defined for $A\in\Herm_0(n)$, simply because $P_k$ is only defined there (a family of matrices $P(t)$, which diagonalizes a family of $A(t)$, can be singular when $t$ approaches to the point where eigenvalues of $A(t)$ coincides). In this way, we only derive the explicit expression of the Stokes matrices on $\Herm_0(n)\subset\Herm(n)$. But one checks that the expression smoothly (actually real analytic) extends to $\Herm(n)$. See also remark \ref{clsdiag}. Thus in the above computation, we can do the diagonalization in stages completely formally, i.e., ignore the issue if $P_k(A)$ is invertible. 

We also remark that it is interesting to study the behaviour of the expression in Theorem \ref{mainthm} under the toric degeneration of Gelfand-Tsetlin systems, see e.g., \cite{NNU}.
\end{rmk}
\begin{rmk}
Apart from Poncar\'{e} rank 1 case, similar results as in Theorem \ref{isomonopro} and Theorem \ref{mainthm} for other linear differential equations have been studied in the literature of Painlev$\rm{\acute{e}}$ transcendents (see \cite{FIKN} for the history), the expression of Stokes matrices of some linear differential equations in terms of the asymptotics of the solutions of the associated nonlinear isomonodromy equations has been a major tool in the analysis of Painlev$\rm{\acute{e}}$ transcendents, see \cite[Part II and Part III]{FIKN} for the case
of Painlev$\rm{\acute{e}}$ II and III. Also, in \cite{GIL}, the expressions of Stokes matrices of certain differential equations with two irregular singularities were given via the asymptotics of the solutions of the associated
isomonodromy equations (which is in this case the $tt^*$
equations of Cecotti and Vafa).
\end{rmk}
\subsection{The diffeomorphism $\Phi_{u_{\rm cat}}(u)$ and the existence part of Theorem \ref{isomonopro}}\label{diffPhi}
For any fixed $u\in\h_{\rm reg}(\mathbb{R})$, the Riemann-Hilbert-Birkhoff map (also known as the dual exponential map) is
\begin{eqnarray}
\nu(u)\colon\Herm(n)\cong\frak u(n)^*\rightarrow \Herm^+(n)\cong{\rm U}(n)^*; \ A\mapsto S_-(u,A) S_+(u,A).
\end{eqnarray}
It follows from \cite[Theorem 2]{Boalch1} that $\nu(u)$ is a diffeomorphism (a real analytic map). Similarly, 
\begin{defi}\label{RHBcat}
The Riemann-Hilbert-Birkhoff map at $u_{\rm cat}$ is
\begin{eqnarray}
\nu(u_{\rm cat})\colon\Herm(n)\rightarrow \Herm^+(n); \ A\mapsto S_-(u_{\rm cat},A) S_+(u_{\rm cat},A).
\end{eqnarray}
\end{defi}
Recall that Theorem \ref{eBoalchthm} states that the map $\nu(u_{\rm cat})$ is a Poisson diffeomorphism.
We leave a proof of Theorem \ref{eBoalchthm} to the end of Section \ref{rRHvsGZ}. As a corollary, we get the existence part of Theorem \ref{isomonopro}. 

\begin{pro}\label{inversesol}
Given any $A\in\Herm(n)$, there exists a unique solution $\Phi(u;A)$ of the isomonodromy equation \eqref{isoeq} on $U_{\rm id}$ with the boundary value $A$ at $u_{\rm cat}$.
\end{pro}
\pf Since $\nu(u)$ and $\nu(u_{\rm cat})$ are diffeomorphisms, for any fixed $A\in\Herm(n)$, there exists a unique function $Y(u)\in\Herm(n)$ of $u\in U_{\rm id}$ such that $S_{\pm}(u, Y(u))=S_\pm(u_{\rm cat},A)$. Because $S_\pm(u_{\rm cat},A)$ is independent of $u$, the function $Y(u)$ is a solution of isomonodromy equation \eqref{isoeq} on $U_{\rm id}$. By Theorem \ref{isomonopro} and Corollary \ref{uucat}, there exists a constant $A'\in \Herm(n)$, as the boundary value of the solution $Y(u)$ at $u_{\rm cat}$ (i.e., as $\frac{u_{k+1}-u_{k}}{u_{k}-u_{k-1}}\rightarrow +\infty$ and $u_2-u_1\rightarrow 0$), such that  
\[S_\pm(u_{\rm cat},A)=S_{\pm}(u, Y(u))=S_\pm(u_{\rm cat},A').\]
Since $\nu(u_{\rm cat})$ is a diffeomorphism, 
it follows that $A=A'$, i.e., the boundary value of the solution $Y(u)$ at $u_{\rm cat}$ is $A$.
\qed
\begin{defi}\label{phiu}
We introduce the map 
\[\Phi_{u_{\rm cat}}(u):\Herm(n)\rightarrow \Herm(n)~;~A\mapsto \Phi(u;A),\]
where $\Phi(u;A)$ is the solution of \eqref{isoeq} with the boundary value $A$ at $u_{\rm cat}$. Then by Corollary \ref{uucat} it is just the diffeomorphism such that \[\nu(u)\circ \Phi_{u_{\rm cat}}(u)=\nu(u_{\rm cat}).\]
\end{defi}

\subsection{The leading terms of Stokes matrices via the Gelfand-Tsetlin systems}\label{leadterm}

This subsection is a proof of Proposition \ref{StokesinGZ}. For any $B\in\Herm(n)$, denote by the unitary matrix
\begin{equation}\label{gu}
g(u;B):=\overrightarrow{\prod_{k=1,...,n-1}}(z_k)^{\frac{\delta_k(B)}{2\pi\I }}.
\end{equation}
Here recall that $z_k$ are the coordinates \begin{align*}
z_0=u_1+\cdots+u_n, \ z_1=u_2-u_1, \ z_2=\frac{u_3-u_2}{u_2-u_1}, \ ...... , \ z_{n-1}=\frac{u_n-u_{n-1}}{u_{n-1}-u_{n-2}}.
\end{align*}
Since the Stokes matrices do not depend on $z_0$, we can assume $z_0=0$ and thus $u\in\frak t_{\rm reg}(\mathbb{R}).$
\begin{pro}\label{leadingterm}
For any fixed $A\in\Herm(n)$, we have that as $u\in U_{\rm id}$ and $u\rightarrow u_{\rm cat}$, i.e., $z_{1}\rightarrow 0+, z_k\rightarrow+\infty$,
\begin{align}\label{relimitGZ}
S_\pm(u,A)=S_\pm\left(u_{\rm cat}, \hspace{2mm} g(u;-A)\cdot A\cdot g(u;-A)^{-1}+\sum_{k=2}^{n-1}\mathcal{O}(z_k^{-1})\right),
\end{align}
\end{pro}
\pf 
Suppose we are given the solution $\Phi(u;\Phi_0)$ of the isomonodromy equation \eqref{isoeq} on $U_{\rm id}$ with the boundary value $\Phi_0\in\Herm(n)$. Then it follows from the definition of Stokes matrices at $u_{\rm cat}$ that
\begin{align}\label{leadinginv}
S_\pm(u_{\rm cat},\Phi_0)=S_\pm(u, \Phi(u;\Phi_0)).
\end{align}
If we set \[A:=\Phi(u;\Phi_0)=g(u;\Phi_0)\cdot \Phi_0\cdot g(u;\Phi_0)^{-1}+\sum_{k=2}^{n-1}\mathcal{O}(z_k^{-1})\] for any $u\in U_{\rm id}$ close to $u_{\rm cat}$, then following Proposition \ref{reverseasym}
\begin{equation}\label{PhiA}
    \Phi_0= g(u;-A)\cdot A\cdot g(u;-A)^{-1}+\sum_{k=2}^{n-1}\mathcal{O}(z_k^{-1}).
\end{equation}
Therefore, setting $A=\Phi(u;\Phi_0)$ and applying \eqref{PhiA} to the identity \eqref{leadinginv} leading to \eqref{relimitGZ}.
\qed 

\vspace{2mm}
{\bf Proof of Proposition \ref{StokesinGZ} }: On the one hand, following Theorem \ref{mainthm} we have explicit expression of $S_\pm(u_{\rm cat}, \Phi_0)$ for any $\Phi_0\in \Herm(n)$. On the other hand, following Proposition \ref{leadingterm}, we have
\begin{align}\nonumber
    S_\pm\left(u,~A\right)
&=S_\pm\left(u_{\rm cat},~g(u;-A)\cdot A\cdot g(u;-A)^{-1}+\sum_{k=2}^{n-1}\mathcal{O}(z_k^{-1})\right)\\ \label{transinvers}
&= S_\pm\left(u_{\rm cat},~g(u;-A)\cdot A\cdot g(u;-A)^{-1}\right)+\sum_{k=2}^{n-1}\mathcal{O}(z_k^{-1}).
\end{align}
Here the second identity uses the explicit expression of the Stokes matrices at $u_{\rm cat}$.
Therefore, we only need to find the expression of the eignvalues $\lambda^{(k)}_i$ and minors $\Delta^{1,...,k-1,k}_{1,...,k-1,k+1}$ of $g(u;-A)\cdot A\cdot g(u;-A)^{-1}$ in terms of those of $A$. The conjugation action by $g(u;-A)$ on $A$ is better to be understood in terms of the Thimm action.

First, under the Gelfand-Tsetlin action and angle coordinates, the Thimm action (see Section \ref{GTcoor}) of an element \[\theta^{(k)}={\rm diag}(e^{\I \theta^{(k)}_1},...,e^{\I \theta^{(k)}_k})\in T(k), \text{ for } k=1,....,n-1\] on $A\in\Herm(n)$ is described by
\begin{align}
\lambda^{(i)}_j(\theta^{(k)}\bullet A)&=\lambda^{(i)}_j(A), \ \ 1\le j\le i\le n,\\ \label{torusGT1} \psi^{(i)}_j(\theta^{(k)}\bullet A)&=\psi^{(i)}_j(A)+\delta_{ki}\theta^{(k)}_{j}, \ \ 1\le j\le i\le n-1.
\end{align}
Second, the conjugation by a diagonal element \[d={\rm diag}(d_1,...,d_n)\in T(n)\subset U(n)\] is described by
\begin{equation}\label{torusGT}
\lambda^{(i)}_j({\rm Ad}_dA)= \lambda^{(i)}_j(A), \hspace{5mm} \alpha^{(i)}_j({\rm Ad}_dA)=d_i\cdot \alpha^{(i)}_j \cdot d_{i+1}^{-1}.
\end{equation}

For any $A$ and $u$, let us introduce an element in the product of torus
\begin{equation}\label{toriele}
(u_2-u_1)^{\frac{-\lambda^{(1)}(A)}{2\pi\I }}\times \overrightarrow{\underset{k=2,...,n-1}{\prod} }\left(\frac{u_{k+1}-u_{k}}{u_{k}-u_{k-1}}\right)^{\frac{-\lambda^{(k)}(A)}{2\pi\I }}\in T(1)\times \cdots \times T(n-1)\end{equation}
where $\lambda^{(k)}(A):={\rm diag}(\lambda^{(k)}_1,...,\lambda^{(k)}_{k})$. Let us introduce the diagonal matrix 
\begin{equation}\label{diagele}D(u;A)={\rm diag}(1,(u_2-u_1)^{-{A_{22}}/{2\pi\I }},...,(u_n-u_{n-1})^{-{A_{nn}}/{2\pi\I }})\in T(n).
\end{equation}
Then one checks 
\begin{align*}
&{\rm Ad}_{g(u;-A)}A\\
=&{\rm Ad}_{D(u;A)}\left((u_2-u_1)^{\frac{-\lambda^{(1)}(A)}{2\pi\I }}\times \overrightarrow{\underset{k=2,...,n-1}{\prod} }\left(\frac{u_{k+1}-u_{k}}{u_{k}-u_{k-1}}\right)^{\frac{-\lambda^{(k)}(A)}{2\pi\I }}\bullet A\right).
\end{align*}
Here $\bullet$ denotes the Thimm action of the element in \eqref{toriele} on $A$. Therefore, by \eqref{torusGT} and \eqref{torusGT1} we have
\begin{align*}
\lambda^{(k)}_i({\rm Ad}_{g(u;-A)}A)&=\lambda^{(k)}_i(A),\\
a^{(k)}_i({\rm Ad}_{g(u;-A)}A)&=a^{(k)}_i(A)\cdot (u_k-u_{k-1})^{\frac{-A_{kk}}{2\pi\I }}(u_{k+1}-u_{k})^{\frac{A_{k+1,k+1}}{2\pi\I }} \left(\frac{u_{k}-u_{k-1}}{u_{k+1}-u_k}\right)^{\frac{\lambda^{(k)}_i}{2\pi\I }}\\
&=a^{(k)}_i(A)\cdot (u_k-u_{k-1})^{\frac{\lambda^{(k)}_i(A)-A_{kk}}{2\pi\I }} (u_k-u_{k-1})^{\frac{A_{k+1,k+1}-\lambda^{(k)}_i(A)}{2\pi\I }},
\end{align*}
where the $(u_k-u_{k-1})^{\frac{-A_{kk}}{2\pi\I }}(u_{k+1}-u_{k})^{\frac{A_{k+1,k+1}}{2\pi\I }}$ term in the second identity comes from the conjugation action of the diagonal matrix $D(u;A)$. 

In the end, the proposition follows from the real analyticity and the explicit expression of $S_\pm(u_{\rm cat}, \Phi_0)$ in Theorem \ref{mainthm} (provided replacing the Gelfand-Tsetlin coordinates of $\Phi_0={\rm Ad}_{g(u;-A)}A$ by those of $A$). \qed

\vspace{2mm}
One can prove similar results for the other entries of $S_\pm(u,A)$. Proposition \ref{StokesinGZ}  gives the expression of the leading terms of Stokes matrices as $u\rightarrow u_{\rm cat}$ from $U_{\rm id}$, in terms of the Gelfand-Tsetlin action and angle coordinates. Furthermore, we see that the leading terms include a fast spin on the corresponding Liouville torus of the Gelfand-Tsetlin integrable systems, which is cancelled out as far as the regularized limit is considered.  

\subsection{Proof of Theorem \ref{introcor}: regularized limits of Stokes matrices at a caterpillar point}\label{pf:introthm3}
\begin{proof}[Proof of Theorem \ref{introcor}]
A manipulation of Gauss decomposition and monodromy relation \eqref{Ckmonodromy} shows the following equivariant property of Riemann-Hilbert-Birkhoff map $\nu(u_{\rm cat}): {\rm Herm}(n)\rightarrow {\rm Herm}^+(n)$ at $u_{\rm cat}$
(see Proposition \ref{ConnGZ} for a stronger equivariant property of $\nu(u_{\rm cat})$ with respect to the Thimm action),
\begin{align*}
&\nu\left(u_{\rm cat},g(u;-A)\cdot A\cdot g(u;-A)^{-1}\right)
\\
=&D(u;A)\cdot \left((u_2-u_1)^{\frac{-\lambda^{(k)}(A)}{2\pi\I }}\times \overrightarrow{\underset{k=2,...,n-1}{\prod} }\left(\frac{u_{k+1}-u_{k}}{u_{k}-u_{k-1}}\right)^{\frac{-\lambda^{(k)}(A)}{2\pi\I }}\bullet \nu(u_{\rm cat},A) \right)\cdot D(u;A)^{-1}.
\end{align*}
Here we take the notation from \eqref{toriele} and \eqref{diagele} and in particular $\bullet$ denotes the Thimm action. It also follows from Proposition \ref{ConnGZ} that
\[\lambda^{(k)}_i(A)= {\rm log}\left(\lambda^{(k)}_i\left(\nu(u_{\rm cat},A)\right)\right), \text{ for } 1\le i\le k\le n,\]
here recall that the eigenvalues are ordered. Therefore, we get
\begin{equation}\label{deltaequiva}
\nu\left(u_{\rm cat},g(u;-A)A g(u;-A)^{-1}\right)=G\left(u;-\nu(u_{\rm cat},A)\right)\cdot \nu(u_{\rm cat},A) \cdot G\left(u;-\nu(u_{\rm cat},A)\right)^{-1},
\end{equation}
where for any positive definite Hermitian matrix $B\in\Herm^+(n)$, if $B=B_+^\dagger B_+$ for an upper triangular matrix, then
\[G(u;B):=\overrightarrow{\prod_{k=1,...,n-1}}(z_k)^{\frac{{\rm log}(\delta_k(B_-)\delta_k(B_+))}{2\pi\I }}.\]
Note that 
\begin{align*}
    S_\pm\left(u,~A\right)
=&S_\pm\left(u_{\rm cat},~g(u;-A)\cdot A\cdot g(u;-A)^{-1}+\sum_{k=2}^{n-1}\mathcal{O}(z_k^{-1})\right)\\
=& S_\pm\left(u_{\rm cat},~g(u;-A)\cdot A\cdot g(u;-A)^{-1}\right)+\sum_{k=2}^{n-1}\mathcal{O}(z_k^{-1}).
\end{align*}
Therefore, the identities \eqref{transinvers} and \eqref{deltaequiva} give
\begin{align*}
&\nu(u, A)
\\
=&
S_-(u,~A)S_+(u,~A)\\
=&S_-\left(u_{\rm cat},~g(u;-A)\cdot A\cdot g(u;-A)^{-1}\right)S_+\left(u_{\rm cat},~g(u;-A)\cdot A\cdot g(u;-A)^{-1}\right)+\sum_{k=2}^{n-1}\mathcal{O}(z_k^{-1})\\
=&\nu\left(u_{\rm cat},g(u;-A)\cdot A\cdot g(u;-A)^{-1}\right)+\sum_{k=2}^{n-1}\mathcal{O}(z_k^{-1})\\
=&G\left(u;-\nu(u_{\rm cat},A)\right)\cdot \nu(u_{\rm cat},A) \cdot G\left(u;-\nu(u_{\rm cat},A)\right)^{-1}+\sum_{k=2}^{n-1}\mathcal{O}(z_k^{-1}).
\end{align*}
Following the definition of $G(u;B)$, the above identity implies
\begin{equation*}
G(u;\nu(u,A))\cdot \nu\left(u, A\right) \cdot G(u;\nu(u,A))^{-1}= \nu(u_{\rm cat},A)+\sum_{k=2}^{n-1}\mathcal{O}(z_k^{-1}),
\end{equation*}
which proves the theorem provided writing the Riemann-Hilbert-Birkhoff maps $\nu(u,A)$, $\nu(u_{\rm cat},A)$ and the gauge transform $G(u;\nu(u,A) )$ by the Stokes matrices via the definitions. \end{proof}

\begin{ex}\label{rank2ex}[2 by 2 cases]
Let us consider the rank two case, that is
\begin{equation*}
\frac{dF}{dz}=\left(\left(\begin{array}{cc}
    \I u_1 & 0  \\
    0 & \I u_2
  \end{array}\right)
-\frac{1}{2\pi\I z}{\left(
  \begin{array}{cc}
    t_1 & a  \\
    \bar{a} & t_2
  \end{array}
\right)}\right)\cdot F.
\quad
\end{equation*} 
Following \rm[\cite{BJL} Proposition 8], the Stokes matrices (with respect to the chosen branch of ${\rm log}(z)$) are
\[S_-(u,A)=\left(
  \begin{array}{cc}
    e^{\frac{t_1}{2}} &  0  \\
    \frac{ \bar{a}\cdot  e^{\frac{t_1-t_2}{4}} ({u_2-u_1})^{\frac{t_1-t_2}{2\pi \I }}}{\Gamma(1-\frac{\lambda_1-t_1}{2\pi \I })\Gamma(1-\frac{\lambda_2-t_1}{2\pi \I })}  & e^{\frac{t_2}{2}}
  \end{array}\right), \ \ S_+(u,A)=\left(
  \begin{array}{cc}
    e^{\frac{t_1}{2}} &  \frac{ a\cdot e^{\frac{t_1+t_2}{4}} ({u_2-u_1})^{\frac{t_2-t_1}{2\pi \I }}}{\Gamma(1+\frac{\lambda_1-t_1}{2\pi \I })\Gamma(1+\frac{\lambda_2+t_1}{2\pi \I })}   \\
    0 & e^{\frac{t_2}{2}}
  \end{array}\right).\]
Here $\lambda_1,\lambda_2$ are eigenvalues of $\left(
  \begin{array}{cc}
    t_1 & a  \\
    \bar{a} & t_2
  \end{array}
\right)$. By the definition of $\delta_k$ given in \eqref{delta1}, we have that $\delta_1(S_\pm(u,A))$ is the diagonal part of $S_\pm(u,A)$. Then ${\frac{{\rm log}(\delta_1(S_-)\delta_1(S_+))}{2\pi\I }}={{\rm diag}(\frac{t_1}{2\pi\I },\frac{t_2}{2\pi\I })}$, and $z_1=u_2-u_1$, thus
\begin{align}\label{2by2exp}
{\left(\begin{array}{cc}
    (u_2-u_1)^{\frac{t_1}{2\pi\I }} & 0  \\
    0 & (u_2-u_1)^{\frac{t_2}{2\pi\I }}
  \end{array}\right)} S_-(u,A)S_+(u,A)\left(\begin{array}{cc}
    (u_2-u_1)^{\frac{-t_1}{2\pi\I }} & 0  \\
    0 & (u_2-u_1)^{\frac{-t_2}{2\pi\I }}
  \end{array}\right)
\end{align}
approaches to $S_-(u_{\rm cat},A)S_+(u_{\rm cat},A)$, as $u_2-u_1\rightarrow 0$, where
 \[S_-(u_{\rm cat},A)^\dagger=S_+(u_{\rm cat},A)=\left(
  \begin{array}{cc}
    e^{\frac{t_1}{2}} &  \frac{ a\cdot e^{\frac{t_1+t_2}{4}}}{\Gamma(1+\frac{\lambda_1-t_1}{2\pi \I })\Gamma(1+\frac{\lambda_2-t_1}{2\pi \I })}   \\
    0 & e^{\frac{t_2}{2}}
  \end{array}\right).\] 
Actually, as we have seen in the above computation or in the proof of Proposition \eqref{strongerpro}, the expression \eqref{2by2exp} equals to $S_-(u_{\rm cat},A)S_+(u_{\rm cat},A)$ (not necessary to take the limit as $z_1=u_2-u_1\rightarrow 0$). 
\end{ex}

\subsection{The explicit expression of the other entries of Stokes matrices}\label{classicalRLL}
Following Theorem \ref{eBoalchthm} (whose proof is given in Section \ref{rRHvsGZ}), the Riemann-Hilbert-Birkhoff map $\nu(u_{\rm cat})$ at the caterpillar point is a Poisson map. That is, in terms of the r-matrix formulation of the Poisson brackets on the dual Poisson Lie group, Theorem \ref{eBoalchthm} can be rewritten via the classical RLL formulation (see {\cite[Formula (235)]{AM1}}):
\begin{cor}\label{sclrll}
The	matrix valued function $\nu(u_{\rm cat})=S_-(u_{\rm cat},A) S_+(u_{\rm cat},A)\in \Herm^+(n)$ satisfies
\begin{equation}\label{StokesPoisson} 
		\{\nu^1 \otimes \nu^2 \} = r_+ \nu^1 \nu^2 + \nu^1 \nu^2 r_- - \nu^1 r_+ \nu^2 - \nu^2 r_- \nu^1.
	\end{equation}
	Here
	\begin{align*}
		r_+ & = \frac{1}{2}  \sum_{i = 1}^n E_{ii} \otimes E_{ii} + \sum_{1 \le
			i < j \le n} E_{ij} \otimes E_{ji} \in \mathrm{End} (\mathbb{C}^n) \otimes
		\mathrm{End} (\mathbb{C}^n), \\
		r_- & = - \frac{1}{2}  \sum_{i = 1}^n E_{ii} \otimes E_{ii} - \sum_{1
			\le j < i \le n} E_{ij} \otimes E_{ji} \in \mathrm{End} (\mathbb{C}^n)
		\otimes \mathrm{End} (\mathbb{C}^n), 
	\end{align*}
	and the tensor notation is used: both sides are $\mathrm{End} (\mathbb{C}^n)
	\otimes \mathrm{End} (\mathbb{C}^n)$ valued functions on
	$A\in\mathfrak{g}\mathfrak{l}_n$, $\nu^1 := \nu(u_{\rm cat}) \otimes 1$, $\nu^2 := 1
	\otimes \nu(u_{\rm cat})$, and the $ij, kl$ coefficient of the matrix $\{\nu^1, \nu^2 \}$ is
	defined as the Poisson bracket $\{\nu_{ij}, \nu_{kl} \}$ of the functions on the
	canonical linear Poisson space $\mathfrak{g}\mathfrak{l}_n^{\ast} \cong
	\mathfrak{g}\mathfrak{l}_n$.
\end{cor}
The other entries of $S_\pm(u_{\rm cat},A)$ are uniquely determined by the sub-diagonal ones through the identity \eqref{StokesPoisson}, and one can actually find their explicit expressions, as well as their leading terms via the Gelfand-Tsetlin systems like in Section \ref{leadterm}. Since the explicit expressions of the other entries are not used in this paper, we omit the computation.

\begin{rmk}
The quantum analog of Corollary \ref{sclrll}, or equivalently Theorem \ref{Boalchthm}, is just Theorem \ref{introthm1}. In particular, the quantum Stokes matrices of quantum confluent hypergeometric equation \eqref{introqeq} satisfy the quantum RLL relation in the Faddeev-Reshetikhin-Takhtajan \cite{FRT} realization of quantum groups, and then Corollary and Theorem \ref{Boalchthm} follows by taking a semiclassical limit \cite{TLXu}.
\end{rmk}

\subsection{Regularized limits of Stokes matrices with degenerate irregular terms}\label{dirreg}
Although the results in the previous sections are obtained under the assumption that $u\in\h_{\rm reg}(\mathbb{R})$, their generalization to the degenerate case, i.e., $u\in\h(\mathbb{R})\setminus\h_{\rm reg}(\mathbb{R})$, is direct and is given in this subsection.

Given any fixed partition $\underline{d}$ of $n$, i.e., a set of integers $\{d_j\}_{j=1,...,m}$ such that $n=d_1+\cdots+d_m$, let us consider the subspace of $\h(\mathbb{R})$
\[U^{\underline{d}}_{\rm id}=\{u={\rm diag}(\underbrace{u_1,...,u_1}_{d_1}, \underbrace{u_2,...,u_2}_{d_2},....,\underbrace{u_m,...,u_m}_{d_m})\in\h(\mathbb{R})~|~ u_i<u_j \ \ \text{if} \ \ i< j\}.\] Any $n\times n$ matrix $A$ can be seen as a block matrix $A=(A_{ij})_{1\le i,j\le m}$ according to the tuple $\underline{d}$, where each $A_{ij}$ denotes the $d_i\times d_j$ block/submatrix formed by the corresponding $d_i$ rows and $d_j$ columns of $A$. Then for each $k=0,...,m-1$, we denote by $\delta_k^{\underline{d}}(A)$ the blocked matrix
 \begin{equation}\label{bdelta}\delta_k^{\underline{d}}(A)_{ij}=\left\{
          \begin{array}{lr}
             A_{ij},   & if \ \ 1\le i, j\le k, \ \text{or} \ i=j  \\
           0, & \text{otherwise}.
             \end{array}
\right.\end{equation}
\begin{rmk}
The identity \eqref{bdelta} depends on the partition $d_1,...,d_m$, which generalizes the definition of $\delta_k(A)$ given in \eqref{delta1} corresponding to the case $d_1=\cdots=d_n=1$. Although this more general notation can lead to confusion with \eqref{delta1}, the context typically eliminates any ambiguity. 
\end{rmk}

Accordingly, we consider the $n\times n$ system of partial differential equations for a function $F(z,u)\in GL(n)$ of $(z,u)\in \mathbb{C}\times U^{\underline{d}}_{\rm id}$,
\begin{align}\label{bisoStokeseq1}
\frac{\partial F}{\partial z}&=\left(\I u-\frac{1}{2\pi\I }\frac{\Phi(u)}{z}\right)\cdot F,\\
\label{bisoStokeseq2}
\frac{\partial F}{\partial u_k}&=\left(\I E(\underline{d})_k z-\frac{1}{2\pi\I } {\rm ad}^{-1}_u{\rm ad}_{E(\underline{d})_k}\Phi(u)\right)\cdot F, \ \text{for all} \ k=1,...,m,
\end{align}
where the residue $\Phi(u): U^{\underline{d}}_{\rm id}\rightarrow {\rm Herm}(n)$ is a solution of the blocked isomonodromy deformation equation
\begin{equation}\label{bisoeq}
\frac{\partial \Phi}{\partial u_k} =\frac{1}{2\pi\I }[\Phi,{\rm ad}^{-1}_u{\rm ad}_{E(\underline{d})_k}\Phi], \ \text{for all} \ k=1,...,m.\end{equation} 
Here
\begin{equation}\label{bEk}
    E(\underline{d})_k={\rm diag}(0,...,0,{\rm Id}_{d_k},0,...,0)
\end{equation}
is the block diagonal matrix whose $k-$th block diagonal entry is the rank $d_k$ identity matrix ${\rm Id}_{d_k}$. Note that ${\rm ad}_{E(\underline{d})_k}\Phi$ takes values in the space ${\frak {gl}}_n^{od}(\underline{d})$ of block off diagonal matrices with respect to the tuple $(\underline{d})$, and that ${\rm ad}_u$
is invertible when restricted to the subspace ${\frak {gl}}_n^{od}(\underline{d})\subset {\frak {gl}}_n$. Note that \eqref{bisoeq} is the compatibility condition of the systems \eqref{bisoStokeseq1} and \eqref{bisoStokeseq2}.

For any fixed $u\in\h_{\rm reg}(\mathbb{R})$, since $\Phi(u)/2\pi\I $ is skew-Hermitian matrix, the ordinary differential equation \eqref{introisoStokeseq1} is nonresonant, and thus has a unique formal solution $\hat{F}$ around
$z = \infty$. Similar to Theorem \ref{uniformresum}, the Borel-Laplace transform of the unique formal solution $\hat{F}(z)$ produces actual fundamental solutions $F_+(z), F_-(z)$ of equation \eqref{introisoStokeseq1} with prescribed asymptotics within the sectors $S(-\pi,\pi)$ and $S(-2\pi,0)$ respectively. Then the transition matrices between the two actual solutions are a pair of Stokes matrices $S_\pm(u,\Phi(u))\in{\rm GL}(n)$. Similar to Section \ref{secSC}, one also defines the connection matrix $C(u,\Phi(u))$ of equation \eqref{introisoStokeseq1}.

Furthermore, for any solution $\Phi(u)$ of the equation \eqref{bisoeq}, the canonical solutions $F_\pm(z;u,A)$ of \eqref{introisoStokeseq1} automatically satisfy the equation \eqref{introisoStokeseq2}. Thus the Stokes matrices $S_\pm(\Phi(u),u)$ of \eqref{bisoStokeseq1}
are locally constants (independent of $u\in U^{\underline{d}}_{\rm id}$).

Now in a similar way, we can prove results analog to Proposition \ref{strongerpro}, and an analog of Theorem \ref{isomonopro}.
\begin{thm}
For any solution $\Phi(u)$ of the equation \eqref{bisoeq} on $U^{\underline{d}}_{\rm id}$, there exists a unique constant $\Phi_0\in\Herm(n)$ such that as the real numbers $\frac{u_{k+1}-u_{k}}{u_{k}-u_{k-1}}\rightarrow +\infty$ for all $k=2,...,m-1$ and $u_2-u_1\rightarrow 0$,
\begin{equation*}
\Phi(u)={\rm Ad}{\left((u_2-u_1)^{\frac{\delta^{\underline{d}}_1(\Phi_0)}{2\pi\I }}\overrightarrow{\underset{k=2,...,m-1}{\prod} }\left(\frac{u_{k+1}-u_{k}}{u_{k}-u_{k-1}}\right)^{\frac{\delta^{\underline{d}}_k(\Phi_0)}{2\pi\I }}\right)}\Phi_0 +\sum_{k=2}^{m-1}\mathcal{O}\left(\left(\frac{u_{k+1}-u_{k}}{u_{k}-u_{k-1}}\right)^{-1}\right).
\end{equation*}
Furthermore, given any $\Phi_0\in\Herm(n)$ there exists a unique real analytic solution $\Phi(u)$ of \eqref{introisoeq} with the prescribed asymptotics \eqref{firstasy}.
\end{thm}
Similarly, we can prove a result analog to Theorem \ref{relautomorphism} for the systems \eqref{bisoStokeseq1} and \eqref{bisoStokeseq2}. Then the same argument as in Section \ref{pf:introthm3} leads to the following analog of Theorem \ref{introcor} for the degenerate $u$ case.
\begin{thm}\label{dreglimit}
For any fixed $A\in\Herm(n)$, the limit of the matrix valued function of $u\in U^{\underline{d}}_{\rm id}$
\begin{equation}
{\rm Ad}{\left((\frac{1}{u_{2}-u_{1}})^{\frac{{\rm log}(\delta^{\underline{d}}_1(S_-)\delta^{\underline{d}}_1(S_+))}{2\pi\I }}\cdot\overrightarrow{\underset{k=2,...,m-1}{\prod} }(\frac{u_{k}-u_{k-1}}{u_{k+1}-u_{k}})^{\frac{{\rm log}(\delta^{\underline{d}}_k(S_-)\delta^{\underline{d}}_k(S_+))}{2\pi\I }}\right)}  \left(S_-(u,A)S_+(u,A)\right), 
\end{equation}
as $\frac{u_{k+1}-u_{k}}{u_{k}-u_{k-1}}\rightarrow +\infty$ for all $k=2,...,m-1$ and $u_2-u_1\rightarrow 0$, equals to
\begin{equation}
\left(\overrightarrow{\prod_{k=1,...,m}}C\left(E(\underline{d})_k,\delta^{\underline{d}}_k(A)\right)\right)\cdot e^A\cdot \left(\overrightarrow{\prod_{k=1,...,m}}C\left(E(\underline{d})_k,\delta^{\underline{d}}_k(A)\right)\right)^{-1}.
\end{equation}
Here each $C\left(E(\underline{d})_k,\delta^{\underline{d}}_k(A)\right)$ denotes the connection matrix of the system
\begin{equation}\label{eq:relStokes}
\frac{dF}{dz}=\left(\I E(\underline{d})_k-\frac{1}{2\pi \I }\frac{\delta^{\underline{d}}_{k}(A)}{z}\right)F,
\end{equation}
where recall $\delta^{\underline{d}}_k(A)$ and $E(\underline{d})_k$ are defined in \eqref{bdelta} and \eqref{bEk}.
\end{thm}

This theorem will be used to derive the expression of quantum Stokes matrices at a caterpillar point $u_{\rm cat}$, see Section \ref{sec:qStokescat}.

\section{Applications in Poisson geometry}\label{appPoisson}
In this section, we prove Theorem \ref{eBoalchthm}, and give a new proof of Theorem \ref{Boalchthm}. Section \ref{rRHvsGZ} shows that the Riemann-Hilbert-Birkhoff maps at caterpillar points intertwine Gelfand-Tsetlin systems and their multiplicative analogs. Then Section \ref{endsection} shows the Poisson geometric nature of taking closure of Stokes matrices, and provides a new proof to Theorem \ref{Boalchthm}.

\subsection{Riemann-Hilbert-Birkhoff map at the caterpillar point is compatible with Gelfand-Tsetlin}\label{rRHvsGZ}
Recall that in Section \ref{GTcoor}, we have introduced the Gelfand-Tsetlin maps and Thimm torus actions. In this subsection, we first introduce their multiplicative analogs \cite{FR} on $\Herm^+(n)$, and then prove that the map $\nu(u_{\rm cat})$ is compatible with them.

{\bf Logarithmic Gelfand-Tsetlin maps.} 
Let $\Herm^+(n)\subset \Herm(n)$ denote the subset of positive definite
Hermitian matrices, and define a logarithmic Gelfand-Tsetlin map 
\begin{equation}\label{eq:logmomentmap}
\mu\colon \Herm^+(n)\to \mathbb{R}^{\frac{n(n+1)}{2}},
\end{equation}
taking $A$ to the collection of numbers
$\mu^{(k)}_i(A)=\log(\lambda^{(k)}_i(A))$. Here recall that $\lambda^{(k)}_i(A)'s$ are the ordered set of
eigenvalues of the upper left $k$-th submatrix $A^{(k)}$ of $A\in \Herm(n)$. Then $\mu$ is a continuous
map from $\Herm^+(n)$ onto the Gelfand-Tsetlin cone $\mathcal{C}(n)$.

{\bf Thimm torus actions.}
Let $\mathcal{C}_0(n)\subset \mathcal{C}(n)$ denote the subset where all of the eigenvalue inequalities
\eqref{eq:cone} are strict. Let $\Herm_0^+(n)$ denote the intersection of $\Herm_0(n)$ and $\Herm^+(n)$, i.e., $\Herm_0^+(n)=\mu^{-1}(\mathcal{C}_0(n))$. Then the actions of Thimm torus on 
$\Herm_0(n)$ defined in \eqref{eq:taction} restrict to a torus action on $\Herm_0^+(n)$. The action preserves the logarithmic Gelfand-Tsetlin map $\mu$. 

Recall from Definition \ref{RHBcat} that the Riemann-Hilbert-Birkhoff map at $u_{\rm cat}$ is
\[\nu(u_{\rm cat})\colon\Herm(n)\rightarrow \Herm^+(n); \ A\mapsto S_-(u_{\rm cat},A) S_+(u_{\rm cat},A).\]
\begin{pro}\label{ConnGZ}
The Riemann-Hilbert-Birkhoff map 
\[\nu(u_{\rm cat}): {\rm Herm}(n)\rightarrow {\rm Herm}^+(n)\] 
is a diffeomorphism compatible with the Gelfand-Tsetlin systems. That is
\begin{itemize}
\item[(a).] $\nu(u_{\rm cat})$ intertwines the
  Gelfand-Tsetlin maps: $\mu\circ \nu(u_{\rm cat})=\lambda$;
\item[(b).] $\nu(u_{\rm cat})$ intertwines the Thimm's torus actions on $\Herm_0(n)$ and $\Herm_0^+(n)$. 
\end{itemize}
\end{pro}
The proof of the proposition relies on the following linear algebra fact.
\begin{lem}\label{lem:GWiso}
For each $0<k\le n$, let $\Psi^{(k)}:\Herm(k)\rightarrow {\rm U}(k)$ be a smooth map satisfying the conditions
\begin{enumerate}
\item[(1).]\label{it:a}  $\Psi^{(k)}$ is a ${\rm U}(k-1)$-equivariant map, i.e., $\Psi^{(k)}(gAg^{-1})={\rm Ad}_g\Psi^{(k)}(A)$, for any $g\in {\rm U}(k-1)\subset U(k)$;
\item[(2).]\label{it:b} for any $A\in \Herm(k)$, there is a block decomposition of $\Psi^{(k)}(A) e^A \Psi^{(k)}(A)^{-1}$ taking the form \[\Psi^{(k)}(A) e^A \Psi^{(k)}(A)^{-1}=\left(\begin{array}{cc}
{\rm Id_{k-1}}&0\\
{B^{(k)}}^\dagger & 1
\end{array} \right)\left(\begin{array}{cc}
e^{A^{(k-1)}} & 0\\
0& \star
\end{array} \right)\left(\begin{array}{cc}
{\rm Id_{k-1}} & B^{(k)}\\
0& 1
\end{array} \right),\]
where $B^{(k)}$ is a column with $k-1$ elements, and ${B_+^{(k)}}^\dagger$ the conjugate transpose.
\end{enumerate}
Let us extend $\Psi^{(k)}:\Herm(k)\rightarrow {\rm U}(k)$ as a map from $\Herm(n)$ to ${\rm U}(n)$, using the projection of $\Herm(n)$ onto $\Herm(k)$ and the natural inclusion ${\rm U}(k)\subset {\rm U}(n)$.
Let $\Psi:=\Psi^{(1)}\cdot\cdot\cdot \Psi^{(n)}$ be the map from $\Herm(n)$ to ${\rm U}(n)$ given by the pointwise multiplication, then the map
\[ \Gamma_{\Psi}:={\rm Ad}_\Psi\circ {\rm exp}\colon \Herm(n)\to
  \Herm^+(n); \ A\mapsto \Psi(A)e^A\Psi(A)^{-1}\]
is a diffeomorphism compatible with the Gelfand-Tsetlin systems.
\end{lem}
\pf We will prove this theorem inductively on $n$. When $n=1$, the result is obvious. 
For the inductive step $n>1$, we first assume that the map $\Psi_{n-1}:=\Psi^{(1)}\cdot\cdot\cdot \Psi^{(n-1)}\colon\Herm(n-1)\rightarrow\Herm^+(n-1)$ is such that ${\rm Ad}_{\Psi_{n-1}}\circ {\rm exp}$ satisfies the conditions $(a)$ and $(b)$ (provided replacing $\nu(u_{\rm cat})$ by $\Gamma_\Psi$). 

Now we try to prove that the map ${\rm Ad}_{\Psi_{n-1} \Psi^{(n)}} \circ {\rm exp}$ also satisfy $(a)$ and $(b)$. First, using the identity in the assumption $(2)$
$$\Psi^{(n)}(A) e^A \Psi^{(n)}(A)^{-1}=\left(\begin{array}{cc}
{\rm Id}_{n-1}&0\\
{B^{(n)}}^\dagger & 1
\end{array} \right)\left(\begin{array}{cc}
e^{A^{(n-1)}} & 0\\
0& \star
\end{array} \right)\left(\begin{array}{cc}
{\rm Id}_{n-1} & B^{(n)}\\
0& 1
\end{array} \right),$$
we have
\begin{align*}
{\rm Ad}_{\Psi_{n-1}(A)\Psi^{(n)}(A)}e^A=\left(\begin{array}{cc}
\Psi_{n-1} e^{A^{(n-1)}}\Psi_{n-1}^{-1} & \Psi_{n-1} e^{A^{(n-1)}}B^{(n)}\\
{B^{(n)}}^\dagger e^{A^{(n-1)}}\Psi_{n-1}^{-1}& \star
\end{array} \right).
\end{align*}
Hence by the assumption about the map ${\rm Ad}_{\Psi_{n-1}}\circ {\rm exp}$, we observe that ${\rm Ad}_{\Psi_{n-1} \Psi^{(n)}}\circ{\rm exp}$ intertwines the
Gelfand-Tsetlin maps.

For the Thimm's torus action, by the assumption for $n-1$, the map ${\rm Ad}_{\Psi_{n-1}}\circ{\rm exp}$
intertwines the $T(k)\subset U(k)$ actions on 
$\Herm_0(n-1)$ and $\Herm^+_0(n-1)$ for $1\le k<n-1$. That is
\begin{equation}\label{inductiveaction}\Psi_{n-1}(t\bullet A)e^{t\bullet A}\Psi_{n-1}(t\bullet A)^{-1}=t\bullet \left(\Psi_{n-1}(A) e^A\Psi_{n-1}(A)^{-1}\right).
\end{equation}
Here recall that the $k$-torus 
$T(k)\subset {\rm U}(k)\subset U(n)$ acts on 
$\Herm_0(n)$ as in \eqref{eq:taction}.
Furthermore, using the assumption $(1)$ on the equivariance of $\Psi^{(n)}$, we have 
$${\rm Ad}_{\Psi_{n-1}(A)\Psi^{(n)}(A)}e^A={\rm Ad}_{\Psi^{(n)}(\Psi_{n-1}(A)A\Psi_{n-1}(A)^{-1})}e^{\Psi_{n-1}(A)A\Psi_{n-1}(A)^{-1}}.$$
Together with \eqref{inductiveaction} and the definition of the torus action, we obtain
\begin{align*}
&\Gamma_\Psi(t\bullet A)\\
=&{\rm Ad}_{\Psi_{n-1}(t\bullet A)\Psi^{(n)}(t\bullet A)}e^{t\bullet A }
\\=&{\rm Ad}_{\Psi^{(n)}(t \bullet (\Psi_{n-1}(A)A\Psi_{n-1}(A)^{-1}))}e^{t\bullet (\Psi_{n-1}(A)A\Psi_{n-1}(A)^{-1})}\\
=&t\bullet\left({\rm Ad}_{\Psi^{(n)}(\Psi_{n-1}(A)A\Psi_{n-1}(A)^{-1})}e^{\Psi_{n-1}(A)A\Psi_{n-1}(A)^{-1}}\right)\\
=&t\bullet \Gamma_\Psi(A),
\end{align*}
for any $t\in T(k)$ and $1\le k\le n$. It finishes the proof.\qed

\vspace{2mm}

{\bf Proof of Proposition \ref{ConnGZ}.} By the definition of $\nu(u_{\rm cat})$, we only need to prove that the connection map 
\[C^{(k)}:\Herm(k)\rightarrow {\rm U}(k)~;~B\mapsto C(E_k,B)\] 
of the $k\times k$ system \begin{equation}\label{krankeq}
\frac{dF}{dz}=\left(\I E_k-\frac{1}{2\pi \I }\frac{B}{z}\right)F,
\end{equation} 
(or equivalently the system \eqref{eq:relativeStokes}) satisfies the assumption $(1)$ and $(2)$ in Lemma \ref{lem:GWiso}. To see that, let $S^{(k)}_\pm(B)\in{\rm GL}(k)$ denote the two Stokes matrices of \eqref{eq:relativeStokes}. By Definition \ref{defiStokes} and the asymptotics of the canonical solutions $F_\pm$ of \eqref{krankeq}, we have (here to derive the second formula,
the change in choice of ${\rm log}(z)$ is accounted for)
\begin{align*} F_+(z)F_-^{-1}(z)=e^{\I E_kz}z^{-\frac{\delta_{k-1}(B)}{2\pi \I }}e^{\frac{-\delta_{k-1}(B)}{2}} S^{(k)}_+z^{\frac{\delta_{k-1}(B)}{2\pi \I }}e^{-\I E_kz}\sim {\rm Id}_k, \ \ \text{as} \ z\rightarrow\infty  \ \text{ in } {\rm Sect}_-,\\
F_-(ze^{-2\pi\I })F_+^{-1}(z)=e^{\I E_kz}z^{-\frac{\delta_{k-1}(B)}{2\pi \I }}S^{(k)}_-e^{\frac{-\delta_{k-1}(B)}{2}}z^{\frac{\delta_{k-1}(B)}{2\pi \I }}e^{-\I E_kz}\sim {\rm Id}_k, \ \ \text{as} \ z\rightarrow\infty  \ \text{ in } {\rm Sect}_+,
\end{align*}
where $\delta_{k-1}(B)$ is the projection of $B$ to the centralizer of $E_{k}$ in ${\Herm}(k)$. It follows that the Stokes matrices take the form
\begin{eqnarray}\label{StokesForm}
S^{(k)}_-(B)=\left(\begin{array}{cc}
e^{\frac{B^{(k-1)}}{2}}&0\\
{b_-^{(k)}} & \ast
\end{array} \right), \ \ \ S^{(k)}_+(B)=\left(\begin{array}{cc}
e^{\frac{B^{(k-1)}}{2}} & b^{(k)}\\
0& \ast
\end{array} \right),
\end{eqnarray}
where $b^{(k)}$ is a column vector. Furthermore, the monodromy relation \eqref{monodromyrelation} gives rise to
\begin{equation}\label{Ckmonodromy}C^{(k)}(B)e^{B}C^{(k)}(B)^{-1}=S_-^{(k)}(B) S_+^{(k)}(B), \hspace{3mm} for \ B\in\Herm(k).
\end{equation}

It follows from \eqref{StokesForm} and \eqref{Ckmonodromy} that the connection map $C^{(k)}:\Herm(k)\rightarrow {\rm U}(k)$ of the system \eqref{krankeq} satisfies the assumption $(1)$ and $(2)$ in the beginning of this subsection. Since \eqref{krankeq} is seen as a subsystem of \eqref{eq:relativeStokes}, the proof follows from Lemma \ref{lem:GWiso} and the definition of $\nu(u_{\rm cat})$.\qed 

\begin{rmk}
Now we see that the chosen minus sign in the coefficients of linear systems \eqref{krankeq} or \eqref{eq:relativeStokes} is to ensure the monodromy relation \eqref{Ckmonodromy}, which is further to ensure the compatibility with the chosen Gelfand-Tsetlin chain $\frak u(1)\subset \cdots \subset \frak u(n-1) \subset \frak u(n)$ as the upper left corners embeddings.
\end{rmk}
\begin{rmk}
There exists a family of integrable systems $\mathcal{F}(u)$ on $\frak u(n)\cong\Herm(n)$ defined by Mishchenko and Fomenko, parametrized by $\widetilde{\frak t_{\rm reg}}(\mathbb{R})$, such that $\mathcal{F}(u_{\rm cat})$ coincides with the Gelfand-Tsetlin system. More generally, we expect that there exists a family of mutiplicative analog $m\mathcal{F}(u)$ of $\mathcal{F}(u)$ on $U(n)^*\cong\Herm^+(n)$ such that the Riemann-Hilbert-Birkhoff map $\nu(u)$ intertwines $\mathcal{F}(u)$ and $m\mathcal{F}(u)$ for the same $u\in\widetilde{\frak t_{\rm reg}}(\mathbb{R})$.
\end{rmk}

{\bf Proof of Theorem \ref{eBoalchthm}.} It follows from \cite{AM} that any smooth map from $\Herm_0(n)\subset \Herm(n)$ to $\Herm_0^+(n)$, intertwining the Gelfand-Tsetlin systems, uniquely extends to a Poisson diffeomorphism from $\Herm(n)$ to $\Herm^+(n)$. Following Proposition \ref{ConnGZ}, the map $\nu(u_{\rm cat})$ intertwines the Gelfand-Tsetlin system on the open dense subset $\Herm_0(n)$. Since $\nu(u_{\rm cat})$ is already real analytic map defined on $\Herm(n)$, it is thus a Poisson diffeomorphism. \qed

\subsection{Riemann-Hilbert-Birkhoff maps are Poisson}\label{endsection}

\begin{pro}\label{poissonnature}
For any $u\in U_{\rm id}$, the diffeomorphism given in Definition \ref{phiu} \[\Phi_{u_{\rm cat}}(u):\Herm(n)\cong {\frak u}(n)^*\rightarrow \Herm(n)\cong {\frak u}(n)^*; A\mapsto \Phi(u;A)\]
is a Poisson isomorphism. 
\end{pro}
\pf On the one hand, due to the Hamiltonian description of the isomonodromy equation, the solution $\Phi(u;A)$ can be seen as a Hamiltonian flow. Thus for any fixed $u, u'\in U_{\rm id}$, the map
\[\Phi_{u_{\rm cat}}(u')\circ \Phi_{u_{\rm cat}}(u)^{-1}:\Herm(n)\rightarrow \Herm(n); \Phi(u;A)\mapsto \Phi(u';A)\]
is Poisson. On the other hand, for any $u\in U_{\rm id}$ the map
\[G(u):\Herm(n)\rightarrow \Herm(n); A\mapsto 
g(u;-A) \cdot A\cdot g(u;-A)^{-1},\]
corresponding to a time $u$ flow of the Gelfand-Tsetlin system, is a Poisson isomorphism. Here $g(u;-A)$ is defined in \eqref{gu}. It follows that the map
\[ G(u')\circ \Phi_{u_{\rm cat}}(u')\circ \Phi_{u_{\rm cat}}(u)^{-1}:\Herm(n)\rightarrow \Herm(n)\]
is a Poisson isomorphism.

Now by Theorem \ref{isomonopro} and Definition \ref{solasy}, we see that 
\[(G(u')\circ \Phi_{u_{\rm cat}}(u'))(A)=G(u')(\Phi(u';A))\rightarrow A, \ \ \text{as} \ u'\rightarrow u_{\rm cat}.\]
Therefore, for fixed $u$ the Poisson maps $G(u')\circ \Phi_{u_{\rm cat}}(u')\circ \Phi_{u_{\rm cat}}(u)^{-1}$ has the inverse map $\Phi^{-1}_{u_{\rm cat}}(u)$ of $\Phi_{u_{\rm cat}}(u)$ as a limit when $u'\rightarrow u_{\rm cat}$. It finishes the proof.
\qed

\vspace{2mm}
Given the relation $\nu(u_{\rm cat})=\nu(u)\circ \Phi_{u_{\rm cat}}(u)$ (see Definition \ref{phiu}), this proposition states that the closure of Stokes matrices preserves the Poisson geometry nature. For example, it follows that the map $\nu(u_{\rm cat})$ at the boundary point $u_{\rm cat}$ is Poisson. 

\vspace{2mm}
{\bf A new proof of Theorem \ref{Boalchthm}:} we have seen the relation $\nu(u_{\rm cat})=\nu(u)\circ \Phi_{u_{\rm cat}}(u)$ (see Definition \ref{phiu}). Following from Theorem \ref {eBoalchthm} and Proposition \ref{poissonnature}, the diffeomorphisms $\nu(u_{\rm cat})$ and $\Phi_{u_{\rm cat}}(u)$ are both Poisson. Therefore, the Riemann-Hilbert map $\nu(u)$ is Poisson. \qed

\section{Regularized limit of Stokes matrices in the De Concini-Procesi space}\label{closureStokes}
The infinite point $u_{\rm cat}$ is actually a special boundary point in the De Concini-Procesi space $\widetilde{\frak t_{\rm reg}}(\mathbb{R})$. In this section, we generalize the result in Section \ref{exStokesviaiso} from the special boundary point $u_{\rm cat}$ to arbitrary boundary points on $\widetilde{\frak t_{\rm reg}}(\mathbb{R})$. First, in Section \ref{subsect-dCP}, we recall the definition of the space $\widetilde{\frak t_{\rm reg}}(\mathbb{R})$. In Sections \ref{coorchart} and \ref{paraboundary}, we recall its stratification and coordinate charts, and introduce a parameterization of a boundary point by planar rooted coloring tree. Then in Section \ref{arybv}, we introduce the boundary value of solutions of the isomonodromy equation at an arbitrary boundary point of $\widetilde{\frak t_{\rm reg}}(\mathbb{R})$. Later in Section \ref{connboundary}, we introduce the (regularized limits of) connection and Stokes matrices at the boundary point, which depend on some discrete choices. In Section \ref{connembedd}, we show that how the connection and Stokes matrices explicitly depend on the discrete choices, and introduce the wall-crossing formula. After these preparations, in Section \ref{Limitiso}, we generalize the results in Section \ref{exStokesviaiso} from the caterpillar point $u_{\rm cat}$ to an arbitrary boundary point. 

\subsection{The De Concini-Procesi space}\label{subsect-dCP}
Let us take simple Lie algebra $\frak{sl}_n $ with the Cartan subalgebra $\frkt $, set of roots $\Pi \subset \frkt^*$, positive roots $ \Pi_+ $, and simple roots $ \{ \alpha_i\} $.

We let $ \CG $ denote the \emph{minimal building set} associated to the set of roots.  To define $ \CG $, let $ \CG' $ denote the set of all non-zero subspaces of $\frkt^* $ which are spanned by a subset of $ \Pi $.   Let $ V \in \CG' $.  We say that $ V  = V_1 \oplus \cdots \oplus V_k $ is a \emph{decomposition} of $ V $ if  $ V_1, \dots, V_k \in \CG'$, and if whenever $ \alpha  \in \Pi $ and $ \alpha \in V $, then $ \alpha \in V_i $ for some $ i$.  From Section 2.1 of \cite{dCP}, every element of $ \CG' $ admits a unique decomposition. Then we define $ \CG $ as the set of indecomposable elements of $ \CG' $. The set $ \CG $ has the following description. There is an action of the Weyl group $ W $ on $\frkt $ that preserves $ \Pi $. Thus, we get actions of $ W $ on $ \CG$ and $\CG' $. If $ J \subseteq I $ is a non-empty, connected subset of the Dynkin diagram $ I $ of $ \g $, we can form $ V_J = {\rm span}(\alpha_j : j \in J) $.  Then $ V_J \in \CG $.  In fact, every $ V \in \CG $ is of the form $ w(V_J) $ for some $ w \in W $ and $ J $ as above.

Let $\frkt_{\rm reg} = \{ \chi \in \frkt : \alpha(\chi) \ne 0, \text{ for all } \alpha \in \Pi \} $.  For any $ V \in \CG$, we have a map $ \frkt_{\rm reg} \rightarrow \mathbb{P}(\g/V^\perp) $.

\begin{defi}
The De Concini-Procesi space $\widetilde{\frkt_{\rm reg}}\subset \frkt\times \prod_{V \in \CG} \mathbb{P}(\frkt/V^\perp) $ is the closure of the image of the map $ \frkt_{\rm reg} \rightarrow  \frkt\times \prod_{V \in \CG} \mathbb{P}(\frkt/V^\perp)$.  
\end{defi}
Since the root system $\Pi$ is defined over $\mathbb{R}$, the variety is defined over $\mathbb{R}$
and so it make sense to consider the real points $\widetilde{\frkt_{\rm reg}}(\mathbb{R})$. The space $\widetilde{\frkt_{\rm reg}}(\mathbb{R})$ is called the De Concini-Procesi space of $\frak t_{\rm reg}(\mathbb{R})$, i.e., the space of $n\times n$ diagonal matrices $u={\rm diag}(u_1,...,u_n)$ with distinct real eigenvalues and such that $\sum_{k=1}^nu_k=0$. 
\begin{rmk}
As explained in the introduction, by definition, the Stokes matrices $S_\pm(u,A)$ are invariant under the translation action on $u\in\frak{h}_{\rm reg}(\mathbb{R})$. Thus let $\mathbb{R}$ act on $\frak{h}_{\rm reg}(\mathbb{R})$ by translation, then for any fixed $A$, $S_\pm(u,A)$ are parameterized by $\frak t_{\rm reg}(\mathbb{R})\cong \frak{h}_{\rm reg}(\mathbb{R})/\mathbb{R}$. In this section, we assume $u\in \frak t_{\rm reg}(\mathbb{R})$.
\end{rmk}

\subsection{Coordinate charts on $\widetilde{\frak t_{\rm reg}}(\mathbb{R})$}\label{coorchart}

The space $\widetilde{\frak t_{\rm reg}}(\mathbb{R})$ has a stratification, with the strata indexed by rooted trees with $n$ colored leaves.
Let $RT$ be such a tree, then the corresponding
stratum $\mathcal{M}_{RT}$ is the product of $\mathcal{M}_{0,d(I)}$ over all internal vertices $I$ of $RT$ with $d(I)$ the index
of $I$. In particular, $0$-dimensional strata correspond to binary rooted trees with $n$ ordered leaves, while $1$-dimensional strata correspond to almost binary trees (with exactly
one $4$-valent internal vertex). In the following, let us introduce coordinates chart on $\widetilde{\frak t_{\rm reg}}(\mathbb{R})$, and explain more on the stratification.

We denote by $RT(u)$ a planar rooted tree $RT$ with $n$ leaves colored by the components $u_1, . . . , u_n$ of $u$. We say that $RT$ is compatible with $\sigma\in S_n$ if all internal vertices of the tree are in the lower half plane, all leaves are on the horizontal line $y = 0$ and are colored by $u_{\sigma(1)},...,u_{\sigma(n)}$ from left to right.

To any planar binary rooted tree $BRT$ compatible with $\sigma$, one can assign a set of $n-1$ coordinates $z_I$, indexed by internal vertices $I$ of $BRT$, in an appropriate neighborhood $U_{BRT_\sigma}$ of the corresponding 0-dimensional stratum. The coordinate ring of the open
chart $U_{BRT_\sigma}\subset \widetilde{\frak t_{\rm reg}}(\mathbb{R})$ is generated by the following set of coordinate functions $z_I$ on $U_{BRT_\sigma}$ indexed by inner vertices $I$
of the tree $BRT$,
\begin{align}
z_{I}=\left\{
          \begin{array}{lr}
            u_{r(I)}- u_{l(I)},   & \text{ if $I$ is the root vertex},\\
            \frac{u_{r(I)}-u_{l(I)}}{u_{r(I')}-u_{l(I')}},   & \text{ if $I$ is any other vertex},
             \end{array}
\right.
\end{align}
where $I'$ is the preceding vertex of $I$ in $BRT$, i.e. $I':= {\rm max}\{J \in BRT~|~ J < I\}$ in the partial ordering $<$ of the vertices of $BRT_\sigma$ with the root being the minimal element, and for any vertex $I$, $l(I)\in [1, . . . , n]$ is such that 
$\sigma(l(I))$ is the maximal index of the $u_i’s$ in the left branch at $I$, and analogously, $r(I)\in [1,...,n]$ is such that $\sigma(r(I))$ the minimal 
index of the $u_i's$ in the right branch at $I$.

Then any evaluation of $z_I's$ represents a point in $\widetilde{\frak t_{\rm reg}}(\mathbb{R})$. In particular, the point $z_I=0$ for all inner vertices $I$ is the corresponding $0$-dimensional stratum (origin point of the chart $U_{BRT_\sigma}$). And a point, with $k$ of the $n-1$ coordinates $z_I$ are zero, is in the $n-1-k$ dimensional stratum. 
Actually, the stratum corresponding to a tree $RT$ (not necessary binary) lies in the closure of the one corresponding to another
tree $RT'$ if and only if $RT'$ is obtained from $RT$ by contracting some edges. In the coordinate charts, the stratum $\mathcal{M}_{RT'}$ corresponding to a rooted tree $RT'$ in the local
coordinates determined by a binary rooted tree $BRT$ can be described as follows.

\begin{pro}
The stratum $\mathcal{M}_{RT'}$ has a nonempty intersection with the coordinate chart $U_{{BRT_\sigma}}$ if and only if $RT'$ is
obtained from $BRT$ by contracting some edges. In the latter case, $\mathcal{M}_{RT'}$ is a subset of $U_{BRT_\sigma}$ defined as follows: $z_I\ne 0$ if the (unique) edge of $BRT$ which ends at $I$ is contracted in $RT'$, and $z_I = 0$ otherwise.
\end{pro} 

\subsection{Parameterization of boundary points of $\widetilde{\frak t_{\rm reg}}(\mathbb{R})$}\label{paraboundary}
From Section \ref{coorchart}, we have seen a cover of $\widetilde{\frak t_{\rm reg}}(\mathbb{R})$ by coordinate charts $U_{BRT_\sigma}$ associated to planar binary rooted trees. And boundary points of $\widetilde{\frak t_{\rm reg}}(\mathbb{R})\setminus \frak t_{\rm reg}(\mathbb{R})$ in $U_{BRT_\sigma}$ are parameterized by coordinates $z_I's$ with some of them being zero.
In the following, let us introduce another way to parameterize the boundary points, that is particularly easy to work with for our purpose.
\begin{defi}
We denote by $RT_\sigma(u)$ a planar rooted tree $RT$ (not necessary binary) with $n$ leaves colored by $u$ and compatible with $\sigma\in S_n$.
\end{defi}
\begin{figure}[H]
\begin{center}
  \begin{tikzpicture}[scale=0.7]
  \draw
  (0,0)--(0,1)--(2.5,4) node[above]{$u_3$}
  (0,1)--(-2.5,4) node[above]{$u_5$}
  (1.25,2.5)--(0,4) node[above]{$u_6$}
  (1.75,3.1)--(1,4) node[above]{$u_2$}
  (1.75,3.1)--(1.75,4) node[above]{$u_1$}
  (-1.75,3.1)--(-1,4) node[above]{$u_4$};
  \end{tikzpicture}
  \caption{A planar binary rooted tree with $6$ leaves colored by $u_1,...,u_6$.}
 \end{center}
 \end{figure}

Then $RT_\sigma(u)$ represents a point in $\widetilde{\frak t_{\rm reg}}(\mathbb{R})$ as follows. If the planar rooted tree $RT_\sigma(u)$ can be obtained from a binary tree $BRT$ (compatible with $\sigma$) by contracting some edges, the coloring $u=(u_1,...,u_n)$ of the leaves of $RT_\sigma(u)$ determines a point $u_0$ in the chart $U_{BRT_\sigma}\subset \widetilde{\frak t_{\rm reg}}(\mathbb{R})$: in the coordinates $\{z_I\}$ associated to $BRT$ and $\sigma$, the point $u_0$ has the coordinates 
\begin{align}\label{xcoor1}
z_I&= \frac{u_{r(I)}-u_{l(I)}}{u_{r(I')}-u_{l(I')}}, \ \text{ if $I$ and $I'$ are contracted in $RT_\sigma(u)$}, \\ \label{xcoor2}
z_I&=0, \ \ \ \ \ \ \text{otherwise}.
\end{align}
where $I'$ is the preceding vertex of $I$ in the binary tree $BRT$. We call $RT_\sigma(u)$ a representative of $u_0$. 

For any $\sigma\in S_n$, we take the connected component $U_{\sigma}$ of $\frak t_{\rm reg}(\mathbb{R})$, and denote by $\overline{U_{\sigma}}$ its closure in $\widetilde{\frkt_{\rm reg}}(\mathbb{R})$. Note that $\overline{U_{\sigma}}$ is consisting of the points which can be represented by certain planar rooted tree $RT_{\sigma}(u)$ compatible with $\sigma$.
From the expression \eqref{xcoor1}, we see that on the one hand, for different coloring $u$, $RT_\sigma(u)$ can represent a same point. In particular, if we have a binary tree $RT_\sigma$, then for all $u\in U_{\sigma}$, $RT_\sigma(u)$ represents a same point. On the other hand, if the ordering $\sigma$ is obtained from $\sigma'$ by reversing the order
of the descendants of any internal vertex of $RT_\sigma(u)$, then $RT_\sigma(u)$ and $RT_{\sigma'}(u)$ represent a same point.

\begin{ex}
Let us take the following planar binary tree $BRT_{\sigma}$ with coloring
\begin{figure}[H]
\begin{center}
  \begin{tikzpicture}[scale=1.2]
  \draw
  (0,0)--(0,1)--(3,4) node[above]{$u_{\sigma(n)}$}
  (0,1)--(-3,4) node[above]{$u_{\sigma(1)}$}
  (-2.5,3.5)--(-2,4) node[above]{$u_{\sigma(2)}$}
  (-2,3)--(-1,4) node[above]{$u_{\sigma(3)}$}
  (-1.25,2.25)--(0.5,4) node[above]{$\cdots $}
  (-0.5,1.5)--(2,4) node[above]{$u_{\sigma(n-1)}$};
  
  \end{tikzpicture}
  \caption{A planar rooted tree with coloring}
 \end{center}
 \end{figure}
Let us denote the vertices of the above tree $BRT_\sigma$ in the partial ordering by $I^{(1)}>I^{(2)}>\cdots >I^{(n)}$. And let 
$z_{I^{(1)}}, ..., z_{I^{(n)}}$ be the corresponding coordinates, then (following the definition in \cite[Page 16]{Sp}) the point $u^\sigma_{\rm cat}\in U_{BRT_{\sigma}}$ with coordinates $z_{I^{(k)}}=0$ for all $k=1,...,n$ is called a caterpillar point. For $\sigma={\rm id}\in S_n$, we simply denote $u^{\rm id}_{\rm cat}$ by $u_{\rm cat}$.
\end{ex}

Now associated to a planar labelled rooted tree $RT_\sigma(u)$, let us introduce another set of coordinates on the connected component $U_{\sigma}$.
\begin{defi}\label{def:tcoor}
For any vertex $I\in RT_\sigma(u)$ with $E(I)$ outgoing edges,
we denote its branches by $I_1,..., I_{E(I)}$ (counting from left to right). We define $r_I$ and $\{t_{I_j}\}_{j=1,...,E(I)}$ as the following collection of $1+E(I)$ variables assigned to the vertex $I$ and the $E(I)$ branches of $I$
\begin{itemize}
\item $r_I=u_{min(I_{E(I)})}-u_{max(I_{1})}\in\mathbb{R}_{>0}$;
    \item for $I$ an internal vertex,
\begin{equation}\label{tcoor}
t_{I_1}=0, \ t_{I_{E(I)}}=1 \hspace{3mm} \text{and} \hspace{3mm} t_{I_j}=\frac{u_{min(I_{j})}-u_{max(I_{1})}}{u_{min(I_{E(I)})}-u_{max(I_{1})}}, \ \text{for} \ 1<j<E(I),
\end{equation}
here $max(I_1)\in [1,...,n]$ is such that 
$\sigma(max(I_1))$ is the maximal index of the $u_i’s$ in the first branch $I_1$ of $I$, and analogously, $min(I_{j})\in [1,...,n]$ is such that $\sigma(min(I_{j}))$ the minimal 
index of the $u_i's$ in the $j$-th branch $I_{j}$ of $I$;
\item for $I$ the root vertex,  \begin{equation}\label{tcoor1}
t_{I_1}=0, \hspace{3mm} \text{and} \hspace{3mm} t_{I_j}={u_{r(I_{j})}-u_{l(I_{1})}}, \ \text{ for } \ 1<j\le E(I).
\end{equation}
\end{itemize}

\end{defi}
Note that the collection of variables $\{t_{I_j}\}_{j=1,...,E(I)}, \{r_I\}$, for all vertices $I$ of $RT_\sigma(u)$, define a new coordinate system on $U_\sigma$. And the change of coordinates from the system $\{t_{I_j}, r_I\}$ to the coordinate system $\{u_i\}$ of $U_\sigma$ is 
\begin{equation}\label{urtchange}
u_{min(I_{j})}-u_{max(I_{1})}=r_I\cdot t_{I_j} \text{ for } 1\le j\le E(I).
\end{equation}

If the planar labelled rooted tree $RT_\sigma(u)$ with coloring has $k$ internal vertex, then it represents a point $u_0\in\overline{U_{\sigma}}$ in the codimension $k$ strata. The numbers of $\{r_I\}$ and $\{t_{I_j}\}$ variables are exactly $k$ and $n-k-1$.
Actually, by \eqref{xcoor1}-\eqref{xcoor2} and \eqref{tcoor}, the variables $\{t_{I_j}\}$ for all $I$ already specify the point $u_0$.

To be more precise, under the coordinate transform \eqref{urtchange}, the points $u(\{t_{I_j}\}, \{r_I\})\in U_\sigma$ are functions of the variables $\{r_I\}$ and $\{t_{I_j}\}$. Although the components $u_1,...,u_n$ of $u$ depend on the variables $\{r_I\}$ and $\{t_{I_j}\}$, while for fixed $\{t_{I_j}\}$, by \eqref{xcoor1}-\eqref{xcoor2} and \eqref{tcoor}, the tree $RT_\sigma(u)$ with coloring $u(\{t_{I_j}\}, \{r_I\})$ represents the same point $u_0$ for all positive real numbers $\{r_I\}$.
It is useful to think of $u_0$ as the limit of the point $u(\{t_{I_{j}}\}, \{r_I\})\in U_\sigma$ as ${r_I}/{r_{I'}}\rightarrow 0+$ for all $I$ and its preceding vertex $I'$ in $RT_\sigma(u)$. And the tree $RT_\sigma(u)$ with coloring $u$ encodes the way to take limit of the components $u_1,...,u_n$ of $u\in U_\sigma$.

\begin{ex}\label{catex}
Let us take the following binary tree with coloring, that represents the caterpillar point $u_{\rm cat}$. Denote its vertices in the partial ordering by $I^{(1)}>I^{(2)}>\cdots >I^{(n)}$. Its internal vertex $I^{(k)}$ for $k=1,...,n-1$ has two branches $I^{(k)}_1, I^{(k)}_2$. By Definition \ref{def:tcoor} the numbers associated to the two branches are just $t_{I^{(k)}_1}=0$ and $t_{I^{(k)}_2}=1$, and the number $r_{I^{(k)}}$ associated to the vertex $I^{(k)}$ is $u_{k+1}-u_k$. And the point $u_{\rm cat}$ is then the limit $\frac{r_{I^{(k-1)}}}{r_{I^{(k)}}}=\frac{u_{k}-u_{k-1}}{u_{k+1}-u_{k}}\rightarrow 0+$ for all $k$.
\begin{figure}[H]
\begin{center}
  \begin{tikzpicture}[scale=1.2]
  \draw
  (0,0)--(0,1)--(3,4) node[above]{$u_{n}$}
  (0,1)--(-3,4) node[above]{$u_{1}$}
  (-2.5,3.5)--(-2,4) node[above]{$u_{2}$}
  (-2,3)--(-1,4) node[above]{$u_{3}$}
  (-1.25,2.25)--(0.5,4) node[above]{$\cdots $}
  (-0.5,1.5)--(2,4) node[above]{$u_{n-1}$};
  
  \end{tikzpicture}
 \caption{A planar rooted tree with coloring that represents the caterpillar point $u_{\rm cat}$}
 \end{center}
 \end{figure}
 
\end{ex}

\begin{ex}\label{ex3}
Let us consider planar tree with coloring in the following figure. Denote its three internal vertices in the partial ordering by $I^{(l)}>I$, $I^{(r)} >I$, where $I^{(l)}$ is the vertex that has $k$ branches $I^{(l)}_1, ..., I^{(l)}_k$, while $I^{(r)}$ has $n-k$ branches $I^{(r)}_1, ..., I^{(r)}_{n-k}$. By Definition \ref{def:tcoor} the numbers associated to the two branches of $I$ are $t_{I_1}=0$ and $t_{I_2}=1$, and the numbers associated to the branches of $I^{(l)}$ and $I^{(r)}$ are $t_{I^{(l)}_i}=\frac{u_i-u_1}{u_k-u_1}$ (for $i=1,...,k$) and $t_{I^{(r)}_j}=\frac{u_j-u_{k+1}}{u_n-u_{k+1}}$ (for $j=k+1,...,n$) respectively. The numbers $r_{I}=u_{k+1}-u_k$, $r_{I^{(l)}}=u_k-u_1$ and $r_{I^{(r)}}=u_n-u_{k+1}$. Note that the boundary point $u_0$ represented by the tree is specified by the $t$ variables, and is independent of the $r$ variables. And the coloring tree represents (or the corresponding boundary point $u_0$) the limit 
\begin{equation*}
 \frac{r_{I^{(l)}}}{r_{I}}=\frac{u_k-u_{1}}{u_{n}-u_{1}}\rightarrow 0+, \ \  \frac{r_{I^{(r)}}}{r_{I}}=\frac{u_n-u_{k+1}}{u_{n}-u_{1}}\rightarrow 0+, \ \ r_I=u_n-u_1\rightarrow 0+,
\end{equation*}
with the given fixed ratios $\frac{u_i-u_1}{u_k-u_1}$ and $\frac{u_j-u_{k+1}}{u_n-u_{k+1}}$ for $i=1,...,k$ and $j=k+1,...,n$. 
\begin{figure}[H]\label{figure3}
\begin{center}
  \begin{tikzpicture}[scale=0.7]
  \draw
  (0,0)--(0,1)--(4,4) node[above]{$u_n$}
  (0,1)--(-4,4) node[above]{$u_1$}
  (-2,2.5)--(-3,4) node[above]{$u_2$}
   (-2,2.5)--(-2,4) node[above]{...}
  (2,2.5)--(2,4) node[above]{$...$}
  (2,2.5)--(1,4) node[above]{$u_{k+1}$}
  (2,2.5)--(3,4) node[above]{$...$}
  (-2,2.5)--(-1,4) node[above]{$u_k$};
  \end{tikzpicture}
   \caption{The planar rooted tree for Example \ref{ex3}}
 \end{center}
 \end{figure} 
 \end{ex}

\subsection{Boundary value of solutions of the isomonodromy equation at an arbitrary boundary point}\label{arybv}
\begin{defi}
Suppose $u_0\in \widetilde{\frak t_{\rm reg}}(\mathbb{R})\setminus {\frak t_{\rm reg}}(\mathbb{R})$ is a boundary point represented by a planar rooted colored tree $RT_\sigma(u)$. For any $A\in\Herm(n)$ and any internal vertex $I$ of $RT_\sigma(u)$, we denote by $\delta_I(A)\in\Herm(n)$ the matrix whose $(i,j)$ entry is
\begin{equation}\label{takeI}
\delta_I(A)_{ij}=\left\{
          \begin{array}{lr}
            A_{ij},   & \text{if } i,j\in Branch(I), \ \text{or} \ i=j \\
           0, & \text{otherwise},
             \end{array}
\right. 
\end{equation}
where $Branch(I)$ is the subset of $\{1,...,n\}$ formed by all index of the coloring $u_k's$ of leaves at $I$.
\end{defi}
Then the following theorem describes the asymptotics of solutions of the isomonodromy equations as $u\rightarrow u_0$, that is an analog of Theorem \ref{isomonopro}.
\begin{thm}\label{isomonopro2}
For any point $u_0$ represented by a planar rooted tree $RT_\sigma(u)$ and any solution $\Phi(u)$ of the isomonodromy equation \eqref{isoeq} on $U_{\sigma}$, there exists a constant $A(u_0)\in\Herm(n)$ such that as $r_I/r_{I'}\rightarrow 0+$ for all $I$ and its preceding vertex $I'$ in $RT_\sigma(u)$,
\begin{equation}\label{gboundaryval}
\left(\overrightarrow{\prod_{I\in RT_\sigma(u)} }\left(\frac{r_I}{r_{I'}}\right)^{\frac{\delta_I(\Phi(u))}{2\pi\I }}\right)\cdot \Phi(u)\cdot \left(\overrightarrow{\prod_{I\in RT_\sigma(u)} }\left(\frac{r_I}{r_{I'}}\right)^{\frac{\delta_I(\Phi(u))}{2\pi\I }}\right)^{-1} \rightarrow A(u_0),
\end{equation}
where the product $\overrightarrow{\prod}$ is taken over all internal vertices $I$ of $RT_\sigma(u)$, and is taken with $I$ to the right of ${J}$ if $I<J$ in the partial order. Here $\Phi(u)^I$ is defined as \eqref{takeI} associated to the vertex $I$. Furthermore, given any $A(u_0)\in\Herm(n)$, there exists a unique real analytic solution $\Phi(u)$ of \eqref{introisoeq} on $U_\sigma$ such that the limit \eqref{gboundaryval} holds. 
\end{thm}
\begin{defi}
We call $\A(u_0)$ the boundary value of $\Phi(u)$ at $u_0\in\overline{U_{\sigma}}$, and denote by $\Phi^\sigma(u;A(u_0))$ the solution on $U_\sigma$ with the given boundary value $A(u_0)$ at $u_0$. 
\end{defi}
In the following, let us study how to express the Stokes matrices $S_\pm\left((u, \Phi^\sigma(u;\Phi(u_0))\right)$ via the boundary value $A(u_0)$.

\subsection{Connection matrices and Stokes matrices parameterized by $\widetilde{\frak t_{\rm reg}}(\mathbb{R})\setminus \frak t_{\rm reg}(\mathbb{R})$ }\label{connboundary}
Suppose $u_0\in \overline{U_{\sigma}}$ is a boundary point represented by a planar rooted colored tree $RT_\sigma(u)$. Recall that for any vertex $I\in T_\sigma(u)$ with $E(I)$ outgoing edges, we have assigned a collection of numbers $\{t_{I_j}\}_{j=1,...,E(I)}$ by \eqref{tcoor} and \eqref{tcoor1}.  Then with respect to the permutation $\sigma\in S_n$ and associated to each $I$, let us introduce an $n\times n$ diagonal matrix $t_\sigma^{I}={\rm diag}((t_\sigma^{I})_{1},...,(t_\sigma^{I})_{n})$ with
\begin{align}\label{ucoor1}
(t_\sigma^{I})_{k}&=t_{I_j}, \ \ \ \text{if the coloring} \ u_{\sigma(k)} \text{ is in the $j$-th branch $I_j$ of $I$} \ \ \text{for certain} \ 1\le j\le E(I),\\ \label{ucoor2}
(t_\sigma^{I})_{k}&=0, \ \ \ \ \text{otherwise},
\end{align}
Besides, for any $A\in\Herm(n)$ we define $\delta_I(A)\in\Herm(n)$ as in \eqref{takeI}. 

Now associated to $RT_\sigma(u)$ and any $A\in\Herm(n)$, there is a finite set of $n\times n$ linear systems labelled by the vertices of $RT_\sigma(u)$: for any vertex $I\in RT_\sigma(u)$, the linear system is
\begin{eqnarray}\label{Isystem}
\frac{dF}{dz}=\left(\I t_\sigma^{I}-\frac{1}{2\pi\I }\frac{\delta_I(A)}{z}\right)F.
\end{eqnarray}
The system has only two anti-Stokes directions, the two halves of the imaginary axis. We will choose the right half plane ${\rm Sect}_+:=\{z\in\mathbb{C}~|~ {\rm Re}(z)>0\}$ as the initial Stokes sector, and take the branch of ${\rm log}(z)$ which is real on $\mathbb{R}_{>0}$. Let us denote by $C(t_\sigma^{I},\delta_I(A))$ and $S_\pm(t_\sigma^{I},\delta_I(A))$ the associated connection and Stokes matrices. 
Let $U(n)^{I}\subset U(n)$ be the stabilizer subgroup of $t_\sigma^{I}\subset \Herm(n)$ under the conjugation action of $U(n)$ on $\Herm(n)$. Then the connection map \[C(t_\sigma^{I}):\Herm(n)\rightarrow U(n);~A\mapsto C(t_\sigma^{I},\delta_I(A))\]
is equivariant with respect to the group $U(n)^{I}$, and (similar to the $u_{\rm cat}$ case given in \eqref{StokesForm}) the Stokes matrices $S_\pm(t_\sigma^{I},\delta_I(A))$ are blocked according to the index set of the coloring $u_i's$ in the branches $I_1,...,I_{E(I)}$ at the vertex $I$. Furthermore, we have the monodromy relation
\begin{eqnarray}\label{monoI}
C(t_\sigma^{I},\delta_I(A))e^{\delta_I(A)}C(t_\sigma^{I},\delta_I(A))^{-1}=S_-(t_\sigma^{I},\delta_I(A))S_+(t_\sigma^{I},\delta_I(A)).
\end{eqnarray}

\begin{defi}\label{ConnectionTree}
For any $A\in\Herm(n)$, the {\it connection matrix} $C^\sigma(u_0,A)$ at a boundary point $u_0\in \overline{U_{\sigma}}$ (with respect to the connected component $U_\sigma$) is
\begin{eqnarray}\label{connbound}
C^\sigma(u_0,A):=\overrightarrow{\prod_{I\in RT_\sigma(u)}}C(t_\sigma^{I},\delta_I(A)),
\end{eqnarray}
where $RT_\sigma(u)$ is any representative of $u_0$, the product is taken over all vertices $I$ of $RT_\sigma(u)$, and is taken with $C(t_\sigma^{I},\delta_I(A))$ to the right of $C(t_\sigma^{I'},\delta_{I'}(A))$ if $I'<I$ in the partial ordering. 
\end{defi}

Let us check the definition is independent of the different choices of the representatives $RT_\sigma(u)$ of $u_0$. By comparing identities \eqref{xcoor1} and \eqref{xcoor2} with \eqref{tcoor}, we see that the irregular term $t_\sigma^{I}$ in the equation \eqref{Isystem} only depends on the variables $\{t_{I_j}\}$. Since the $\{t_{I_j}\}$ variables for all $I$ specify the point $u_0$, the right hand side of \eqref{connbound} only depends on the point $u_0$ and the chosen compatible permutation $\sigma\in S_n$. That is the right hand side of \eqref{connbound} 
is independent of the different choices of the coloring $u$, as long as $RT_\sigma(u)$ is compatible with $\sigma$ and represents the same point $u_0$. Recall that by Definition \eqref{def:tcoor} we think of the boundary point $u_0$ as the limit of the point $u(\{t_{I_{j}}\}, \{r_I\})\in U_\sigma$ as ${r_I}/{r_{I'}}\rightarrow 0+$ for all $I$ and its preceding vertex $I'$ in $RT_\sigma(u)$. Then as we will see, the connection matrix $C^\sigma(u_0,A)$ is the regularized limit of the (ordinary) connection matrices $C(u(\{t_{I_{j}}\}, \{r_I\}),A)$ as ${r_I}/{r_{I'}}\rightarrow 0+$.

In particular, if $u_0$ is in the $0$-dimensional strata, then any coloring $u\in U_{\sigma}$ of the binary tree $RT_\sigma(u)$ will define the same $C^\sigma(u_0,A).$ The definition of $C^\sigma(u_0,A)$ also depends on the choice of $\sigma$. That is for $u_0\in \overline{U_\sigma}\cap \overline{U_{\sigma'}}$, i.e., there exists different $\sigma,\sigma'\in S_n$ such that $u_0$ can be represented by certain $RT_\sigma(u)$ and $RT_{\sigma'}(u)$, the connection matrices $C^\sigma(u_0,A)$ and $C^{\sigma'}(u_0,A)$ will be in general different, see Section \ref{connembedd}. 

Now since each factor $C(t_\sigma^{I},\delta_I(A))$ in the product is unitary, we have
\begin{pro}
The connection matrix $C^\sigma(u_0,A)$ is unitary.
\end{pro}

Given any $\sigma\in S_n$, denote by $P_\sigma\in{\rm GL}_n$ the corresponding permutation matrix.
\begin{defi}\label{StokesTree}
For any $u_0\in \overline{U_\sigma}$ and $A\in\Herm(n)$, the {\it Stokes matrices} $S^\sigma_\pm(u_0,A)$ at $u_0$ with respect to $U_\sigma$ are the unique triangular matrices, with ${P_{\sigma}}^{-1} e^{\frac{[A]}{2}}P_\sigma$ as diagonal part, determined by the identity (Gauss decomposition) \begin{eqnarray}
{C^\sigma(u_0,A)}e^A{C^\sigma(u_0,A)}^{-1}=P_\sigma S_{+}^\sigma(u_0,A)S_{-}^\sigma(u_0,A){P_\sigma}^{-1}.
\end{eqnarray}
Here $[A]$ is the diagonal part of $A$.
\end{defi}

\begin{ex}
If we take $u_0=u_{\rm cat}$ and the compatible $\sigma={\rm id}\in S_n$, then Definitions \ref{ConnectionTree} and \ref{StokesTree} coincide with the definition of connection and Stokes matrices at $u_{\rm cat}$ introduced in Section \ref{exviaasy}.
\end{ex}

We stress that just like the connection matrices, the Stokes matrices, at a boundary point $u_0$, not only depend on the boundary point itself, but also depend on the discrete choices of the $\sigma\subset S_n$ such that $u_0$ is in the closure of the connected component $U_\sigma$. 
\subsection{Different choices of planar embeddings and wall-crossing formula}\label{connembedd}
In this subsection, we study how connection and Stokes matrices associated to planar trees depend on the choices of planar embeddings. 

Suppose that a given point $u_0$ lies at the intersection of $\overline{U_{\sigma}}$ and $\overline{U_{\sigma'}}$, for different $\sigma, \sigma'\in S_n$. Then there are two equivalent representatives of $u_0$, i.e., planar embeddings $RT_\sigma(u)$ and $RT_{\sigma'}(u)$ represent the same rooted tree coloring $u$ compatible with $\sigma$ and $\sigma'$ respectively. Here equivalence means that $RT_{\sigma'}(u)$ is obtained from $RT_\sigma(u)$ by reversing the order of the descendants of some internal vertex of $RT_\sigma(u)$. Let us assume that $\sigma$ and $\sigma'$ differ by reversing only one cycle. Then $RT_{\sigma'}(u)$ is given by reversing the order of the descendants of some vertex $I$ of $RT_\sigma(u)$. 

Assume that $C(t_\sigma^J,A^J)$ and $C(t_{\sigma'}^J,A^J)$ are the connection matrices associated to the vertex $J$ of $RT_\sigma(u)$ and $RT_{\sigma'}(u)$ respectively. Then by the assumption and the definition of the connection matrices, we have 
\begin{eqnarray*}
C(t_{\sigma'}^J,A^J)=\left\{
          \begin{array}{lr}
             C(t_{\sigma}^J,A^J),   & if \ J\ne I,  \\
           C_-(t_{\sigma}^J,A^J), & J=I,
             \end{array}
\right. 
\end{eqnarray*}
where $C_-(t_\sigma^{I},\delta_I(A))$ is the connection matrix of the system \eqref{Isystem} but with respect to the Stokes sector ${\rm Sect}_-$. To be more precise, similar to Definition \ref{connectionmatrix}, the connection
matrix $C_-(t_\sigma^{I},\delta_I(A))$ of the system \eqref{Stokeseq} (with $(u,A)$ replacing by $(t_\sigma^{I},\delta_I(A))$ is determined by 
$F_0(z)=F_-(z)\cdot C_-(t_\sigma^{I},\delta_I(A))$, where $F_-(z)$ is the canonical solution on ${\rm Sect}_-$. Furthermore, since a negative half cycle (i.e., in clockwise direction) around $0$ is a positive half cycle around $\infty$, we get the monodromy relation
\begin{equation*}
C_-(t_\sigma^{I},\delta_I(A))=\left(S_+(t_\sigma^{I},\delta_I(A))S_-(t_\sigma^{I},\delta_I(A))\right)^{-\frac{1}{2}}S_+(t_\sigma^{I},\delta_I(A))C(t_\sigma^{I},\delta_I(A)),
\end{equation*}
where $S_\pm(t_\sigma^{I},\delta_I(A))$ are the Stokes matrices of the system \eqref{Isystem}. In summary, we have the following wall-crossing formula of the connections matrices as $u$ crosses the common face of $U_\sigma\cup U_{\sigma'}$.
\begin{pro}\label{diffembedd}
Given any point $u_0$ in $\overline{U_{\sigma}}\cap \overline{U_{\sigma'}}$, where $\sigma$ differs with $\sigma'$ by reversing some segments, let $RT_\sigma(u)$ and $RT_{\sigma'}(u)$ be its two representatives. Then for any $A\in\Herm(n)$, $C^{\sigma'}(u_0,A)$ equals to $C^\sigma(u_0,A)=\overrightarrow{\prod}_{J\in RT_{\sigma}(u)}C(t_\sigma^J,A)$ after replacing the corresponding factors $C(t_\sigma^J,A)$ by $C_-(t_\sigma^J,A)$ in the product. 
\end{pro}
Following Definition \ref{StokesTree}, we also get a wall-crossing formula of the Stokes matrices as $u$ crosses the common face of $U_\sigma\cup U_{\sigma'}$.

\begin{ex}[Caterpillar points]\label{cactusaction}
For each $1\le i\le n$, let $\tau_i\in S_n$ be the permutation reversing the segment $[1,...,i]$. Then associated to $\tau_i$, there is a representative $RT_{\tau_i}(u)$ of the caterpillar point $u_{\rm cat}$ as in the following figure. Note that for all $\tau_i$, $i=1,...,n$, the trees $RT_{\tau_i}(u)$ represent the same point.
 
\begin{figure}[H]
\begin{center}
  \begin{tikzpicture}[scale=1.2]
  \draw
  (0,0)--(0,1)--(3,4) node[above]{$u_{n}$}
  (0,1)--(-3,4) node[above]{$u_{i}$}
  (-2.5,3.5)--(-2,4) node[above]{$\cdots$}
  (-2,3)--(-1,4) node[above]{$u_1$}
  (-1.5,2.5)--(0,4) node[above]{$u_{i+1}$}
(-1,2)--(1,4) node[above]{$\cdots$ }
  (-0.5,1.5)--(2,4) node[above]{$u_{n-1}$};
  
  \end{tikzpicture}
  \caption{A caterpillar point with a planar embedding given by $\tau_i$}
 \end{center}
 \end{figure}
For any $A\in\Herm(n)$, we denote by $S_+^{\tau_i}(u_{\rm cat}, A)$ the Stokes matrix associated to the planar tree. For $\tau_1={\rm id}\in S_n$, it is just the Stokes matrix $S_{+}(u_{\rm cat}, A)$. For any other $\tau_i$, $S_+^{\tau_i}(u_{\rm cat}, A)$ is the Stokes matrix at the caterpillar point but with respect to a different choice of planar embedding. Thus we can use Proposition \ref{diffembedd} to express $S_+^{\tau_i}(u_{\rm cat},A)$ by $S_+(u_{\rm cat},A)$.

\begin{pro}[Wall-crossing formula at $u_{\rm cat}$]\label{procac}
Under the change of the planar embedding, the associated Stokes matrices at $u_{\rm cat}$ change as \begin{eqnarray}
S_+(u_{\rm cat},A)=\left(\begin{array}{cc}
S_{i+} & B \\
0 & C
\end{array} \right)\rightarrow S_+^{\tau_i}(u_{\rm cat},A)=\left(\begin{array}{cc}
P_{i}S_{i+}^\dagger P_i^{-1} & P_i(S_{i+} S_{i+}^\dagger)^{-\frac{1}{2}}S_{i+} B \\
0 & C
\end{array} \right),
\end{eqnarray}
where $S_{i+}$ are the left-top $i$-th principal submatrices of $S_+(u_{\rm cat},A)$, and $S_{i+}^\dagger$ is the complex conjugate of $S_{i+}$, and $P_i$ is the $i\times i$ permutation matrix associated to $\tau_i$.
\end{pro}
\pf By Definition \ref{StokesTree}, the Stokes matrix $S^{\tau_i}_+(u_{\rm cat}, A)$ associated to $\tau_i$, is given by
\begin{eqnarray}\label{monoreltau}
{C^{\tau_i}(u_{\rm cat}, A)}e^A{C^{\tau_i}(u_{\rm cat}, A)}^{-1}=P_i S_{-}^{\tau_i}(u_{\rm cat}, A) S_{+}^{\tau_i}(u_{\rm cat}, A){P_i}^{-1},
\end{eqnarray}
where \[C^{\tau_i}(u_{\rm cat}, A)=C_-(E_1,\delta_1(A))\cdot\cdot\cdot C_-(E_i,\delta_1(A_i))C(E_{i+1},\delta_{i+1}(A))\cdot\cdot\cdot C^{(n)}.\] Here for $1\le k\le i$, $C_-(E_k,\delta_k(A))$ is the connection matrix of 
$\frac{dF}{dz}=\left(-\I E_k-\frac{1}{2\pi\I }\frac{\delta_k(A)}{z}\right)F,$ and recall that
$C(E_k,\delta_k(A))$, $S_\pm(E_k,\delta_k(A))$ are the connection and Stokes matrices of $\frac{dF}{dz}=\left(\I E_k-\frac{1}{2\pi\I }\frac{\delta_k(A)}{z}\right)F.$ 
\begin{lem} For any $A\in\Herm(n)$, we have
\begin{eqnarray}\label{monorel2}
C_-(E_1,\delta_1(A))\cdot\cdot\cdot C_-(E_i,\delta_1(A_i))=(S_{i+}S_{i-})^{-\frac{1}{2}}S_{i+} C(E_1,\delta_1(A))\cdot\cdot\cdot C(E_i,\delta_1(A_i)).
\end{eqnarray}
\end{lem}
\pf Consider a rank $i$ linear system $\frac{dF(z)}{dz} = \left(\I u-\frac{1}{2\pi\I }\frac{\Phi(u)}{z}\right)\cdot F(z)$. Since a negative half cycle (i.e., in clockwise direction) around $0$ is a positive half cycle around $\infty$, we get the monodromy relation
\begin{eqnarray*}
C_-(u,\Phi(u))=(S_{+}(u,\Phi(u))S_{-}(u,\Phi(u)))^{-\frac{1}{2}}S_{+}(u,\Phi(u))C(u,\Phi(u)).
\end{eqnarray*} 
If we take $\Phi(u)$ the solution of the corresponding rank $i$ isomonodromy equation with the boundary value $A$ at $u_{\rm cat}$, then the above relation implies the identity \eqref{monorel2}. It finishes the proof of Proposition \ref{procac}.
\qed 

\vspace{2mm}
The lemma implies 
\[C^{\tau_i}(u_{\rm cat}, A)=(S_{i+}S_{i-})^{-\frac{1}{2}}S_{i+}\cdot C(u_{\rm cat}, A).\]
Then by a direct blocked matrix manipulation, the proposition follows from the defining relations \eqref{monorelation} and \eqref{monoreltau} of $S_{+}(u_{\rm cat}, A)$ and $S_{+}^{\tau_i}(u_{\rm cat}, A)$, and the relation $S_{i-}(u_{\rm cat}, A)=S_{i+}(u_{\rm cat}, A)^{\dagger}$. \qed
\end{ex}

\begin{rmk}
Proposition \ref{procac} can be used to interpret the cactus group actions on the Gelfand-Tsetlin cones \cite{BK} from the perspective of WKB approximation. See Section \ref{tropicalisomo}.
\end{rmk}

\subsection{Regularized limit of Stokes matrices }\label{Limitiso}
In this subsection, we generalize the results in Section 3, particularly the Theorem \ref{relautomorphism}, from the caterpillar points to any boundary points.

\begin{thm}[The expression of Stokes matrices at the boundary via the isomonodromy deformation]\label{opentoboundary}
Let us denote by $\Phi(u;A(u_0))$ the solution of \eqref{isoeq} with the boundary value $A(u_0)$ at $u_0\in \overline{U_{\sigma}}$ as in Theorem \ref{isomonopro2}, then we have
\begin{eqnarray*}
S_\pm(u,\Phi^\sigma(u;A(u_0))=S_\pm^{\sigma}(u_0,A(u_0)).
\end{eqnarray*}
\end{thm}

For $u_0=u_{\rm cat}$ and $\sigma={\rm id}$ (i.e., $U_\sigma=U_{\rm id}$), they recover Theorem \ref{isomonopro} and Theorem \ref{relautomorphism}. For a general boundary point $u_0$, in order to prove them, we can repeat the recursive procedures in Sections \ref{asyisoeq}--\ref{exviaasy}. The only difference is that now we need to block the system \eqref{introeq} in a different way, i.e., study the set of intermediate linear systems determined by $u_0$ and the associated isomonodromy equations restricted to each variable $r_I$ (instead of the variable $u_i$ as in the caterpillar point case). Since the proofs are rather same, the complete proofs will be omitted here. 

And the analog of Theorem \ref{introcor} is
\begin{thm}[Regularized limits of Stokes matrices at the boundary]\label{reglimitbound}
For any $A\in\Herm(n)$ and a boundary point $u_0\in \overline{U_{\sigma}}$ represented by a planar rooted colored tree $RT_\sigma(u)$, we have
\begin{eqnarray*}
{\rm Ad}\left(\overrightarrow{\prod_{I\in RT_\sigma(u)} }\left(\frac{r_{I'}}{r_{I}}\right)^{\frac{{\rm log}(\delta_I(S_-)\delta_I(S_+))}{2\pi\I }}\right) (S_-(u, A)S_+(u,A)) \rightarrow S_-^{\sigma}(u_0,A)S_+^{\sigma}(u_0,A),
\end{eqnarray*}
as $r_I/r_{I'}\rightarrow 0+$ for all vertex $I$ and its preceding vertex $I'$ in $RT_\sigma(u)$.
Here $\delta_I(S_-(u,A))\delta_I(S_+(u,A))$ is a positive definition Hermitian matrix, and ${\rm log}$ takes its logarithm.
\end{thm}
It can be proved in the same way to the proof of Theorem \ref{introcor} given in Section \ref{pf:introthm3}, except that we now need the following equivariance, generalizing the identity \ref{deltaequiva} at $u_{\rm cat}$,
\begin{align*}
&S_-^{\sigma}\left(u_0,
{\rm Ad}\left(\overrightarrow{\prod_{I\in RT_\sigma(u)} }\left(\frac{r_I}{r_{I'}}\right)^{\frac{\delta_I(A)}{2\pi\I }}\right) A\right)\cdot S_+^{\sigma}\left(u_0,
{\rm Ad}\left(\overrightarrow{\prod_{I\in RT_\sigma(u)} }\left(\frac{r_I}{r_{I'}}\right)^{\frac{\delta_I(A)}{2\pi\I }}\right)A\right)\\
&={\rm Ad}\left(\overrightarrow{\prod_{I\in RT_\sigma(u)} }\left(\frac{r_I}{r_{I'}}\right)^{\frac{{\rm log}(\delta_I(S_-)\delta_I(S_+))}{2\pi\I }}\right) (S^\sigma_-(u_0, A)S^\sigma_+(u_0,A)).
\end{align*}
The above identity follows from a manipulation of the equivariance of the connection maps $C(t_\sigma^{I})$ with respect to $U(n)^{I}$ for each vertex $I$, and the (blocked Gauss decomposition) monodromy relation \eqref{monoI}. 
\begin{rmk}
Theorem \ref{opentoboundary} can be understood as a branching rule of the Stokes matrices $S(u,A)$ as $u$ degenerates to $u_0$ according to the associated planar tree $T_\sigma(u)$. It relates the Stokes matrices at $u\in {\h_{\rm reg}}(\mathbb{R})$ to the analytic data of a set of linear systems of lower ranks via the solutions of isomonodromy equation with prescribed asymptotics, as some of components $u_i$ of $u$ collapse (in a comparable speed).
\end{rmk}

As can be seen in this paper, the application of the explicit expression, of Stokes matrices at a caterpillar point, relies on the fact that various involved structures we are interested in are preserved under the closure of Stokes matrices. For example, as for the Poisson structures, let us introduce

\begin{defi}\label{RHTree}
The Riemann-Hilbert-Birkhoff map at a boundary point $u_0\in\overline{U_{\sigma}}$ (with respect to the connected component $U_\sigma$) is
\begin{eqnarray}
\nu^{\sigma}(u_0):\Herm(n)\rightarrow \Herm^+(n); \ A\mapsto P_\sigma S_{-}^{\sigma}(u_0,A) {P_\sigma}^{-1} e^{[A]}P_\sigma S_{+}^{\sigma}(u_0,A){P_\sigma}^{-1}.
\end{eqnarray}
\end{defi}

Then, we can similarly generalize the discussion for $u_{\rm cat}$ to any other boundary point (with a chosen planar embedding), and then prove that the map $\nu^{\sigma}(u_0)$ is a Poisson isomorphism. In this way, we prove that the Poisson geometric nature of the Stokes matrices, i.e., Theorem \ref{Boalchthm}, is preserved under taking the regularized limits from $\h_{\rm reg}(\mathbb{R})$ to any boundary point in $\widetilde{\frak t_{\rm reg}}(\mathbb{R})$. That is

\begin{thm}
The Riemann-Hilbert-Birkhoff map $\nu^{\sigma}(u_0):\Herm(n)\rightarrow \Herm^+(n)$ is a Poisson Diffeomorphism.
\end{thm}

\section{Quantum Stokes matrices and their explicit expression at caterpillar points}\label{qStokes}
This section gives the quantum analog of the results in Section \ref{exStokesviaiso}, i.e., the expression of the regularized limits, as well as the leading terms, of quantum Stokes matrices as $u\rightarrow u_{\rm cat}$, in terms of the Gelfand-Tsetlin basis. In Section \ref{qbeginsection}, we introduce the quantum Stokes matrices of the linear system \eqref{introeq}. In Sections \ref{sec:qStokescat}, we recall the regularized limits of quantum Stokes matrices at the caterpillar point $u_{\rm cat}$. In Section \ref{GZbasis}, we introduce the Gelfand-Tsetlin basis and quantum minors. In Sections \ref{proofofqS} and \ref{qGZsystem}, we obtain the explicit expressions and the leading terms of quantum Stokes matrices as $u\rightarrow u_{\rm cat}$ in terms of the Gelfand-Tsetlin basis. In Section \ref{reglimit} we interpret the regularized limit of quantum Stokes matrices via the viewpoint of isomonodromy deformation. In the end, in Section \ref{reglimit2}, we prove that the quantum Stokes matrices at $u_{\rm cat}$ give rise to representation of quantum groups.

\subsection{Quantum Stokes matrices}\label{qbeginsection}
In this subsection, we recall the quantum Stokes matrices of the linear system \eqref{introqeq} associated to a representation $L(\lambda)$
\begin{equation}\label{texteq}
\frac{dF_h}{dz}=h\left(\I u+\frac{1}{2\pi \I }\frac{T}{z}\right)\cdot F_h.
\end{equation}
First, since $h$ is a real number, the system is nonresonant and has a unique formal fundamental solution, see \cite{Xu2} for more details.
\begin{pro}\label{uniformal}
For any nonzero real number $h$ and $u\in\h_{\rm reg}(\mathbb{R})$, the ordinary differential equation \eqref{texteq} has a unique formal fundamental solution taking the form \begin{eqnarray}\label{formalsum}
\widehat{F_h}(z)=\widehat{H_h}(z) e^{{h\I uz}}z^{h[T]}, \ \ \ {\it for} \ \widehat{H_h}=1+H_1(h)z^{-1}+H_2(h)z^{-2}+\cdot\cdot\cdot, \end{eqnarray}
where each coefficient $H_m(h)\in{\rm End}(L(\lambda))\otimes{\rm End}(\mathbb{C}^n)$, and $[T]$ denotes the diagonal part of $T$, i.e., $\delta T=\sum_{k} e_{kk}\otimes E_{kk}.$
\end{pro}
\begin{proof}
Plugging \eqref{formalsum} into the equation \eqref{Stokeseq} gives rise to the equation for $\widehat{H_h},$
\begin{eqnarray}
\frac{1}{h}\frac{d\widehat{H_h}}{dz}+\hat{H_h}\cdot \Big(\I u+
\frac{1}{2\pi\I}\frac{[T]}{z}\Big)=\Big( \I u +
\frac{1}{2\pi\I}\frac{T}{z}\Big)\cdot \widehat{H_h}.\end{eqnarray}
Comparing the coefficients of $z^{-m-1}$, we see that $H_m$ satisfies 
\begin{eqnarray}\label{simHm}
[2\pi u, H_{m+1}]=(\frac{2\pi\I m}{h}+T)\cdot H_{m}-H_{m}\cdot   \delta T.\end{eqnarray}
Set $\{E_{kl}\}_{1\le k,l\le n}$ the standard basis of ${\rm End}(\mathbb{C})$. Then \[T=\sum_{k,l} e_{kl}\otimes E_{kl}, \ \text{ and } u= \sum_i 1\otimes u_iE_{ii}.\] Plugging $H_m(h)=\sum_{k,l} H_{m, kl}(h)\otimes E_{kl}$, with $H_{m,kl}(h)\in{\rm End}(L(\lambda))$, into the equation \eqref{simHm} gives rise to 
\begin{align} \nonumber
&\sum_{k,l}(u_k-u_l) H_{m+1, kl}(h)\otimes E_{kl}\\ \label{recuHh}
=&\sum_{k,l}\frac{2\pi m}{h}H_{m,kl}(h)\otimes E_{kl}+
\sum_{k,l,j} e_{kj} H_{m, jl}(h) \otimes E_{kl}-\sum_{k,l}  H_{m, kl}(h) e_{ll}\otimes E_{kl}.
\end{align}
Here $e_{kl}'s$ are understood as elements in ${\rm End}(L(\lambda))$ via the given representation. That is for $k\ne l$ 
\begin{equation}\label{Hknelh}
(u_k-u_l) H_{m+1, kl}=\frac{2\pi m}{h}H_{m,kl}(h)+\sum_{j=1}^n e_{kj} H_{m, jl}(h)- H_{m, kl}(h) e_{ll} \ \in {\rm End}(L(\lambda)),
\end{equation}
and for $k= l$ (replacing $m$ by $m+1$ in \eqref{recuHh}), 
\begin{eqnarray}\label{Hkelh}
0=\sum_{j\ne k} e_{kj} H_{m+1, jk}(h)+\frac{2\pi(m+1)}{h} H_{m+1, kk}(h)+[e_{kk}, H_{m+1, kk}(h)]\ \in {\rm End}(L(\lambda)).
\end{eqnarray}

Suppose $H_m(h)$ is given, let us check that the above recursive relation have a unique solution $H_{m+1}(h)$. First note that, since $u_k\ne u_l$ for $k\ne l$, the identity \eqref{Hknelh} uniquely defines the "off-diangonal" part $H_{m+1, kl}(h)$ ($k\ne l$) of $H_{m+1}(h)$ from $H_m(h)$. Furthermore, since $h$ is real, we have $\frac{2\pi(m+1)}{2}{\rm Id}+{\rm ad}_{e_{kk}}$ is invertible on ${\rm End}(L(\lambda))$ for any integer $m+1$. Thus, the condition \eqref{Hkelh} uniquely defines the "diagonal" part $H_{m+1, kk}(h)$ of $H_{m+1}(h)$ from the off diagonal part. \end{proof}

Note that under a choice of basis in $L(\lambda)$, the system \eqref{texteq} becomes a special case of the equation \eqref{Stokeseq} with rank $n\times{\rm dim}(L(\lambda))$. Then one can follow the standard resummation procedure as in Section \ref{beginsection} to study the Stokes phenomenon. In particular, the Borel-Laplace transforms of the formal solutions $\widehat{F_h}$ given in Proposition \ref{uniformal} produce actual solutions with the prescribed asymptotics in the corresponding sectoral regions. And again, let us choose the branch of ${\rm log}(z)$, which is real on the positive real axis, with a cut along $i\mathbb{R}_{\ge 0}$. Analog to Theorem \ref{uniformresum}, we have that 
\begin{thm}\cite{Xu2}\label{quantumsol}
For any $u\in\h_{\rm reg}(\mathbb{R})$, on ${\rm Sect}_\pm$ there
is a unique (therefore canonical) holomorphic fundamental solution $F_{h\pm}(z;u,A)\in {\rm End}(L(\lambda))\otimes {\rm End}(\mathbb{C}^n)$ of \eqref{texteq} such that $F_{h\pm}\cdot e^{-\I huz} z^{\frac{-h[T]}{2\pi \I }}$ can be analytically continued to $S(-\pi,\pi)$ and $S(-2\pi,0)$ respectively, and
\begin{eqnarray*}
\lim_{z\rightarrow\infty}F_{h+}(z;u,A)\cdot e^{-\I h uz}\cdot z^{\frac{-h[T]}{2\pi\I }}&=&1, \ \ \ as \ \ \ z\in {S}(-\pi,\pi),
\\
\lim_{z\rightarrow\infty}F_{h-}(z;u,A)\cdot e^{-\I h uz}\cdot z^{\frac{-h[T]}{2\pi\I }}&=&1, \ \ \ as \ \ \ z\in {S}(-2\pi,0),
\end{eqnarray*}
where $[T]={\rm diag}(E_{11},...,E_{nn})\in {\rm End}(L(\lambda))\otimes {\rm End}(\mathbb{C}^n)$. 
\end{thm}

\begin{defi}\cite{Xu2}
The {\it quantum Stokes matrices} of \eqref{texteq} (with respect
to ${\rm Sect}_+$ and the chosen branch of ${\rm log}(z)$) are the elements $S_{h\pm}(u)\in {\rm End}(L(\lambda))\otimes {\rm End}(\mathbb{C}^n)$ determined by
\begin{eqnarray}\label{defi:qStokes}
F_{h+}(z)=F_{h-}(z)\cdot e^{\frac{h[T]}{2}} S_{h+}(u), \  \ \ \ \ 
F_{h-}(ze^{-2\pi\I })=F_{h+}(z)\cdot S_{h-}(u)e^{\frac{-h[T]}{2}},
\end{eqnarray}
where the first (resp. second) identity is understood to hold in ${\rm Sect}_-$
(resp. ${\rm Sect}_+$) after $ F_{h+}$ (resp. $F_{h-}$) has been analytically continued anticlockwise around $z=\infty$. 
\end{defi}
Let us assume $u\in U_{\rm id}\subset\h_{\rm reg}(\mathbb{R})$, then just like the classical case, the asymptotics in Theorem \ref{quantumsol} ensures that $S_{h+}$ is a upper triangular matrix, and $S_{h-}$ is lower triangular with entries in ${\rm End}(L(\lambda))$, see e.g., \cite{Xu2}. 

\begin{rmk}
Following \cite{Xu2}, the (formal solution and quantum Stokes matrices of) linear system \eqref{texteq} is interpreted as a quantization of the (ones of) linear system \eqref{introeq} in the framework of deformation quantization. The study of the quantum Stokes matrices is generalized to arbitrary order pole cases: in \cite{Xu7}, a quantum analog of meromorphic linear systems of ODEs with pole of order $k$, as well as its quantum Stokes matrices, is introduced. As for $k=2$, it becomes the equation \eqref{texteq}. The quantum Stokes matrices at pole of order $k$ is then interpreted as a quantization of the space of the classical Stokes matrices. 
\end{rmk}

\subsection{Quantum Stokes matrices at a caterpillar point}\label{sec:qStokescat} 
Since \eqref{texteq} can be seen as a special form of the general system \eqref{Stokeseq}, we have
\[S_{h\pm}(u)=S_\pm(hu,-hT),\]
where the right hand side is understood as the (classical) Stokes matrices of the equation \eqref{Stokeseq} with $u$ and $A$ replaced by $hu$ and $-hT$. Here the irregular term $\I hu= {\rm diag}(\I hu_1,...,\I hu_n)\in {\rm End}(L(\lambda))\otimes {\rm End}(\mathbb{C}^n)$ is degenerate, i.e., has repeated eigenvalues.
One can still study the regularized limits of $S_{h\pm}(u)$ as the components $u_i$ of $u\in\h_{\rm reg}(\mathbb{R})$ collapse in a comparable speed, and introduce the quantum Stokes matrices for any boundary point in $\widetilde{\frak t_{\rm reg}}(\mathbb{R})\setminus {\h_{\rm reg}}(\mathbb{R})$. In particular, the construction in Section \ref{dirreg} enables us to define the quantum Stokes matrices at $u_{\rm cat}$ as follows.

First, the system \eqref{texteq} has rank $n\times {\rm dim}(L(\lambda))$. In terms of the notation in Section \ref{dirreg}, let us take the partition $\underline{d}$ of $n\times {\rm dim}(L(\lambda))$ with $d_1=\cdots=d_n={\rm dim}(L(\lambda))$, and assume that the irregular term $u$ in the Stokes matrices $S_\pm(hu,-hT)$ of \eqref{texteq} lives in $U^{\underline{d}}_{\rm id}$. Following Theorem \ref{dreglimit}, the matrix $S_-(hu,-hT)\cdot S_+(hu,-hT)$ has a regularized limit as $u\rightarrow u_{\rm cat}$ from $u\in U^{\underline{d}}_{\rm id}$. 

\begin{defi}\label{qScat}
The quantum Stokes matrices at $u_{\rm cat}$, with respect to the choice of $U^{\underline{d}}_{\rm id}$, are the upper and lower $n\times n$ triangular matrices $S_{h\pm}(u_{\rm cat})$ (having the same diagonal part), with entries valued in ${\rm End}(L(\lambda))$, such that (the blocked Gauss decomposition)
\[S_{h-}(u_{\rm cat})S_{h+}(u_{\rm cat})\]
equals to the regularized limit of the function $S_-(hu,-hT)\cdot S_+(hu,-hT)$ on $U^{\underline{d}}_{\rm id}$, i.e., the limit of 
\begin{equation}\label{gauge}
{\rm Ad}{\left((\frac{1}{u_{2}-u_{1}})^{\frac{{\rm log}(\delta^{\underline{d}}_1(S_{h-})\delta^{\underline{d}}_1(S_{h+}))}{2\pi\I }}\overrightarrow{\underset{k=2,...,m-1}{\prod} }(\frac{u_{k}-u_{k-1}}{u_{k+1}-u_{k}})^{\frac{{\rm log}(\delta^{\underline{d}}_k(S_{h-})\delta^{\underline{d}}_k(S_{h+}))}{2\pi\I }}\right)}  \left(S_{h-}(u) S_{h+}(u)\right),
\end{equation}
as $\frac{u_{k+1}-u_{k}}{u_{k}-u_{k-1}}\rightarrow +\infty$ for all $k=2,...,m-1$ and $u_2-u_1\rightarrow 0$,
where $\delta^{\underline{d}}_k(S_{h\pm}(u))$ takes the ${\rm End}(\lambda)$ valued entries of the matrices $S_{h\pm}(u)$ as in \eqref{bdelta}.
\end{defi}

\subsection{Quantum minors and Gelfand-Tsetlin basis}\label{GZbasis}
To write down the explicit expression of $S_{h\pm}(u_{\rm cat})$, as a non-commutative version of the formula in Theorem \ref{mainthm}, let us introduce

\begin{defi}
    The matrix $T(x)=T-x \operatorname{Id}$ is the characteristic matrix of $T$, where $x$ is an indeterminate commuting with all generators $e_{ij}$, for $1\leq i,j \leq n$.
\end{defi}

\begin{defi}
For $1\leq m \leq n$, given two sequences, $a=\{a_1,\cdots,a_m\}$ and $b=\{b_1,\cdots,b_m\}$, with elements in $\{1,2,\cdots,n-1,n\}$, the corresponding quantum
minor of the matrix $T(x)$ is defined by 
\begin{align*}
\Delta^{a_1,...,a_m}_{b_1,...,b_m}(T(x)):=\sum_{\sigma\in S_m}(-1)^\sigma T(x_1)_{a_{\sigma(1)},b_1}\cdots T(x_m)_{a_{\sigma(m)},b_m}\in U({\frak {gl}}_n)[x],
\end{align*}
where $x_k:=x+k-1$ for $k=1,...,m$, and $(-1)^{\sigma}$ means the signature of the permutation $\sigma$ in the symmetry group $S_m$ of $m$ elements.
\end{defi}

Since $\zeta^{(k)}_i$ does not commute with all $e_{ab}$, for $1\leq a, b \leq n$. So $\Delta^{a_1,...,a_m}_{b_1,...,b_m}\big(T(\zeta^{(k)}_i)\big)$ is not well defined. However,

\begin{defi}\label{lrminor}
Suppose that $\Delta^{a_1,...,a_m}_{b_1,...,b_m}\left(T\left(x\right)\right)=\sum_{i=0}^{n-1}r_i x^i$ with coefficient $r_i\in U\left({\mathfrak {gl}}_n\right)$. For any element $\zeta$ (that may not commute with $e_{ij}$ for all $1\le i,j\le n$), we define $\left(\Delta_L\right)^{a_1,...,a_m}_{b_1,...,b_m}\big(T(\zeta)\big)$ and $\left(\Delta_R\right)^{a_1,...,a_m}_{b_1,...,b_m}\big(T(\zeta)\big)$ in the following way
\begin{align*}
&\left(\Delta_L\right)^{a_1,...,a_m}_{b_1,...,b_m}\big(T(\zeta)\big)=\zeta^{n-1} r_{n-1}+\zeta^{n-2} r_{n-2}+\cdots +\zeta r_{1}+r_0,\\
&\left(\Delta_R\right)^{a_1,...,a_m}_{b_1,...,b_m}\big(T(\zeta)\big)=r_{n-1} \zeta^{n-1}+ r_{n-2} \zeta^{n-2}+\cdots +r_{1} \zeta+r_0.
\end{align*}
\end{defi}

For any $1\le k\le n$, let 
\begin{eqnarray}
M_k(\zeta):=\Delta^{1,...,k}_{1,...,k}\left(T(\zeta)\right)
\end{eqnarray}
be the upper left $k\times k$ quantum-minor of $T$. It is known that the subalgebra, generated in $U({\frak {gl}}_n)$ by the coefficients in all $M_k(\zeta)$ for all $1\le k\le n$, is a maximal commutative subalgebra and is called the Gelfand-Tsetlin subalgebra. 

\begin{defi}\label{qroots}
Let $\{\zeta^{(k)}_i\}_{1\le i\le k\le n}$ denote the roots of $M_k(\zeta) = 0$ for all $k=1,...,n$ (in an appropriate splitting extension).
\end{defi}

The action of the Gelfand-Tsetlin subalgebra on a highest weight representation has simple spectrum and the corresponding orthonomral eigenbasis is called a Gelfand-Testlin basis. The action of the quantum minors and the roots $\{\zeta^{(k)}_i\}_{1\le i\le k\le n}$ on the representation can be explicitly expressed under the basis, see Proposition \ref{GZaction}. Before giving the explicit action, let us introduce the Gelfand-Tsetlin basis in a more conventional way. 

Recall that $\{e_{ij}\}_{i,j=1,...,n}$ is the standard basis of the Lie algebra ${\gl}_n$. Denote by ${\gl}_{n-1}$ the subalgebra spanned by the elements  $\{e_{ij}\}_{i,j=1,...,n-1}$. Finite dimensional irreducible representations of ${\gl}_n$ are parameterized by the highest weight, i.e., $n$-tuples of numbers $\lambda=(\lambda_1,...,\lambda_n)$ with
\[\lambda_i-\lambda_{i+1}\in\mathbb{Z}_+, \hspace{5mm}, \ \forall \ i=1,...,n-1.\]
We denote by $L(\lambda)$ the corresponding representation. It has a lowest vector $\xi_0$ such that $E_{ii}\xi =\lambda_i \xi$ for $i=1,...,n$, and $E_{ij}\xi=0$ for $1\le j<i\le n$. Then the simple branching rule for the reduction from ${\gl}_n$ to ${\gl}_{n-1}$ states that the restriction of $L(\lambda)$ to the subalgebra ${\gl}_{n-1}$
is isomorphic to the
direct sum of pairwise inequivalent irreducible representations
\[L(\lambda)|_{{\gl}_{n-1}}\cong \underset{\lambda^{(n-1)}}{\bigoplus} L'(\lambda^{(n-1)}),\]
where the summation is over the highest weights $\lambda^{(n-1)}$ of ${\gl}_{n-1}$ satisfying the interlacing conditions
\begin{eqnarray}\label{interlacing}
\lambda_i^{(n)}-\lambda^{(n-1)}_i\in\mathbb{Z}_{\ge 0}, \hspace{5mm} \lambda^{(n-1)}_i-\lambda^{(n)}_{i+1}\in\mathbb{Z}_{\ge 0}, \hspace{5mm} \forall \ i=1,...,n-1.
\end{eqnarray}
Thus a chain of subalgebras
\[ {\gl}_1\subset\cdots\subset {\gl}_{n-1}\subset {\gl}_n\]
produces a decomposition of $L(\lambda)$ into one dimensional subspaces, and the one dimensional subspaces are parameterized by the Gelfand-Tsetlin patterns. Such a pattern $\Lambda$ (for fixed $\lambda^{(n)}:=\lambda$) is a collection of numbers $\{\lambda^{(i)}_j(\Lambda)\}_{1\le i\le j\le n-1}$ satisfying the interlacing conditions
\begin{equation}\label{interineq}
\lambda_j^{(i)}(\Lambda)-\lambda^{(i-1)}_j(\Lambda)\in\mathbb{Z}_{\ge 0}, \hspace{5mm} \lambda^{(i-1)}_j(\Lambda)-\lambda^{(i)}_{j+1}(\Lambda)\in\mathbb{Z}_{\ge 0}, \hspace{5mm} \forall \ i=1,...,n-1.
\end{equation}
\begin{defi}
We denote by $P_{GT}(\lambda;\mathbb{Z})$ the set of Gelfand-Tsetlin patterns in $L(\lambda)$, seen as the set of integer points in the real Gelfand-Tsetlin polytope $P_{GT}(\lambda;\mathbb{R})$. 
\end{defi}
Given the decomposition of $L(\lambda)$ into the one dimensional subspaces, one gets a basis by choosing a nonzero vector from each subspace. In particular, there exists a basis $\xi_\Lambda$ of $L(\lambda)$, called a Gelfand-Tsetlin basis, parameterized by all patterns $\Lambda\in P_{GT}(\lambda;\mathbb{Z})$. The basis is denoted by $\xi_\Lambda(u_{\rm cat})$ in the introduction, for simplicity let us drop the symbol $u_{\rm cat}$. The structure of the basis obtained in this way is summarized in the following proposition (see e.g., \cite{Molev}). 
\begin{pro}\label{GZaction}
There exists an orthonormal basis $\xi_\Lambda $ of $L(\lambda)$, called the Gelfand-Tsetlin basis, parameterized by all patterns $\Lambda\in P_{GZ}(\lambda;\mathbb{Z})$, such that for any $1\le i\le k\le n$, the actions of $e_{kk}$, $\zeta^{(k)}_i$ and $\alpha^{(k)}_i$ on the basis ${\xi_\Lambda}$ of $L(\lambda)$ are given by
\begin{align}\label{eigenvalues0}
e_{kk}\cdot \xi_\Lambda &=\left(\sum_{i=1}^k\lambda^{(k)}_i(\Lambda)-\sum_{i=1}^{k-1}\lambda^{(k-1)}_i(\Lambda)\right)\xi_\Lambda , \\
\label{eigenvalues}
\zeta^{(k)}_i\cdot \xi_\Lambda &=\left(\lambda^{(k)}_i(\Lambda)-i+1\right)\xi_\Lambda ,
\end{align}
and
\begin{align}\nonumber
&(\Delta_L)^{1,...,k}_{1,...,k-1,k+1}\left(T(\zeta^{(k)}_i)\right) \cdot \xi_\Lambda \\ \label{ashift}
=&(-1)^{k+i}\sqrt{-\frac{\prod_{l=1,l\ne i}^{k}(\zeta^{(k)}_i-\zeta^{(k)}_l)\prod_{l=1}^{k+1}(\zeta^{(k)}_i-\zeta^{(k+1)}_l-1)\prod_{l=1}^{k-1}(\zeta^{(k)}_i-\zeta^{(k-1)}_l)}{\prod_{l=1,l\ne i}^{k}(\zeta^{(k)}_i-\zeta^{(k)}_l-1)}}\cdot \xi_{\Lambda+\delta^{k}_i} 
\end{align}
where the pattern $\Lambda+\delta^{(k)}_i$ is obtained from $\Lambda$ by replacing $\lambda^{(k)}_i$ by $\lambda^{(k)}_i+ 1$. It is supposed that $\xi_\Lambda$ is zero if $\Lambda$ is not a pattern.
\end{pro}
\pf 
We refer to \cite[Sectioin 2]{Molev} for the existence of the orthogonal Gelfand-Tsetlin basis $\{\eta_\Lambda\}$, and the norms $|\eta_\Lambda|$ of the elements in the basis \cite[Proposition 2.4]{Molev} (we remark that the basis $\{\xi_\Lambda\}$ is the orthonormal basis in the original work of Gelfand and Tsetlin \cite{GT}, and is different from the orthogonal one $\{\eta_\Lambda\}$ from \cite{Molev}). We then perform the explicit normalization to get the orthonormal basis $\{\xi_\Lambda\}=\eta_\Lambda/|\eta_\Lambda|$ from $\{\eta_\Lambda\}$. Then the actions of $e_{kk}$ and $\zeta^{(k)}_i$ on the orthogonal basis $\{\eta_\Lambda\}$ given in \cite[Theorem 2.3]{Molev} imply the identities \eqref{eigenvalues0} and \eqref{eigenvalues}. Furthermore, the identity \eqref{ashift} follows from the action of $(\Delta_L)^{1,...,k}_{1,...,k-1,k+1}$ on the orthogonal basis $\{\eta_\Lambda\}$ given in \cite[Section 2.5]{Molev}, provided the normalization, from the orthogonal to the orthonormal basis, is accounted for. We also refer the reader to \cite{LX} for a detailed proof of \eqref{ashift}. \qed

\subsection{The explicit expression of quantum Stokes matrices at a caterpillar point}\label{proofofqS}

\begin{thm}\label{explicitS}
For any $1\le k\le n-1$, the $(k,k+1)$-entry of $S_{h+}(u_{\rm cat})$, as an element in ${\rm End}(L(\lambda))$ is given by
\begin{align*}
(S_{h+})_{k,k+1}&=
(\I h)^{\frac{h\small{(e_{kk}-e_{k+1,k+1}-1})}{2\pi \I }} e^{-\frac{he_{kk}}{2}}\\
\times & \sum_{i=1}^k\left(\frac{\prod_{l=1,l\ne i}^{k}\Gamma(h\frac{\zeta^{(k)}_i-\zeta^{(k)}_l}{2\pi \I })}{\prod_{l=1}^{k+1}\Gamma(1+h\frac{\zeta^{(k)}_i-\zeta^{(k+1)}_l-1}{2\pi \I })}\frac{\prod_{l=1,l\ne i}^{k}\Gamma(1+h\frac{\zeta^{(k)}_i-\zeta^{(k)}_l-1}{2\pi \I })}{\prod_{l=1}^{k-1}\Gamma(1+h\frac{\zeta^{(k)}_i-\zeta^{(k-1)}_l}{2\pi \I })}\right)\cdot  (\Delta_L)^{1,...,k}_{1,...,k-1,k+1}\left(\frac{h}{2\pi\I }T(\zeta^{(k)}_i)\right),\end{align*}
and the $(k+1,k)$-entry of $S_{h-}(u_{\rm cat})$ is given by
\begin{align*}
(S_{h-})_{k+1,k}&=
(\Delta_R)_{1,...,k}^{1,...,k-1,k+1}\left(\frac{-h}{2\pi\I }T(\zeta^{(k)}_i)\right) \\
\times &
\left(\sum_{i=1}^k \frac{\prod_{l=1,l\ne i}^{k}\Gamma(-h\frac{\zeta^{(k)}_i-\zeta^{(k)}_l}{2\pi \I })}{\prod_{l=1}^{k+1}\Gamma(1-h\frac{\zeta^{(k)}_i-\zeta^{(k+1)}_l+1}{2\pi \I })}\frac{\prod_{l=1,l\ne i}^{k}\Gamma(1-h\frac{\zeta^{(k)}_i-\zeta^{(k)}_l-1}{2\pi \I })}{\prod_{l=1}^{k-1}\Gamma(1-h\frac{\zeta^{(k)}_i-\zeta^{(k-1)}_l}{2\pi \I })}\right)(-\I h)^{\frac{-h\small{(e_{kk}-e_{k+1,k+1}-1})}{2\pi \I }} e^{\frac{-he_{kk}}{2}} .
\end{align*}
Here for any constant $c$, $\Delta^{1,...,k}_{1,...,k-1,k+1}\left(cT(\zeta^{(k)}_i)\right):=c^k\Delta^{1,...,k}_{1,...,k-1,k+1}\left(T(\zeta^{(k)}_i)\right)$.
\end{thm}

\pf 
The above explicit expressions can be derived in the same way as the ones in Theorem \ref{mainthm}, except that in the non-commutative setting we should use Lemma \ref{qGZ} to exchange the orders in the product of the ${\rm End}(L(\lambda))$ valued functions $\zeta^{(k)}_i$ and $\alpha^{(k)}_i$.  

First, following Definition \ref{qScat} and Theorem \ref{dreglimit}, we have
\begin{eqnarray}\label{qconnection}
S_{h-}(u_{\rm cat})S_{h+}(u_{\rm cat})=\left(\overrightarrow{\prod_{k=2,...,n}}C_h(E_k,\delta_k(T))\right)\cdot e^{hT}\cdot \left(\overrightarrow{\prod_{k=2,...,n}}C_h(E_k,\delta_k(T))\right)^{-1},
\end{eqnarray}
where each $C_h(E_{k+1},\delta_{k+1}(T))$ denotes the connection matrix of
\begin{eqnarray}\label{eq:qrelativeStokes}
\frac{dF}{dz}=h\left(\I E_{k+1}+\frac{1}{2\pi \I }\frac{\delta_{k+1}(T)}{z}\right)\cdot F.
\end{eqnarray}
Here (by abuse of notation) $E_{k+1}\in {\rm End}(\mathbb{C}^{n})\otimes {\rm End}(L(\lambda))$ denotes the matrix whose $(k+1,k+1)$ entry is $1\in {\rm End}(L(\lambda))$ and other entries are zero. Therefore, to get the explicit expressions of $S_{h\pm}(u_{\rm cat})$, we need to compute each $C_h(E_{k+1},\delta_{k+1}(T))$.

Now the key observation is that, just like the classical case, the equation \eqref{eq:qrelativeStokes} has a fundamental solution given by confluent hypergeometric function $_kF_k$. 

\vspace{2mm}

{\bf Diagonalization in stages.} Just as in the classical case, to simplify the equation \eqref{eq:qrelativeStokes}, let us diagonalize the upper left $k$-th submatrix of its coefficient matrix. 
Recall that in the classical case, the diagonalization in stages are only taken on the open dense subset $\Herm(n)_0$ of $\Herm(n)$. In the quantum case, the analog of the space $\Herm(n)_0$ is a subspace of the representation $L(\lambda)$ spanned by the basis elements $\xi_\Lambda$ such that the inequalities in \eqref{interineq} are strict, i.e., 
\[L(\lambda)_0:={\rm span}\{\xi_\Lambda~|~\lambda_j^{(i)}(\Lambda)-\lambda^{(i-1)}_j(\Lambda)\in\mathbb{Z}_{>0}, \ \lambda^{(i-1)}_j-\lambda^{(i)}_{j+1}\in\mathbb{Z}_{>0}, \ \forall \ i, j \}\subset L(\lambda).\]
However, just as in the classical setting, see Remark \ref{clsdiag} and \ref{smoothex}, the Stokes and connection matrices of \eqref{eq:qrelativeStokes} have no singularities. It is only the diagonalization method of computation (the matrix $P_k$) that introduce the singularities. Therefore, to compute the Stokes and connection matrices of the equation \eqref{eq:qrelativeStokes}, it is enough to do the diagonalization in stages formally. Thus for simplicity, in below we do the computation in a formal setting, i.e., we formally introduce $\mathcal{P}_k$ and $\mathcal{Q}_k$ as below ignoring the pole issue in the denominator. 

Recall that $\zeta^{(k)}_1,...,\zeta^{(k)}_k$ of $M_k(\zeta)$ denote the roots of the quantum minor $\Delta^{1,...,k}_{1,...,k}\left(T(\zeta-\frac{h}{2}(k-1))\right)$, which act diagonally on the Gelfand-Zeitlin basis in $L(\lambda)$. 
For any $2\le k\le n$, let us take the $n\times n$ matrix $\mathcal{Q}_k$,
\begin{align*}
(\mathcal{Q}_k)_{ij}&=\frac{(-1)^{k+j}}{\sqrt{\prod_{l=1,l\ne i}^k(\zeta^{(k)}_i-\zeta^{(k)}_l)\prod_{l=1}^{k-1}(\zeta^{(k)}_i-\zeta^{(k-1)}_l)}}(\Delta_L)^{1,...,\hat{j},...,k}_{1,...,k-1}\left(T(\zeta^{(k)}_i)\right), \ \text{if} \ 1\le i,j\le k\\ 
\nonumber (\mathcal{Q}_k)_{ii}&=1, \ \ \text{if} \ i>k, \\ 
\nonumber (\mathcal{Q}_k)_{ij}&=0, \ \ \text{otherwise},
\end{align*}
and the $n\times n$ matrix $\mathcal{P}_k$ with entries in ${\rm End}(L(\lambda))$,
\begin{align*}(\mathcal{P}_k)_{ij}&=(\Delta_R)_{1,...,\hat{i},...,k}^{1,...,k-1}\left(T(\zeta^{(k)}_j)\right)\cdot \frac{(-1)^{k+i}}{\sqrt{\prod_{l=1,l\ne j}^k(\zeta^{(k)}_j-\zeta^{(k)}_l)\prod_{l=1}^{k-1}(\zeta^{(k)}_j-\zeta^{(k-1)}_l)}}, \ \text{if} \ 1\le i,j\le k\\ 
\nonumber (\mathcal{P}_k)_{ii}&=1, \ \ \text{if} \ i>k, \\ 
\nonumber (\mathcal{P}_k)_{ij}&=0, \ \ \text{otherwise}.
\end{align*}
The following lemmas (the analog of the classical linear algebra facts in Section \ref{diagonalization}) follow from the Laplace expansion of quantum minors, see e.g., \cite{LX}.
\begin{lem}\label{qGZ}
We have \[\mathcal{P}_k\cdot\mathcal{Q}_k=\mathcal{Q}_k\cdot \mathcal{P}_k={\rm Id}_n,\] and
the matrix $\delta_{k+1}(T_{k}):=\mathcal{Q}_k\cdot \delta_{k+1}(T)\cdot \mathcal{P}_k$ takes the form
\[\delta_{k+1}(T_{k})=\begin{pmatrix}
    \zeta^{(k)}_1+(k-1) & & & &\alpha^{(k)}_1 & \star \\
    & \ddots & && \vdots & \star\\
    & &  \zeta^{(k)}_k+(k-1) && \alpha^{(k)}_k & \star \\
    &&&&\\
    \beta_1^{(k)}&\cdots &\beta_k^{(k)}&& E_{k+1} & 0\\
    \star & \star & \star && \star & {\rm diag}(E_{k+2,k+2},...,E_{nn})
    
    \end{pmatrix}.
\]
Here 
\begin{eqnarray}\label{a}
\alpha^{(k)}_j(T)=\frac{(-1)^{k+j}}{\sqrt{\prod_{l=1,l\ne i}^k(\zeta^{(k)}_i-\zeta^{(k)}_l)\prod_{l=1}^{k-1}(\zeta^{(k)}_i-\zeta^{(k-1)}_l)}}(\Delta_L)^{1,...,k}_{1,...,k-1,k+1}\left(T(\zeta^{(k)}_i)\right),
\end{eqnarray}
and for any $1\le k\le n-1$, and $1\le i,j\le k$, we have the commutators
\begin{eqnarray}\label{commutator}
[\zeta^{(k)}_i,\alpha^{(k)}_j]=\delta_{ij}\alpha^{(k)}_j.
\end{eqnarray}
\end{lem}    
\begin{lem}
For any $1\le k\le n-1$, we have
\begin{equation}
\mathcal{P}_k=\mathcal{L}^{(k+1)}\mathcal{P}_{k+1}^{-1} 
\end{equation}
where the $n\times n$ matrix $\mathcal{L}^{(k+1)}$ is 
\begin{eqnarray}
&&\mathcal{L}^{(k+1)}_{ij}=\alpha^{(k)}_i\cdot \frac{1}{\mathcal{N}_j^{(k+1)}(\zeta^{(k+1)}_j-\zeta^{(k)}_i)}, \ \ \ 1\le i\le k, \ \ 1\le j\le k+1,\\ \label{L2}
&&\mathcal{L}^{(k+1)}_{k+1,j}=\frac{1}{\mathcal{N}_j^{(k+1)}}, \ \ \ \ \ 1\le j\le k+1.\\
&&\nonumber \mathcal{L}^{(k+1)}_{ii}=1, \ \ if \ i>k+1, \ \ \ \mathcal{L}^{(k+1)}_{ij}=0, \ otherwise.
\end{eqnarray}
Here the normalizer is
\begin{equation}
\mathcal{N}_j^{(k+1)}(T)=\sqrt{\frac{\prod_{v=1, v\ne j}^{k+1}(\zeta^{(k+1)}_j-\zeta^{(k+1)}_v) }{\prod_{v=1}^{k}(\zeta^{(k+1)}_j-\zeta^{(k)}_v)}}.
\end{equation}
\end{lem} 

{\bf Confluent hypergeometric functions in representation spaces.} Using Lemma \ref{qGZ} to "formally" diagonalize the upper left $k$-th submatrix of the coefficient of \eqref{eq:qrelativeStokes}, we are left with the equation
\begin{eqnarray}\label{qrelequation}
\frac{dF}{dz}=h\left(\I E_{k+1}+\frac{1}{2\pi \I }\frac{\delta_{k+1}(T_{k})}{z}\right)\cdot F.
\end{eqnarray}
A solution of this equation, as a function valued in ${\rm End}(\mathbb{C}^n)\otimes {\rm End}(L(\lambda)_0)$, is described using the generalized confluent hypergeometric functions as follows. First recall that the confluent hypergeometric functions associated to any $\alpha_j\in\mathbb{C}$ and $\beta_j\in\mathbb{C}\setminus \{0, -1, -2, ...\}, 1\le j\le m$, are
\begin{eqnarray}\label{qgchf}
_kF_k(\alpha_1,...,\alpha_k,\beta_1,...,\beta_k;z)=\sum_{n=0}^\infty \frac{(\alpha_1)_n\cdots(\alpha_k)_n}{(\beta_1)_n\cdots(\beta_k)_n}\frac{z^n}{n!},
\end{eqnarray}
where $(\alpha)_0=1$, $(\alpha)_n=\alpha\cdots (\alpha+n-1),$ $n\ge 1$. Since $\{\zeta^{(k)}_i\}_{1\le i\le k\le n}$ are commutative elements in ${\rm End}(L(\lambda))$, we can introduce the following confluent hypergeometric functions valued in ${\rm End}(L(\lambda))$,
\begin{align*}
H(z)_{ij}=\scriptsize{\frac{1}{h(\zeta^{(k+1)}_j-\zeta^{(k)}_i)}}\cdot {_kF_k}(\alpha_{ij,1},...,\alpha_{ij,k},\beta_{ij,1},...,\widehat{\beta_{ij,j}},...,\beta_{ij,k+1};\I z), \ \ \ 1\le i\le k, 1\le j\le k+1,\\
H(z)_{k+1,j}= {_kF_k}(\alpha_{k+1j,1},...,\alpha_{k+1j,k},\beta_{k+1j,1},...,\widehat{\beta_{k+1j,j}},...,\beta_{k+1j,k+1};\I z), \ \ i=k+1, \ 1\le j\le k+1,
\end{align*}
with the variables $\{\alpha_{ij,l}\}$ and $\{\beta_{ij,l}\}$ given by
\begin{eqnarray*}
&&\alpha_{ij,i}=\frac{h}{2\pi \I }\left(\zeta^{(k+1)}_j-\zeta^{(k)}_i\right), \ \ \ 1\le i\le k, \ 1\le j\le k+1,\\
&&\alpha_{(k+1,j), i}=1+\frac{h}{2\pi \I }\left(\zeta^{(k+1)}_j-\zeta^{(k)}_{i}\right), \ \ \ 1\le j\le k+1,\\
&&\alpha_{ij, l}=1-\frac{h}{2\pi \I }\left(\zeta^{(k+1)}_j-\zeta^{(k)}_l\right), \ \ \ \ l\ne i, \  1\le l\le k, \ 1\le i, j\le k+1, \\
&&\beta_{ij, l}=1+\frac{h}{2\pi \I }\left(\zeta^{(k+1)}_j-\zeta^{(k+1)}_l\right), \ \ \ \ l\ne j, \ 1\le l\le k+1, \ \ 1\le i, j\le k+1.
\end{eqnarray*}

\begin{lem}\label{qsol}
The equation \eqref{qrelequation} has a solution $F_h(z)\in {\rm End}(\mathbb{C}^{n})\otimes {\rm End}(L(\lambda))$ taking the form
\begin{eqnarray}
F_h(z)=\left(\begin{array}{cc}
{\rm diag}(ha^{(k)}_1, ..., ha^{(k)}_k,1) \ \ & 0\\
0 \ \ & {\rm Id}_{n-k-1}
\end{array} \right) \left(\begin{array}{cc}
\left(H_{ij}(z)\right)_{i,j=1}^{k+1} \ \ & 0\\
0 \ \ & {\rm Id}_{n-k-1}
\end{array} \right)
\cdot z^{\frac{h[T_{k+1}]}{2\pi\I }},
\end{eqnarray}
where $[T_{k+1}]$ is the diagonal part of $T_{k+1}$, i.e.,
\[{ [T_{k+1}]}=\left(\begin{array}{cc}
{\rm diag}({\zeta^{(k+1)}_1+k},...,{\zeta^{(k+1)}_{k+1}+k}) \ \ & 0\\
0 \ \ & {\rm diag}({E_{k+2,k+2}},...,{E_{n,n}})
\end{array} \right).\]
\end{lem}
\begin{proof}
We can plug the expression into the equation \eqref{qrelequation}. Then the first $k$ rows of the equation can be verified directly. And the $k+1$-th row can be verified using the following identity for $m$ a parameter
\begin{align*}
    &\frac{\prod_{i=1}^{k+1}\big(m+\frac{h \zeta^{\left({k+1}\right)}_j}{2\pi \mathrm{i} }-\frac{h \zeta^{\left({k+1}\right)}_i}{2\pi \mathrm{i} }\big)}{\prod_{i=1}^{k}\big(m+\frac{h \zeta^{\left({k+1}\right)}_j}{2\pi \mathrm{i} }-\frac{h \zeta^{\left(k\right)}_i}{2\pi \mathrm{i} }\big)}\\
    =&m+\frac{h \zeta^{\left(n\right)}_j}{2\pi \mathrm{i} }-\frac{h\left(e_{{k+1},{k+1}}-k\right)}{2\pi \mathrm{i} }-\sum_{l=1}^{k}\frac{h \beta^{\left(k\right)}_l}{2\pi \mathrm{i} } \frac{1}{m+\frac{h \zeta^{\left({k+1}\right)}_j}{2\pi \mathrm{i} }-\frac{h}{2\pi \mathrm{i} }\big(\zeta^{\left(k\right)}_l-1\big)}\frac{h \alpha^{\left(k\right)}_l}{2\pi \mathrm{i} }.
\end{align*}  
The identity in turn can be verified by comparing the actions of its left and right hand sides on the Gelfand-Tsetlin basis. See \cite{LX} for more details.
\end{proof}
We denote by $C_h(E_{k+1},\delta_{k+1}(T_k))$ the connection matrix of the equation \eqref{qrelequation}.
Then the Stokes matrices $S_{h+}(E_{k+1},\delta_{k+1}(T_k))$ and {\it the normalized connection matrix} \[\widetilde{C_h}(E_{k+1}, \delta_{k+1}(T)):=C_h(E_{k+1},\delta_{k+1}(T_k))\cdot \mathcal{L}^{(k+1)}\] of the linear system \eqref{qrelequation} are described by
\begin{pro}\label{qexStokes}
$(1)$ The entries $s_{ij}$ of $S_{h+}(E_{k+1},\delta_{k+1}(T_k))$ are
\begin{eqnarray*}
s_{j,k+1}=\frac{(\I h)^{\frac{h(\zeta^{(k)}_j+k-1-e_{k+1,k+1}-1)}{2\pi\I }} e^{\frac{he_{k+1,k+1}}{2}}\prod_{l=1,l\ne j}^{k}\Gamma(1+h\frac{\zeta^{(k)}_j-\zeta^{(k)}_l-1}{2\pi \I })}{\prod_{l=1}^{k+1}\Gamma(1+h\frac{\zeta^{(k)}_j-\zeta^{(k+1)}_l-1}{2\pi \I })}\cdot {\alpha^{(k)}_j}, \ \ \ \text{for} \ j=1,...,k;
\end{eqnarray*}
and
\[s_{ii}=e^{\frac{-e_{ii}}{2}}, \ \ \ i=1,...,n; \hspace{5mm} s_{ij}=0, \ \ otherwise.\]

$(2)$. The entries $c_{ij}$ of the matrix $\widetilde{C_h}(E_{k+1}, \delta_{k+1}(T))$ are
\begin{eqnarray*} c_{ij}=- \frac{(\I h)^{\frac{h(\zeta^{(k)}_i-\zeta^{(k+1)}_j)}{2\pi\I }}e^{\frac{h(\zeta^{(k+1)}_j-2\zeta^{(k)}_i)}{2}}\prod\limits_{l=1}^{k+1}\Gamma(1+h\frac{\zeta^{(k+1)}_j-\zeta^{(k+1)}_l}{2\pi \I })\prod\limits_{l=1}^{k}\Gamma(1+h\frac{\zeta^{(k)}_i-\zeta^{(k)}_l-1}{2\pi \I })}{\prod\limits_{l=1,l\ne i}^{k}\Gamma(1+h\frac{\zeta^{(k+1)}_j-\zeta^{(k)}_l}{2\pi \I })\prod\limits_{l=1,l\ne j}^{k+1}\Gamma(1+h\frac{\zeta^{(k+1)}_l-\zeta^{(k)}_i+1}{2\pi \I })}\cdot \frac{\alpha^{(k)}_i}{2\pi\I {\mathcal{N}^{(k+1)}_j}},
\end{eqnarray*}
for $1\le j\le k+1, 1\le i\le k$;
\begin{eqnarray*} c_{k+1,j}=\frac{(\I h)^{\frac{h(e_{k+1,k+1}-\zeta^{(k+1)}_j-k)}{2\pi\I }}e^{\frac{-he_{k+1,k+1}}{2}} \prod_{l=1}^{k+1}\Gamma(1+h\frac{\zeta^{(k+1)}_j-\zeta^{(k+1)}_l}{2\pi \I })}{{\mathcal{N}^{(k+1)}_j}\prod_{l=1}^{k}\Gamma(1+h\frac{\zeta^{(k+1)}_j-\zeta^{(k)}_l}{2\pi \I })}, \ \ \ \text{for} \ 1\le j\le k+1;
\end{eqnarray*}
and \[\ c_{ii}=1  \ \ \text{for} \ \ k+2\le i\le n; \hspace{5mm} c_{ij}=0 \ \ \text{otherwise}.\]
\end{pro}
Just like in the classical case,
the proposition follows from the known asymptotics of the functions $_kF_k$. In the quantum case, we only need to use the commutative relation \eqref{commutator} to commute $a^{(k)}_j$ and the arguments $\zeta^{(k)}_j$ in the computation. We refer the reader to \cite{LX} for the full details.

\vspace{2mm}
{\bf Subdiagonal entries of quantum Stokes matrices.} By the equivariant property, when restricts to $L(\lambda)_0\subset L(\lambda)$, the connection matrices of the equations \eqref{eq:qrelativeStokes} and \eqref{qrelequation} are related by
\[C_h(E_{k+1},\delta_{k+1}(T))=\mathcal{P}_k\cdot C_h(E_{k+1},\delta_{k+1}(T_k))\cdot \mathcal{Q}_k, \]
where $\mathcal{P}_k$ is given in Lemma \ref{qGZ}. By the relation $\widetilde{C_h}(E_{k+1}, \delta_{k+1}(T))=C_h(E_{k+1},\delta_{k+1}(T_k))\cdot \mathcal{L}^{(k+1)}$ and $\mathcal{Q}_k\mathcal{P}_{k+1}=\mathcal{L}^{(k+1)}$, we get \begin{eqnarray}\label{crc}
C_h(E_{k+1},\delta_{k+1}(T))=\mathcal{P}_k\cdot \widetilde{C_h}(E_{k+1}, \delta_{k+1}(T))\cdot \mathcal{Q}_{k+1}.
\end{eqnarray} 
Plugging \eqref{crc} for all $k=1,...,n-1$ into \eqref{qconnection} and using $\mathcal{P}_k \mathcal{Q}_k={\rm Id}$ lead to
\begin{eqnarray*}
S_{h-}(u_{\rm cat})S_{h+}(u_{\rm cat})=\left(\overrightarrow{\prod_{k=2,...,n}}\widetilde{C_h}(E_k, \delta_{k}(T))\right)\cdot e^{hT_n}\cdot \left(\overrightarrow{\prod_{k=2,...,n}}\widetilde{C_h}(E_{k}, \delta_{k}(T))\right)^{-1}.
\end{eqnarray*}
Here recall $T_n={\rm diag}\left(\zeta^{(n)}_1+(n-1),...,\zeta^{(n)}_n+(n-1)\right)\in {\rm End}(L(\lambda))\otimes {\rm End}(\mathbb{C}^n)$.
Then just as in the classical case, a manipulation of (blocked) Gauss decomposition, as well as the explicit formula of  $\widetilde{C_h}(E_{k+1}, \delta_{k+1}(T))$ and $S_{h+}(E_{k+1},\delta_{k+1}(T_k))$ given in Proposition \ref{qexStokes}, gives rise to the formula of $S_{h+}(E_{k+1},\delta_{k+1}(T_k))_{k,k+1}$ in Theorem \ref{explicitS}. \qed

\subsection{The leading terms of quantum Stokes matrices via the Gelfand-Tsetlin basis}\label{qGZsystem}

Given any $L(\lambda)$, let us think of the associated quantum Stokes matrices as blocked matrices. Let us take a proper norm on the vector space ${\rm End}(\lambda)$, then up to a slight modification, Proposition \ref{leadingterm} and Proposition \ref{StokesinGZ} can be applied to the quantum/blocked cases. That is if we introduce the matrix \[g(u;T)=\overrightarrow{\prod_{k=1,...,n-1}}(z_k)^{\frac{\delta_k(T)}{2\pi\I }}\in {\rm End}(L(\lambda))\otimes {\rm End}(\mathbb{C}^n),\]
(here recall that $z_k$ are the coordinates in \eqref{newcoor}) then
\begin{pro}\label{qStokesGZ}
For any fixed none zero real number $h$, we have as $u\in U_{\rm id}$ and $u\rightarrow u_{\rm cat}$, 
\begin{eqnarray}\label{qrelimitGZ}
S_{h\pm}(u)=S_{h\pm}\left(u_{\rm cat}, \hspace{2mm} g(hu;hT)\cdot hT\cdot g(hu;hT)^{-1}+\sum_{k=2}^{n-1}\mathcal{O}(z_k^{-1})\right).
\end{eqnarray}
Furthermore, the subdiagonal entries of $S_{h\pm}(u)$ satisfy
\begin{eqnarray}\label{hnonuni}
S_{h+}(u)_{k,k+1}={l}^{(+)}_{k,k+1}(u)+o(l^{(\pm)}_{k,k+1}(u)), \hspace{3mm} \text{as} \ u\rightarrow u_{\rm cat} \ \text{ from } U_{\rm id},\\
S_{h-}(u)_{k+1,k}={l}^{(-)}_{k+1,k}(u)+o(l^{(-)}_{k+1,k}(u)), \hspace{3mm} \text{as} \ u\rightarrow u_{\rm cat} \ \text{ from } U_{\rm id},
\end{eqnarray}
where the leading terms ${l}^{(\pm)}_{k,k+1}(u)$ are given by the expressions $S_{h+}(u_{\rm cat})_{k,k+1}$ and $S_{h-}(u_{\rm cat})_{k+1,k}$ in Theorem \ref{explicitS} provided we replace respectively $(\Delta_L)^{1,...,k}_{1,...,k-1,k+1}\left(\frac{h}{2\pi\I }T(\zeta^{(k)}_i)\right)$ and $(\Delta_R)_{1,...,k}^{1,...,k-1,k+1}\left(\frac{-h}{2\pi\I }T(\zeta^{(k)}_i)\right)$ by 
\begin{equation}\label{modify0}
(u_k-u_{k-1})^{\frac{he_{kk}-h\zeta^{(k)}_i}{2\pi\I }}\left({u_{k+1}-u_k}\right)^{\frac{h\zeta^{(k)}_i-he_{k+1,k+1}}{2\pi\I }}\cdot (\Delta_L)^{1,...,k}_{1,...,k-1,k+1}\left(\frac{h}{2\pi\I }T(\zeta^{(k)}_i)\right)
\end{equation}
and  
\begin{equation}\label{modify1}
(\Delta_R)_{1,...,k}^{1,...,k-1,k+1}\left(\frac{-h}{2\pi\I }T(\zeta^{(k)}_i)\right) \cdot (u_k-u_{k-1})^{\frac{h\zeta^{(k)}_i-he_{kk}}{2\pi\I }}\left({u_{k+1}-u_k}\right)^{\frac{he_{k+1,k+1}-h\zeta^{(k)}_i}{2\pi\I }}.
\end{equation}
\end{pro}
By Proposition \ref{GZaction} and \ref{qStokesGZ}, the leading terms of entries of quantum Stokes matrices, as elements in ${\rm End}(L(\lambda))$, can be written down explicitly in terms of the Gelfand-Tsetlin basis.

\subsection{Regularized limits of quantum Stokes matrices from the isomonodromy deformation}\label{reglimit}
In this subsection, let us explain the regularized limit of quantum Stokes matrices given in Proposition \ref{qStokesGZ} from the viewpoint of isomonodromy deformation. It follows, from the general theory of isomonodromic deformation of meromorphic linear system of ordinary differential equations \cite{JMMS, JMU} that 
\begin{pro}
As a function of $u\in \h_{\rm reg}(\mathbb{R})$, the entries of the Stokes matrix $S_{h\pm}(u)$ satisfy 
\begin{eqnarray}\label{Sisoeq}
\frac{\partial S_{h\pm}(u)_{kl}}{\partial u_i}=h\sum_{j=1,...,n, \text{ and } j\ne i}\frac{[e_{ij}e_{ji}, S_{h\pm}(u)_{kl}]}{u_i-u_j}\in {\rm End}(L(\lambda)).
\end{eqnarray}
Here $e_{ij}$ is seen as an element in ${\rm End}(L(\lambda))$ via the given representation of $\frak{gl}_n$ on $L(\lambda).$
\end{pro}

Now let us consider the equation for a function $W(u)\in {\rm End}(L(\lambda))$
\begin{eqnarray}\label{Geq}
\frac{\partial W}{\partial u_i}=h\sum_{j\ne i}\frac{e_{ij}e_{ji}}{u_i-u_j} \cdot W, \hspace{5mm} i=1,...,n. 
\end{eqnarray}
The equation \eqref{Geq} was called the Casimir equation \cite{TL}. In terms of the variables $z_i$ introduced in \eqref{newcoor}, the equation takes the form
\[\frac{\partial W}{\partial z_k}=\left(\frac{V_k}{z_k}+Reg\right)\cdot W, \ \text{ for } k=1,...,n-1\]
where $Reg$ stands for some ${\rm End}(L(\lambda))$ valued rational function of $z_k$ regular at $z_k=\infty$, and $V_k:=\sum_{1\le j\le k-1} e_{kj}e_{jk}\in {\rm End}(L(\lambda))$. 

The elements $V_k$ commute with each other, and are the analog of the Jucys-Murphy elements in the group algebra of $S_n$. Then let $W(u)\in {\rm End}(L(\lambda))$ be the solution of the equation \eqref{Geq}
with the prescribed asymptotics
\begin{eqnarray}\label{Gasy}
W(u)\cdot \overrightarrow{\prod_{k=1,...,n-1} }(z_k)^{h{V_k}} \rightarrow 1, \hspace{3mm} \text{ as } z_k\rightarrow\infty.
\end{eqnarray}
The existence of $W(u)$ follows from the general theory that a formal solution of differential equations with regular singularities is in fact analytic. 

It follows from \eqref{Sisoeq} and \eqref{Geq} that for any $k,l$ the element $W(u)^{-1} \cdot S_{\pm}(u;h)_{kl} \cdot W(u)\in{\rm End}(L(\lambda))$ is a constant on $u\in U_{\rm id}\subset \frak t_{\rm reg}(\mathbb{R})$. By the asymptotics \eqref{Gasy} and \eqref{modify0}, \eqref{modify1} at the caterpillar point $u_{\rm cat}$, we have
\begin{pro}\label{SWS}
The Stokes matrices $S_+(u_{\rm cat};h)\in{\rm End}(L(\lambda))\otimes {\rm End}(\mathbb{C}^n)$ at the caterpillar point $u_{\rm cat}$ satisfy
\begin{eqnarray}\label{GStokes}
S_{h\pm}(u_{\rm cat})_{kl}:=W(u)^{-1} \cdot S_{\pm}(u;h)_{kl} \cdot W(u), \hspace{3mm} k,l=1,...,n.
\end{eqnarray}
\end{pro}
Then we notice that letting $u\rightarrow u_{\rm cat}$ from $U_{\rm id}$ in the identity \eqref{GStokes} (for $l=k+1$), the asymptotics \eqref{Gasy} of $W(u)$ matches up with the singular part in the expressions \eqref{hnonuni} and \eqref{modify0}. Just like the classical case, it interprets the leading asymptotics of the quantum Stokes matrices via the isomonodromy approach.

\subsection{Representation of quantum groups from quantum Stokes matrices at $u_{\rm cat}$}\label{reglimit2}
Now we can generalize Theorem \ref{introthm1} to the caterpillar point $u_{\rm cat}$. 
\begin{thm}\label{quantumgpucat}
For any fixed none zero real number $h$, the map (with $q=e^{h/2}$)
\begin{equation*}
\begin{split}
\mathcal{S}_q(u_{\rm cat}): U_q(\frak{gl}_n)&\rightarrow {\rm End}(L(\lambda))~;\\
e_i&\mapsto \frac{S_+(u_{\rm cat})_{i,i}^{-1}\cdot S_+(u_{\rm cat})_{i,i+1}}{q^{-1}-q}, \\
f_i&\mapsto \frac{S_-(u_{\rm cat})_{i+1,i}\cdot S_-(u_{\rm cat})_{i,i}^{-1}}{q-q^{-1}}, \\
q^{h_i}&\mapsto S_+(u_{\rm cat})_{i,i}
\end{split}
\end{equation*}
defines a representation of the Drinfeld-Jimbo quantum group $U_q(\frak{gl}_n)$ on the vector space $L(\lambda)$. 
\end{thm}
\begin{proof}
It follows from the fact that the gauge transformation transformation $W(u)$ in Proposition \ref{SWS} preserves the defining relations of quantum group. For example, for any $u\in\h_{\rm reg}(\mathbb{R})$ by Theorem \ref{introthm1} we have
\begin{align}\nonumber
    &S_\pm(u)_{k,k+1}^2 S_\pm(u)_{k-1,k} - (q+q^{-1})S_\pm(u)_{k,k+1} S_\pm(u)_{k-1,k} S_\pm(u)_{k,k+1} \\ \label{qser}
    +& S_\pm(u)_{k-1,k} S_\pm(u)_{k,k+1}^2=0.
\end{align}
Since the conjugation by $W(u)$ preserves the above relation (and the other defining relations), the corresponding entries of $S_\pm(u_{\rm cat})$ satisfy the same relation \eqref{qser} (and the other defining relations).
\end{proof}

The above theorem can also be checked directly from the identities \eqref{hnonuni}-\eqref{modify1}. For example, using \eqref{hnonuni} and \eqref{modify0}, one checks $S_\pm(u_{\rm cat})$ satisfy the same relation \eqref{qser} (after taking the regularized limits as $u\rightarrow u_{\rm cat}$).

\begin{rmk}
More generally, Theorem \ref{introthm1} can be generalized from $\h_{\rm reg}(\mathbb{R})$ to the De Concini-Procesi space. Therefore, the regularized limit of quantum Stokes matrices preserves the underlying quantum group structures.
\end{rmk}

\section{Crystals and cactus group actions arising from the WKB approximation}\label{tropicalisomo}
In the previous section, we introduce the relation between the regularized limits of quantum Stokes matrices at $u_{\rm cat}$ and Gelfand-Tsetlin basis. In this section, we deepen the relation between Stokes phenomenon and representation theory by showing that the WKB approximation of quantum Stokes matrices at a caterpillar point gives rise to $\frak{gl}_n$-crystals. Section \ref{sec:WKB} computes explicitly the WKB approximation. Then Section \ref{sec:WKBdatum} introduces the notion of WKB datum of quantum Stokes matrices. In the end, Section \ref{cactuswall} gives a realization of the cactus group action on crystals via the wall-crossing formula of Stokes matrices at a caterpillar point.

\subsection{WKB approximation of quantum Stokes matrices at caterpillar points}\label{sec:WKB}
Since the derivative in equation \eqref{introeq} is multiplied by a small parameter $1/h$, we will call the leading term, as $h\rightarrow +\infty$, of Stokes matrices of \eqref{introeq} as the WKB approximation. In this subsection, we study the WKB approximation of Stokes matrices at $u_{\rm cat}$ given in \eqref{explicitS}.

First, for any Gelfand-Tsetlin pattern $\Lambda$ of $L(\lambda)$ given in Section \ref{GZbasis}, set
\begin{align}\label{x}
x^{(k)}_j(\Lambda)&:=-\lambda^{(k)}_j+\lambda^{(k-1)}_{j-1}-\lambda^{(k)}_{j-1}+\lambda^{(k+1)}_{j}, \hspace{3mm} 1\le j\le k+1,\\
y^{(k)}_j(\Lambda)&:=\lambda^{(k)}_j-\lambda^{(k-1)}_j+\lambda^{(k)}_{j+1}-\lambda^{(k+1)}_{j+1}, \hspace{3mm} 0\le j\le k,
\end{align}
and
\begin{align}\label{capx}
X^{(k)}_j(\Lambda)&:=\sum_{i=1}^jx^{(k)}_i(\Lambda), \hspace{3mm} 1\le j\le k+1,\\
Y^{(k)}_j(\Lambda)&:=\sum_{i=j}^ky^{(k)}_i(\Lambda), \hspace{3mm} 0\le j\le k.
\end{align}
Note that $Y^{(k)}_0(\Lambda)=Y^{(k)}_j(\Lambda)-X^{(k)}_j(\Lambda)$ for any $j=0,...,k+1$. 
Furthermore, we define 
\begin{align}
wt_k(\Lambda)&=\sum_{i=1}^k\lambda^{(k)}_i(\Lambda)-\sum_{i=1}^{k-1}\lambda^{(k-1)}_i(\Lambda),\\
\label{maximal}
\varepsilon_k(\Lambda)&={\rm max}\{X_1^{(k)}(\Lambda), X_2^{(k)}(\Lambda),...,X_k^{(k)}(\Lambda)\},\\
\phi_i(\Lambda)&={\rm max}\{Y_1^{(k)}(\Lambda), Y_2^{(k)}(\Lambda),...,Y_k^{(k)}(\Lambda)\},
\label{minimal}
\end{align}
and define the functions $l_1(\Lambda)<\cdots<l_{m_k}(\Lambda)$ of $\Lambda\in P_{GT}(\lambda;\mathbb{Z})$ be those ordered labels such that \begin{eqnarray}\label{l}
X^{(k)}_{l_j}(\Lambda)=\varepsilon_k(\Lambda).
\end{eqnarray}
Let us denote by 
\begin{equation}\label{PkGZ}
    P^k_{GT}(\lambda;\mathbb{Z}):=\{\Lambda\in P_{GT}(\lambda;\mathbb{Z})~|~\Lambda+\delta^{(k)}_{l_1}\in P_{GT}(\lambda;\mathbb{Z}) \ \text{and} \ \Lambda+\delta^{(k)}_{l_i}\notin P_{GT}(\lambda;\mathbb{Z}) \ \text{for} \ i=2,...,m_k\}.
\end{equation}
In particular, all $\Lambda$ satisfying $m_k(\Lambda)=1$ and $\Lambda+\delta^{(k)}_{l_i}\notin P_{GT}(\lambda;\mathbb{Z})$ are in $P^k_{GT}(\lambda;\mathbb{Z})$. From \eqref{x}--\eqref{l}, the function $m_k(\Lambda)$ is defined for any point $\Lambda$ in the real polytope $P_{GT}(\lambda;\mathbb{R})$. Thus $P^k_{GT}(\lambda;\mathbb{R}):=P^k_{GT}(\lambda;\mathbb{Z})\otimes_\mathbb{Z}\mathbb{R}$ is an open dense part of $P_{GT}(\lambda;\mathbb{R})$, whose complements are cut out by various equalities between $\lambda^{(j)}_i(\Lambda)$ for $j=k-1,k,k+1$. By this reason,
the elements in the subset $P^k_{GT}(\lambda;\mathbb{Z})$ are called generic (it gives precise meaning of the word generic in Theorem \ref{WKBthm}).

Let $\{\xi_\Lambda \}$ be the Gelfand-Tsetlin basis given in Section \ref{GZbasis}. Then we have
\begin{pro}\label{qleading}
For $k=1,...,n-1$ and for any $\Lambda\in P^k_{GT}(\lambda;\mathbb{Z})$, there exists real valued functions $\theta_{k}(\xi_\Lambda)\in [0,2\pi)$ of the patterns $\Lambda$ (independent of $q=e^{h/2}$) such that
\begin{equation}\label{qexpansion}
S_{h+}(u_{\rm cat})_{k,k+1}\cdot \xi_\Lambda \sim q^{-wt_{k+1}(\Lambda)-\varepsilon_k(\Lambda)+\I \theta_{k}} \xi_{\Lambda+\delta^{(k)}_{l_1(\Lambda)}} , \hspace{3mm} \text{as} \ \ q\rightarrow \infty.
\end{equation}
Here recall the pattern $\Lambda+\delta^{(k)}_{l_1}$ is obtained from $\Lambda$ by replacing $\lambda^{(k)}_{l_1}$ by $\lambda^{(k)}_{l_1}+ 1$. 
\end{pro}
\pf By the asymptotics of gamma function
\[{\rm In}(\Gamma(1+z))\sim z{\rm In}(z)-z+\frac{1}{2}{\rm In}(z)+\frac{1}{2}{\rm In}(2\pi)+\frac{1}{12z}+o\left(\frac{1}{z}\right), \ \ \text{as} \ |z|\rightarrow\infty, \ \ |{\rm arg}(z)|<\pi,\] 
for $r$ a real number and $h\rightarrow+\infty$, we have
\begin{equation}\label{gamma}
{\rm In}\left(\Gamma\left(1+\frac{rh}{2\pi \I }\right)\right)\sim \frac{rh{\rm In}(h)}{2\pi \I }+\frac{rh}{2\pi \I }{\rm In}\left(\frac{|r|}{2\pi}\right)- \frac{|r|h}{4}-\frac{rh}{2\pi \I }+\frac{1}{2}{\rm In}\left(\frac{rh}{\I }\right). 
\end{equation}
Here we use ${\rm In}(\frac{\pm1}{\I })=\mp\frac{\pi\I }{2}$ to separate the real and imaginary part of the $h$ linear terms. By \eqref{eigenvalues}, \[\zeta^{(k)}_i\cdot \xi_\Lambda =(\lambda^{(k)}_i(\Lambda)-i+1)\cdot\xi_\Lambda ,\] 
by abuse of notation, we will take $\zeta^{(k)}_i$ as the number $\lambda^{(k)}_i(\Lambda)-i+1$ when the vector $\xi_\Lambda$ is specified.
Then by \eqref{gamma} and the interlacing inequalities between $\lambda^{(i)}_j$ for $i=k-1,k,k+1$,  
\begin{eqnarray*}\nonumber
&&{\rm In}\left(\frac{\prod_{l=1,l\ne i}^{k}\Gamma(h\frac{\zeta^{(k)}_i-\zeta^{(k)}_l}{2\pi \I })}{\prod_{l=1}^{k+1}\Gamma(1+h\frac{\zeta^{(k)}_i-\zeta^{(k+1)}_l-1}{2\pi \I })}\frac{\prod_{l=1,l\ne i}^{k}\Gamma(1+h\frac{\zeta^{(k)}_i-\zeta^{(k)}_l-1}{2\pi \I })}{\prod_{l=1}^{k-1}\Gamma(1+h\frac{\zeta^{(k)}_i-\zeta^{(k-1)}_l}{2\pi \I })}\right)\\
&\sim&\frac{h{\rm In}(h)}{2\pi \I }A^{(k)}_i+\I h\theta^{(k)}_i+\frac{h}{4}B^{(k)}_i+C^{(k)}_i, \hspace{5mm} \text{as} \ h\rightarrow+\infty,
\end{eqnarray*}
where
\begin{eqnarray*}
&&A^{(k)}_i=\sum_{l=1}^{k-1}\zeta^{(k-1)}_l+\sum_{l=1}^{k+1}\zeta^{(k+1)}_l-2\sum_{l=1}^{k}\zeta^{(k)}_l+2,\\
&&B^{(k)}_i=-2\sum_{l=1}^{i-1}(\zeta^{(k)}_i-\zeta^{(k)}_l)+2\sum_{l=i+1}^{k}(\zeta^{(k)}_i-\zeta^{(k)}_l)+\sum_{l=1}^{i-1}(\zeta^{(k)}_i-\zeta^{(k-1)}_l)\\&&\hspace{12mm}-\sum_{l=i}^{k-1}(\zeta^{(k)}_i-\zeta^{(k-1)}_l)+\sum_{l=1}^{i}(\zeta^{(k)}_i-\zeta^{(k+1)}_l)-\sum_{l=i+1}^{k+1}(\zeta^{(k)}_i-\zeta^{(k+1)}_l),\\ \label{constantc}
&&C^{(k)}_i=\frac{1}{2}{\rm In}\left(\frac{(2\pi\I)^{2k-2}}{h^{2k}}\cdot\frac{-\prod_{l=1,l\ne i}^{k}(\zeta^{(k)}_i-\zeta^{(k)}_l-1)}{\prod_{l=1,l\ne i}^{k}(\zeta^{(k)}_i-\zeta^{(k)}_l)\prod_{l=1}^{k+1}(\zeta^{(k)}_i-\zeta^{(k+1)}_l-1)\prod_{l=1}^{k-1}(\zeta^{(k)}_i-\zeta^{(k-1)}_l)}\right),\end{eqnarray*}
and (since $\theta^{(k)}_i$'s are real, they do not affect the leading asymptotics, but for completeness let us still list them)
\begin{eqnarray*}
&&\theta^{(k)}_i=\frac{1}{2\pi}A^{(k)}_i-\sum_{l=1}^{i-1}\frac{\zeta^{(k)}_i-\zeta^{(k)}_l}{2\pi}{\rm In}(\frac{\zeta^{(k)}_l-\zeta^{(k)}_i}{2\pi})+\sum_{l=i+1}^{k}(\frac{\zeta^{(k)}_l-\zeta^{(k)}_i}{2\pi}){\rm In}(\frac{\zeta^{(k)}_i-\zeta^{(k)}_l}{2\pi})\\ &&\hspace{12mm}-\sum_{l=1}^{i-1}\frac{\zeta^{(k)}_l+1-\zeta^{(k)}_i}{2\pi}{\rm In}(\frac{\zeta^{(k)}_i-\zeta^{(k)}_l-1}{2\pi})+\sum_{l=i+1}^{k}(\frac{\zeta^{(k)}_l-\zeta^{(k)}_i-1}{2\pi}){\rm In}(\frac{\zeta^{(k)}_i+1-\zeta^{(k)}_l}{2\pi})\\ &&\hspace{12mm}+\sum_{l=1}^{i-1}(\frac{\zeta^{(k-1)}_l-\zeta^{(k)}_i-1}{2\pi}){\rm In}(\frac{\zeta^{(k)}_i+1-\zeta^{(k-1)}_l}{2\pi})-\sum_{l=i}^{k-1}(\frac{\zeta^{(k)}_i-\zeta^{(k-1)}_l-1}{2\pi}){\rm In}(\frac{\zeta^{(k-1)}_l+1-\zeta^{(k)}_i}{2\pi})\\&&\hspace{12mm}+\sum_{l=1}^{i}(\frac{\zeta^{(k+1)}_l-\zeta^{(k)}_i}{2\pi}){\rm In}(\frac{\zeta^{(k)}_i-\zeta^{(k+1)}_l}{2\pi})-\sum_{l=i+1}^{k+1}(\frac{\zeta^{(k)}_i-\zeta^{(k+1)}_l}{2\pi}){\rm In}(\frac{\zeta^{(k+1)}_l-\zeta^{(k)}_i}{2\pi}).\end{eqnarray*}
Finally, the action of $S_{h+}(u_{\rm cat})_{k,k+1}$ on $\xi_\Lambda $ is given by
\begin{align*}\nonumber
&2\pi\I \cdot (\I h)^{\frac{h\small{(e_{kk}-e_{k+1,k+1}-1)}}{2\pi \I }} e^{\frac{-he_{kk}}{2}}\\
&\times \sum_{i=1}^k \frac{\prod_{l=1,l\ne i}^{k}\Gamma(h\frac{\zeta^{(k)}_i-\zeta^{(k)}_l}{2\pi \I })}{\prod_{l=1}^{k+1}\Gamma(1+h\frac{\zeta^{(k)}_i-\zeta^{(k+1)}_l-1}{2\pi \I })}\frac{\prod_{l=1,l\ne i}^{k}\Gamma(1+h\frac{\zeta^{(k)}_i-\zeta^{(k)}_l-1}{2\pi \I })}{\prod_{l=1}^{k-1}\Gamma(1+h\frac{\zeta^{(k)}_i-\zeta^{(k-1)}_l}{2\pi \I })}\Delta^{1,...,k}_{1,...,k-1,k+1}\left(\frac{h}{2\pi\I }T(\zeta^{(k)}_i)\right)\cdot\xi_\Lambda 
\\
\sim& \sum_{i=1}^k\left(( h)^{\frac{h\small{(e_{kk}-e_{k+1,k+1}-1)}}{2\pi \I }}e^{\frac{h{\rm In}(h)}{2\pi \I }A^{(k)}_i}\right)\left(e^{\frac{-h\small{(e_{kk}+e_{k+1,k+1}+1)}}{4}}e^{\I h\theta^{(k)}_i+\frac{h}{4}B^{(k)}_i}\right)\cdot 2\pi\I e^{C^{(k)}_i}\Delta^{1,...,k}_{1,...,k-1,k+1}\cdot\xi_\Lambda \\
=&\sum_{i=1}^k q^{-wt_{k+1}(\Lambda)-X^{(k)}_i(\Lambda)+2\I h \theta^{(k)}_i} \xi_{\Lambda+\delta^{(k)}_i} .
\end{align*}
Here in the last equality, we use the identities \eqref{eigenvalues0}, \eqref{eigenvalues}, \eqref{x} and \eqref{capx} to get 
\begin{align*}
(h)^{\frac{h\small{(e_{kk}-e_{k+1,k+1}-1)}}{2\pi \I }}e^{\frac{h{\rm In}(h)}{2\pi \I }A^{(k)}_i}&=1,\\
e^{\frac{-h\small{(e_{k+1,k+1}+e_{kk}+1)}}{4}}e^{\frac{h}{4}B^{(k)}_i+\I h\theta^{(k)}_i}\cdot \xi_{\Lambda+\delta^{(k)}_i} &=q^{-wt_{k+1}(\Lambda)-X^{(k)}_i(\Lambda)+2\I h \theta^{(k)}_i}\xi_{\Lambda+\delta^{(k)}_i} ,
\end{align*}
and use the
expression \eqref{ashift} of $\Delta^{1,...,k}_{1,...,k-1,k+1}\left(T(\zeta^{(k)}_i)\right) \cdot \xi_\Lambda$ and the expression of $C^{(k)}_i$ to get
\begin{align*}2\pi\I \cdot e^{C^{(k)}_i}\Delta^{1,...,k}_{1,...,k-1,k+1}\left(\frac{h}{2\pi\I }T(\zeta^{(k)}_i)\right)\cdot\xi_\Lambda &=\left(\frac{h}{2\pi\I}\right)^{k}e^{C^{(k)}_i}\Delta^{1,...,k}_{1,...,k-1,k+1}\left(T(\zeta^{(k)}_i)\right)\cdot\xi_\Lambda \\
&=\xi_{\Lambda+\delta^{(k)}_i} .
\end{align*}
Then the proposition follows from the definitions \eqref{maximal} and \eqref{l} of $\varepsilon_k(\Lambda)$ and $l_{i}(\Lambda)$ with $i=1,...,m_k$, and the assumption $m_k(\Lambda)=1.$ The constant $\theta_k$ in \eqref{qexpansion} is just $2\theta^{(k)}_{l_1}$.
\qed
\subsection{WKB operators}\label{sec:WKBdatum}
In this subsection, we introduce a combinatorial structure to encode the WKB leading terms (of entries) of quantum Stokes matrices at a caterpillar point. 

Recall that the subset $P_{GT}^k(\lambda;\mathbb{Z})\subset P_{GT}(\lambda;\mathbb{Z})$ is defined in \eqref{PkGZ}. Let us denote by $E^k_{GT}\subset E_{GT}(\lambda)$ the subset consisting of the basis elements parametrized by $P_{GT}^k(\lambda;\mathbb{Z})$, then the formula \eqref{qexpansion} shows that taking the $q$ leading term of $S_{h+}(u_{\rm cat})_{k,k+1}$ naturally induces a map
\begin{eqnarray*}
\widetilde{e_k}(u_{\rm cat}): \ E^k_{GT}\rightarrow E_{GT}~;~
\widetilde{e_k}(u_{\rm cat})(\xi_{\Lambda} )=\xi_{\Lambda+\delta^{(k)}_{l_1}} .
\end{eqnarray*}
In the following, for simplicity, if there is no ambiguity we will write $\widetilde{e_k}$ for $\widetilde{e_k}(u_{\rm cat})$. Equivalently, as we identify $E_{GT}(\lambda)$ with $P_{GT}(\lambda;\mathbb{Z})$ by mapping $\xi_\Lambda$ to $\Lambda$, $\widetilde{e_k}$ can be seen as a map
\begin{eqnarray*}
\widetilde{e_k}(u_{\rm cat}): \ P_{GT}^k(\lambda;\mathbb{Z})\rightarrow P_{GT}(\lambda;\mathbb{Z})~;~
\widetilde{e_k}({\Lambda})={\Lambda+\delta^{(k)}_{l_1}}.
\end{eqnarray*}
The map has a canonical extension to the whole $P_{GT}(\lambda;\mathbb{Z})$ as follows. 

First, the expression of $\widetilde{e_k}$ is universal, i.e., doesn't depend on the choice of $\lambda$. Thus for any positive integer $N$, if we set $N\lambda=(N\lambda^{(n)}_1, ..., N\lambda^{(n)}_n)$, then $\widetilde{e_k}$ can be equivalently seen as 
a map from $P^k_{GT}(\lambda;\frac{\mathbb{Z}}{N}):=\{\Lambda\in P^k_{GT}(\lambda;\mathbb{R})~|~\lambda^{(i)}_j(\Lambda)\in \frac{\mathbb{Z}}{N}\}$ to $P_{GT}(\lambda;\frac{\mathbb{Z}}{N}):=\{\Lambda\in P_{GT}(\lambda;\mathbb{R})~|~\lambda^{(i)}_j(\Lambda)\in \frac{\mathbb{Z}}{N}\}$
\begin{eqnarray}\label{dsystem}
\widetilde{e_k}^{\frac{1}{N}}( \Lambda):=\Lambda+\frac{1}{N}\delta^{(k)}_{l},
\end{eqnarray}
where $\lambda^{(i)}_j(\Lambda+\frac{1}{N}\delta^{(k)}_{l})=\lambda^{(i)}_j+\frac{\delta_{ik}\delta_{jl}}{N}$. 
Now as $N\rightarrow\infty$, the discrete "dynamical system" \eqref{dsystem} approaches to a unique continuous system in the inner part of the whole real polytope $P_{GT}(\lambda;\mathbb{R})$:
\begin{eqnarray}\label{csystem}
\widetilde{e_k}^t( \Lambda):=\Lambda+t_1\delta^{(k)}_{r_1}+\cdots +t_m\delta^{(k)}_{r_m},
\end{eqnarray}
where the time
\[t\in \Big[0, \hspace{2mm} \sum_{i=1}^k({\rm min}(\lambda^{(k+1)}_i(\Lambda),\lambda^{(k-1)}_{i-1}(\Lambda))-\lambda^{(k)}_i(\Lambda))\Big],\] 
and $t_1,...,t_{m}$, $r_1,...,r_m$ are determined by
\begin{eqnarray*}
&&{r_1}=l_1(\Lambda+s\delta^{(k)}_{r_1}), \hspace{3mm} \text{for all} \ 0\le s< t_1,\\
&&r_i=l_1(\Lambda+t_1\delta^{(k)}_{r_1}+\cdots t_{i-1}\delta^{(k)}_{r_{i-1}}+s\delta^{(k)}_{r_i}), \hspace{3mm} \text{for all} \ 0\le s<t_i, \ for \ i=2,...,m,\\
&&t=t_1+\cdots + t_m.
\end{eqnarray*}
Here the function $l_1(\Lambda)$ is defined for any point $\Lambda$ in $P_{GT}(\lambda;\mathbb{R})$ just as \eqref{x}--\eqref{l}. 

This 
continuation is unique, thus canonically determines an extension of $\widetilde{e_k}$ to the complements of the "generic part" $P^k_{GT}(\lambda;\mathbb{Z})$ in $P_{GT}(\lambda;\mathbb{Z})$. It gives rise to
\begin{defi}\label{CWKBe}
For each $k$, the WKB operator $\widetilde{e_k}$ from $E_{GT}\cong P_{GT}(\lambda;\mathbb{Z})$ to $E_{GT}\cup \{0\}$ is given by
\begin{eqnarray}\label{WKBoperator}
\widetilde{e_k}\cdot \xi_{\Lambda}(u_{\rm cat}):=\xi_{\Lambda+\delta^{(k)}_{l_1(\Lambda)}}(u_{\rm cat}), \hspace{5mm} \forall\Lambda\in P_{GT}(\lambda;\mathbb{Z})
\end{eqnarray}
Here recall that $l_1(\Lambda)$ is the integer given by \eqref{l}, and it is supposed that $\widetilde{e_k}\cdot \xi_\Lambda$ is zero if $\Lambda+\delta^{(k)}_{l_1(\Lambda)}$ doesn't belong to $P_{GT}(\lambda;{\mathbb{Z}})$.
\end{defi}
The above computation and discussion carry to the WKB approximation of the entries $S_{h-}(u_{\rm cat})_{k+1,k}$ of lower triangular Stokes matrix $S_{h-}(u_{\rm cat})$. It produces the same set $E_{GT}(\lambda)$, and induces operators $\widetilde{f_k}$ on $E_{GT}(\lambda)$
\begin{eqnarray}\label{CWKBf}
\widetilde{f_k}\cdot \xi_{\Lambda}=\xi_{\Lambda-\delta^{(k)}_{l_{m_k}(\Lambda)}}, \hspace{5mm} \forall\Lambda\in P_{GT}(\lambda;\mathbb{Z}).
\end{eqnarray}
One checks that $\widetilde{e_k}$ and $\widetilde{f_k}$ satisfy that
for all $\Lambda, \Lambda'\in P_{GT}(\lambda;{\mathbb{Z}})$, \begin{eqnarray*}
&&\varepsilon_k(\Lambda)={\rm max}\{j:\widetilde{e_k}^j(\xi_\Lambda)\ne 0\},\\
&&\phi_k(\Lambda)={\rm max}\{j:\widetilde{f_k}^j(\xi_\Lambda)\ne 0\},
\end{eqnarray*}
and
\[\widetilde{e_k}\cdot \xi_\Lambda=\xi_\Lambda' \hspace{3mm} if \ and \ only \ if \hspace{3mm} \widetilde{f_k}\cdot \xi_\Lambda'=\xi_\Lambda.\]

\begin{defi}
We call $(E_{GT}(\lambda),\widetilde{e_k},\widetilde{f_k},\varepsilon_k,\phi_k)$ the WKB datum of the Stokes matrices $S_{h\pm}(u_{\rm cat})$ at the caterpillar point associated to the representation $L(\lambda)$.
\end{defi}

\subsection{WKB datum are crystals}\label{WKBcry}
In Section \ref{sec:crystal} we recall the notion of crystals and their tensor products. In Section \ref{glncry} we prove that the WKB datum of quantum Stokes matrices at $u_{\rm cat}$ are $\frak{gl}_n$-crystals. Furthermore, in Section \ref{sec:tensor} we also realize the tensor products of $\frak{gl}_n$-crystals by WKB analysis.
\subsubsection{Crystals and tensor products}\label{sec:crystal}
Let $\g$ be a semisimple Lie algebra with a Cartan datum $(A,\Delta_+=\{\alpha_i\}_{i\in I},\Delta^\vee_+=\{\alpha_i^\vee\}_{i\in I}, P, P^\vee)$ be a Cartan datum, where $P\subset \h^*$ denotes the weight lattice, $I$ denotes the set of vertices of its Dynkin diagram, $\alpha_i\in I$ denote its simple
roots, and $\alpha^\vee_i$ the simple coroots.

\begin{defi}\label{WKBdata}
A $\g$-crystal is a finite set $B$ along with maps
\begin{eqnarray*}
wt&:& B\rightarrow P,\\
\tilde{e}_i,\tilde{f}_i&:& B\rightarrow B\cup \{0\}, \hspace{2mm} i\in I,\\
\varepsilon_i,\phi_i&:& B\rightarrow \mathbb{Z}\cup \{-\infty\}, \hspace{2mm} i\in I,
\end{eqnarray*}
satisfying for all $b,b'\in B$, and $i\in I$,
\begin{itemize}
    \item $\tilde{f}_i(b)=b'$ if and if $b=\tilde{e}_i(b')$, in which case
    \[wt(b')=wt(b)-\alpha_i, \hspace{3mm} \varepsilon_i(b')=\varepsilon_i(b)+1, \hspace{3mm} \phi_i(b')=\phi_i(b)-1.\]
    \item $\phi_i(b)=\varepsilon_i(b)+\langle wt(b),\alpha_i^\vee\rangle$, and if $\phi_i(b)=\varepsilon_i(b)=-\infty$, then $\tilde{e}_i(b)=\tilde{f}_i(b)=0.$
\end{itemize}
The map $wt$ is called the weight map, $\tilde{e}_i$ and $\tilde{f}_i$ are called Kashiwara operators or crystal operators.
\end{defi}

\begin{defi}
Let $B_1$ and $B_2$ be two crystals. The tensor product $B_1\otimes B_2$ is the crystal with the underlying set $B_1\times B_2$ (the Cartesian product) and structure maps
 \begin{eqnarray*}
wt(b_1,b_2)&=&wt(b_1)+wt(b_2),\\
\tilde{e}_i(b_1,b_2)&=&\left\{
          \begin{array}{lr}
             (\tilde{e}_i(b_1),b_2),   & if \ \varepsilon_i(b_1)>\phi(b_2)  \\
           (b_1,\tilde{e}_i(b_2)), & otherwise
             \end{array}\right. \\
\tilde{f}_i(b_1,b_2)&=&\left\{
          \begin{array}{lr}
             (\tilde{f}_i(b_1),b_2),   & if \ \varepsilon_i(b_1)\ge \phi(b_2)  \\
           (b_1,\tilde{f}_i(b_2)), & otherwise.
             \end{array}
\right. 
\end{eqnarray*}
\end{defi}
\subsubsection{${\gl}_n$-crystals from WKB approximation}\label{glncry}
Let us take the Cartan datum of type $A_{n-1}$, where $I=\{1, 2, ..., n-1\}$, the weight lattice $P=\mathbb{C}\{v_1,...,v_n\}$, and $\alpha_i = v_i-v_{i+1}$, $\alpha_i^\vee$ is given by $\langle v_i, \alpha_j^\vee\rangle=\delta_{ij}-\delta_{i,j+1}.$ Let us consider the WKB datum, given in Definition \ref{WKBdata}, of the quantum Stokes matrices at the caterpillar point $u_{\rm cat}$ associated to the representation $L(\lambda)$. Let us define a weight map
\[wt:E_{GT}(\lambda)\rightarrow P~;~ \xi_\Lambda\rightarrow \sum_{k=1}^nwt_k(\Lambda)v_k=\sum_{k=1}^n(\sum_{i=1}^k\lambda^{(k)}_i(\Lambda)-\sum_{i=1}^{k-1}\lambda^{(k-1)}_i(\Lambda))v_k.\]
The following theorem gives a realization of $\frak{gl}_n$-crystals from the WKB approximation.
\begin{thm}\label{WKBcrystal}
The WKB datum $(E_{GT}(\lambda),\widetilde{e_k}(u_{\rm cat}),\widetilde{f_k}(u_{\rm cat}),\varepsilon_k,\phi_k)$ with the weight map $wt$ is a ${\gl}_n$-crystal. 
\end{thm}
\pf By the formula \eqref{WKBoperator} and \eqref{CWKBf}, the WKB operators on the finite set $E_{GT}(\lambda)$ are:
\begin{eqnarray}\label{eoperator}
 \widetilde{e_k}\cdot \xi_\Lambda=\left\{
          \begin{array}{lr}
             \xi_{\Lambda+\delta^{(k)}_l},   & if \ \ \varepsilon_k(\Lambda)>0,   \\
            0, & if \ \ \varepsilon_k(\Lambda)=0,
             \end{array}
\right.
\text{where} \ l ={\rm min}\{j=1,...,k~|~X_j^{(k)}(\Lambda)=\varepsilon_k(\Lambda)\},\\\label{foperator}
 \widetilde{f_k}\cdot \xi_\Lambda=\left\{
          \begin{array}{lr}
             \xi_{\Lambda-\delta^{(k)}_l},   & if \ \ \phi_k(\Lambda)>0,   \\
            0, & if \ \ \phi_k(\Lambda)=0,
             \end{array}
\right.
\text{where} \ l ={\rm max}\{j=1,...,k~|~Y_j^{(k)}(\Lambda)=\phi_k(\Lambda)\}.\end{eqnarray}
Then one can verify directly the conditions in Definition \ref{WKBdata} via the explicit expressions. 

To avoid the direct but lengthy computation, one can also compare the WKB datum with the known ${\gl}_n$-crystal structures in literature: the explicit realization \eqref{eoperator} and \eqref{foperator} of the WKB operators $\widetilde{e_k}$, $\widetilde{f_k}$ coincide with the known ${\gl}_n$-crystal operators realized on the Gelfand-Tsetlin basis, see e.g., \cite{HKc}; or equivalently, under the natural bijection between semistandard Young tableaux and Geland-Testlin patterns, the WKB datum $(E_{GT},\widetilde{e_k},\widetilde{f_k},\varepsilon_k,\phi_k)$ coincides with the $\frak{gl}_n$-crystal structure on semistandard Young tableaux, see e.g., \cite{HKbook}.  \qed

\subsubsection{Tensor products from WKB approximation}\label{sec:tensor}
Recall that given any representation $L(\lambda)$, the quantum Stokes matrices at $u_{\rm cat}$ produce operators $S_{h\pm}(u_{\rm cat})_{ij}\in {\rm End}(L(\lambda)).$ Now given two representations $L(\lambda_1)$ and $L(\lambda_2)$, let us consider the actions of 
\begin{align}\label{copro}
S_{h+}(u_{\rm cat})_{k,k}\otimes S_{h+}(u_{\rm cat})_{k,k+1}+S_{h+}(u_{\rm cat})_{k,k+1}\otimes S_{h+}(u_{\rm cat})_{k+1,k+1}\\ 
S_{h-}(u_{\rm cat})_{k+1,k}\otimes S_{h-}(u_{\rm cat})_{k,k}+S_{h-}(u_{\rm cat})_{k+1,k+1}\otimes S_{h-}(u_{\rm cat})_{k+1,k},
\end{align} 
for $k=1,...,n-1$, on the tensor product $L(\lambda_1)\otimes L(\lambda_2)$. Let $\{\xi_{\Lambda_1}\}$ and $\{\xi_{\Lambda_2}\}$ be the basis of $L(\lambda_1)$ and $L(\lambda_2)$ respectively. Similar to Proposition \ref{qleading}, one can compute the WKB leading term of the operators in \eqref{copro} under the basis $\xi_{\Lambda_1}\otimes \xi_{\Lambda_2}$. (Since the diagonal elements $S_{h+}(u_{\rm cat})_{k,k}$ of quantum Stokes matrices have a rather simple expression, the computation is direct and is omitted here.) Then following the same argument as previous sections, one verifies
\begin{pro}
The WKB approximation of the operators in \eqref{copro} induces the crystal operators $\widetilde{e_k}$ and $\widetilde{f_k}$ on the tensor product $E(u_{\rm cat};\lambda_1)\otimes E(u_{\rm cat};\lambda_2)$ of $\frak{gl}_n-$crystals. 
\end{pro}

\subsection{Proof of Theorem \ref{WKBthm}}
The discussions in this subsection, particularly Proposition \ref{qleading}, show that for each $k=1,...,n-1$, there exists canonical operators $\widetilde{e_k}(u_{\rm cat})$ acting on the finite set $E_{GT}(\lambda)$ such that for any generic element $\xi(u_{\rm cat})\in E_{GT}(\lambda)$, there exist real valued functions $c_{k1}(\xi)$ and $\psi_{k1}(h,u,\xi)$ such that
\begin{eqnarray*}
\mathop{\rm lim}\limits_{h\rightarrow +\infty}\left(S_{h+}(u_{\rm cat})_{k,k+1}\cdot e^{c_{k1}(\xi) h+\I  \psi_{k1}(h,u,\xi)} \xi(u_{\rm cat})\right)=\widetilde{e_k}(\xi(u_{\rm cat})).
\end{eqnarray*}
Then following Proposition \ref{qStokesGZ}, the leading terms ${l}^{(+)}_{k,k+1}(u)$ of the subdiagonal entries of $S_{h+}(u_{\rm cat})$ as $u\rightarrow u_{\rm cat}$ only differ from $S_{h+}(u_{\rm cat})_{k,k+1}$ by some extra terms, appearing in the replacement of $\Delta^{1,...,k}_{1,...,k-1,k+1}\left(\frac{h}{2\pi\I }T(\zeta^{(k)}_i)\right)$ in \eqref{modify0}. As $h\rightarrow +\infty$, the extra terms are fast spin and thus there exist real valued functions $\gamma_{k1}(u,h,\xi)$ such that
\[S_{h+}(u_{\rm cat})_{k,k+1}\cdot e^{\I  \gamma_{k1}(h,u,\xi)}\xi(u_{\rm cat})= {l}^{(+)}_{k,k+1}(u)\cdot \xi(u_{\rm cat}),\]
that is (by Proposition \ref{qStokesGZ})
\[S_{h+}(u_{\rm cat})_{k,k+1}\cdot e^{\I  \gamma_{k1}(h,u,\xi)}\xi(u_{\rm cat})= \mathop{\rm lim}\limits_{u\rightarrow u_{\rm cat}}  S_{h+}(u)_{k,k+1}\cdot \xi(u_{\rm cat}).\]
Therefore, if we let $\theta_{k1}(h,u,\xi):=\psi_{k1}+\gamma_{k1}$, then
\begin{eqnarray*}
\mathop{\rm lim}\limits_{h\rightarrow +\infty}\left(\mathop{\rm lim}\limits_{u\rightarrow u_{\rm cat}} S_{h+}(u)_{k,k+1}\cdot e^{c_{k1}(\xi) h+\I  \theta_{k1}(h,u,\xi)} \xi(u_{\rm cat})\right)=\widetilde{e_k}(\xi(u_{\rm cat})).
\end{eqnarray*}
Similar results hold for $S_{h-}$. 

Furthermore, by Theorem \ref{WKBcrystal} the WKB datum $(E(u_{\rm cat};\lambda)=E_{GT},\widetilde{e_k}(u_{\rm cat}),\widetilde{f_k}(u_{\rm cat}),\varepsilon_k,\phi_k)$ with the weight map $wt$ is a ${\gl}_n$-crystal. It finishes the proof of Theorem \ref{WKBthm}.

\subsection{The cactus group action on the crystals arsing from the wall-crossing formula}\label{cactuswall}

\begin{defi} The Cactus group $Cact_n$ is a group (with a unit) generated by elements $\sigma_{ij}$, $1\le i<j\le n$, subject to the set of relations
\begin{itemize}
    \item $\sigma_{ij}^2=1,$ if $1\le i< j\le n$,
    \item $\sigma_{ij}\sigma_{kl}=\sigma_{kl}\sigma_{ij}$, if $j<k$,
    \item $\sigma_{ij}\sigma_{kl}\sigma_{ij}=\sigma_{i+j-l,i+j-k}$, if $i\le k<l\le j$.
\end{itemize}
\end{defi}
Let us set $\sigma_i:=\sigma_{1,i+1}$, $1\le i\le n-1$. It is clear that $\sigma_i^2=1$, and the elements $\sigma_1,...,\sigma_{n-1}$ generate the Cactus group. The cactus group acts on $\frak{gl}_n$-crystals by the Sch\"utzenberger involution, see e.g., \cite{HKRW}. In particular, if we take the realization of a crystal by the Gelfand-Teitlin basis, the action can be described by the results in \cite{BK, BK2} as follows.
Given the representation $L(\lambda)$ and the Gelfand-Teitlin basis $E_{GT}(\lambda)$,
\begin{thm}\cite{BK,BK2}\label{cactusact}
There is a $Cact_n$ group action on the set of $E_{GT}(\lambda)$ generated by
\begin{equation}
\sigma_i=t_1t_2t_1t_3t_2t_1,...,t_it_{i-1}...t_1, \ \forall i=1,..., n-1,
\end{equation}
where each $t_j$ for $j=1,...,n-1$ is an operator on $E_{GT}(\lambda)$: for any basis element $\xi_\Lambda$, the action $t_j(\xi_\Lambda)$ of $t_j$ on $\xi_\Lambda$ is another basis element $\xi_{\Lambda'}$ with the pattern $\Lambda'$ uniquely determined by 
\begin{align}
\lambda^{(k)}_i(\Lambda')&=\lambda^{(k)}_i(\Lambda), \ \ \ \text{for} \ k\ne j,\\
\lambda^{(j)}_i(\Lambda')&={\rm Max}(\lambda^{(j+1)}_i(\Lambda),\lambda^{(j-1)}_{i-1}(\Lambda))-{\rm Max}(\lambda^{(j)}_i(\Lambda)-\lambda^{(j+1)}_{i+1}(\Lambda),\lambda^{(j)}_i(\Lambda)-\lambda^{(j-1)}_{i}(\Lambda)).
\end{align}
Here we presuppose that $\lambda^{(j)}_0=-\infty$ and $\lambda^{(j)}_{j+1}=+\infty$.
\end{thm}

Recall that the regularized limits $S_{h\pm}(u_{\rm cat})$ simply encode the leading terms of $S_{h\pm}(u)$ as $u\rightarrow u_{\rm cat}$ from the connected component $U_{\rm id}$. For each $1\le i\le n$, let $\tau_i\in S_n$ be the permutation reversing the segment $[1,...,i]$. Then as $u\rightarrow u_{\rm cat}$ from the connected component $U_{\tau_i}=\{u\in \h_{\rm reg}(\mathbb{R})~|~u_{\tau_i(1)}<\cdots <u_{\tau_i(n)}\}$, the regularized limits are given by different $S_{h\pm}^{\tau_i}(u_{\rm cat})\in{\rm End}(L(\lambda))\otimes {\rm End}(\mathbb{C}^n)$. Similar to the Wall-crossing formula at $u_{\rm cat}$ (in the classical case) in Proposition \ref{procac}, we have 
\begin{pro}\label{tauaction}
The off-diagonal elements of $S_{h\pm}^{\tau_i}(u_{\rm cat})$ are
\begin{equation}
S_{h+}^{\tau_i}(u_{\rm cat})_{k,k+1}  =\left\{
          \begin{array}{lr}
             S_{h+}^{\tau_i}(u_{\rm cat})_{i+1-k,i+2-k},   & \text{if} \ \ 1\le k \le i,   \\
           S_{h+}^{\tau_i}(u_{\rm cat})_{k,k+1}, & \text{if} \ \ i< k \le n.
             \end{array}
\right.  
\end{equation}
and 
\begin{equation}
S_{h-}^{\tau_i}(u_{\rm cat})_{k+1,k}  =\left\{
          \begin{array}{lr}
             S_{h-}^{\tau_i}(u_{\rm cat})_{i+2-k,i+1-k},   & \text{if} \ \ 1\le k \le i,   \\
           S_{h-}^{\tau_i}(u_{\rm cat})_{k,k+1}, & \text{if} \ \ i< k \le n.
             \end{array}
\right.  
\end{equation}
\end{pro}

Similar to Theorem \ref{WKBthm}, the WKB approximation of $S_{h\pm}^{\tau_i}(u_{\rm cat})$ also leads to a $\frak{gl}_n$-crystal \[\left(E_{GT}(\lambda),\{\widetilde{e_k}^{\tau_i}(u_{\rm cat})\}_{k}, \{\widetilde{f_k}^{\tau_i}(u_{\rm cat})\}_k\right),\] 
and by Proposition \ref{tauaction},
\begin{align}
    \widetilde{e_k}^{\tau_i}(u_{\rm cat})=\left\{
          \begin{array}{lr}
             \widetilde{e_{i+1-k}}(u_{\rm cat}),   & \text{if} \ \ 1\le k \le i,   \\
           \widetilde{e_k}(u_{\rm cat}), & \text{if} \ \ i< k \le n.
             \end{array}
\right.  \\
 \widetilde{f_k}^{\tau_i}(u_{\rm cat})=\left\{
          \begin{array}{lr}
             \widetilde{f_{i+1-k}}(u_{\rm cat}),   & \text{if} \ \ 1\le k \le i,   \\
           \widetilde{f_k}(u_{\rm cat}), & \text{if} \ \ i< k \le n.
             \end{array}
\right.  
\end{align}

An analog of Theorem \ref{WKBthm} is then
\begin{thm}
For each $k=1,...,n-1$, there exists canonical operators $\widetilde{e_k}^{\tau_i}(u_{\rm cat})$ and $\widetilde{f_k}^{\tau_i}(u_{\rm cat})$ acting on the finite set $E_{GT}(\lambda)$ such that for any generic element $\xi(u_{\rm cat})\in E_{GT}(\lambda)$, there exist real valued functions $c'_{ki}(\xi)$ and $\theta'_{ki}(h,u,\xi)$ with $i=1,2$ such that
\begin{eqnarray*}
\mathop{\rm lim}\limits_{h\rightarrow +\infty}\left(\mathop{\rm lim}\limits_{u\rightarrow u_{\rm cat} \text{ from } U_{\tau_i}} S^{\tau_i}_{h+}(u)_{k,k+1}\cdot e^{c'_{k1}(\xi) h+\I  \theta'_{k1}(h,u,\xi)} \xi(u_{\rm cat})\right)=\widetilde{e_k}^{\tau_i}(\xi(u_{\rm cat})),\\
\mathop{\rm lim}\limits_{h\rightarrow +\infty}\left(\mathop{\rm lim}\limits_{u\rightarrow u_{\rm cat} \text{ from } U_{\tau_i}}  S^{\tau_i}_{h-}(u)_{k+1,k}\cdot e^{c'_{k2}(\xi) h+\I  \theta'_{k2}(h,u,\xi)} \xi(u_{\rm cat})\right)=\widetilde{f_k}^{\tau_i}(\xi(u_{\rm cat})).
\end{eqnarray*}
Furthermore, the set $E_{GT}(\lambda)$ equipped with the operators $\widetilde{e_k}(u_{\rm cat})$ and $\widetilde{f_k}(u_{\rm cat})$ is a ${\gl}_n$-crystal. Here $u\rightarrow u_{\rm cat}$ from $U_{\tau_i}$ means taking the limit $\frac{u_{\tau_i(k+1)}-u_{\tau_i(k)}}{u_{\tau_i(k)}-u_{\tau_i(k-1)}}\rightarrow +\infty$. 
\end{thm}

Now since the actions of $\widetilde{e_k}$, $\widetilde{e_k}^{\tau_i}$, and the action of the generators $\sigma_i$ on $E_{GT}(\lambda)$ are explicitly given. A straightforward but lengthy computation, using the Cauchy inequalities, verifies that
\begin{pro}\label{knowncactus}
The action of the generators $\{\sigma_{i}\}_{i=1,...,n-1}$ on $E_{GT}(\lambda)$ of $Cact_n$ given in Theorem \ref{cactusact} satisfy
\begin{eqnarray}
\sigma_i\circ \widetilde{e_k}=  \widetilde{e_k}^{\tau_i}\circ \sigma_i, \hspace{3mm} \text{ for all } k=1,...,n-1.
\end{eqnarray}
\end{pro}
It gives a proof of Theorem \ref{cactusact0}.

\Addresses
\end{document}